\documentclass[a4paper]{book}
\usepackage{amscd,amsmath,amssymb,amsthm,verbatim}
\newtheorem{Theorem}{Theorem}[section]

\newtheorem{Problem}{Problem}[section]
\newtheorem{Remark}{Remark}[]
\newtheorem{Definition}{Definition}[section]
\newtheorem{Proposition}[Theorem]{Proposition}
\newtheorem{Lemma}[Theorem]{Lemma}
\newtheorem{Corollary}[Theorem]{Corollary}
\begin{document}
\title{Navier-Stokes Equation on the Rectangle\thanks{Supported by FCT (Portuguese Foundation for Science and Tecnology).}}
\author{S\'{e}rgio S. Rodrigues\thanks{SISSA-ISAS, via Beirut 2-4, Trieste 34014, Italy\,$\And$\,University of Aveiro, Portugal; e-mail:
srodrigs@sissa.it}}
\date{April 11, 2005}
\maketitle
\numberwithin{equation}{section}
\tableofcontents
\chapter{Introduction}
Following part of the work iniciated by A. Agrachev and A. Sarychev in
\cite{pas} we study controllability, by means of low modes forcing, of incompressible 2D Navier-Stokes (NS) Equations on the two dimensional rectangle with tangent velocity on the boundary.\par
In the present paper we deal with the 2D NS system
\begin{align}
u_t + (u\cdot\nabla)u+\nabla p &= \nu\Delta u + F(x_1,x_2)+v(t,x_1,x_2)\label{nse1}\\
\nabla\cdot u&=0\label{divo1}\quad\text{on}\quad R\\
u\cdot\mathbf n&= 0\quad\text{on}\quad\partial R\label{tgb1};\\
\nabla^\bot\cdot u&=0\quad\text{on}\quad\partial R\label{w0b1}
\end{align}
Where $R:=\{(x_1,x_2)\in\mathbb R^2\mid a_1< x_1< a_2;\,b_1< x_2
< b_2\}$ and $\nabla^\bot:=\begin{pmatrix}-\frac{\partial}{\partial x_2}\\ \;\;\frac{\partial}{\partial x_1}\end{pmatrix}$ and $\mathbf n$ is the unit normal to the boundary. In the equation \eqref{nse1} $u$ is the \emph{velocity of the fluid ``particle''}; $p$ is the \emph{pressure}; the only nonlinear term of the equation ---  $(u\cdot\nabla)u$ --- is called the \emph{inertial term}; $\nu\Delta u$ is called the \emph{viscosity term}, $\nu>0$ is the \emph{coeficient of viscosity}; $F$ is an \emph{external force} and; $v$ is a control at our disposal. We are interested in the case where $v$ is a degenerate forcing, i.e., $v$ is a finite sum of harmonics --- $v=\sum_{k\in\mathcal K^1}v_k(t)E_k$, where $E_k$ are eigenfunctions of the Stokes operator. So the components $v_k(t),\;k\in\mathcal K^1,\,t\in[0,\,T]$ are our controls which are measurable essentially bounded functions. We shall study Galerkin approximations, say, big enough to contain a set of modes we want to observe.\par
A natural way to study the NSE is to study its evolution on subspaces of Sobolev spaces; such subspaces depend on the boundary conditions.\par
We shall note by $L^2(R)$ the space of Lebesgue measurable square integrable real functions defined on $R$ and by $\mathbf L^2(R)$ the product space $L^2(R)^2$. Similarly $H^1:=\{f\in L^2(R)\mid\,\frac{\partial f}{\partial x_j}\in L^2(R),\,j=1,2 \}$ and, $\mathbf H^1(R):=H^1(R)^2$.\par   
The most studied boundary conditions are the full Dirichlet ($u=0$ on the boundary, also called no-slip) and periodic ($u_1$ and $u_2$ are periodic). For the full Dirichlet conditions the study is done in the spaces
\begin{eqnarray*}
H_D&:=&\{u\in \mathbf L^2(R)\mid \nabla\cdot u=0\,\And\, u\cdot\mathbf n=0\,\text{on}\,\partial R \}\\
V_D&:=&\{u\in \mathbf H^1(R)\mid \nabla\cdot u=0\,\And\, u=0\,\text{on}\,\partial R \}
\end{eqnarray*}
and for the periodic boundary conditions the study is done in         
\begin{eqnarray*}
H_P&:=&\{u\in\mathbf L_{per}^2(R)\mid \nabla\cdot u=0 \}\\
V_P&:=&\{u\in\mathbf H_{per}^1(R)\mid \nabla\cdot u=0 \}.
\end{eqnarray*}
In the periodic case sometimes the special case of zero average space --- $\int_R u\,dx=0$ --- is considered, in this case we have to take this new condition into account in the definition of the spaces $H_P$ and $V_P$. For more information on the boundary conditions above see for example \cite{temams} or \cite{temnsf}.\par
For the boundary conditions \eqref{tgb1}, \eqref{w0b1} the spaces 
\begin{eqnarray}
H&:=&\{u\in \mathbf L^2(R)\mid \nabla\cdot u=0\,\And\, u\cdot\mathbf n=0\,\text{on}\,\partial R \};\notag\\
V&:=&\text{closure of}\;\mathcal D_1(R)\;\text{on}\;\mathbf H^1(R);\label{spaces}\\
D(A)&:=&\{u\in\mathbf H^2(R)\mid\,\nabla\cdot u=0\,\And\, (u\cdot\mathbf n=0\wedge\nabla^\bot\cdot u=0)\,\text{on}\,\partial R\};\notag
\end{eqnarray}
where $\mathcal D_1(R):=\{u\in C^\infty(\overline R)\mid\,\nabla\cdot u=0\And (u\cdot\mathbf n=0\,\wedge\,\nabla^\bot\cdot u=0)\,\text{on}\,\partial R\}$,
are those where we shall consider the evolution of the NSE on.\par
In Chapter \ref{Ch:EUC} we prove the existence of weak and strong solutions for \eqref{nse1}-\eqref{w0b1}, as well as its uniqueness and continuous dependence in the initial data. For the case of weak solutions we proceed as in \cite{temams} for the case of no-slip boundary conditions. To obtain strong solutions we just have to ask for some regularity on the initial data.\par
In Chapter \ref{Ch:Gal}, for strong solutions, we prove controllability, by means of low modes forcing, of Galerkin approximations of the infinite-dimensional system associated with \eqref{nse1}-\eqref{w0b1}.\\
In Chapter \ref{Ch:COP} we prove the so called \emph{controllability in observed component}, again by means of low modes forcing. In other words we prove that the projection, onto any finite dimensional subspace (the space spanned by the modes we want to observe), of the attainable set from any point is surjective.\par
Finally, we end with Chapter \ref{L2AC} and with the prove of $\mathbf L^2$-{\em Approximate Controllability} the is a straightforward corollary of some tools we have presented in Chapter \ref{Ch:COP}.\\ \ \par
The author is gratefull to A. Agrachev and A. Sarychev for the inspiring and helpfull discussions on the subject and, for the sugestions in the improvement of the text.\par
The author would like to thank FCT (Portuguese Foundation for Science and Tecnology) for financial support and, SISSA-ISAS (International School for Advanced Studies) for hospitality. 

\chapter{Existence, Continuity and Uniqueness.}\label{Ch:EUC}

\section{The Spaces.}\label{SHV}
\subsection{Recollection of Auxiliary Material on the Spaces $\mathbf L^2(R)$ and $\mathbf H^1(R)$.}
Recall that since $L^2(R)$ is a Hilbert space for the scalar product
\begin{equation}\label{scp(,)1}
(u,\,v)_1:=\int_R uv\,dx,
\end{equation}
then $\mathbf L^2(R)$ is a Hilbert space for the product topology and the scalar product is
\begin{equation}\label{scp(,)}
(u,\,v):=\int_R u\cdot v\,dx.
\end{equation}
We note that
$$
(u,\,v)=\int_R u\cdot v\,dx=\sum_{i=1}^2\int_R u_iv_i\,dx=(u_1,\,v_1)_1+(u_2,\,v_2)_1.
$$
The norms associated with the previous scalar products shall be represented by
\begin{eqnarray}
|u|_1^2&:=&(u,\,u)_1\label{n||1}\\
|u|^2&:=&(u,\,u)\label{n||}
\end{eqnarray}
We note that
$$
|u|^2:=|u_1|_1^2+|u_2|_1^2
$$
so, the norm $|\cdot|$ is the product norm $|\cdot|_1\times|\cdot|_1$.
Similarly, the Sobolev space $H^1(R)$ is a Hilbert space for the scalar product\begin{equation}\label{scp((,))01}
((u,\,v))_{01}:=(u,\,v)_1+\sum_{i=1}^2(\frac{\partial u}{\partial x_i},\,\frac{\partial v}{\partial x_i})_1=(u,\,v)_1+(\nabla u,\,\nabla v).
\end{equation}
then $\mathbf H^1(R)$ is a Hilbert space for the product topology and the scalar product is
\begin{equation}\label{scp((,))0}
((u,\,v))_0:=\sum_{j=1}^2((u_j,\,v_j))_{01}=(u,\,v)+\sum_{i=1}^2 (\nabla u_i,\,\nabla v_i)
\end{equation}
The norms associated with the previous scalar products shall be represented by
\begin{eqnarray}
\|u\|_{01}^2&:=&((u,\,u))_{01}\label{n||||01}\\
\|u\|_0^2&:=&((u,\,u))_0\label{n||||0}
\end{eqnarray}
We note that
$$
\|u\|_0^2:=\|u_1\|_{01}^2+\|u_2\|_{01}^2
$$
so, the norm $\|\cdot\|_0$ is the product norm $\|\cdot\|_{01}\times\|\cdot\|_{01}$.
\subsection{Some Properties of the Spaces $V$ and $H$.}

\begin{Lemma}\label{HVclsmooth}
$H$ coincides with the closure of $\mathcal D_1(R)$ in $\mathbf L^2(R)$.
\end{Lemma}
\begin{proof}
It is well known from the study of the NSE with full Dirichlet boundary conditions that $H$ is the closure of $\mathcal D(R)$ in $\mathbf L^2(R)$, where $\mathcal D(R)$ is the set of solenoidal (or divergence free) smooth functions with support in $R$. Since
$$
\mathcal D(R)\subset\mathcal D_1(R)\subset H
$$
we conclude that $H$ coincides with the closure of $\mathcal D_1(R)$ in $\mathbf L^2(R)$.
\end{proof}
We use the same argument, since $\mathcal D_1\subset V\subset H$ to conclude that
\begin{Corollary}\label{VdenseH}
$H$ is the closure of $V$ in $\mathbf L^2$.\;\footnote{In other words $V$ is dense in the closed space $H$}
\end{Corollary}
In the study of NSE some classical imbedding and compactness theorems are frequently used, we start by presenting some of them which we shall need (for our particular equation).
\begin{Proposition}\label{imb.H1-Lq(1)}
For all $u\in H^1(R)$ and for all $q,\,1\leq q<\infty$ we have
$$
\|u\|_{L^q(R)}\leq C_1(q)\|u\|_{01}.\quad
$$
\end{Proposition}

\begin{Corollary}\label{imb.H1-Lq}
For all $u\in \mathbf H^1(R)$ and for all $q,\,1\leq q<\infty$ we have
$$
\|u\|_{\mathbf L^q(R)}\leq C_1(q)\|u\|_0.
$$
\end{Corollary}
\begin{Proposition}\label{comp.H1-L2(1)}
The imbedding
\begin{align*}
H^1(R)\,&\to\,L^2(R)\\
u\,&\mapsto\, u
\end{align*}
is compact.
\end{Proposition}
\begin{Corollary}\label{comp.H1-L2}
The imbedding
\begin{align*}
\mathbf H^1(R)\,&\to\, \mathbf L^2(R)\\
u\,&\mapsto\, u
\end{align*}
is compact.
\end{Corollary}
\begin{Corollary}\label{VctcpH}
The inclusion
\begin{align}
i:V\,&\to\,H\label{inc.V-H}\\
u\,&\mapsto\, u\notag
\end{align}
is continuous and compact.
\end{Corollary}
\begin{proof}
By Corollary \ref{imb.H1-Lq}, inclusion \eqref{inc.V-H} is continuous and, by Corollary \ref{comp.H1-L2} it is compact.
\end{proof}
 \subsection{Poincar\'{e} Inequality. Equivalent Norms.}
We have the following Lemma:\;\footnote{See \cite{temsp68} subsections II.1.4 and III.2.2}
\begin{Lemma}\ \\
\begin{itemize}
\item For any seminorm $p$ in $\mathbf H^1(R)$ satisfying ``$(p(a)=0\wedge a\in\mathbb R^2)\Rightarrow a=0$'' we have
$$
|u|\leq c\Bigl(\|\nabla u\|_{\mathbf L^2(R)^2}+p(u)\Bigr),\quad\forall u\in \mathbf H^1(R)
$$
and,
\item The seminorm
$$
p(u):=\int_\Gamma |u\cdot\mathbf n|\,d\Gamma
$$
satisfies the required condition.
\end{itemize}
\end{Lemma}
\begin{Remark}
The product space $\mathbf L^2(R)^2$ is a Hilbert space for the scalar product
$$
(U,\,V):=((U_1,\,V_1))_0+((U_2,\,V_2))_0,\quad U=\begin{pmatrix}U_1\\U_2\end{pmatrix},\,V=\begin{pmatrix}V_1\\V_2\end{pmatrix}\in\mathbf L^2(R)^2
$$
and, $\nabla\begin{pmatrix}V_1\\V_2\end{pmatrix}:=\begin{pmatrix}\nabla V_1\\\nabla V_2\end{pmatrix}\quad\forall V\in\mathbf L^2(R)^2$.
\end{Remark}
A consequence of this lemma is a Poincar\'e-like inequality:
\begin{equation}
\forall u\in V[\,|u|\leq c\|\nabla u\|_{\mathbf L^2(R)^2}\,].\label{poinc.in}
\end{equation}
\begin{Corollary}\label{eq.norms.V}
The norms $\|u\|:=\|\nabla u\|_{\mathbf L^2(R)^2}$ and $\|\cdot\|_0$ are equivalent in $V$.\footnote{Note that $\|\cdot\|_0=|\cdot|+\|\cdot\|$.}
\end{Corollary}
\section{The Duals of $V$ and $H$.}\label{S.VH}
From now we shall consider the space $V$ endowed with the norm
$\|\cdot\|$. From Corollaries \ref{VdenseH}, \ref{VctcpH} and
\ref{eq.norms.V} we obtain that the
inclusion
\eqref{inc.V-H}(considering $V$ endowed with
$\|\cdot\|$) is dense continuous and compact.\par Since $H$ is an
Hilbert Space due to the Riesz Representation Theorem we can
identify $H$ with its dual $H^\prime$. In this way we arrive to
the inclusions
\begin{equation}\label{VHH'V'}
V\subset H\equiv H^\prime\subset V^\prime
\end{equation}
where both inclusions are linear, dense and continuous. Indeed we already know that the first inclusion has these properties. In the second inclusion we consider the map $i^\prime$, dual to the map $i$ defined in \eqref{inc.V-H},
\begin{align*}
i^\prime:H^\prime &\to V^\prime\\
<i(v),\,h>\,&=\,<v,\,i^\prime(h)>,
\end{align*}
that is linear and continuous --- we consider $V^\prime$ endowed with the classical norm $\|f\|:=\sup_{\|v\|=1}|<f,\,v>|$. The density of $H^\prime$ in $V^\prime$ is a corollary of the following Lemma, which one can find in \cite{bre}, section (II.6),
\begin{Lemma}\label{NR}
Let $G:D(G)\subset E\,\to\, F$ be a (nonbounded) closed operator with $\overline{D(G)}=E$. Then
\begin{eqnarray*}
(i)& &N(G)=R(G^\prime)^\bot;\\
(ii)& &N(G^\prime)=R(G)^\bot;\\
(iii)& &N(G)^\bot\supset\overline{R(G^\prime)};\\
(iv)& &N(G^\prime)^\bot=\overline{R(G)}.
\end{eqnarray*}
\end{Lemma}
Where $\overline{X}$ means closure of $X$; $N(G)$ stays for the Kernel of $G$ and $R(G)$ for the image of $G$. Applying this Lemma to our inclusions the injectivity of the first inclusion implies the density of the second --- $(i)$. Moreover we also see that, since the first inclusion is dense, it follows the injectivity of the second $(ii)$.
\begin{Remark}
The domain $D(i)$ of $i$ is $V$. By the continuity of $i$ we conclude its closedness. Recall that a operator $G:E\to F$ is said closed if for any sequence $(u_n)$ in $D(G)$ such that $\begin{cases} \;\;u_n\to u&\text{in}\;E\\ Au_n\to f&\text{in}\;F\end{cases}$ we have  $\begin{cases} \;\;u&\in\,D(G)\\ Au&=\,f\end{cases}$.
\end{Remark}

\section{Fourier Series. A Basis in $\mathcal D_1(R)$ for $H$ and $V$.}\label{Sfourier}

We start by noting that if in the case of the general rectangle $R:=[a_1,\,a_2]\times[b_1,\,b_2]\quad a_1<a_2,\,b_1<b_2$ we make the change of variables
$$
z_1:=x_1-a_1,\quad z_2:=x_2-b_1,
$$
since the differential operators $\nabla,\,\Delta$ are invariant under tranlations, the Navier-Stokes Equation does not change in the variables $(z_1,\,z_2)$. To simplify the exposition from now we shall deal with the rectanlge
$$
R:=\{(x_1,x_2)\in\mathbb R^2\mid 0\leq x_1\leq a;\,0\leq x_2\leq b\}.
$$
First we note that under the condition $u\cdot\mathbf n=0$ on $\partial R$, the condition $\nabla^\bot\cdot u=0$ on $\partial R$ equivals $\frac{\partial u_1}{\partial x_2}=\frac{\partial u_2}{\partial x_1}=0$ on $\partial R$, because on the left and right faces we have $\frac{\partial u_1}{\partial x_2}=0$ and, on the top and botton faces $\frac{\partial u_2}{\partial x_1}=0$. It is well known that $\{\sin(\frac{n\pi x}{L})\mid\,n\in\mathbb N_0\}$ is a basis for the functions in $L^2([0,\,L])$ vanishing on $0$ and $L$ and, $\{\cos(\frac{n\pi x}{L})\mid\,n\in\mathbb N\}$ is a basis for the functions in $L^2([0,\,L])$ with vanishing velocity on the boundary points $0$ and $L$.\footnote{Here $\mathbb N$ is the set of natural numbers and $\mathbb N_0:=\mathbb N\setminus\{0\}$} So
$$
\Bigl\{\sin\biggl(\frac{n\pi x}{a}\biggr)\cos\biggl(\frac{m\pi x}{b}\biggr)\mid\,n\in\mathbb N_0,\,m\in\mathbb N\Bigr\}
$$
is a basis for $L^2(R)$ vanishing on the right and left faces of the rectangle and,
$$
\Bigl\{\cos\biggl(\frac{n\pi x}{a}\biggr)\sin\biggl(\frac{m\pi x}{b}\biggr)\mid\,n\in\mathbb N,\,m\in\mathbb N_0\Bigr\}
$$
is a basis for $L^2(R)$ vanishing on the botton and top faces of the rectangle.\par
Now consider a function $u\in H$, this function, being an element of $\mathbf L^2(R)$, can be written as
$$
u=\begin{pmatrix} u_1\\u_2\end{pmatrix}=
\begin{pmatrix}
\sum_{\begin{subarray}{l}n\in\mathbb N_0\\ m\in\mathbb N\end{subarray}}u_{1k}\sin\bigl(\frac{n\pi x_1}{a}\bigr)\cos\bigl(\frac{m\pi x_2}{b}\bigr)\\
\sum_{\begin{subarray}{l}m\in\mathbb N_0\\ n\in\mathbb N\end{subarray}}u_{2k}\cos\bigl(\frac{n\pi x_1}{a}\bigr)\sin\bigl(\frac{m\pi x_2}{b}\bigr)
\end{pmatrix}.
$$
Since $u$ is solenoidal --- $\nabla\cdot u=0$, we conclude that
$$
u=\begin{pmatrix} u_1\\u_2\end{pmatrix}=
\begin{pmatrix}
\sum_{k\in\mathbb N_0^2}u_{1k}\sin\bigl(\frac{k_1\pi x_1}{a}\bigr)\cos\bigl(\frac{k_2\pi x_2}{b}\bigr)\\
\sum_{k\in\mathbb N_0^2}u_{2k}\cos\bigl(\frac{k_1\pi x_1}{a}\bigr)\sin\bigl(\frac{k_2\pi x_2}{b}\bigr)
\end{pmatrix}.
$$
and that for each $k\in\mathbb N_0^2$ we have
$$
u_{1k}\frac{k_1\pi}{a}+u_{2k}\frac{k_2\pi}{b}=0.
$$
This means that $\frac{u_{1k}b}{-k_2\pi}=\frac{u_{2k}a}{k_1\pi}=:u_k$ and, we arrive to
$$
u=\sum_{k\in\mathbb N_0^2}u_kW_k,\quad W_k:=
\begin{pmatrix}
\frac{-k_2\pi}{b}\sin\bigl(\frac{k_1\pi x_1}{a}\bigr)\cos\bigl(\frac{k_2\pi x_2}{b}\bigr)\\
\frac{k_1\pi}{a}\cos\bigl(\frac{k_1\pi x_1}{a}\bigr)\sin\bigl(\frac{k_2\pi x_2}{b}\bigr)
\end{pmatrix}.
$$
We put
\begin{equation}\label{basisW}
\mathcal W:=\{W_k\mid\,k\in\mathbb N_0^2\}.
\end{equation}
\subsection{Fourier Characterization of $H$ and $V$.}
We start by computing the scalar product $(W_k,\,W_z)$ between two elements of $\mathcal W$:
\begin{align*}
(W_k,\,W_z)&=\int_R W_k\cdot W_z\,dx\\
&=\int_R \frac{k_2z_2\pi^2}{b^2}\sin\bigl(\frac{k_1\pi x_1}{a}\bigr)\cos\bigl(\frac{k_2\pi x_2}{b}\bigr)\sin\bigl(\frac{z_1\pi x_1}{a}\bigr)\cos\bigl(\frac{z_2\pi x_2}{b}\bigr)\,dx\\
&+\int_R\frac{k_1z_1\pi^2}{a^2}\cos\bigl(\frac{k_1\pi x_1}{a}\bigr)\sin\bigl(\frac{k_2\pi x_2}{b}\bigr)\cos\bigl(\frac{z_1\pi x_1}{a}\bigr)\sin\bigl(\frac{z_2\pi x_2}{b}\bigr)\,dx\\
&=\begin{cases}
0&\text{if}\;k\neq z\\
-\bar k\frac{ab}{4}&\text{if}\;k=z.\quad\footnotemark
\end{cases}
\end{align*}
\footnotetext{This follows from the fact that the family of sines is orthogonal in $L^2([0,\,\pi])$, as well is the family of cosines.}
where $\bar k:=-\pi^2\Big(\frac{k_1^2}{a^2}+\frac{k_2^2}{b^2}\Big)$. Hence the square of the $\mathbf L^2(R)$ norm of each element in the orthogonal family $\mathcal W$ equals $-\bar k\frac{ab}{4}$.\par
Given $u,\,v\in H$ the scalar product between them is then
$$
(u,\,v)=\sum_{k\in\mathbb N_0^2}-\bar k\frac{ab}{4}u_kv_k,
$$
where
$$
u=\sum_{k\in\mathbb N_0^2}u_kW_k;\quad v=\sum_{k\in\mathbb N_0^2}v_kW_k.
$$
The norm of $u$ in $H$ results
\begin{equation}\label{F.norm.H}
|u|^2=\frac{ab}{4}\sum_{k\in\mathbb N_0^2}-\bar ku_k^2
\end{equation}
and $u\in H$ means that
$$
\sum_{k\in\mathbb N_0^2}-\bar ku_k^2<+\infty.
$$\par
If $u,\,v\in V$ we want to compute the scalar product
$$
((u,\,v))=((\nabla u_1,\,\nabla v_1))+((\nabla u_2,\,\nabla v_2)).
$$
First we note that
$$
\nabla u_1=
\begin{pmatrix}
\sum_{k\in\mathbb N_0^2}u_k\bigl(-\frac{k_2\pi}{b}\frac{k_1\pi}{a}\bigr)\cos\bigl(\frac{k_1\pi x_1}{a}\bigr)\cos\bigl(\frac{k_2\pi x_2}{b}\bigr)\\
\sum_{k\in\mathbb N_0^2}u_k\bigl(-\frac{k_2\pi}{b}\frac{-k_2\pi}{b}\bigr)\sin\bigl(\frac{k_1\pi x_1}{a}\bigr)\sin\bigl(\frac{k_2\pi x_2}{b}\bigr)
\end{pmatrix},
$$
$$
\nabla u_2=
\begin{pmatrix}
\sum_{k\in\mathbb N_0^2}u_k\bigl(\frac{k_1\pi}{a}\frac{-k_1\pi}{a}\bigr)\sin\bigl(\frac{k_1\pi x_1}{a}\bigr)\sin\bigl(\frac{k_2\pi x_2}{b}\bigr)\\
\sum_{k\in\mathbb N_0^2}u_k\bigl(\frac{k_1\pi}{a}\frac{k_2\pi}{b}\bigr)\cos\bigl(\frac{k_1\pi x_1}{a}\bigr)\cos\bigl(\frac{k_2\pi x_2}{b}\bigr)
\end{pmatrix}
$$
and we obtain analogous expressions for $\nabla v_1$ and $\nabla v_2$. Hence, doing the computations, we arrive to
$$
((u,\,v))=\sum_{k\in\mathbb N_0^2}\bar k^2\frac{ab}{4}u_kv_k.
$$
The norm of $u$ in $V$ results
\begin{equation}\label{F.norm.V}
\|u\|^2=\frac{ab}{4}\sum_{k\in\mathbb N_0^2}\bar k^2u_k^2
\end{equation}
and $u\in V$ means that
$$
\sum_{k\in\mathbb N_0^2}\bar k^2u_k^2<+\infty.
$$
From caracterizations \eqref{F.norm.H} and \eqref{F.norm.V} we can see easily, as referred in Corollary \ref{eq.norms.V}, that the norms $\|\cdot\|_0$ and $\|\cdot\|$ are equivalent in $V$. Indeed, $-\bar k$ goes to $+\infty$ when $\max\{k_1,\,k_2\}$ does\footnote{For $M:=\max\{k_1,\,k_2\}$ and $m:=max\{a^2,\,b^2\}$, we have $-\bar k=\pi^2(\frac{k_1^2}{a^2}+\frac{k_2^2}{b^2})\geq \pi^2\frac{M^2}{m}$ that goes to $\infty$ when $M$ does.}  so, there is $N\in\mathbb N_0$ and a constant $C\geq 1$ such that if $\max\{k_1,\,k_2\}>N$ we have $-\bar k\geq 1$ and,
\begin{equation}\label{const.k.k2}
\sum_{\max{\{k_1,\,k_2\}}\leq N}-\bar k\leq C\sum_{\max{\{k_1,\,k_2\}}\leq N}\bar k^2.
\end{equation}
Therefore $|u|\leq C\|u\|$ and then $\|\cdot\|_0$ and $\|\cdot\|$ are equivalent.\par
When $u\in V$ we can not guarantee, in general, neither $\Delta u\in H$ nor $\Delta u\in\mathbf L^2$. Anyway if we compute $(\Delta u,\,v)$ for $u,\,v\in V$ we obtain
\begin{equation}\label{delta.(())}
(\Delta u,\,v)=\sum_{k\in\mathbb N_0^2}\bar k^2u_k^2<+\infty=-((u,\,v)).
\end{equation}
Hence the result of $(\Delta u,\,v)$ is finite and is properly defined even if $\Delta u\notin\mathbf L^2$ and so the operation $(\Delta u,\,v)$ seems to have no sense in itself. Below we will give a sense to it.
\section{The Operators $A$ and $B$.}
\subsection{The Operator $A$.}
First we note that a consequence of the identifications \eqref{VHH'V'}, the scalar product in $H$ of $f\in H$ and $u\in V$ is the same as the scalar product of $f$ and $u$ in the duality between $V$ and $V^\prime$
\begin{equation}\label{scalH=scalVV'}
\forall f\in H\;\forall u\in V[\,< f,\,u>=(f,\,u)\,].
\end{equation}
For each $u\in V$, the form
\begin{align}
V\,&\to\,\mathbb R\\
v\,&\mapsto\, ((u,\,v))\notag
\end{align}
is linear and continuous on $V$. Therefore there is an element of $V^\prime$ we shall denote by $Au$ satisfying
\begin{equation}\label{defA}
< Au,\,v>=((u,\,v)),\quad \forall v\in V
\end{equation}
and $A$ is clearly linear and continuous; we can also see that $A$ is an isomorphism between $V$ and $V^\prime$. Indeed the form $a(u,\,v):=((u,\,v))$ is bilinear, continuous and coersive on $V\times V$ so, by Lax-Milgram Theorem we can conclude that $A$ is an isomorphism between $V$ and $V^\prime$.\par
By equation \eqref{delta.(())} we see that
\begin{equation}\label{<A>=(delta)}
< Au,\,v>=(-\Delta u,\,v)\quad\forall u,\,v\in V.
\end{equation}
\subsubsection{Further Properties of $A$.}
As we have said above for $v\in V$, $Av$ can be not in $\mathbf L^2$. Now we define a subset $D(A)$ of $V$, we shall call the domain of $A$, such that
$$
A:D(A)\to H.
$$
The domain of $A$ is defined by
\begin{align*}
D(A)&:=\{v\in\mathbf H^2\mid\,v\in V\And \nabla^\bot\cdot v=0\,\text{on}\,\Gamma\}\\&=\{v\in\mathbf H^2\mid\,\nabla\cdot v=0\And( v\cdot\mathbf n=0\,\wedge\,\nabla^\bot\cdot v=0\,\text{on}\,\Gamma)\}\\
&=\,\text{closure of}\,\mathcal D_1(R)\,\text{in}\,\mathbf H^2.
\end{align*}
Where $\Gamma:=\partial R$ and $\mathbf H^2$ is the Sobolev space $\mathbf H^2=(W^{2,2})^2$ where
$$
W^{2,2}:=\{u\in L^2\mid\,\frac{\partial^{|\alpha|}u}{\partial x^\alpha}\in L^2,\quad\forall |\alpha|\leq 2\}
$$
which is a Hilbert space for the scalar product
$$
(u,\,v)_2:=\sum_{\begin{subarray}{l}i\in\{1,2\}\\|\alpha|\leq 2\end{subarray}}(\frac{\partial^{|\alpha|}u_i}{\partial x^\alpha},\,\frac{\partial^{|\alpha|}v_i}{\partial x^\alpha})_1=\sum_{|\alpha|\leq 2}(\frac{\partial^{|\alpha|}u}{\partial x^\alpha},\,\frac{\partial^{|\alpha|}v}{\partial x^\alpha}).
$$
The norm corresponding to this scalar product is defined by
$$
|u|_2^2=|u|^2+\|u\|^2+|u|_{[2]}^2,
$$
with
$$
|u|_{[2]-}^2=\sum_{|\alpha|= 2}|\frac{\partial^{|\alpha|}u}{\partial x^\alpha}|^2=\sum_{|\alpha|= 2}\Bigl(|\partial_1^2 u|^2+|\partial_2^2 u|^2+|\partial_{1,2}^2 u|^2\Bigr)\quad\footnotemark
$$
\footnotetext{$\partial_j^2$ stays for $\frac{\partial^2}{\partial x_j^2}$ and, $\partial_{1,2}^2$ stays for $\frac{\partial^2}{\partial x_1\partial x_2}$.}
We can see that the norm $|\cdot|_2$ is equivalent to that $|\cdot|_{2+}$ coming from the scalar product defined by
$$
(u,\,v)_{2+}=(u,\,v)_{2}+(\partial_{1,2}^2u,\,\partial_{1,2}^2v).
$$
Indeed $|u|_2^2\leq|u|_{2+}^2\leq2|u|_2^2$.\par
Now we look for a Fourier caracterization of the elements on $D(A)$: For that we need to write down the expressions for the second order derivatives of $u$:
\begin{eqnarray*}
\partial_{1,2}^2 u&=&\sum_{k\in\mathbb N_0^2}u_k
\begin{pmatrix}
-\frac{k_2\pi}{b}\Bigl(-\frac{k_2\pi}{b}\Bigr)\frac{k_1\pi}{a}\cos\bigl(\frac{k_1\pi x_1}{a}\bigr)\sin\bigl(\frac{k_2\pi x_2}{b}\bigr)\\
\frac{k_1\pi}{a}\Bigl(-\frac{k_1\pi}{a}\Bigr)\frac{k_2\pi}{a}\sin\bigl(\frac{k_1\pi x_1}{a}\bigr)\cos\bigl(\frac{k_2\pi x_2}{b}\bigr)
\end{pmatrix};\\
\partial_1^2 u&=&\sum_{k\in\mathbb N_0^2}u_k
\begin{pmatrix}
\frac{k_2\pi}{b}\Bigl(\frac{k_1\pi}{a}\Bigr)^2\sin\bigl(\frac{k_1\pi x_1}{a}\bigr)\cos\bigl(\frac{k_2\pi x_2}{b}\bigr)\\
-\frac{k_1\pi}{a}\Bigl(\frac{k_1\pi}{a}\Bigr)^2\cos\bigl(\frac{k_1\pi x_1}{a}\bigr)\sin\bigl(\frac{k_2\pi x_2}{b}\bigr)
\end{pmatrix};\\
\partial_2^2 u&=&\sum_{k\in\mathbb N_0^2}u_k
\begin{pmatrix}
\frac{k_2\pi}{b}\Bigl(\frac{k_2\pi}{b}\Bigr)^2\sin\bigl(\frac{k_1\pi x_1}{a}\bigr)\cos\bigl(\frac{k_2\pi x_2}{b}\bigr)\\
-\frac{k_1\pi}{a}\Bigl(\frac{k_2\pi}{b}\Bigr)^2\cos\bigl(\frac{k_1\pi x_1}{a}\bigr)\sin\bigl(\frac{k_2\pi x_2}{b}\bigr)
\end{pmatrix}.
\end{eqnarray*}
Computing $(u,\,v)_{2+}$ for $u,\,v\in D(A)$ we obtain
\begin{align*}
(u,\,v)_{2+}&=\frac{ab}{4}\sum_{k\in\mathbb N_0^2}-\bar ku_kv_k+\frac{ab}{4}\sum_{k\in\mathbb N_0^2}\bar k^2u_kv_k\\
&+\frac{ab}{4}\sum_{k\in\mathbb N_0^2}\pi^6\Biggl[2\Bigl(\frac{k_2^2k_1}{b^2a}\Bigr)^2+2\Bigl(\frac{k_1^2k_2}{a^2b}\Bigr)^2+\Bigl(\frac{k_2k_1^2}{ba^2}\Bigr)^2\\
&\qquad\qquad\qquad\qquad\qquad+\Bigl(\frac{k_2^2k_1}{b^2a}\Bigr)^2+\Bigl(\frac{k_1^3}{a^3}\Bigr)^2+\Bigl(\frac{k_2^3}{b^3}\Bigr)^2\Biggr]u_kv_k\\
&=\frac{ab}{4}\sum_{k\in\mathbb N_0^2}(-\bar k+\bar k^2-\bar k^3)u_kv_k.
\end{align*}
Hence for an element $u$ of $D(A)$:
$$
|u|_{2+}^2=\frac{ab}{4}\sum_{k\in\mathbb N_0^2}(-\bar k+\bar k^2-\bar k^3)u_k^2.
$$
Using the constants $C$ and $N$ of the equation \eqref{const.k.k2} and choosing $C_1$ such that
\begin{equation}\label{const.k2.k3}
\sum_{\max\{k_1,\,k_2\}\leq N}\bar k^2\leq C_1\sum_{\max\{k_1,\,k_2\}\leq N}-\bar k^3,
\end{equation}
we conclude that
$$
\sum_{k\in\mathbb N_0^2}-\bar k^3\leq\sum_{k\in\mathbb N_0^2}(-\bar k+\bar k^2-\bar k^3)\leq ((C+1)C_1+1)\sum_{k\in\mathbb N_0^2}-\bar k^3.
$$
We define
$$
(u,\,v)_{[2]}:=(u,\,v)_{[2]-}+(\partial_{1,2}^2u,\,\partial_{1,2}^2v).
$$
Then all the norms $|u|_{2+};\; |u|_2;\; |u|_{[2]-}$ and $|u|_{[2]}$ are equivalent. Moreover for the last one we have the nice representations:
\begin{align*}
(u,\,v)_{[2]}&=\frac{ab}{4}\sum_{k\in\mathbb N_0^2}-\bar k^3u_kv_k;\\
|u|_{[2]}^2&=\frac{ab}{4}\sum_{k\in\mathbb N_0^2}-\bar k^3u_k^2
\end{align*}
Now we compute
$$
\nabla(\nabla\cdot u)+\nabla^\bot(\nabla^\bot\cdot u),\quad u=\begin{pmatrix} u_1\\u_2\end{pmatrix}:
$$
\begin{align*}
\nabla(\nabla\cdot u)+\nabla^\bot(\nabla^\bot\cdot u)=&\nabla(\partial_1u_1+\partial_2u_2)+\nabla^\bot(-\partial_2u_1+\partial_1u_2)\\
=&\begin{pmatrix}\partial_1^2u_1+\partial_{1,2}^2u_2 \\\partial_{1,2}^2u_1+\partial_2^2u_2\end{pmatrix}+
\begin{pmatrix}\partial_2^2u_1-\partial_{1,2}^2u_2 \\-\partial_{1,2}^2u_1+\partial_1^2u_2\end{pmatrix}\\
=&\begin{pmatrix}\partial_1^2u_1+\partial_2^2u_1 \\\partial_1^2u_2+\partial_2^2u_2\end{pmatrix}
\end{align*}
so
\begin{equation}\label{D.n.nbot}
\Delta u=\nabla(\nabla\cdot u)+\nabla^\bot(\nabla^\bot\cdot u)
\end{equation}
and, for $u$ solenoidal we have $\Delta u=\nabla^\bot(\nabla^\bot\cdot u)$. Moreover for each $u\in D(A)$ we have
\begin{align}
\nabla\cdot\Delta u&=\nabla\cdot\nabla^\bot(\nabla^\bot\cdot u)=0;\notag\\
\mathbf n\cdot\Delta u&=\mathbf n\cdot\nabla^\bot(\nabla^\bot\cdot u)=0\,\text{on}\,\Gamma;\notag\\
\Delta u&\in H.\label{uDA.AuH}\quad\footnotemark
\end{align}
\footnotetext{Since $u\in\mathbf H^2,\;\Delta u\in\mathbf L^2$.}
The operator $A$ coincides with $-\Delta$ in $D(A)$ because, due to \eqref{scalH=scalVV'} and to the definition of $A$ we have
$$
< Au,\,v>=((u,\,v))=(-\Delta u,\,v)=<-\Delta u,\,v>,\quad\forall v\in V.
$$
\subsection{The Operator $B$.}
We define the form
\begin{align}
b:(\mathbf H^1(R))^3&\to\mathbb R\notag\\
(u,\,v,\,w)&\mapsto\sum_{i,j=1}^2\int_R u_i(\partial_iv_j)w_j\,dx.\label{def.b}\
\end{align}

\begin{Lemma}
The form $b$ defined in \eqref{def.b} is trilinear and continuous on the product space $(\mathbf H^1(R))^3$.
\end{Lemma}
\begin{proof}
The trilinearity is clear: $b(u,\,v,\,w)=((u\cdot\nabla)v,\,w)$. For the continuity we start to compute $|\int_R u_i(\partial_iv_j)w_j\,dx|$: Since $u,\,v,\,w\in\mathbf H^1(R)$, by proposition \ref{imb.H1-Lq(1)}, we have that
$$
u_i,\,w_j\in L^4(R)\And \partial_iv_j\in L^2(R).
$$
By H$\ddot{\text{}o}$lder Inequality we obtain that $u_i(\partial_iv_j)w_j\in L^1(R)$ and that
$$
|\int_R  u_i(\partial_iv_j)w_j\,dx|\leq \int_R |u_i(\partial_iv_j)w_j|\,dx\leq |u_i|_{L^4}|\partial_iv_j|_{L^2}|w_j|_{L^4}.
$$
Therefore
$$
|b(u,\,v,\,w)|\leq |u|_{\mathbf L^4}|\nabla v|_{(\mathbf L^2)^2}|w|_{\mathbf L^4}.
$$
So, by Corollaries \ref{imb.H1-Lq} and \ref{eq.norms.V}
$$
|b(u,\,v,\,w)|\leq C\|u\|\|v\|\|w\|.
$$
\end{proof}
The form $b$, being continuous in $(\mathbf H^1(R))^3$, is continuous in $V^3$ and, for each pair $(u,\,v)\in V^2$ we define the operator $B(u,\,v)\in V^\prime$ by
\begin{align}
B(u,\,v):V&\to\mathbb R\\
w&\mapsto <B(u,\,v),\,w>=b(u,\,v,\,w)
\end{align}
and we set
$$
B(u):=B(u,\,u)\in V^\prime\quad\forall u\in V.
$$

\section{Classical and Weak Formulations.}
Classicaly the existence problem for \eqref{nse1}--\eqref{w0b1}, renaming $\tilde F:=F+v$  amounts to finding a vector function
$$
u:\overline R\times[0,\,T]\to\mathbb R^2
$$
and a scalar function
$$
p:\overline R\times[0,\,T]\to\mathbb R,
$$
such that
\begin{align}
u_t + (u\cdot\nabla)u+\nabla p &= \nu\Delta u + \tilde F\quad&&\text{in}\quad R\times]0,\,T[;\label{c.s.nse}\\
\nabla\cdot u&=0\quad&&\text{in}\quad R\times]0,\,T[;\label{c.s.divo}\\
u\cdot\mathbf n&= 0\quad&&\text{on}\quad\partial R\times]0,\,T[;\label{c.s.tgb}\\
\nabla^\bot\cdot u&=0\quad&&\text{on}\quad\partial R\times]0,\,T[;\label{c.s.w0b}\\
u(x,\,0)&=u_0(x)\quad&&\text{in}\quad R.\label{c.s.inicond}
\end{align}
Where $\tilde F$ and $u_0$ are given and defined in $R\times[0\,T]$ and $R$ respectively. $u_0$ is the position at time 0 of the system.
\begin{Lemma}\label{D.n.n}
For $u,\,v\in V$ we have $(\Delta u,\,v)=-((u,\,v))$.
\end{Lemma}
\begin{proof}
First we note that $(\Delta u,\,v)=-(A^{\frac{1}{2}}u,\,A^{\frac{1}{2}}v)$ where $A^{\frac{1}{2}}$ is the linear map defined by
$$
A^{\frac{1}{2}}u=A^{\frac{1}{2}}\sum_{k\in\mathbb N_0^2}u_kW_k:=\sum_{k\in\mathbb N_0^2}(-\bar k)^{\frac{1}{2}}u_kW_k.
$$
It is clear that $A^{\frac{1}{2}}$ maps $V$ onto $H$ continuously. Indeed for $u\in V$
$$
|A^{\frac{1}{2}}u|^2=\frac{ab}{4}\sum_{k\in\mathbb N_0^2}-\bar k(-\bar k)u_k^2=\|u\|^2<\infty
$$
and, the continuity follows from the linearity and from $|A^{\frac{1}{2}}(u-v)|^2=\|u-v\|^2$.\\
Therefore it is enough to prove the Lemma for $u,\,v\in\mathcal D_1(R)$, by continuity it will be true for $u,\,v\in V$.
\begin{align*}
(\Delta u,\,v)&=\sum_{i=1}^2\int_R\Delta u_iv_i\,dx=-\sum_{i=1}^2\int_R\nabla u_i\cdot\nabla v_i\,dx+\sum_{i=1}^2\int_R\nabla\cdot(\nabla u_iv_i)\,dx\\
&=-((u,\,v))+\sum_{i=1}^2\int_\Gamma v_i\nabla u_i\cdot\mathbf n\,d\Gamma.
\end{align*}
Now we compute
$$
\int_\Gamma v_1\nabla u_1\cdot\mathbf n\,d\Gamma=\int_\Gamma v_1\partial_1u_1\mathbf n_1\,d\Gamma+\int_\Gamma v_1\partial_2u_1\mathbf n_2\,d\Gamma.
$$
The first term vanishes because $v_1$ vanishes where $\mathbf n_1$ does not\,\footnote{$\mathbf n_1$ vanishes on the top and botton faces of the rectangle and $v_1$ vanishes on the left and right ones.} so we have
\begin{align*}
\int_\Gamma v_1\nabla u_1\cdot\mathbf n\,d\Gamma&=0+\int_{f_b\cup f_t} v_1\partial_2u_1\mathbf n_2\,d\Gamma\\
&=\int_{f_b\cup f_t}\partial_1u_2v_1\mathbf n_2\,d\Gamma=0.\qquad[\text{by}\;\nabla^\bot\cdot u=0\;\text{on}\;\Gamma],
\end{align*}
where $f_b,\;f_t$ and (below) $f_l,\;f_r$ stay for the botton, top, left and right faces of $R$.\\
Analogously we conclude that
$$
\int_\Gamma v_2\nabla u_2\cdot\mathbf n\,d\Gamma=
\int_{f_l\cup f_r} v_2\partial_2u_1\mathbf n_1\,d\Gamma=0.
$$
Therefore
$$
(\Delta u,\,v)=-((u,\,v))\;\text{for}\;u,\,v\in \mathcal D_1(R).
$$
\end{proof}
\begin{Lemma}\label{D.nbot.nbot}
For $u,\,v\in V$ we have $(\Delta u,\,v)=-(\nabla^\bot\cdot u,\,\nabla^\bot\cdot v)$.
\end{Lemma}
\begin{proof}
Again it is enough to prove the Lemma for $u,\,v\in\mathcal D_1(R)$. By \eqref{D.n.nbot}
\begin{align*}
\int_R\Delta u\cdot v\,dx&=\int_R\nabla^\bot(\nabla^\bot\cdot u)\cdot v\,dx\\
&=-\int_R(\nabla^\bot\cdot u)(\nabla^\bot\cdot v)\,dx+\int_R\nabla^\bot\cdot\Bigl((\nabla^\bot\cdot u) v\Bigr)\,dx\\
&=-\int_R(\nabla^\bot\cdot u)(\nabla^\bot\cdot v)\,dx+\int_R\nabla\cdot\Bigl((\nabla^\bot\cdot u) v\Bigr)^\circ\,dx\\
&=-\int_R(\nabla^\bot\cdot u)(\nabla^\bot\cdot v)\,dx+\int_\Gamma\Bigl((\nabla^\bot\cdot u) v\Bigr)^\circ\cdot\mathbf n\,d\Gamma\\
&=-\int_R(\nabla^\bot\cdot u)(\nabla^\bot\cdot v)\,dx.
\end{align*}
Where $\begin{pmatrix}u_1\\u_2\end{pmatrix}^\circ:=\begin{pmatrix}u_2\\-u_1\end{pmatrix}$.
\end{proof}
\begin{Corollary}
The norms $u\mapsto\|u\|$ and $u\mapsto|\nabla^\bot\cdot u|$ coincide in $\mathcal D_1(R)$ and then in $V$.
\end{Corollary}
If $u$ and $p$ are classical solutions of \eqref{c.s.nse}--\eqref{c.s.inicond} --- $u\in C^2(\overline R\times[0,\,T]),\,p\in C^1(\overline R\times[0,\,T])$. Then clearly $u\in L^2(0,\,T,\,V)$ and if we fix $v\in\mathcal D_1(R)$, take the scalar product and use Lemma \ref{D.n.n} obtain
\begin{equation}
\frac{d}{dt}(u,\,v)+\nu((u,\,v))+b(u,\,u,\,v)=(\tilde F,\,v).
\end{equation}
By continuity the previous expression holds for each $v\in V$.\footnote{$v_k\to v$ in $\|\cdot\|\Rightarrow\, v_k\to v$ in $|\cdot|$. $(u,\,\cdot)$ continuous in $|\cdot|$ and $((u,\,\cdot))$ continuous in $\|\cdot\|$} So is natural to define the following weak formulation of problem \eqref{c.s.nse}--\eqref{c.s.inicond}:
\begin{Problem}\label{w.p.leray}
Given
\begin{eqnarray}
& &\tilde F\in L^2(0,\,T,\,V^\prime)\label{w.s.F}\\
&\And&\notag\\
& &u_0\in H\label{w.s.u0}\\
&to\, find&\notag\\
& &u\in L^2(0,\,T,\,V)\label{w.s.u}\\
&\text{satisfying}&\text{(in the distribution sense)}\notag\\
& &\frac{d}{dt}(u,\,v)+\nu((u,\,v))+b(u,\,u,\,v)=<\tilde F,\,v>,\quad\forall v\in V,\qquad\label{w.s.nse}\\
&\text{and}&\notag\\
& &u(0)=u_0\label{w.s.inicond}.
\end{eqnarray}
\end{Problem}
\begin{Remark}
Note that \eqref{w.s.u} is not sufficient to give sense to \eqref{w.s.inicond}. But we will show that if we have in addition \eqref{w.s.nse} then $u$ coincides almost everywhere with a continuous function giving a meaning to \eqref{w.s.inicond}.
\end{Remark}
Now we present a proplem equivalent to Problem \ref{w.p.leray}:
\begin{Problem}\label{w.p.AB}
Given
\begin{eqnarray}
& &\tilde F\in L^2(0,\,T,\,V^\prime),\label{AB.F}\\
&\And&\notag\\
& &u_0\in H,\label{AB.u0}\\
&to\, find&\notag\\
& &u\in L^2(0,\,T,\,V),\qquad u^\prime\in L^1(0,\,T,\,V^\prime)\label{AB.uu'}\\
&\text{satisfying}&\notag\\
& &u^\prime+\nu Au+Bu=\tilde F\quad\text{on}\quad]0,\,T[,\label{AB.nse}\\
&\text{and}&\notag\\
& &u(0)=u_0\label{AB.inicond}.
\end{eqnarray}
\end{Problem}
To verify this equivalence we will need the following lemmas
\begin{Lemma}\label{primit}
Let $X$ be a Banach space with dual $X^\prime$ and let $u$ and $g$  be two functions belonging to $L^1(a,\,b,\,X)$. Then the following three conditions are equivalent
\begin{enumerate}
\item $u$ is a.e. equal to a primitive function of $g$,
$$
u(t)=\xi+\int_0^t g(s)\,ds,\quad\xi\in X,\,\text{a.e.}\,t\in[a,\,b];
$$\label{primit1}
\item For each test function $\phi\in\mathcal D(]a,\,b[)$,
$$
\int_a^b u(t)\phi^\prime(t)\,dt=-\int_a^b g(t)\phi(t)\,dt\qquad\bigg(\phi^\prime=\frac{d\phi}{dt}\biggr);
$$\label{primit2}
\item For each $\eta\in X^\prime$,
$$
\frac{d}{dt}<u,\,\eta>=<g,\,\eta>
$$
in the scalar distribution sense, on $]a,\,b[$\label{primit3}
\end{enumerate}
In particular if \ref{primit1}--\ref{primit3} is satisfied, $u$ is a.e. equal to a continuous function $[a,\,b]\to X$.
\end{Lemma}
The proof of Lemma \ref{primit} can be found in \cite{temams} section 3.1.
\begin{Lemma}
If $u\in L^2(0,\,T,\,V)$, the function $Bu$ defined by
$$
<Bu(t),\,v>=b(u(t),\,u(t),\,v),\quad\forall v\in V,\,\text{a.e.}\,t\in [0,\,T],
$$
belongs to $L^1(0,\,T,\,V^\prime)$.
\end{Lemma}
\begin{proof}
For almost all $t,$ $Bu(t)$ is an element of $V^\prime$. The measurability of the function
$$
[0,\,T]\ni t\mapsto Bu(t)\in V^\prime,
$$
is easy to check. Indeed it follows from the measurability of $t\mapsto u(t)\in V$ and from the continuity of $u\mapsto B(u)\in V^\prime$. By the continuity and trilinearity of $b$ on $V$ we have
\begin{equation}\label{bw<w2}
\|Bw\|_{V^\prime}\leq C\|w\|^2,\quad\forall w\in V,\quad C\in\mathbb R^+.
\end{equation}
Hence
$$
\int_0^T\|Bu(t)\|_{V^\prime}\,dt\leq C\int_0^T\|u(t)\|^2\,dt<+\infty.
$$
\end{proof}
Now let $u$ satisfy both \eqref{w.s.u} and \eqref{w.s.nse} then, due to the identity \eqref{scalH=scalVV'}, to definition of $A$ (given in \eqref{defA}) and to the previous lemma we can write \eqref{w.s.nse} as
$$
\frac{d}{dt}<u,\,v>=<\tilde F-\nu Au - Bu,\,v>,\quad\forall v\in V.
$$
By the linearity and continuity of $A:V\to V^\prime$ we have that
\begin{equation}\label{AuL2V'}
Au\in L^2(0,\,T;\,V^\prime)
\end{equation}
so, $\tilde F-\nu Au - Bu\in L^1(0,\,T;\,V^\prime)$. Therefore, by Lemma \ref{primit}\;\footnote{Note that $V$, being a Hilbert space, is reflexive.} and
$$
\begin{cases}
&u^\prime\in L^1(0,\,T;\,V^\prime)\\
&u^\prime=f-\nu Au-Bu,
\end{cases}
$$
$u$ is a.e. equal to a continuous function $[0,\,T]\to V^\prime$. Therefore \eqref{w.s.inicond} is meaningful and any solution of Problem \ref{w.p.leray} is a solution of Problem \ref{w.p.AB}. Since any solution of Problem \ref{w.p.AB} is clearly a solution of Problem \ref{w.p.leray} we conclude that these two problems are equivalent.
\section{Existence}
\subsection{Fourier Transform. Fractional Derivatives.}
We want to prove the existence of a solution for Problem \ref{w.p.leray}, for that we will need another compactness theorem envolving fractional derivatives of a function.
Given a function $f$ from $\mathbb R$ into a Hilbert space $X_1$ we denote its Fourier Transform by
$$
\hat f=\int_{-\infty}^{+\infty} e^{-2i\pi t\tau}f(t)\,dt.
$$
The derivative (in $t$) of order $\gamma$ of $f$ is the Fourier Transform of $(2i\pi\tau)^\gamma\hat f$ or
$$
\widehat{D_t^\gamma f(\tau)}=(2\pi\tau)^\gamma \hat f(\tau).
$$
Now let $X_0,\,X,\,X_1$ be Hilbert spaces such that
$$
X_0\subset X\subset X_1
$$
where, the injections are continuous and the first injection is in addiction compact.
Given $\gamma>0$, we define the space
$$
\mathcal H^\gamma(\mathbb R,\,X_0,\,X_1)=\{v\in L^2(\mathbb R,\,X_0)\mid\,D_t^\gamma v\in L^2(\mathbb R,\,X_1)\}.
$$
$\mathcal H^\gamma(\mathbb R,\,X_0,\,X_1)$ in a Hilbert space for the norm $\|\cdot\|_{\mathcal H^\gamma(\mathbb R,\,X_0,\,X_1)}$ defined by
$$
\|v\|_{\mathcal H^\gamma(\mathbb R,\,X_0,\,X_1)}^2=\|v\|_{L^2(\mathbb R,\,X_0)}^2+\||\tau|^\gamma\hat v\|_{L^2(\mathbb R,\,X_1)}.
$$
For any set $K\subset\mathbb R$ we define the space
$$
\mathcal H_K^\gamma(\mathbb R,\,X_0,\,X_1)=\{u\in\mathcal H^\gamma(\mathbb R,\,X_0,\,X_1)\mid\,\mathrm{supp}\, u\subset K\}.
$$
Now we can state the compactness theorem we will need:
\begin{Theorem}\label{calHinL2}
Let $X_0\subset X\subset X_1$ be Hilbert spaces with both inclusions being continuous and the first one being also compact. Then for any bounded set $K$ and any $\gamma>0$ the injection of $\mathcal H_K^\gamma(\mathbb R,\,X_0,\,X_1)$ into $L^2(\mathbb R,\,X)$ is compact.
\end{Theorem}
The proof can be found in \cite{temams} section 3.2.3.
\subsection{The Existence Theorem.}
In the proof of the Existence Theorem \ref{ex.w.sol} we will need the following Lemma:
\begin{Lemma}\label{bskew}
Fixing the first variable in $V$, the form $b$ defined in \eqref{def.b} results skew-symmetric in the last two variables, i.e.,
$$
\forall u\in V\forall v,\,w\in\mathbf H^1[\,b(u,\,v,\,w)=-b(u,\,w,\,v)\,].
$$
\end{Lemma}
Due to the trilinearity of $b$ the previous Lemma is equivalent to the following corollary
\begin{Corollary}\label{b(uvv)=0}
Fixing the first variable in $V$, we have
$$
\forall u\in V\forall v\in\mathbf H^1[\,b(u,\,v,\,v)=0\,].
$$
\end{Corollary}
\begin{proof}
We prove the statement for $v\in C^\infty(\overline R)$. Then by the continuity  of $b$ it holds for $v\in H^1$.\\
For $u\in V$ and $v\in C^\infty(\overline R)$ we have
\begin{align*}
b(u,\,v,\,v)&=\sum_{i,j=1}^2\int_R u_i\partial_i v_jv_j\,dx=\sum_{i,j=1}^2\int_R u_i\partial_i\frac{(v_j)^2}{2}\,dx\\
&=\sum_{i,j=1}^2\int_R \partial_i\bigl(u_i\frac{(v_j)^2}{2}\bigr)\,dx-\sum_{i,j=1}^2\int_R \partial_i u_i\frac{(v_j)^2}{2}\,dx\\
&=\sum_{j=1}^2\int_R \nabla\cdot\bigl(u\frac{(v_j)^2}{2}\bigr)\,dx-\sum_{j=1}^2\int_R (\nabla\cdot u)\frac{(v_j)^2}{2}\,dx\\
&=\sum_{j=1}^2\int_\Gamma u\frac{(v_j)^2}{2}\cdot\mathbf n\,d\Gamma-0=\sum_{j=1}^2\int_\Gamma \frac{(v_j)^2}{2}u\cdot\mathbf n\,d\Gamma=0.
\end{align*}
\end{proof}
\begin{Theorem}\label{ex.w.sol}
Given $\tilde F$ and $u_0$ satisfying \eqref{AB.F} and \eqref{AB.u0}.
There is at least one
function $u$ satysfying
\eqref{AB.uu'}-\eqref{AB.inicond}.
\end{Theorem}
\begin{Remark}
The proof that follows is completely analogous to that for the case of full Dirichlet conditions that can be found in \cite{temams}.
\end{Remark}
\begin{proof}
We start by defining, for each $m\in\mathbb N_0$, an approximate solution $u^m$ of \eqref{w.s.nse} as follows:
\begin{align}
&u^m:=\sum_{\max\{i_1,\,i_2\}\leq m} u_i^m(t)W_i\label{app.w.s.um}\\
((u^m)^\prime(t),\,W_j)&+\nu((u^m(t),\,W_j))+b(u^m(t),\,u^m(t),\,W_j)=<\tilde F(t),\,W_j>,\label{app.w.s.nse}\\
&\qquad t\in[0,\,T],\quad \max\{j_1,\,j_2\}\leq m,\notag\\
&u^m(0)=u_0^m.\label{app.w.s.inicond}
\end{align}
Where $u_0^m$ is the orthogonal projection of $u_0$ onto $span\{W_i\mid\,\max\{i_1,\,i_2\}\leq m\}$.\\
From \eqref{app.w.s.nse} we obtain the nonlinear system of differential equations in the functions $u_i^m,\quad \max\{i_1,\,i_2\}\leq m$:
\begin{multline}
\sum_{\max\{i_1,\,i_2\}\leq m} (u_i^m)^\prime(t)(W_i,\,W_j)+\nu\sum_{\max\{i_1,\,i_2\}\leq m} u_i^m((W_i,\,W_j))\\
+\sum_{\begin{subarray}{l}\max\{i_1,\,i_2\}\leq m\\\max\{l_1,\,l_2\}\leq m\end{subarray}} u_i^m u_l^m b(W_i,\,W_l,\,W_j)=<\tilde F(t),\,W_j>,
\end{multline}
that reduces to the ODE's system\footnote{In section \ref{Sfourier} we have seen that the family $\mathcal W$ is orthogonal in both $H$ and $V$.}
\begin{multline}\label{app.ODE}
(u_j^m)^\prime(t)\bigl(-\bar\jmath\frac{ab}{4}\bigr)+\nu u_j^m \bar\jmath^2\frac{ab}{4}\\
+\sum_{\begin{subarray}{l}\max\{i_1,\,i_2\}\leq m\\\max\{l_1,\,l_2\}\leq m\end{subarray}} u_i^m u_l^m b(W_i,\,W_l,\,W_j)=<\tilde F(t),\,W_j>;\\ \max\{j_1,\,j_2\}\leq m.
\end{multline}
Now we note that \eqref{app.w.s.inicond} is the same as the $m$ scalar conditions
\begin{equation}\label{ic.app.ODE}
u_j^m(0)=\text{the projection of $u_0$ onto $span\{W_j\}$}=u_{0j}.\quad\footnotemark
\end{equation}
\footnotetext{If we write $u_0=\sum_{i\in\mathbb N_0^2} u_{0i}W_i$.}
 Therefore system \eqref{app.ODE} with initial condition \eqref{ic.app.ODE} has a maximal solution defined in $[0,\,t_{max}[$ and, we will see that $t_{max}=T$. Indeed, if $t_{max}< T$ then $|u^m(t)|$ should tend to $+\infty$ as $t$ goes to $t_{max}$. Below we show that $|u^m(t)|$ remains bounded in $[0,\,T]$ so we have that the maximal  solution is defined in $[0,\,T]$.\\
{\sc A priori estimates}: We multiply \eqref{app.w.s.nse} by $u_j^m,\quad \max\{j_1,\,j_2\}\leq m$ and add the equations obtained. Taking Corollary \ref{b(uvv)=0} into account, we arrive to
$$
((u^m)^\prime(t),\,u^m(t))+\nu\|u^m(t)\|^2=<\tilde F(t),\,u^m(t)>.
$$
Hence
\begin{multline*}
\frac{d}{dt}|u^m(t)|^2+2\nu\|u^m(t)\|^2=2<\tilde F(t),\,u^m(t)>\\
\leq 2\|\tilde F(t)\|_{V^\prime}\|u^m(t)\|\leq\nu\|u^m(t)\|^2+\frac{1}{\nu}\|\tilde F(t)\|_{V^\prime}^2\quad\footnotemark
\end{multline*}
\footnotetext{By expanded $(\sqrt{\nu}a-\frac{1}{\sqrt{\nu}}b)^2\geq 0$ with $a=\|u^m(t)\|$ and $b=\|\tilde F(t)\|_{V^\prime}$.}
so,
\begin{equation}\label{est.eq}
\frac{d}{dt}|u^m(t)|^2+\nu\|u^m(t)\|^2\leq\frac{1}{\nu}\|\tilde F(t)\|_{V^\prime}^2.
\end{equation}
In particular $\frac{d}{dt}|u^m(t)|^2\leq\frac{1}{\nu}\|\tilde F(t)\|_{V^\prime}^2$ and, integrating on $[0,\,s]$ we obtain
\begin{equation}\label{|u|finite}
|u^m(s)|^2\leq|u_0^m|^2+\frac{1}{\nu}\int_0^s \|\tilde F(t)\|_{V^\prime}^2\,dt\leq|u_0|^2+\frac{1}{\nu}\int_0^T \|\tilde F(t)\|_{V^\prime}^2\,dt.
\end{equation}
Therefore
\begin{equation}\label{um.bdd.LinfH}
the\; sequence\; (u^m)\; remains\; in\; a\; bounded\; set\; of\; L^\infty(0,\,T,\,H).
\end{equation}
Now we integrate \eqref{est.eq} over $[0,\,T]$ and obtain:
\begin{align*}
|u^m(T)|^2+\nu\int_0^T\|u^m(t)\|^2\,dt&\leq|u_0^m|^2+\frac{1}{\nu}\int_0^T \|\tilde F(t)\|_{V^\prime}^2\,dt\\&\leq|u_0|^2+\frac{1}{\nu}\int_0^T \|\tilde F(t)\|_{V^\prime}^2\,dt.
\end{align*}
Therefore
\begin{equation}\label{um.bdd.L2V}
the\; sequence\; (u^m)\; remains\; in\; a\; bounded\; set\; of\; L^2(0,\,T,\,V).
\end{equation}
Now we extend $u^m$ to the entire real line putting
$$
\tilde{u}^m(t):=\begin{cases} u^m(t)&\text{if}\,t\in[0,\,T]\\0&\text{if}\,t\notin[0,\,T]\end{cases}
$$
and the Fourier transform of $\tilde{u}^m$ will be denoted by $\hat{u}^m$.\\
{\sc Application of Theorem \ref{calHinL2}.} We start by computing the integral
\begin{equation}\label{int.fract}
\int_{-\infty}^{+\infty}|\tau|^{2\gamma}|\hat{u}^m(\tau)|^2\,d\tau.
\end{equation}
Equation \eqref{app.w.s.nse} with $u^m$ replaced by $\tilde u^m$ results in
\begin{multline}\label{tilde.nse}
\frac{d}{dt}(\tilde u^m(t),\,W_j)=<\tilde f^m,\,W_j>+(u_0^m,\,W_j)\delta_0-(u^m(T),\,W_j)\delta_T;\\ \max\{j_1,\,j_2\}\leq m
\end{multline}
where $\delta_0,\,\delta_T$ are the Dirac distributions at $0$ and $T$ and,
\begin{eqnarray*}
f^m&=&\tilde F-\nu Au^m-Bu^m,\\
\tilde f^m(t)&:=&\begin{cases}f^m(t)&\text{on}[0,\,T]\\0&\text{outside}[0,\,T] \end{cases}.
\end{eqnarray*}
Using the Fourier Transform \eqref{tilde.nse} becomes
\begin{multline}\label{hat.nse}
2\pi i\tau(\hat u^m(\tau),\,W_j)=<\hat f^m(\tau),\,W_j>\\+(u_0^m,\,W_j)-(u^m(T),\,W_j)\exp(-2\pi iT\tau).
\end{multline}
Multiplying \eqref{hat.nse} by $\hat g_i^m(\tau)$ and adding the obtained $m$ equations we arrive to
\begin{multline}\label{hat.nse2}
2\pi i\tau|\hat u^m(\tau)|^2=<\hat f^m(\tau),\,\hat u^m(\tau)>\\+(u_0^m,\,\hat u^m(\tau))-(u^m(T),\,\hat u^m(\tau))\exp(-2\pi iT\tau).
\end{multline}
By \eqref{bw<w2} we have
\begin{align}
\|f^m(t)\|_{V^\prime}&\leq\|\tilde F(t)\|_{V^\prime}+\nu\|u^m(t)\|+C\|u^m(t)\|^2;\; \text{so}\notag\\
\int_0^T\|f^m(t)\|_{V^\prime}\,dt&\leq\int_0^T\|\tilde F(t)\|_{V^\prime}+\nu\|u^m(t)\|+C\|u^m(t)\|^2\,dt\label{intfmbdd}
\end{align}
by \eqref{AB.F} and \eqref{um.bdd.L2V} the integral \eqref{intfmbdd} remains bounded. Hence \footnote{It is known that $|f(\tau)|\leq C\int_{\mathbb R}|f(t)|\,dt$. See for example \cite{stric}.}
$$
\sup_{\tau\in\mathbb R}\|\hat f^m(\tau)\|_{V^\prime}\leq\,\text{const},\quad\forall m\in\mathbb N_0.
$$
By \eqref{|u|finite} both $|u^m(0)|$ and $|u^m(T)|$ are finite. Then by \eqref{hat.nse2}
\begin{equation}\label{hatubdd}
|\tau||\hat u^m(\tau)|^2\leq C_1\|\hat u^m(\tau)\|+C_2|\hat u^m(\tau)|\leq D\|\hat u^m(\tau)\|.\quad\footnotemark
\end{equation}
\footnotetext{$C_1,\,C_2$ and $D$ are constants.}\par
Fix $\gamma<\frac{1}{4}$ and define the real function $Q(x):=\frac{x^{2\gamma}+x}{1+x},\quad x\in[0,\,+\infty[$. $Q$ is continuous and bounded,\footnote{$\lim_{x\to+\infty}Q(x)=1,\,Q(0)=0$} then we can find a constant $D_1\in\mathbb R^+$ such that for all $\tau\in\mathbb R$:
\begin{align*}
&Q(|\tau|)\leq D_1\\
&\Leftrightarrow|\tau|^{2\gamma}\leq D_1\frac{1+|\tau|}{1+|\tau|^{1-2\gamma}}.
\end{align*}
Therefore the integral \eqref{int.fract} is bounded by
\begin{align*}
&D_1\int_{-\infty}^{+\infty} \frac{1+|\tau|}{1+|\tau|^{1-2\gamma}}|\hat u^m(\tau)|^2\,d\tau\\
\text{by \eqref{hatubdd}}&\\
\leq &D_2\int_{-\infty}^{+\infty} \frac{\|\hat u^m(\tau)\|}{1+|\tau|^{1-2\gamma}}\,d\tau + D_3\int_{-\infty}^{+\infty} \|\hat u^m(\tau)\|^2\,d\tau.
\end{align*}
The integral $D_3\int_{-\infty}^{+\infty} \|\hat u^m(\tau)\|^2\,d\tau$, by the Parseval Equality and \eqref{um.bdd.L2V}, there is a constant $D_4$ such that
\begin{equation}\label{int2bdd}
D_3\int_{-\infty}^{+\infty} \|\hat u^m(\tau)\|^2\,d\tau\leq D_4,\qquad\forall m
\end{equation}
For the integral $D_2\int_{-\infty}^{+\infty} \frac{\|\hat u^m(\tau)\|}{1+|\tau|^{1-2\gamma}}\,d\tau$ we apply Schwartz Inequality and Parseval Equality to obtain
$$
D_2\int_{-\infty}^{+\infty} \frac{\|\hat u^m(\tau)\|}{1+|\tau|^{1-2\gamma}}\,d\tau\leq D_2\Biggl(\int_{-\infty}^{+\infty} \frac{1}{(1+|\tau|^{1-2\gamma})^2}\,d\tau\Biggr)^\frac{1}{2}\Biggl(\int_{0}^{T} \|u^m(t)\|^2\,dt\Biggr)^\frac{1}{2}\quad\footnotemark
$$
\footnotetext{Since $\gamma<\frac{1}{4},\quad \frac{1}{1+|\tau|^{1-2\gamma}}\in L^2(\mathbb R)$. Indeed $ \int\frac{1}{(1+|\tau|^{1-2\gamma})^2}=2\int_0^{+\infty}\frac{1}{(1+\tau^{1-2\gamma})^2}$. Putting $x=1+\tau^{1-2\gamma}$ we see that the last integral equals $2\int_1^{+\infty}\frac{1}{x^2}\frac{1}{1-2\gamma}(x-1)^\frac{2\gamma}{1-2\gamma}\leq\frac{2}{1-2\gamma}\int_1^{+\infty}\frac{1}{x^{2-\frac{2\gamma}{1-2\gamma}}}$ and the last integral converges if $2-\frac{2\gamma}{1-2\gamma}>1\leftrightarrow \gamma<\frac{1}{4}$.}
and, this product is finite and bounded as $m\to+\infty$, i.e., there is a constant $D_5$ such that
\begin{equation}\label{int1bdd}
D_2\int_{-\infty}^{+\infty} \frac{\|\hat u^m(\tau)\|}{1+|\tau|^{1-2\gamma}}\,d\tau\leq D_5,\qquad\forall m.
\end{equation}
By \eqref{int2bdd} and \eqref{int1bdd} we conclude that the integral \eqref{int.fract} is finite:
$$
\int_{-\infty}^{+\infty}|\tau|^{2\gamma}|\hat{u}^m(\tau)|^2\,d\tau\leq D_4+D_5=:E.
$$
The finiteness of \eqref{int.fract} with \eqref{um.bdd.L2V} implies that
\begin{equation}\label{um.bdd.Hgamma}
the\; sequence\; (u^m)\; remains\; in\; a\; bounded\; set\; of\; \mathcal H^\gamma(\mathbb R,\,V,\,H).
\end{equation}
{\sc The limit.} We now need the following lemmas that can be found in \cite{bre} section III.6:
\begin{Lemma}\label{wconvL2}
Let $E$ be a reflexive Banach space and let $(x_n)$ a bounded sequence in $E$. Then there is a subsequence $(x_{\sigma(n)})$ of $(x_n)$ such that $x_{\sigma(n)}\rightharpoonup x$, for some $x\in E$.\quad\footnotemark
\end{Lemma}
\footnotetext{We shall use the symbols: $\rightharpoonup$ for weak convergence; $\rightharpoonup_\ast$ for weak-star convergence and, $\rightarrow$ for strong convergence.}
\begin{Lemma}\label{wstconvLinf}
Let $E$ be a separable Banach space and let $(f_n)$ a bounded sequence in $E^\prime$. Then there is a subsequence $(f_{\sigma(n)})$ of $(f_n)$ such that $f_{\sigma(n)}\rightharpoonup_\ast f$, for some $f\in E^\prime$.
\end{Lemma}
Lemma \ref{wstconvLinf} and \eqref{um.bdd.LinfH} implies the existence of a subsequence $(u^{\sigma(m)})$ of $u^m$ and $u\in L^\infty(0,\,T,\,H)$ such that
$$
u^{\sigma(m)}\rightharpoonup_\ast u,\,\text{in}\,L^\infty(0,\,T,\,H).\quad\footnotemark
$$
\footnotetext{If wanted, without lack of generality we may assume $u^m\rightharpoonup_\ast u,\,\text{in}\,L^\infty(0,\,T,\,H)$.}
Analogously, Lemma \ref{wconvL2} and \eqref{um.bdd.L2V} implies the existence of a subsequence $(u^{\alpha(\sigma(m))})$ of $u^\sigma(m)$ and $v\in L^2(0,\,T,\,V)$ such that
$$
u^{\alpha(\sigma(m))}\rightharpoonup v,\,\text{in}\,L^2(0,\,T,\,V).
$$
The sequence $(u^m)$ is in the space $\mathcal H_{[0,\,T]}^\gamma(\mathbb R,\,V,\,H)$ which injection in $L^2(\mathbb R,\,H)$ is compact due to Theorem \ref{calHinL2}.\footnote{With $V=X_0,\And H=X=X_1$.} Then \eqref{um.bdd.Hgamma} implies the existence of a subsequence $\beta(\alpha(\sigma(m)))$ of $\alpha(\sigma(m))$ and $w\in L^2(0,\,T,\,H)$ satisfying
$$
u^{\beta(\alpha(\sigma(m)))}\rightarrow w,\,\text{in}\,L^2(0,\,T,\,H).
$$
We put $\eta:=\beta\circ\alpha\circ\sigma$ and we obtain
$$
u=v=w\in L^2(0,\,T,\,H)\cap L^2(0,\,T,\,V)\cap L^\infty(0,\,T,\,H)
$$
and
\begin{align}
&u^{\eta(m)}\rightharpoonup_\ast u,\,\text{in}\,L^\infty(0,\,T,\,H);\label{um.cwst.LinfH}\\
&u^{\eta(m)}\rightharpoonup u,\,\text{in}\,L^2(0,\,T,\,V);\label{um.cw.L2V}\\
&u^{\eta(m)}\rightarrow u,\,\text{in}\,L^2(0,\,T,\,H).\label{um.cs.L2H}
\end{align}
Indeed $L^\infty(0,\,T,\,H)=(L^1(0,\,T,\,H))^\prime$ and $L^2(0,\,T,\,H)\subset L^1(0,\,T,\,H)$. So for each $f\in L^2(0,\,T,\,H)$:
$$
\begin{cases}
u^{\eta(m)}(f)&\to u(f)\\
u^{\eta(m)}(f)&\to v(f).
\end{cases}
$$
Since both $w,\,v\in L^2(0,\,T,\,H)=(L^2(0,\,T,\,H))^\prime$, we can rewrite the previous expressions as
$$
\forall f\in(L^2(0,\,T,\,H))^\prime
\begin{cases}
f(u^{\eta(m)})&\to f(u)\\
f(u^{\eta(m)})&\to f(v).
\end{cases}
$$
Then $u=v\in L^2(0,\,T,\,H)$.\\
On the other side $(L^2(0,\,T,\,H))^\prime\subseteq(L^2(0,\,T,\,V))^\prime$ because a converging sequence in $L^2(0,\,T,\,V)$ converges in $L^2(0,\,T,\,H)$ too.\\
So for each $f\in(L^2(0,\,T,\,H))^\prime$ we have
$$
\begin{cases}
f(u^{\eta(m)})&\to f(v)\\
f(u^{\eta(m)})&\to f(w).
\end{cases}
$$
But since both $w,\,v\in L^2(0,\,T,\,H)$ we have that $v=w$.\footnote{Note that since both $L^2(0,\,T,\,V)$ and $L^2(0,\,T,\,H)$ are both Hilbert spaces then each of them coincide with its dual. But we can not do these identifications here if we want to compare their duals, because the spaces $H$ and $V$ the norms are different. If we do the identifications we will obtain $(L^2(0,\,T,\,H))^\prime\supseteq(L^2(0,\,T,\,V))^\prime$ and that is not true.}\\
{\sc Ending.} To complete the proof we will need the following Lemma:
\begin{Lemma}\label{for.lim}
If a sequence $(u^n)$ satisfies $u^n\rightharpoonup u$ in $L^2(0,\,T,\,V)$ and $u^n\rightarrow u$ in $L^2(0,\,T,\,H)$, then for any vector function $w$ with components in $C^1(\overline{R}\times[0,\,T])$,
$$
\int_0^T b(u^n(t),\,u^n(t),\,w(t))\,dt\rightarrow\int_0^T b(u(t),\,u(t),\,w(t))\,dt.
$$
\end{Lemma}
\begin{proof}
From
$$
b(u^n,\,u^n,\,w)-b(u,\,u,\,w)=b(u^n-u,\,u^n,\,w)+b(u,\,u^n-u,\,w)
$$
we obtain
\begin{align*}
&\Bigl|\int_0^T b(u^n(t),\,u^n(t),\,w(t))-b(u(t),\,u(t),\,w(t))\,dt\Bigr|\\
\leq&\int_0^T |b(u^n(t)-u(t),\,u^n(t),\,w(t))|\,dt+\int_0^T|b(u(t),\,u^n(t)-u(t),\,w(t))|\,dt\\
(\text{using} &\text{\eqref{|b|<||and||||} below, and Lemma \ref{bskew},})\\ 
\leq& C\int_0^T \Bigl[|u^n(t)-u(t)|^\frac{1}{2}\|u^n(t)-u(t)\|^\frac{1}{2}|u^n(t)|^\frac{1}{2}\|u^n(t)\|^\frac{1}{2}\|w\|\\
&\qquad+|u(t)|^\frac{1}{2}\|u(t)\|^\frac{1}{2}|u^n(t)-u(t)|^\frac{1}{2}\|u^n(t)-u(t)\|^\frac{1}{2}\|w\|\Bigr]\,dt\\
\leq&  C_1\Bigl(\int_0^T |u^n(t)-u(t)|\|u^n(t)-u(t)\|\,dt\Bigr)^\frac{1}{2}\|u^n(t)\|_{L^2(0,\,T,\,V)}\\
&\qquad+C_1\|u(t)\|_{L^2(0,\,T,\,V)}\Bigl(\int_0^T|u^n(t)-u(t)|\|u^n(t)-u(t)\|\,dt\Bigr)^\frac{1}{2};
\end{align*}
but $u^n$ is bounded in $L^2(0,\,T,\,V)$, because $u^n\rightharpoonup u$ in $L^2(0,\,T,\,V)$, then
\begin{align*}
&\Bigl|\int_0^T b(u^n(t),\,u^n(t),\,w(t))-b(u(t),\,u(t),\,w(t))\,dt\Bigr|\\
\leq& C_2\|u^n(t)-u(t)\|_{L^2(0,\,T,\,H)}^\frac{1}{2}.
\end{align*}
Hence $\bigl|\int_0^T b(u^n(t),\,u^n(t),\,w(t))-b(u(t),\,u(t),\,w(t))\,dt\bigr|$ goes to $0$ as $n$ goes to $+\infty$.
\end{proof}
Multiplying \eqref{app.w.s.nse} by a function $\phi\in C^1([0,\,T])$ and, integrating
we obtain:
\begin{align*}
&\int_0^T((u^m)^\prime(t),\,W_j\phi(t))\,dt+\int_0^T \nu((u^m(t),\,W_j\phi(t)))\,dt\\+&\int_0^T b(u^m(t),\,u^m(t),\,W_j\phi(t))\,dt=\int_0^T <\tilde F(t),\,W_j\phi(t)>\,dt\\
\Longleftrightarrow&-\int_0^T(u^m(t),\,W_j\phi^\prime(t))\,dt+\nu\int_0^T ((u^m(t),\,W_j\phi(t)))\,dt\\+&\int_0^T b(u^m(t),\,u^m(t),\,W_j\phi(t))\,dt=(u^m_0,\,W_j)\phi(t)(0)-(u^m(T),\,W_j)\phi(t)(T)\\+&\int_0^T <\tilde F(t),\,W_j\phi(t)>\,dt.\end{align*}
Due to Lemma \ref{for.lim} we can take the limite in the nonlinear term and, by \eqref{um.cw.L2V} we can take the limite in the linear terms.\footnote{Note that \eqref{um.cw.L2V} implies that $u^{\eta(m)}\rightharpoonup u$ in $L^2(0,\,T,\,H)$.} Hence
\begin{multline}
-\int_0^T(u(t),\,W_j\phi^\prime(t))\,dt+\nu\int_0^T ((u(t),\,W_j\phi(t)))\,dt\\+\int_0^T b(u(t),\,u(t),\,W_j\phi(t))\,dt\\=(u_0,\,W_j)\phi(0)-(u(T),\,W_j)\phi(T)+\int_0^T <\tilde F(t),\,W_j\phi(t)>\,dt\label{dist.nse.Wj}
\end{multline}
Equation \eqref{dist.nse.Wj} being true for all $W_j$ by linearity will be true for any finite combination of functions in $\mathcal W$ and by continuity will be true for all $v\in V$. Then taking a test function $\phi\in\mathcal D(]0,\,T[)$ in \eqref{dist.nse.Wj} we conclude that
\begin{align}
&-\int_0^T(u(t),\,v\phi^\prime(t))\,dt+\nu\int_0^T ((u(t),\,v\phi(t)))\,dt+\int_0^T b(u(t),\,u(t),\,v\phi(t))\,dt\notag\\=&(u_0,\,v)\phi(0)-(u(T),\,v)\phi(T)+<\tilde F(t),\,v\phi(t)>\,dt\label{dist.nse.v};\\
&\text{then}\notag\\
&-\int_0^T(u(t),\,v\phi^\prime(t))\,dt+\nu\int_0^T ((u(t),\,v\phi(t)))\,dt+\int_0^T b(u(t),\,u(t),\,v\phi(t))\,dt\notag\\=&<\tilde F(t),\,v\phi(t)>\,dt\label{dist.nse.v1}
\end{align}
what means that
$$
Equation\; \eqref{w.s.nse}\; is\; satisfied\; in\; the\; distribution\; sense.
$$
Multiplying \eqref{w.s.nse} by $\phi\in C^1([0,\,T])$ such that $\phi(0)=1,\,\phi(T)=0$ we obtain
\begin{align}
&-\int_0^T(u(t),\,v\phi^\prime(t))\,dt+\nu\int_0^T ((u(t),\,v\phi(t)))\,dt+\int_0^T b(u(t),\,u(t),\,v\phi(t))\,dt\notag\\=&(u(0),\,v)+\int_0^T <\tilde F(t),\,W_j\phi(t)>\,dt.\label{dist.nse}
\end{align}
The equation resulting from \eqref{dist.nse.v} with the same $\phi$ as in \eqref{dist.nse} is
\begin{align}
&-\int_0^T(u(t),\,v\phi^\prime(t))\,dt+\nu\int_0^T ((u(t),\,v\phi(t)))\,dt+\int_0^T b(u(t),\,u(t),\,v\phi(t))\,dt\notag\\=&(u_0,\,v)+<\tilde F(t),\,v\phi(t)>\,dt\label{dist.nse1}
\end{align}
From \eqref{dist.nse} and \eqref{dist.nse1} we conclude that \eqref{w.s.inicond} is satisfied (and, we have finished the proof of theorem \ref{ex.w.sol}).
\end{proof}
\section{Uniqueness.}
\begin{Lemma}
\begin{equation}\label{|b|<||and||||}
|b(u,\,v,\,w)|\leq c|u|^{\frac{1}{2}}\|u\|^{\frac{1}{2}}\|v\||w|^{\frac{1}{2}}\|w\|^{\frac{1}{2}},\quad\forall u,\,v,\,w\in\mathbf H^1(R).
\end{equation}
If $u$ belongs to $L^2(0,\,T,\,V)\cap L^\infty(0,\,T,\,H)$, then $Bu$ belongs to $L^2(0,\,T,\,V^\prime)$ and
\begin{equation}\label{BuinL2}
\|Bu\|_{L^2(0,\,T,\,V^\prime)}\leq c|u|_{L^\infty(0,\,T,\,H)}\|u\|{L^2(0,\,T,\,V)}.
\end{equation}
\end{Lemma}
See \cite{temsp68} section III 3.2. for \eqref{|b|<||and||||}. Then \eqref{BuinL2} is a corollary of \eqref{|b|<||and||||}. Indeed
\begin{align}
|b(u,\,u,\,w)|=|b(u,\,w,\,u)|&\leq C|u|^{\frac{1}{2}}\|u\|^{\frac{1}{2}}|u|^{\frac{1}{2}}\|u\|^{\frac{1}{2}}\|w\|\notag\\
&=C|u|\|u\|\|w\|,\qquad w,\,u\in V.\label{new.est.b}
\end{align}
In the case $u\in L^2(0,\,T,\,V)\cap L^\infty(0,\,T,\,H)$ we obtain that
\begin{align*}
&\int_0^T \|Bu(t)\|_{V^\prime}^2\,dt\leq C\int_0^T |u|^2\|u\|^2\,dt\\
\leq& C\|u\|_{L^\infty(0,\,T,\,H)}^2\int_0^T \|u(t)\|^2\,dt=C\|u\|_{L^\infty(0,\,T,\,H)}^2\|u\|_{L^2(0,\,T,\,V)}^2<+\infty
\end{align*}
\begin{Theorem}\label{cont.uH}
The solution of Problems \ref{w.p.leray}-\ref{w.p.AB} given by Theorem \ref{ex.w.sol} is unique. Moreover it is a.e. equal to a continuous function from $[0,\,T]$ into $H$ and,
\begin{equation}\label{limt1H}
u(t)\to u(t_1),\quad\text{in}\;H\quad\text{as}\;t\to t_1\quad t_1\in[0,\,T].
\end{equation}
In particular
\begin{align*}
&u(t)\to u_0,\quad\text{in}\;H\quad\text{as}\;t\to0;\\
&u(t)\to u(T),\quad\text{in}\;H\quad\text{as}\;t\to T.
\end{align*}
\end{Theorem}
For the proof we need the following lemma from \cite{temams}:
\begin{Lemma}\label{L:du2=<>}
Let $V,\,H,\,V^\prime$ three Hilbert spaces satisfying $V\subset H\subset V^\prime$  with dense and continuous inclusions. If
$$
u\in L^2(0,\,T,\,V)\And u^\prime\in L^2(0,\,T,\,V^\prime)
$$
then  $u$ is a.e. equal to a continuous function from $[0,\,T]$
into $H$ and
\begin{equation}\label{du2=<>}
\frac{d}{dt}|u|^2=2<u^\prime,\,u>
\end{equation}
holds in the distribution sense on $(0,\,T)$.
\end{Lemma}
\begin{Remark}
Note that \eqref{du2=<>} is meaningful because both
$$
t\mapsto |u(t)|^2\And t\mapsto<u^\prime(t),\,u(t)>
$$
are integrable on $[0,\,T]$.
\end{Remark}
\begin{proof}[Proof of Theorem \ref{cont.uH}]
By \eqref{AB.F}, \eqref{AB.nse}, \eqref{AuL2V'} and \eqref{BuinL2} we have that $u^\prime\in L^2(0,\,T,\,V^\prime)$. Then by \eqref{um.cwst.LinfH}, \eqref{um.cw.L2V}, and Lemma \ref{L:du2=<>} we have that
$$
u\in C([0,\,T],\,H).
$$
Therefore \eqref{limt1H} is satisfied.\\
To prove the uniqueness we consider two solutions $u,\,v$ of the problems and put $w:=u-v$. Then
\begin{align*}
w^\prime=&u^\prime-v^\prime=\tilde F-\nu Au-Bu-(\tilde F-\nu Av-Bv)\\
=&-\nu Aw-Bu+Bv\\
\text{and}\qquad&\\
w(0)=&0.
\end{align*}
Now we take the scalar product with $w$ in the duality between $V$ and $V^\prime$ and obtain
$$
<w^\prime,\,w>+\nu<Aw,\,w>=<Bv,\,w>-<Bu,\,w>.
$$
We note that since $u,\,v\in L^2(0,\,T,\,V)\And u^\prime,\,v^\prime\in L^2(0,\,T,\,V^\prime)$ then $w\in L^2(0,\,T,\,V)\And w^\prime\in L^2(0,\,T,\,V^\prime)$. Hence we can apply \eqref{du2=<>} and obtain
\begin{align*}
\frac{d}{dt}|w(t)|^2+2\nu\|w(t)\|^2&=2b(v(t),\,v(t)\,w(t))-2b(u(t),\,u(t),\,w(t))\\
&=-2b(w(t),\,v(t),\,w(t))=2b(w(t),\,w(t),\,v(t)).\quad\footnotemark
\end{align*}
\footnotetext{$b(w,\,v,\,w)=-b(v,\,v,\,w)+b(u,\,v,\,w)=-b(v,\,v,\,w)+b(u,\,w,\,w)+b(u,\,u,\,w)=-b(v,\,v,\,w)+b(u,\,u,\,w)$.}

Hence using \eqref{new.est.b}
\begin{align*}
\frac{d}{dt}|w(t)|^2+2\nu\|w(t)\|^2&\leq2(C|w(t)|\|w(t)\|\|v(t)\|)\\
&\leq 2\nu\|w(t)\|^2+\frac{C^2}{2\nu}|w(t)|^2\|v(t)\|^2\quad\footnotemark
\end{align*}
\footnotetext{By expanded $(\sqrt{2\nu}a-\frac{C}{\sqrt{2\nu}}b)^2\geq 0$ with $a=\|w(t)\|$ and $b=|w(t)|\|v(t)\|$.}
Threfore
$$
\frac{d}{dt}|w(t)|^2\leq\frac{C^2}{2\nu}|w(t)|^2\|v(t)\|^2.
$$
Since $t\mapsto \frac{d}{dt}|w(t)|^2=2<w^\prime(t),\,w(t)>,\,t\mapsto |w(t)|^2=<w(t),\,w(t)>$ and $t\mapsto \|v(t)\|^2=<v(t),\,v(t)>$ are integrable, by Gronwall Inequality (see \eqref{gronin} below)
\begin{align*}
|w(t)|^2&\leq|w(0)|^2\exp\Biggl(\frac{C^2}{2\nu}\int_0^t\|v(\tau)\|^2\,d\tau\Biggr)\\
&\leq|w(0)|^2\exp\Biggl(\frac{C^2}{2\nu}\int_0^T\|v(\tau)\|^2\,d\tau\Biggr)=0,\qquad t\in[0,\,T]
\end{align*}
(using the fact that $w(0)=0$). Hence
$$
u=v
$$
concluding the uniqueness of the solution.
\end{proof}
\begin{Lemma}[Gronwall Inequality]
Let $g,\,h,\,y,\,\frac{dy}{dt}$ be locally integrable functions sattisfying
\begin{equation}\label{gronin}
\frac{dy}{dt}\leq gy+h\quad\text{for $t\geq t_0$}.
\end{equation}
Then
$$
y(t)\leq y(t_0)\exp\Bigl(\int_{t_0}^t g(\tau)\,d\tau\Bigr)+\int_{t_0}^t h(s)\exp\bigl(-\int_t^s g(\tau)\,d\tau\Bigr)\,ds,\quad t\geq t_0.
$$
[The proof can be found in \cite{temsp68}].
\end{Lemma}
\section{Continuity on Initial Data.}
\begin{Theorem}\label{contindata}
The map
\begin{align*}
\mathbb S:\,H\times L^2(0,\,T,\,V^\prime)\times]0,\,+\infty[&\to C([0,\,T],\,H)\\
(u_0,\,\tilde F,\,\nu)&\mapsto u
\end{align*}
is continuous. Here $u\in C([0,\,T],\,H)$ is the unique solution of Problems \ref{w.p.leray}-\ref{w.p.AB} (see Theorem \ref{cont.uH}).
\end{Theorem}
\begin{proof}
Fix a triple $(u_0,\,\tilde F,\,\nu)\in H\times L^2(0,\,T,\,V^\prime)\times]0,\,+\infty[$ and consider the solution of Theorem \ref{cont.uH} induced by this triple. Then $u$ satisfies
$$
u^\prime + \nu Au + Bu = \tilde F,\qquad u(0)=u_0.
$$
Now consider another triple $(v_0,\,G,\,\eta)\in H\times L^2(0,\,T,\,V^\prime)\times]0,\,+\infty[$. The solution associated with this triple satisfies
$$
v^\prime + \eta Av + Bv = G,\qquad v(0)=v_0.
$$
We put
$$
w:=v-u
$$
and see that $w$ satisfies the equation
$$
w^\prime=G-\tilde F-\eta Aw +(\nu-\eta)Au-Bv+Bu.
$$
Taking the scalar product with $w$ we obtain
\begin{align*}
<w^\prime,\,w>&=<G-\tilde F,\,w>-\eta\|w\|^2+(\nu-\eta)((u,\,w))+b(w,\,w,\,u)\\
\frac{d}{dt}|w|^2&\leq 2\|G-\tilde F\|_{V^\prime}\|w\|-2\eta\|w\|^2+2|\nu-\eta|\|u\|\|w\|+2C|w|\|w\|\|u\|.
\end{align*}
Expanding
\begin{align*}
\Bigl(\sqrt\frac{\eta}{3}\|w\|-\sqrt\frac{3}{\eta}\|G-\tilde F\|_{V^\prime}\Bigr)^2\geq 0;\\
\Bigl(\sqrt\frac{\eta}{3}\|w\|-C\sqrt\frac{3}{\eta}|w|\|u\|\Bigr)^2\geq 0;\\
\Bigl(\sqrt\frac{\eta}{3}\|w\|-|\nu-\eta|\sqrt\frac{3}{\eta}\|u\|\Bigr)^2\geq
0;
\end{align*}
we obtain
\begin{equation}\label{contSest}
\frac{d}{dt}|w|^2+\eta\|w\|^2\leq
\frac{3}{\eta}\|G-\tilde F\|_{V^\prime}^2+|\nu-\eta|^2\frac{3}{\eta}\|u\|^2+C^2\frac{3}{\eta}|w|^2\|u\|^2.
\end{equation}
By \eqref{gronin} we have
\begin{align*}
&\quad|w(t)|^2\\
&\leq |w(0)|^2\exp\Bigl(C^2\frac{3}{\eta}\int_0^t \|u(\tau)\|^2\,d\tau\Bigr)\\&+\int_0^t \Bigl(\frac{3}{\eta}\|G(s)-\tilde F(s)\|_{V^\prime}^2+|\nu-\eta|^2\frac{3}{\eta}\|u(s)\|^2\Bigr)\exp\Bigl(-\frac{3C^2}{\eta}\int_t^s \|u(\tau)\|^2\,d\tau\Bigr)\,ds\\
&\leq \exp\Bigl(C^2\frac{3}{\eta}\int_0^T
\|u(\tau)\|^2\,d\tau\Bigr)\biggl( |w(0)|^2\\&+\int_0^T
\Bigl(\frac{3}{\eta}\|G(s)-\tilde F(s)\|_{V^\prime}^2+|\nu-\eta|^2\frac{3}{\eta}\|u(s)\|^2\Bigr)\,ds\biggr).
\end{align*}
Now if $|\nu-\eta|<\frac{\nu}{2}$, i.e., $\eta\in]\frac{\nu}{2},\,\frac{3\nu}{2}[$, we have
\begin{multline*}
|w(t)|^2\leq \exp\Bigl(C^2\frac{6}{\nu}\int_0^T
\|u(\tau)\|^2\,d\tau\Bigr)\biggl( |w(0)|^2\\+\int_0^T
\Bigl(\frac{6}{\nu}\|G(s)-\tilde F(s)\|_{V^\prime}^2+|\nu-\eta|^2\frac{6}{\nu}\|u(s)\|^2\Bigr)\,ds\biggr).
\end{multline*}
Fix $\varepsilon>0$.\\
Now if we put $E_0:=\exp\Bigl(C^2\frac{6}{\nu}\int_0^T
\|u(\tau)\|^2\,d\tau\Bigr)$ and choose the triple $(v_0,\,G,\,\eta)$ such that
\begin{align*}
|v_0-u_0|=|w(0)|<\frac{\varepsilon}{\sqrt 3}E_0^{-\frac{1}{2}}=:&\alpha_v;\\
\|G-\tilde F\|_{L^2(0,\,T,\,V^\prime)}=\Biggl(\int_0^T\|G(s)-\tilde F(s)\|_{V^\prime}^2\,ds\Biggr)
^\frac{1}{2}<\frac{\varepsilon\sqrt\nu}{3\sqrt 2}E_0^{-\frac{1}{2}}=:&\alpha_g;\\
|\nu-\eta|<\min\{\frac{\nu}{2},\,\delta\}=:&\alpha_\eta,\\\text{where}\quad\delta=
\frac{\varepsilon\sqrt\nu}{3\sqrt
2}\Bigl(\int_0^T\|u(s)\|^2\,ds\Bigr)^{-\frac{1}{2}}E_0^{-\frac{1}{2}},&\qquad
\footnotemark
\end{align*}
\footnotetext{Since $u\in L^2(0,\,T,\,V)$ we have that $\alpha_v$
and $\alpha_\eta$ are finite. We can also see that
$\alpha_v,\,\alpha_g,\,\alpha_\eta$ depend only on the (fixed)
triple $(u_0,\,\tilde F,\,\nu)$.} we obtain
$$
|w(t)|^2<\varepsilon^2\Leftrightarrow |w(t)|<\varepsilon,\qquad t\in[0,\,T]
$$
and, we have the continuity of $\mathbb S$.
\end{proof}

\begin{Theorem}\label{contindata2}
The map
\begin{align*}
\mathbb S_2:\,H\times L^2(0,\,T,\,V^\prime)\times]0,\,+\infty[&\to L^2(0,\,T,\,V)\\
(u_0,\,\tilde F,\,\nu)&\mapsto u
\end{align*}
is continuous. Here $u\in L^2(0,\,T,\,V)$ is the unique solution
of Problems \ref{w.p.leray}-\ref{w.p.AB} (see Theorem
\ref{cont.uH}).
\end{Theorem}
\begin{proof}
As in the proof of Theorem \ref{contindata} we fix
$\varepsilon>0$ and a triple $(u_0,\,\tilde F,\,\nu)\in H\times
L^2(0,\,T,\,V^\prime)\times]0,\,+\infty[$ and consider the
solution of Theorem \ref{cont.uH} induced by this triple. Now
consider another triple $(v_0,\,G,\,\eta)\in H\times
L^2(0,\,T,\,V^\prime)\times]0,\,+\infty[$ and the solution
associated with this triple. Put
$$
w:=v-u
$$
and taking the scalar product with $w$ we arrive to
\eqref{contSest}. Integrating over $[0,\,T]$ we obtain
\begin{multline*}
\frac{1}{\eta}|w(T)|^2-\frac{1}{\eta}|w(0)|^2
+\int_0^T\|w(t)\|^2\,dt\leq\frac{3}{\eta^2}\|G-\tilde F\|_{L^2(0,\,T,\,V^\prime)}^2\\+|\nu-\eta|^2\frac{3}{\eta^2}\|u\|_{L^2(0,\,T,\,V)}^2+\frac{3C^2}{\eta^2}\|w\|_{C([0,\,T],\,H)}\|u\|_
{L^2(0,\,T,\,V)}^2
\end{multline*}
and, if $|\nu-\eta|<\frac{\nu}{2}$ we have
\begin{align*}
\int_0^T\|w(t)\|^2\,dt&\leq\frac{12}{\nu^2}\|G-\tilde F\|_{L^2(0,\,T,\,V^\prime)}^2
+|\nu-\eta|^2\frac{12}{\nu^2}\|u\|_{L^2(0,\,T,\,V)}^2\\
&+\Bigl(\frac{12C^2}{\nu^2}\|u\|_{L^2(0,\,T,\,V)}^2+
\frac{4}{\nu}\Bigr)\|w\|_{C(0,\,T,\,H)}^2.
\end{align*}
Using Theorem \ref{contindata}, there is $\delta>0$ such that if
$(v_0,\,G,\,\eta)$ satisfies
$$
|v_0-u_0|<\delta,\,\|G-\tilde F\|_{L^2(0,\,T,\,V^\prime)}<\delta,\,|\eta-\nu|<\delta,
$$
then
$\|w\|_{C(0,\,T,\,H)}<\bigl(\frac{\varepsilon^2}{3}\bigr)^{\frac{1}{2}}\bigl(\frac{12C^2}{\nu^2}\|u\|_{L^2(0,\,T,\,V)}^2+\frac{4}{\nu}\bigr)^{-\frac{1}{2}}$.\\
Hence for
$$
\delta_1=\min\Bigl\{\frac{\nu}{2},\;\delta,\;\Bigl(\frac{\varepsilon^2}{3}\Bigr)^\frac{1}{2}\Bigl(\frac{12}{\nu^2}\Bigr)^{-\frac{1}{2}}
,\;\Bigl(\frac{\varepsilon^2}{3}\Bigr)^{\frac{1}{2}}\Bigl(\frac{12}{\nu^2}\|u\|_{L^2(0,\,T,\,V)}^2\Bigr)^
{-\frac{1}{2}}\Bigr\},
$$
we have $\|w\|_{L^2(0,\,T,\,V)}<\varepsilon$. Therefore the map
$\mathbb S_2$ is continuous.
\end{proof}

\section{Some More Estimates for the Form $b$.}\label{S:Est.b}
We present here some estimates we shall need later. We
have\footnote{These estimates can be found in \cite{temsp68} section III.3.2.. They can be obtained by interpolation (\cite{lio}), generalized Sobolev inequalities (\cite{pee}) and by a theorem by S.Agmon (\cite{agm}). See \cite{temnsf} for indications how to obtain them.}
$$
|b(u,\,v,\,w)|\leq CK.
$$
where $C$ is a constant and $K$ is one of the following products
\begin{align}
&\|u\|\|v\|\|w\|\quad &&u,\,v,\,w\in\mathbf H^1(R),\label{E.1.1.1}\\
&|u|^\frac{1}{2}\|u\|^\frac{1}{2}\|v\|^\frac{1}{2}\|v\|_2^\frac{1}{2}|w|\quad &&u\in\mathbf H^1(R),\,v\in\mathbf H^2(R),\,w\in\mathbf L^2(R),\label{E.01.12.0}\\
&|u|^\frac{1}{2}\|u\|_2^\frac{1}{2}\|v\||w|\quad &&u\in\mathbf H^2(R),\,v\in\mathbf H^1(R),\,w\in\mathbf L^2(R),\label{E.02.1.0}\\
&|u|\|v\||w|^\frac{1}{2}\|w\|_2^\frac{1}{2}\quad &&u\in\mathbf L^2(R),\,v\in\mathbf H^1(R),\,w\in\mathbf H^2(R),\label{E.0.1.02}\\
&|u|^\frac{1}{2}\|u\|^\frac{1}{2}\|v\||w|^\frac{1}{2}\|w\|^\frac{1}{2}\quad
&&u,\,v,\,w\in\mathbf H^1(R)\label{E.01.1.01}.
\end{align}

\section{Strong Formulation.}

Sometimes we want more regularity for the solutions of problem
\eqref{nse1}--\eqref{w0b1}. Instead of Problems
\eqref{w.p.leray}-\eqref{w.p.AB} where we ask for weak solutions
we consider the following equivalent problems where we look for
solutions more regular than weak.
\begin{Problem}\label{s.p.leray}
Given
\begin{eqnarray}
& &\tilde F\in L^2(0,\,T,\,H)\label{s.p.F}\\
&\And&\notag\\
& &u_0\in V\label{s.p.u0}\\
&to\, find&\notag\\
& &u\in L^2(0,\,T,\,D(A))\cap L^\infty (0,\,T,\,V)\label{s.p.u}\\
&\text{satisfying}&\text{(in the distribution sense)}\notag\\
& &\frac{d}{dt}(u,\,v)+\nu((u,\,v))+b(u,\,u,\,v)=(\tilde F,\,v),\quad\forall v\in V,\qquad
\label{s.p.nse}\\
&\text{and}&\notag\\
& &u(0)=u_0\label{s.p.inicond}.
\end{eqnarray}
\end{Problem}

\begin{Problem}\label{s.p.AB}
Given
\begin{eqnarray}
& &\tilde F\in L^2(0,\,T,\,H),\label{ABs.F}\\
&\And&\notag\\
& &u_0\in V,\label{ABs.u0}\\
&to\, find&\notag\\
& &u\in L^2(0,\,T,\,D(A))\cap L^\infty (0,\,T,\,V),\;\text{and}\notag\\
& &u^\prime\in L^2(0,\,T,\,H)\label{ABs.uu'}\\
&\text{satisfying}&\notag\\
& &u^\prime+\nu Au+Bu=\tilde F\quad\text{on}\quad]0,\,T[,\label{ABs.nse}\\
&\text{and}&\notag\\
& &u(0)=u_0\label{ABs.inicond}.
\end{eqnarray}
\end{Problem}
The equivalence of these problems follows from
\begin{enumerate}
\item A solution of Problem \ref{s.p.AB} is a solution of
Problem \ref{s.p.leray};
\item For a solution $u$ of Problem \ref{s.p.leray} we have
$$
|Bu(t)|\leq
C\|u(t)\|^\frac{3}{2}\|u(t)\|_2^\frac{1}{2}\quad\text{a.e. $t$ by
\eqref{E.01.12.0},\quad so}
$$
$$
\int_0^T|Bu(t)|^4\leq C\int_0^T\|u(t)\|^6\|u(t)\|_2^2\leq
C_1\int_0^T\|u(t)\|_2^2<\infty.
$$
Hence
\begin{enumerate}
\item $Bu\in L^4(0,\,T,\,H)$\label{eq.BL4H}.\\
Since $u\in L^2(0,\,T,\,D(A))$ and $\int_0^T|Au(t)|^2=
\int_0^T\|u(t)\|_2^2<\infty$ we have
\item $Au\in L^2(0,\,T,\,H)$.\label{eq.AL2H}
\end{enumerate}
From \eqref{eq.BL4H} and \eqref{eq.AL2H} we have that
$$
f-\nu Au-Bu\in L^2(0,\,T,\,H)\subseteq L^2(0,\,T,\,V^\prime).
$$
Since $u\in L^2(0,\,T,\,D(A))\subseteq L^2(0,\,T,\,V)$, by Lemma
\eqref{primit} and \eqref{s.p.nse}, we have that
$$
u^\prime=f-\nu Au-Bu \text{ a.e. and } u\in C(0,\,T,\,V^\prime)
\text{``a.e.''}.
$$
Hence
\item $u$ is a solution of Problem \ref{s.p.AB}.
\end{enumerate}
\begin{Lemma}\label{DA.V.H}
The inclusions
$$
D(A)\subseteq V\subseteq H
$$
are both dense continuous and compact.
\end{Lemma}
\begin{proof}
We have already seen in Corollary \ref{inc.V-H} and beginning of
section \ref{S.VH} that the inclusion $V\subseteq H$ has the
required properties. For the first inclusion we have that $V$,
respectively $D(A)$, are the closure of $\mathcal D_1(R)$ in
$\mathbf H^1$, respectively in $\mathbf H^2$. From the density,
continuity and compactness of the inclusion $\mathbf
H^2\subseteq\mathbf H^1$ come the same properties for the
inclusion $D(A)\subseteq V$.
\end{proof}
\section{Existence.}
\begin{Theorem}\label{ex.s.sol}
Given $\tilde F$ and $u_0$ satisfying \eqref{ABs.F} and \eqref{ABs.u0}.
There is at least one function $u$ satisfying
\eqref{ABs.uu'}-\eqref{ABs.inicond}.
\end{Theorem}
\begin{Remark}
The proof that follows is completely analogous to that for the
case of weak solutions and for no-slip boundary conditions that
can be found in \cite{temams}.
\end{Remark}
\begin{proof}{Outlines:}
We define an approximate solution $u^m$, for each $m\in\mathbb N_0$ like in \eqref{app.w.s.um}--\eqref{app.w.s.inicond}, and arrive to the
estimate \eqref{est.eq} and conclusions \eqref{um.bdd.LinfH} and
\eqref{um.bdd.L2V} exactly in the same way. For Theorem we need
only some more estimates: If we multiply \ref{app.w.s.nse} by
$\bar\jmath u_j^m$ and add the obtained equations we arrive to
$$
((u^m(t))^\prime,\,Au^m(t))+\nu((u^m(t),\,Au^m(t)))+Bu^m(t)(Au^m(t))=(\tilde F(t),\,Au^m(t)),
$$
i.e.,
$$
(((u^m(t))^\prime,\,u^m(t)))+\nu(u^m(t),\,u^m(t))_{[2]}+Bu^m(t)(Au^m(t))=(\tilde F(t),\,Au^m(t)),
$$
or
\begin{align*}
\frac{1}{2}\frac{d}{dt}\|u^m(t)\|^2+\nu|u^m(t)|_{[2]}^2&\leq
|Bu^m(t)(Au^m(t))|+|\tilde F(t)||Au^m(t)|\\
\text{(by \eqref{E.01.12.0})}\quad&\leq
C|u^m(t)|^\frac{1}{2}\|u^m(t)\||u^m(t)|_{[2]}
^\frac{3}{2}+|\tilde F(t)||u^m(t)|_{[2]}\\
&\text{(by Young Inequalities)}\footnotemark\\&\leq
\frac{\nu}{4}(|u^m(t)|_{[2]}^\frac{3}{2})^\frac{4}{3}+
C_{\frac{\nu}{4},\frac{4}{3}}(C|u^m(t)|^\frac{1}{2}\|u^m(t)\|)^4\\
&+\frac{\nu}{4}|u^m(t)|_{[2]}^2+C_{\frac{\nu}{4},2}|\tilde F(t)|^2\\
&\leq \frac{\nu}{4}(|u^m(t)|_{[2]}^\frac{3}{2})^\frac{4}{3}+
\frac{3^3}{4\nu^3}(C|u^m(t)|^\frac{1}{2}\|u^m(t)\|)^4\\
&+\frac{\nu}{4}|u^m(t)|_{[2]}^2+\frac{1}{\nu}|\tilde F(t)|^2.
\end{align*}
\footnotetext{Young Inequality: $ab\leq \varepsilon a^p +
C_{\varepsilon,p}^Y b^{p^\prime}$. Where $a,b>0$, $1<p<\infty$,
$\varepsilon>0$, $\frac{1}{p}+\frac{1}{p^\prime}=1$ and
$C_{\varepsilon,p}^Y=\frac{p-1}{\bigl(p^{p^\prime}\bigr)\bigl(\varepsilon^\frac{1}{p-1}\bigr)}$.}

\begin{equation}\label{NE}
\frac{d}{dt}\|u^m\|^2+\nu|u^m|_{[2]}^2\leq
+\frac{3^3}{2\nu^3}C^2|u^m(t)|^2\|u^m(t)\|^4
+\frac{2}{\nu}|\tilde F(t)|^2.
\end{equation}
From equation \eqref{NE} and from Gronwall Inequality we can
derive for $s\in [0,\,T]$ (using \eqref{um.bdd.LinfH} and
\eqref{um.bdd.L2V})
\begin{align*}
&\|u^m(s)\|^2\\
\leq&\exp\Bigl(\int_0^T\frac{3^3}{2\nu^3}C^2|u^m(t)|^2\|u^m(t)\|^2\,dt\Bigr)\biggl(\|u^m_0\|^2+\int_0^T \frac{2}{\nu}|\tilde F(t)|^2\,dt\biggr)\\
\leq&K_1
\end{align*}
for some constant $K_1$ (independent of $m$). Hence
\begin{equation}\label{um.bdd.LinfV}
the\; sequence\; (u^m)\; remains\; in\; a\; bounded\; set\; of\;
L^\infty(0,\,T,\,V).
\end{equation}
Now we integrate \eqref{NE} over $[0,\,T]$ and obtain:
\begin{align*}
&\|u^m(T)\|^2-\|u^m(0)\|^2+\nu\int_0^T|u^m(t)|_{[2]}^2\,dt\\
\leq&\int_0^T\frac{3^3}{2\nu^3}C^2|u^m(t)|^2\|u^m(t)\|^4+\frac{2}{\nu}|\tilde F(t)|^2\,dt\leq K_3
\end{align*}
with $K_3$ being a constant (independent of $m$, using \eqref{um.bdd.LinfH},
\eqref{um.bdd.LinfV} and \eqref{ABs.F}). Therefore
\begin{equation}\label{um.bdd.L2DA}
the\; sequence\; (u^m)\; remains\; in\; a\; bounded\; set\; of\;
L^2(0,\,T,\,D(A)).
\end{equation}
By \ref{um.bdd.LinfV}, \ref{um.bdd.L2DA} and  Lemmas \ref{wconvL2} and \ref{wstconvLinf} we conclude the existence of a subsequence $u^{\sigma(m)}$ such that
\begin{align*}
u^{\sigma(m)}&\rightharpoonup u \;\text{in}\;L^2(0,\,T,\,D(A))\\
u^{\sigma(m)}&\rightharpoonup_\ast u \;\text{in}\;L^\infty(0,\,T,\,V)
\end{align*}
This limit is a
solution for problems \ref{s.p.leray}-\ref{s.p.AB}.
\end{proof}
\section{Uniqueness.}
\begin{Lemma}\label{DA.V.DA'}
The inclusions
$$
D(A)\subseteq V\subseteq D(A)^\prime
$$
are both dense and continuous.
\end{Lemma}
\begin{proof}
We have just seen in Lemma \ref{DA.V.H} that the inclusion
$D(A)\subseteq V$ has the required properties. For the second
inclusion we proceed as in the beginning of section \ref{S.VH}: We
identify $V$ with $V^\prime$.\footnote{Note that here $H$ is not
present. We are doing only one identification.} The continuity of
the second inclusion follows from the continuity of the first
one. The density of the second inclusion follows from Lemma
\ref{NR} and the injectivity of the first one.
\end{proof}
\begin{Theorem}\label{cont.uV}
The solution of problems \ref{s.p.leray}-\ref{s.p.AB} is unique.
Moreover it is a.e. equal to a continuous function from $[0,\,T]$
into $V$.
\end{Theorem}
\begin{proof}
A solution of problems \ref{s.p.leray}-\ref{s.p.AB} is a solution
of problems \ref{w.p.leray}-\ref{w.p.AB} as well, which is unique
by Theorem \ref{cont.uH}. By $u\in L^2(0,\,T,\,D(A))$,
$$
u^\prime\in L^2(0,\,T,\,H)\subseteq
L^2(0,\,T,\,V^\prime)\subseteq L^2(0,\,T,\,D(A)^\prime)
$$
and by lemmas \ref{DA.V.DA'} and \ref{L:du2=<>} we conclude that
$u\in C([0,\,T],V)$.
\end{proof}
\section{Continuity on Initial Data.}
\begin{Theorem}\label{scontindata}
The map
\begin{align*}
\mathbb S_s:\,V\times L^2(0,\,T,\,H)\times]0,\,+\infty[&\to C([0,\,T],\,V)\\
(u_0,\,\tilde F,\,\nu)&\mapsto u
\end{align*}
is continuous. Where $u\in C([0,\,T],\,V)$ is the unique solution
of Problems \ref{s.p.leray}-\ref{s.p.AB} (see Theorem
\ref{cont.uV}).
\end{Theorem}
\begin{proof}
Fix $\varepsilon>0$ and a triple $(u_0,\,\tilde F,\,\nu)\in V\times
L^2(0,\,T,\,H)\times]0,\,+\infty[$ and, consider the solution of
Theorem \ref{cont.uV} induced by this triple. Then $u$ satisfies
$$
u^\prime + \nu Au + Bu = \tilde F,\qquad u(0)=u_0.
$$
Now consider another triple $(v_0,\,G,\,\eta)\in V\times
L^2(0,\,T,\,H)\times]0,\,+\infty[$. The solution associated with
this triple satisfies
$$
v^\prime + \eta Av + Bv = G,\qquad v(0)=v_0.
$$
We put
$$
w:=v-u
$$
and see that $w$ satisfies the equation
$$
w^\prime=G-\tilde F-\eta Aw +(\nu-\eta)Au-Bv+Bu.
$$
Taking the scalar product with $Aw$ we obtain
\begin{equation*}
\frac{d}{dt}\|w\|^2+2\eta|w|_{[2]}^2\leq
2(G-\tilde F,Aw)+2|\nu-\eta|(u,w)_{[2]}+2(Bu-Bv)(Aw).
\end{equation*}
Choosing $\eta$ such that $|\nu-\eta|<\frac{\nu}{2}$:
\begin{equation}\label{estwVDA1}
\frac{d}{dt}\|w\|^2+\nu|w|_{[2]}^2\leq
2(G-\tilde F,Aw)+2|\nu-\eta|(u,w)_{[2]}+2(Bu-Bv)(Aw).
\end{equation}
The last term is the more complicated so we work it a bit:
\begin{align*}
(Bu-Bv)(Aw)&=b(u,\,u,\,Aw)-b(v,\,v,\,Aw)\\
&=b(u-v,\,u,\,Aw)+b(v,\,u,\,Aw)-b(v,\,v,\,Aw)\\
&=-b(w,\,u,\,Aw)-b(v,\,w,\,Aw)\\
&=-b(w,\,u,\,Aw)-b(w,\,w,\,Aw)-b(u,\,w,\,Aw)
\end{align*}
By \eqref{E.01.12.0} we have
\begin{align*}
|b(w,\,u,\,Aw)|&\leq
C|w|^\frac{1}{2}\|w\|^\frac{1}{2}\|u\|^\frac{1}{2}|u|_{[2]}^\frac{1}{2}|w|_{[2]}\\
&\leq C_1\|w\||u|_{[2]}|w|_{[2]}
\end{align*}
and, by \eqref{E.02.1.0}
\begin{align*}
|b(w,\,w,\,Aw)+b(u,\,w,\,Aw)|&\leq
C|w|^\frac{1}{2}|w|_{[2]}^\frac{1}{2}\|w\||w|_{[2]}+C|u|_{[2]}^\frac{1}{2}
|u|^\frac{1}{2}\|w\||w|_{[2]}\\
&\leq
C|w|^\frac{1}{2}\|w\||w|_{[2]}^\frac{3}{2}+C|u|_{[2]}^\frac{1}{2}
|u|^\frac{1}{2}\|w\||w|_{[2]}
\end{align*}
Hence from \eqref{estwVDA1} we have
\begin{align*}
&\frac{d}{dt}\|w\|^2+\nu|w|_{[2]}^2\\
\leq&2|G-\tilde F||w|_{[2]}+2|\nu-\eta||u|_{[2]}|w|_{[2]}+2C_1\|w\||u|_{[2]}|w|_{[2]}\\
+&2C|w|^\frac{1}{2}\|w\||w|_{[2]}^\frac{3}{2}+2C|u|_{[2]}^\frac{1}{2}
|u|^\frac{1}{2}\|w\||w|_{[2]}.
\end{align*}
Applying Young inequalities with suitable exponents and constants
we arrive to
\begin{align*}
&\frac{d}{dt}\|w\|^2+\nu|w|_{[2]}^2\notag\\
\leq&\frac{\nu}{10}|w|_{[2]}^2+\frac{10}{\nu}|G-\tilde F|^2+\frac{\nu}{10}|w|_{[2]}^2+\frac{10}{\nu}|\nu-\eta|^2|u|_{[2]}^2\\
+&\frac{\nu}{10}|w|_{[2]}^2+\frac{10}{\nu}C_1^2\|w\|^2|u|_{[2]}^2\notag\\
+&\frac{\nu}{10}|w|_{[2]}^2+\frac{(10\cdot 3)^3}{4^4\cdot\nu^3}(2C)^4|w|^2\|w\|^4+\frac{\nu}{10}|w|_{[2]}^2+\frac{10}{\nu}C^2|u||u|_{[2]}\|w\|^2.
\end{align*}
Then
\begin{align}
&\frac{d}{dt}\|w\|^2+\frac{\nu}{2}|w|_{[2]}^2\notag\\
\leq&\frac{10}{\nu}|G-\tilde F|^2+\frac{10}{\nu}|\nu-\eta|^2|u|_{[2]}^2+\frac{10}{\nu}C_1^2\|w\|^2|u|_{[2]}^2\notag\\
+&D_1|w|^2\|w\|^4+\frac{10}{\nu}C^2|u||u|_{[2]}\|w\|^2,\label{estwVDA2}
\end{align}
where $D_1=\Bigl(\frac{C}{2}\Bigr)^4\Bigl(\frac{30}{\nu}\Bigr)^3$. By Gronwall Inequality \eqref{gronin} we have for  $t\in [0,\,T]$
\begin{align}
&\|w(t)\|^2\notag\\
\leq& \exp\Bigl(\int_0^T
\frac{10}{\nu}C_1^2|u(s)|_{[2]}^2+D_1|w(s)|^2\|w(s)\|^2+\frac{10}{\nu}C^2|u(s)||u(s)|_{[2]}\,ds\Bigr)\notag\\&\qquad\cdot\biggl(\|w(0)\|^2+\int_0^T
\frac{10}{\nu}|G(s)-\tilde F(s)|^2+\frac{10}{\nu}|\nu-\eta|^2|u(s)|_{[2]}^2\,ds\biggr).\label{grcontSV}
\end{align}
By theorems \ref{contindata} and \ref{contindata2} and from $u\in
L^2(0,\,T,\,D(A))\cap L^\infty(0,\,T,\,V)$ the argument of the
exponential is bounded, say less than a constant $E$, if we choose
the triple $(v_0,\,G,\,\eta)$ such that both $|v_0-u_0|$,
$\|G-\tilde F\|_{L^2(0,\,T,\,V^\prime)}$ and $|\nu-\eta|$ are less than
$\delta$ for some $\delta>0$ sufficiently small.\\
Now put
\begin{multline*}
\delta_1:=\min\Bigl\{\frac{\nu}{2},\;\delta,\;\frac{\varepsilon}{\sqrt3}(\exp E)^{-\frac{1}{2}},\;
\frac{\varepsilon}{\sqrt3}\Bigl(\frac{10}{\nu}\Bigr)^{-\frac{1}{2}}(\exp E)^{-\frac{1}{2}},\;\\
\frac{\varepsilon}{\sqrt3}\Bigl(\frac{10}{\nu}\|u\|_{L^2(0,\,T,\,D(A)}^2\Bigr)^{-\frac{1}{2}}(\exp E)^{-\frac{1}{2}}\Bigr\}.
\end{multline*}
It is clear from \eqref{grcontSV} that if both $|v_0-u_0|$,
$\|G-\tilde F\|_{L^2(0,\,T,\,V^\prime)}$ and $|\nu-\eta|$ are less than
$\delta_1$, we have $\|w(t)\|$ less than $\varepsilon$ and, we
have the continuity of $\mathbb S_s$.
\end{proof}

\begin{Theorem}\label{scontindata2}
The map
\begin{align*}
\mathbb S_{2s}:\,V\times L^2(0,\,T,\,H)\times]0,\,+\infty[&\to L^2(0,\,T,\,D(A))\\
(u_0,\,\tilde F,\,\nu)&\mapsto u
\end{align*}
is continuous. Where $u\in L^2(0,\,T,\,D(A))$ is the unique solution
of Problems \ref{s.p.leray}-\ref{s.p.AB} (see Theorem
\ref{cont.uV}).
\end{Theorem}
\begin{proof}
As in the proof of Theorem \ref{scontindata} we fix
$\varepsilon>0$ and a triple $(u_0,\,\tilde F,\,\nu)\in H\times
L^2(0,\,T,\,V^\prime)\times]0,\,+\infty[$ and consider the
solution of Theorem \ref{cont.uV} induced by this triple. Now
consider another triple $(v_0,\,G,\,\eta)\in H\times
L^2(0,\,T,\,V^\prime)\times]0,\,+\infty[$ and the solution
associated with this triple. Put
$$
w:=v-u
$$
and taking the scalar product with $Aw$ we arrive to
\eqref{estwVDA2}. Integrating over $[0,\,T]$ we obtain

\begin{align}
&\|w(T)\|^2-\|w(0)\|^2+\int_0^T\frac{\nu}{2}\|w\|_{[2]}^2\notag\\
\leq&\frac{10}{\nu}\|G-\tilde F\|_{L^2(0,\,T,\,H)}^2+\frac{10}{\nu}\|u\|_{L^2(0,\,T,\,D(A))}^2|\nu-\eta|^2\notag\\
+&\frac{10}{\nu}C_1^2\|u\|_{L^2(0,\,T,\,D(A))}^2\|w\|_{C([0,\,T],\,V)}^2+TD_1\|w\|_{C([0,\,T],\,H)}^2\|w\|_{C([0,\,T],\,V)}^4\notag\\
+&\frac{10}{\nu}C^2\|u\|_{C([0,\,T],\,H)}^2\|u\|_{L^2(0,\,T,\,D(A))}^2\|w\|_{C([0,\,T],\,V)}^2.\label{estwDA}
\end{align}
By Theorem \ref{scontindata}, for a small enough $\delta>0$ we have
that if both $|v_0-u_0|$, $\|\tilde F-G\|_{L^2(0,\,T,\,H)}$, and
$|\eta-\nu|$ are less than $\delta$, then
$\|w\|_{C([0,\,T],\,V)}<1$. Thus from \eqref{estwDA} we have

\begin{align*}
&\int_0^T\|w\|_{[2]}^2\\
\leq&\frac{20}{\nu^2}\|G-\tilde F\|_{L^2(0,\,T,\,H)}^2+\frac{20}{\nu^2}\|u\|_{L^2(0,\,T,\,D(A))}^2|\nu-\eta|^2\\
+&\Bigl(\frac{4}{\nu}+\frac{20}{\nu^2}C_1^2\|u\|_{L^2(0,\,T,\,D(A))}^2+\frac{2}{\nu}TD_1\|w\|_{C([0,\,T],\,H)}^2\\
&\quad+\frac{20}{\nu^2}C^2\|u\|_{C([0,\,T],\,H)}^2\|u\|_{L^2(0,\,T,\,D(A))}^2\Bigr)\|w\|_{C([0,\,T],\,V)}^2\\
\leq&\frac{20}{\nu^2}\|G-\tilde F\|_{L^2(0,\,T,\,H)}^2+D_2|\nu-\eta|^2+D_3\|w\|_{C([0,\,T],\,V)}^2
\end{align*}
for some constants $D_2,\,D_3$ and all triples $(v_0,\,G,\,\eta)$ satisfying
$$
|v_0-u_0|<\delta,\,\|\tilde F-G\|_{L^2(0,\,T,\,H)}<\delta,\,|\eta-\nu|<\delta.
$$
Then for some $\delta_1$ smaller than $\delta$ we have that if
both $|v_0-u_0|$, $\|\tilde F-G\|_{L^2(0,\,T,\,H)}$, and $|\eta-\nu|$
are less than $\delta_1$, we obtain
$$
\Bigl(\int_0^T\|w\|_{[2]}^2\Bigr)^\frac{1}{2}<\varepsilon.
$$
Therefore the map $\mathbb S_{2s}$ is continuous
\end{proof}

\chapter{Controllability of Galerkin Approximations.}\label{Ch:Gal}

\section{The FCE Procedure.}\label{SecFCE}
In this section we present the FCE Procedure, i.e., a procedure of {\em Factorization}+{\em Convexification}+{\em Extraction}:
\subsection{A Lemma from Linear Algebra.}
A result from Linear Algebra we will need is the following: 
\begin{Lemma}\label{ind.pind}
Fix a linear space $X$ of dimension $N\geq n+m$. Given two families $V:=\{v_i\mid\,i=1,\,\dots\, n\}$ and $W:=\{w_j\mid\,j=1,\,\dots\,m\}$ satisfying: The vectors in $V$ are linearly independent and the vectors in $\pi W=\{\pi w_j\mid\,w_j\in W\}$ are linearly independent as well, where $\pi$ is the projection onto some space $V^c$ transversal to $span(V),\;(X=span(V)\oplus V^c)$. Then the family $V\cup W$ is linearly independent.  
\end{Lemma}

\subsection{Factorization.}\label{sSecFac}
Consider a control-affine system
\begin{equation}\label{affsfac}
\dot{q}=f(q)+\sum_{i=1}^{r}v_i(t)g_i(q)\quad q\in\mathbb R^n,\,v_i\in\mathbb R
\end{equation} 
where $f,\,g_i$ are smooth vector fields and $[g_i,\,g_j]=0\quad i=1,\,\dots,\,r$.\par
In \cite{bas} it is proven that if we decompose the flow of system \eqref{affsfac}
as
\begin{align}
\overrightarrow{exp}\int_0^t (f+gv(\tau))d\tau
:&=\overrightarrow{exp}\int_0^t(f+\sum_{i=1}^{r}v_i(t)g_i(\tau))d\tau\\
\overrightarrow{exp}\int_0^t (Ad\,G_{v(\tau)}^{\tau})f\,d\tau\circ
G_{v(\tau)}^{\tau}&=:\overrightarrow{exp}\int_0^t (\mathrm e^{-gw(\tau)})_\ast f d\tau\circ G_{v(t)}^{\tau}
\end{align}
where $g:=(g_1,\,g_2,\cdots,\,g_r), v:=(v_1,\,v_2,\cdots,\,v_r)^T$ and
$G_{v(t)}^{t}$ stays for the flow
$\overrightarrow{exp}\int_0^t gv(\tau)d\tau=\mathrm e^{gw(t)},\;w(t)=\int_0^tv(\tau)\,d\tau$, then
\begin{equation}
\overline{\mathcal A_{q_0}(f+gv)(t)}=
\overline{\mathcal A_{q_0}((\mathrm e^{-gV})_\ast f)(t)\circ\{G_{v(t)}^t\mid v(t)\in\mathbb R^r\}}\quad\footnotemark
\end{equation}
\footnotetext{Here $V$ is independent of $v$.}
Here $\mathcal A_{x}(y)(t)$ stays for the attainable set at time $t$ from $x$ following the
vector fields $y$.\\
Similarly if we rewrite system \eqref{affsfac} as
$$
\dot{q}=f(q)+\sum_{i=1}^{r}(v_i^1(t)+v_i^2(t))g_i(q)\quad q\in\mathbb R^n,\,v_i^j\in\mathbb R
$$
we arrive to
\begin{equation}\label{decomp}
\overline{\mathcal A_{q_0}(f+gv)(t)}=
\overline{\mathcal A_{q_0}((\mathrm e^{-gV^2})_\ast f_1)(t)\circ\{G_{v^2(t)}^t\mid v^2(t)\in\mathbb R^r\}}
\end{equation}
where $f_1(q):=f(q)+\sum_{i=1}^{r}v_i^1(t)g_i(q)$ and where $v^1,\,v^2$ and $V^2$ are independent.
The system $\dot q\,=\,(\mathrm e^{-gV^2})_\ast f_1(q)$ is called {\em factorized system}.
\begin{Lemma}\label{facclo}
With $(\mathrm e^{-gV^2})_\ast f_1$ and $G_{v^2}^t$ as in equation \eqref{decomp} we have
\begin{multline*}
\overline{\mathcal A_{q_0}((\mathrm e^{-gV^2})_\ast f_1)(t)\circ\{G_{v^2(t)}^t\mid v^2(t)\in\mathbb R^r\}}\\
\supseteq
\overline{\mathcal A_{q_0}((\mathrm e^{-gV^2})_\ast f_1)(t)}\circ\{G_{v^2(t)}^t\mid v^2(t)\in\mathbb R^r\}.
\end{multline*}
\end{Lemma}
\begin{proof}
Let $x\in\overline{\mathcal A_{q_0}((\mathrm e^{-gV^2})_\ast f_1)(t)}\circ\{G_{v^2(t)}^t\mid v^2(t)\in\mathbb
R^r\}$. Then there are $y\in\overline{\mathcal A_{q_0}((\mathrm e^{-gV^2})_\ast f_1)(t)}$ and a
control $u(t)\in\mathbb R^r\;t\in[0,\,T]$ such that $x=y\circ G_{u(t)}^T$. Let
$y_n\longrightarrow y,\quad y_n\in\mathcal A_{q_0}((\mathrm e^{-gV^2})_\ast f_1)(t)$. Hence
$x_n=y_n\circ G_{u(t)}^T$ is a sequence on $\mathcal A_{q_0}((\mathrm e^{-gV^2})_\ast f_1)(t)\circ G_{u(t)}^t$ that
converges to $x$.
\end{proof}
So system \eqref{affsfac} is approximately controllable in time $t$ if
\begin{equation}
\overline{\mathcal A_{q}((\mathrm e^{-gV^2})_\ast f_1)(t)}\circ\{G_{v^2(t)}^t\mid
v^2(t)\in\mathbb R^r\}=\mathbb R^n\quad\forall q\in\mathbb R^n.
\end{equation}
If $g_i$ are constant vector fields --- $g_i(q)=g_i\quad i=1,\,\dots\,r\quad
q\in\mathbb R^n$ --- they commute and the systems
$\dot q=(\mathrm e^{-gX(t)})_\ast f_1(q)$ and
$\dot{q}=f_1(q+gX(t))$ coincide. A corollary of this is
\begin{Corollary}\label{clAfxeV.Rn}
System \eqref{affsfac} (with $g$ constant) is approximately controllable in time $t$ if
$$\overline{\mathcal A_q (f_{1X})(t)}\circ
\{\mathrm e^{gV^2(t)}\}=\mathbb R^n\quad\forall q\in\mathbb R^n.\quad\footnotemark
$$ \footnotetext{Here, for more precision, we should write $\{\mathrm e^{gV^2(t)}\mid V^2(t)\in\mathbb R^r\}$.} Where
\begin{align*}
 f_{1X}(q)&:=f_1(q+gX(t))=f(q+gX(t))+g(q+gX(t))v^1\\
&=f_X(q)+gv^1\quad X(t)\in\mathbb R^r.
\end{align*}
In particular it is approximately controllable in time $t$ if
$$
\overline{\mathcal A_q (f_{1X})(t)}=\mathbb R^n\quad\forall q\in\mathbb R^n.
$$ 
\end{Corollary}

\subsection{Convexification.}\label{sSecCon}
If for some constant vector $\gamma\in\mathbb R^n$, $f(q)+\gamma$ belong to the convex set $Conv\{f_X\mid\,X\in\mathbb R^r\}$, then for every $u^1\in\mathbb R^r$
\begin{multline*}
f(q)+\gamma+gu^1\in Conv\{f_X\mid\,X\in\mathbb R^r\}+span(g)\\
\subseteq Conv\{f_X+gv^1\mid\,X,\,v^1\in\mathbb R^r\}.\footnotemark
\end{multline*}
\footnotetext{Here $span(g)$ means the span of the columns of $g$. $Conv(A)$ stays for convexification of the set $A$.}
This means that we can follow any of the vector fields $f(q)+\gamma+gv^1$ without changing the closure of the attainable set at time $t$ (recall that convexification does not change the closure of attainable set at time $t$ -- see \cite{jur}). In particular system \eqref{affsfac} is approximately controllable at time $t$ if
$$
\overline{\mathcal A_q (f(q)+\gamma+gv^1)(t)}=\mathbb R^n\quad\forall q\in\mathbb R^n.
$$
\subsection{Extraction.}
Let $C$ be a cone and suppose that
$$
f(q)+C\subseteq Conv\{f_X\mid\,X\in\mathbb R^r\}.
$$
Then putting $G:=span(g)$
\begin{align*}
f(q)+C+G &\subseteq f(q)+Conv(C)+G\subseteq Conv\{f_X(q)+G\mid\,X\in\mathbb R^r\}\\
&=Conv\{f_{1X}(q)\mid\,X\in\mathbb R^r\}.
\end{align*}
Now from $Conv(C)+G$ we extract the linear space
$$
G^1:=(G+Conv(C))\cap(G-Conv(C)).
$$
We shall call the directions on $G^1$ \emph{ ``extracted'' directions}. Since clearly $G\subseteq G^1$ because $0\in C$, those directions in $G$ will be called \emph{``old'' directions} and, those in $G^1\setminus G$ \emph{``new'' directions}.\\ 
Adding new directions does not change the closure of attainable sets so, we can say that system \eqref{affsfac} is appoximately controllable in time $t$ if the ``bigger'' system $\dot q=f(q)+g_1v^1$ is, where $v^1\in\mathbb R^{r_1}$, $r_1\;(\geq r)$ is the dimension of $G^1$ and $g_1$ is a matrix whose $r_1$ columns are  vectors spanning $G^1$.
\subsection{Iterating FCEs.}
Iterating FCE Procedures we obtain an increasing sequence
$$
G=:G^0\subseteq G^1 \subseteq\,\dots\,\subseteq G^j \subseteq\, \dots
$$
of subspaces of controlled directions without changing the closure of the attainable set at time $t$. Obviously if for some $p\in\mathbb N$ we have $G^p=\mathbb R^n$, then the controllability in time $t$ is an immediate consequence of Corollary \ref{clAfxeV.Rn} (note that in such a case we can set for $V(t)$ any vector from $\mathbb R^n$).

\section{The Projection onto $H$.}\label{ScProj}
Let $\mathbf L^2(R)_t$ the subspace of $\mathbf L^2(R)$ defined by
$$
\mathbf L^2(R)_t:=\{v\in\mathbf L^2(R)\mid\,v\cdot\mathbf n=0\,\text{on}\,\Gamma\}.
$$
Let $v\in\mathbf L^2(R)_t$. We can write $v$ as
$$
v=\begin{pmatrix}v_1\\v_2\end{pmatrix}=\begin{pmatrix}\sum_{k\in\mathbb N_0\times\mathbb N}v_{1k}\sin\Bigl(\frac{k_1\pi x_1}{a}\Bigr)\cos\Bigl(\frac{k_2\pi x_2}{b}\Bigr)\\\sum_{k\in\mathbb N\times\mathbb N_0}v_{2k}\cos\Bigl(\frac{k_1\pi x_1}{a}\Bigr)\sin\Bigl(\frac{k_2\pi x_2}{b}\Bigr)\end{pmatrix}.
$$
Now we want to write $v$ as a sum of an element $u\in H$ (a solenoidal element) and a gradient of a function $q\in L^2(R)$ --- $v=u+\nabla q$.\\
We put
\begin{align}
&\begin{pmatrix}\sum_{k\in\mathbb N_0\times\mathbb N}v_{1k}\sin\Bigl(\frac{k_1\pi x_1}{a}\Bigr)\cos\Bigl(\frac{k_2\pi x_2}{b}\Bigr)\\\sum_{k\in\mathbb N\times\mathbb N_0}v_{2k}\cos\Bigl(\frac{k_1\pi x_1}{a}\Bigr)\sin\Bigl(\frac{k_2\pi x_2}{b}\Bigr)\end{pmatrix}=u +\nabla q;\label{proj1}\\
&u:=\sum_{k\in\mathbb N_0^2}u_kW_k;\qquad q:=\sum_{k\in\mathbb N^2\setminus\{(0,\,0)\}}q_k\cos\Bigl(\frac{k_1\pi x_1}{a}\Bigr)\cos\Bigl(\frac{k_2\pi x_2}{b}\Bigr).\notag
\end{align}
Then
$$
\nabla q=\begin{pmatrix}\sum_{k\in\mathbb N^2}-\frac{k_1\pi}{a}q_k\sin\Bigl(\frac{k_1\pi x_1}{a}\Bigr)\cos\Bigl(\frac{k_2\pi x_2}{b}\Bigr)\\\sum_{k\in\mathbb N^2}-\frac{k_2\pi}{b}q_k\cos\Bigl(\frac{k_1\pi x_1}{a}\Bigr)\sin\Bigl(\frac{k_2\pi x_2}{b}\Bigr)\end{pmatrix}.
$$
For $k\in\mathbb N^2\setminus\mathbb N_0^2$ we find
$$
\begin{cases}
-q_k\frac{k_1\pi}{a}=v_{1k},\quad&\text{for}\quad k\in\mathbb N_0\times\{0\}\\
-q_k\frac{k_2\pi}{b}=v_{2k},\quad&\text{for}\quad k\in\{0\}\times\mathbb N_0\\

\end{cases}
$$
or,
$$
q_k=\begin{cases}-\frac{a}{k_1\pi}v_{1k},\quad&k\in\mathbb N_0\times\{0\}\\
-\frac{b}{k_2\pi}v_{2k},\quad&k\in\{0\}\times\mathbb N_0\end{cases}.
$$
and, for $k\in \mathbb N_0^2$ we obtain
\begin{align}
-\frac{k_1\pi}{a}q_k=&v_{1k}+\frac{k_2\pi}{b}u_k\label{q2eq1}\\
-\frac{k_2\pi}{b}q_k=&v_{2k}-\frac{k_1\pi}{a}u_k\label{q2eq2}.
\end{align}
Multiplying \eqref{q2eq1} by $-\frac{k_2\pi}{b}$ and \eqref{q2eq2} by $\frac{k_1\pi}{a}$ and then adding the products we obtain
$$
\frac{k_2\pi}{b}v_{1k}-\frac{k_1\pi}{a}v_{2k}=\bar ku_k,
$$
obtaining in this way $u$ from $v$.\\
Similarly we can obtain
$$
-\bigl(\frac{k_2\pi}{b}\bigr)^2q_k-\bigl(\frac{k_1\pi}{a}\bigr)^2q_k=\frac{k_1\pi}{a}v_{1k}+\frac{k_2\pi}{b}v_{2k}
$$
or
$$
\bar kq_k=\frac{k_1\pi}{a}v_{1k}+\frac{k_2\pi}{b}v_{2k}.
$$
Therefore, given
$$
v=\begin{pmatrix}v_1\\v_2\end{pmatrix}=\begin{pmatrix}\sum_{k\in\mathbb N_0\times\mathbb N}v_{1k}\sin\Bigl(\frac{k_1\pi x_1}{a}\Bigr)\cos\Bigl(\frac{k_2\pi x_2}{b}\Bigr)\\\sum_{k\in\mathbb N\times\mathbb N_0}v_{2k}\cos\Bigl(\frac{k_1\pi x_1}{a}\Bigr)\sin\Bigl(\frac{k_2\pi x_2}{b}\Bigr)\end{pmatrix}\in\mathbf L^2(R)_t
$$
we can write it, in a unique way, as a sum $v=u+\nabla q$, where $u\in H$ and $q\in L^2(R)$. So the projection map onto the solenoidal (divergence free) space $H$
\begin{align*}
P^\nabla:\mathbf L^2(R)_t&\to H\\
v&\mapsto u
\end{align*}
is well defined and we have
\begin{align}
u:=P^\nabla v&=\sum_{k\in\mathbb N_0^2}\frac{1}{\bar k}\bigl(\frac{k_2\pi}{b}v_{1k}-\frac{k_1\pi}{a}v_{2k}\bigr)\label{PHv}\\
q&=q_1+q_2\notag\\
&q_1=-\sum_{k\in\mathbb N_0\times\{0\}}v_{1k}\frac{a}{k_1\pi}\cos\Bigl(\frac{k_1\pi x_1}{a}\Bigr)-\sum_{k\in\{0\}\times\mathbb N_0}v_{2k}\frac{b}{k_2\pi}\cos\Bigl(\frac{k_2\pi x_2}{b}\Bigr)\label{q1ofv}\\
&q_2=\sum_{k\in\mathbb N_0^2}\frac{1}{\bar k}\Biggl(v_{1k}\frac{k_1\pi}{a}+v_{2k}\frac{k_2\pi}{b}\Biggr)\cos\Bigl(\frac{k_1\pi x_1}{a}\Bigr)\cos\Bigl(\frac{k_2\pi x_2}{b}\Bigr)\label{q2ofv}.
\end{align}
It remains to show that $u\in H$ and $q_1,\,q_2\in L^2(R)$. That follows from the fact that both $v_1$ and $v_2$ are in $L^2$:
\begin{align*}
|u|^2&=\sum_{k\in\mathbb N_0^2}-\bar k\frac{1}{\bar k^2}\bigl(\frac{k_2\pi}{b}v_{1k}-\frac{k_1\pi}{a}v_{2k}\bigr)^2\frac{ab}{4}\\
&\leq \sum_{k\in\mathbb N_0^2}\bigl(|v_{1k}|+|v_{2k}|\bigr)^2\frac{ab}{4}
\leq \bigl(\|v_{1k}\|_{L^2}+\|v_{2k}\|_{L^2}\bigr)^2;
\end{align*}
\begin{align*}
\|q_1\|_{L^2}^2&\leq\sum_{k\in\mathbb N_0\times\{0\}}\bigl(v_{1k}\frac{a}{k_1\pi}\bigr)^2\frac{a}{2}+\sum_{k\in\{0\}\times\mathbb N_0}\bigl(v_{2k}\frac{b}{k_2\pi}\bigr)^2\frac{b}{2}\\
&\leq\frac{a^2+b^2}{\pi^2}\bigl(\sum_{k\in\mathbb N_0\times\{0\}}\bigl(v_{1k}^2\frac{a}{2}+\sum_{k\in\{0\}\times\mathbb N_0}v_{2k}^2\frac{b}{2}\bigr)\\
&\leq \frac{a^2+b^2}{\pi^2}\bigl(\|v_{1k}\|_{L^2}^2+\|v_{2k}\|_{L^2}^2\bigr);
\end{align*}
\begin{align*}
q_2&=\sum_{k\in\mathbb N_0^2}\frac{1}{\bar k^2}\Biggl(v_{1k}\frac{k_1\pi}{a}+v_{2k}\frac{k_2\pi}{b}\Biggr)^2\frac{ab}{4}\\
&\leq \sum_{k\in\mathbb N_0^2}-\frac{1}{\bar k}\bigl(|v_{1k}|+|v_{2k}|\bigr)^2\frac{ab}{4}\leq -\frac{1}{\,\overline{(1,\,1)}}\bigl(\|v_{1k}\|_{L^2}+\|v_{2k}\|_{L^2}\bigr)^2.
\end{align*}
\section{Spectral Method.}
We want to use {\em Spectral Algorithm}, i.e.,we want to study NS equation in coordinates corresponding to the basis of eigenfunctions $\mathcal W$ defined in \eqref{basisW}.\par
We consider the equation
$$
u_t+(u\cdot\nabla)u+\nabla p=\Delta u+F+v;
$$
and write
$$
u=\sum_{k\in\mathbb N_0^2}u_kW_k.
$$
Then we compute $(u\cdot\nabla)u=\begin{pmatrix}u\cdot\nabla u_1\\u\cdot\nabla u_2\end{pmatrix}:$
\begin{align*}
\nabla u_1&=\begin{pmatrix}\sum_{k\in\mathbb N_0^2}u_k(-\frac{k_2\pi}{b})\frac{k_1\pi}{a}\cos\Bigl(\frac{k_1\pi x_1}{a}\Bigr)\cos\Bigl(\frac{k_2\pi x_2}{b}\Bigr)\\\sum_{k\in\mathbb N_0^2}u_k(-\frac{k_2\pi}{b})(-\frac{k_2\pi}{b})\sin\Bigl(\frac{k_1\pi x_1}{a}\Bigr)\sin\Bigl(\frac{k_2\pi x_2}{b}\Bigr)\end{pmatrix}\\
\nabla u_2&=\begin{pmatrix}\sum_{k\in\mathbb N_0^2}u_k\frac{k_1\pi}{a}(-\frac{k_1\pi}{a})\sin\Bigl(\frac{k_1\pi x_1}{a}\Bigr)\sin\Bigl(\frac{k_2\pi x_2}{b}\Bigr)\\\sum_{k\in\mathbb N_0^2}u_k\frac{k_1\pi}{a}\frac{k_2\pi}{b}\cos\Bigl(\frac{k_1\pi x_1}{a}\Bigr)\cos\Bigl(\frac{k_2\pi x_2}{b}\Bigr)\end{pmatrix}.
\end{align*}
Then
\begin{align*}
&u\cdot\nabla u_1\\
=&\sum_{m,n\in\mathbb N_0^2}u_mu_n(-\frac{n_2\pi}{b})(-\frac{m_2\pi}{b})\frac{m_1\pi}{a}\sin\Bigl(\frac{n_1\pi x_1}{a}\Bigr)\cos\Bigl(\frac{n_2\pi x_2}{b}\Bigr)\\
&\qquad\qquad\qquad\qquad\times\cos\Bigl(\frac{m_1\pi x_1}{a}\Bigr)\cos\Bigl(\frac{m_2\pi x_2}{b}\Bigr)\\
+&\sum_{m,n\in\mathbb N_0^2}u_mu_n\frac{n_1\pi}{a}(-\frac{m_2\pi}{b})(-\frac{m_2\pi}{b})\cos\Bigl(\frac{n_1\pi x_1}{a}\Bigr)\sin\Bigl(\frac{n_2\pi x_2}{b}\Bigr)\\
&\qquad\qquad\qquad\qquad\times\sin\Bigl(\frac{m_1\pi x_1}{a}\Bigr)\sin\Bigl(\frac{m_2\pi x_2}{b}\Bigr)\\
=&\sum_{m,n\in\mathbb N_0^2}u_mu_n\frac{\pi^2}{4ab}m_1n_2\frac{\pi}{b}m_2\Bigl[\sin\Bigl(\frac{(n_1+m_1)\pi x_1}{a}\Bigr)+\sin\Bigl(\frac{(n_1-m_1)\pi x_1}{a}\Bigr)\Bigr]\\& \qquad\qquad\qquad\qquad \times\Bigl[\cos\Bigl(\frac{(n_2+m_2)\pi x_2}{b}\Bigr)+\cos\Bigl(\frac{(n_2-m_2)\pi x_2}{b}\Bigr)\Bigr]\\
+&\sum_{m,n\in\mathbb N_0^2}u_mu_n\frac{\pi^2}{4ab}m_1n_2\frac{\pi}{b}n_2\Bigl[\sin\Bigl(\frac{(n_1+m_1)\pi x_1}{a}\Bigr)+\sin\Bigl(\frac{(n_1-m_1)\pi x_1}{a}\Bigr)\Bigr]\\
& \qquad\qquad\qquad\qquad \times\Bigl[-\cos\Bigl(\frac{(n_2+m_2)\pi x_2}{b}\Bigr)+\cos\Bigl(\frac{(n_2-m_2)\pi x_2}{b}\Bigr)\Bigr].
\end{align*}
We then arrive to
\begin{align}
&u\cdot\nabla u_1\notag\\
=&\sum_{m,n\in\mathbb N_0^2}u_mu_n\frac{\pi^2}{4ab}m_1n_2\frac{\pi}{b}(m_2-n_2)\sin\Bigl(\frac{(n_1+m_1)\pi x_1}{a}\Bigr)\cos\Bigl(\frac{(n_2+m_2)\pi x_2}{b}\Bigr)\notag\\
+&\sum_{m,n\in\mathbb N_0^2}u_mu_n\frac{\pi^2}{4ab}m_1n_2\frac{\pi}{b}(m_2+n_2)\sin\Bigl(\frac{(n_1+m_1)\pi x_1}{a}\Bigr)\cos\Bigl(\frac{(n_2-m_2)\pi x_2}{b}\Bigr)\notag\\
+&\sum_{m,n\in\mathbb N_0^2}u_mu_n\frac{\pi^2}{4ab}m_1n_2\frac{\pi}{b}(m_2-n_2)\sin\Bigl(\frac{(n_1-m_1)\pi x_1}{a}\Bigr)\cos\Bigl(\frac{(n_2+m_2)\pi x_2}{b}\Bigr)\notag\\
+&\sum_{m,n\in\mathbb N_0^2}u_mu_n\frac{\pi^2}{4ab}m_1n_2\frac{\pi}{b}(m_2+n_2)\sin\Bigl(\frac{(n_1-m_1)\pi x_1}{a}\Bigr)\cos\Bigl(\frac{(n_2-m_2)\pi x_2}{b}\Bigr)\label{unu1.1}.
\end{align}
Similarly:
\begin{align*}
&u\cdot\nabla u_2\\
=&\sum_{m,n\in\mathbb N_0^2}u_mu_n(-\frac{n_2\pi}{b})\frac{m_1\pi}{a}(-\frac{m_1\pi}{a})\sin\Bigl(\frac{n_1\pi x_1}{a}\Bigr)\cos\Bigl(\frac{n_2\pi x_2}{b}\Bigr)\\
&\qquad\qquad\qquad\qquad \times\sin\Bigl(\frac{m_1\pi x_1}{a}\Bigr)\sin\Bigl(\frac{m_2\pi x_2}{b}\Bigr)\\
+&\sum_{m,n\in\mathbb N_0^2}u_mu_n\frac{n_1\pi}{a}\frac{m_1\pi}{a}\frac{m_2\pi}{b}\cos\Bigl(\frac{n_1\pi x_1}{a}\Bigr)\sin\Bigl(\frac{n_2\pi x_2}{b}\Bigr)\\
&\qquad\qquad\qquad\qquad \times\cos\Bigl(\frac{m_1\pi x_1}{a}\Bigr)\cos\Bigl(\frac{m_2\pi x_2}{b}\Bigr)\\
=&\sum_{m,n\in\mathbb N_0^2}u_mu_n\frac{\pi^2}{4ab}m_1n_2\frac{\pi}{a}m_1\Bigl[-\cos\Bigl(\frac{(n_1+m_1)\pi x_1}{a}\Bigr)+\cos\Bigl(\frac{(n_1-m_1)\pi x_1}{a}\Bigr)\Bigr]\\
&\qquad\qquad\qquad\qquad \times\Bigl[\sin\Bigl(\frac{(n_2+m_2)\pi x_2}{b}\Bigr)-\sin\Bigl(\frac{(n_2-m_2)\pi x_2}{b}\Bigr)\Bigr]\\
+&\sum_{m,n\in\mathbb N_0^2}u_mu_n\frac{\pi^2}{4ab}m_1n_2\frac{\pi}{a}n_1\Bigl[\cos\Bigl(\frac{(n_1+m_1)\pi x_1}{a}\Bigr)+\cos\Bigl(\frac{(n_1-m_1)\pi x_1}{a}\Bigr)\Bigr]\\
&\qquad\qquad\qquad\qquad \times\Bigl[\sin\Bigl(\frac{(n_2+m_2)\pi x_2}{b}\Bigr)-\cos\Bigl(\frac{(n_2-m_2)\pi x_2}{b}\Bigr)\Bigr]
\end{align*}
and then,
\begin{align}
&u\cdot\nabla u_2\notag\\
=&\sum_{m,n\in\mathbb N_0^2}u_mu_n\frac{\pi^2}{4ab}m_1n_2\frac{\pi}{a}(n_1-m_1)\cos\Bigl(\frac{(n_1+m_1)\pi x_1}{a}\Bigr)\sin\Bigl(\frac{(n_2+m_2)\pi x_2}{b}\Bigr)\notag\\
+&\sum_{m,n\in\mathbb N_0^2}u_mu_n\frac{\pi^2}{4ab}m_1n_2\frac{\pi}{a}(m_1-n_1)\cos\Bigl(\frac{(n_1+m_1)\pi x_1}{a}\Bigr)\sin\Bigl(\frac{(n_2-m_2)\pi x_2}{b}\Bigr)\notag\\
+&\sum_{m,n\in\mathbb N_0^2}u_mu_n\frac{\pi^2}{4ab}m_1n_2\frac{\pi}{a}(m_1+n_1)\cos\Bigl(\frac{(n_1-m_1)\pi x_1}{a}\Bigr)\sin\Bigl(\frac{(n_2+m_2)\pi x_2}{b}\Bigr)\notag\\
+&\sum_{m,n\in\mathbb N_0^2}u_mu_n\frac{\pi^2}{4ab}m_1n_2\frac{\pi}{a}(-m_1-n_1)\cos\Bigl(\frac{(n_1-m_1)\pi x_1}{a}\Bigr)\sin\Bigl(\frac{(n_2-m_2)\pi x_2}{b}\Bigr)\label{unu2.1}.
\end{align}
Grouping the terms in sum \eqref{unu1.1} and \eqref{unu2.1} envolving the product $u_mu_n$ we obtain
\begin{align}
&u\cdot\nabla u_1\notag\\
=&\sum_{\begin{subarray}{l}m,n\in\mathbb N_0^2\\m<n\end{subarray}}u_mu_n\frac{\pi^2}{4ab}m\wedge n\frac{\pi}{b}(m_2-n_2)\sin\Bigl(\frac{(n_1+m_1)\pi x_1}{a}\Bigr)\cos\Bigl(\frac{(n_2+m_2)\pi x_2}{b}\Bigr)\notag\\
+&\sum_{\begin{subarray}{l}m,n\in\mathbb N_0^2\\m<n\end{subarray}}u_mu_n\frac{\pi^2}{4ab}m\vee n\frac{\pi}{b}(m_2+n_2)\sin\Bigl(\frac{(n_1+m_1)\pi x_1}{a}\Bigr)\cos\Bigl(\frac{(n_2-m_2)\pi x_2}{b}\Bigr)\notag\\
+&\sum_{\begin{subarray}{l}m,n\in\mathbb N_0^2\\m<n\end{subarray}}u_mu_n\frac{\pi^2}{4ab}m\vee n\frac{\pi}{b}(m_2-n_2)\sin\Bigl(\frac{(n_1-m_1)\pi x_1}{a}\Bigr)\cos\Bigl(\frac{(n_2+m_2)\pi x_2}{b}\Bigr)\notag\\
+&\sum_{\begin{subarray}{l}m,n\in\mathbb N_0^2\\m<n\end{subarray}}u_mu_n\frac{\pi^2}{4ab}m\wedge n\frac{\pi}{b}(m_2+n_2)\sin\Bigl(\frac{(n_1-m_1)\pi x_1}{a}\Bigr)\cos\Bigl(\frac{(n_2-m_2)\pi x_2}{b}\Bigr)\notag\\
+&\sum_{n\in\mathbb N_0^2} u_n^2\frac{\pi^2}{4ab}n_1n_2\frac{\pi}{b}(2n_2)\sin\Bigl(\frac{2n_1\pi x_1}{a}\Bigr)\quad\footnotemark\label{unu1.2}
\end{align}
\footnotetext{This last sum cooresponds to the sum over the diagonal $\{(m,\,n)\in\mathbb (N_0^2)^2\mid\,m=n\}$.}
and,
\begin{align}
&u\cdot\nabla u_2\notag\\
=&\sum_{\begin{subarray}{l}m,n\in\mathbb N_0^2\\m<n\end{subarray}}u_mu_n\frac{\pi^2}{4ab}m\wedge n\frac{\pi}{a}(n_1-m_1)\cos\Bigl(\frac{(n_1+m_1)\pi x_1}{a}\Bigr)\sin\Bigl(\frac{(n_2+m_2)\pi x_2}{b}\Bigr)\notag\\
+&\sum_{\begin{subarray}{l}m,n\in\mathbb N_0^2\\m<n\end{subarray}}u_mu_n\frac{\pi^2}{4ab}m\vee n\frac{\pi}{a}(m_1-n_1)\cos\Bigl(\frac{(n_1+m_1)\pi x_1}{a}\Bigr)\sin\Bigl(\frac{(n_2-m_2)\pi x_2}{b}\Bigr)\notag\\
+&\sum_{\begin{subarray}{l}m,n\in\mathbb N_0^2\\m<n\end{subarray}}u_mu_n\frac{\pi^2}{4ab}m\vee n\frac{\pi}{a}(m_1+n_1)\cos\Bigl(\frac{(n_1-m_1)\pi x_1}{a}\Bigr)\sin\Bigl(\frac{(n_2+m_2)\pi x_2}{b}\Bigr)\notag\\
+&\sum_{\begin{subarray}{l}m,n\in\mathbb N_0^2\\m<n\end{subarray}}u_mu_n\frac{\pi^2}{4ab}m\wedge n\frac{\pi}{a}(-m_1-n_1)\cos\Bigl(\frac{(n_1-m_1)\pi x_1}{a}\Bigr)\sin\Bigl(\frac{(n_2-m_2)\pi x_2}{b}\Bigr)\notag\\
+&\sum_{n\in\mathbb N_0^2} u_n^2\frac{\pi^2}{4ab}n_1n_2\frac{\pi}{a}(2n_1)\sin\Bigl(\frac{2n_2\pi x_2}{b}\Bigr)\label{unu2.2}.
\end{align}
Where the order under the sum sign, in \eqref{unu1.2} and \eqref{unu2.2}, is the lexycographical one and, by definition:
$$
m\vee n:=m_1n_2+n_1m_2;\qquad\qquad m\wedge n:=m_1n_2-n_1m_2.
$$
Put
$$
C_{m,n}^\vee=\frac{\pi^2}{4ab}m\vee n\And C_{m,n}^\wedge=\frac{\pi^2}{4ab}m\wedge n.
$$
We rewrite\footnote{$m<n$ implies $m_1\leq n_1$.}
\begin{align}
&u\cdot\nabla u_1\notag\\
=&\sum_{\begin{subarray}{l}m,n\in\mathbb N_0^2\\m<n\end{subarray}}u_mu_nC_{m,n}^\wedge\frac{\pi}{b}(m_2-n_2)\sin\Bigl(\frac{|n_1+m_1|\pi x_1}{a}\Bigr)\cos\Bigl(\frac{|n_2+m_2|\pi x_2}{b}\Bigr)\notag\\
+&\sum_{\begin{subarray}{l}m,n\in\mathbb N_0^2\\m<n\end{subarray}}u_mu_nC_{m,n}^\vee\frac{\pi}{b}(m_2+n_2)\sin\Bigl(\frac{|n_1+m_1|\pi x_1}{a}\Bigr)\cos\Bigl(\frac{|n_2-m_2|\pi x_2}{b}\Bigr)\notag\\
+&\sum_{\begin{subarray}{l}m,n\in\mathbb N_0^2\\m<n\end{subarray}}u_mu_nC_{m,n}^\vee\frac{\pi}{b}(m_2-n_2)\sin\Bigl(\frac{|n_1-m_1|\pi x_1}{a}\Bigr)\cos\Bigl(\frac{|n_2+m_2|\pi x_2}{b}\Bigr)\notag\\
+&\sum_{\begin{subarray}{l}m,n\in\mathbb N_0^2\\m<n\end{subarray}}u_mu_nC_{m,n}^\wedge\frac{\pi}{b}(m_2+n_2)\sin\Bigl(\frac{|n_1-m_1|\pi x_1}{a}\Bigr)\cos\Bigl(\frac{|n_2-m_2|\pi x_2}{b}\Bigr)\notag\\
+&\sum_{n\in\mathbb N_0^2} u_n^2\frac{\pi^2}{4ab}n_1n_2\frac{\pi}{b}(2n_2)\sin\Bigl(\frac{2n_1\pi x_1}{a}\Bigr)\label{unu1.3};
\end{align}
and
\begin{align}
&u\cdot\nabla u_2\notag\\
=&\sum_{\begin{subarray}{l}m,n\in\mathbb N_0^2\\m<n\end{subarray}}u_mu_nC_{m,n}^\wedge\frac{\pi}{a}(n_1-m_1)\cos\Bigl(\frac{|n_1+m_1|\pi x_1}{a}\Bigr)\sin\Bigl(\frac{|n_2+m_2|\pi x_2}{b}\Bigr)\notag\\
+&\sum_{\begin{subarray}{l}m,n\in\mathbb N_0^2\\m<n\\m_2<n_2\end{subarray}}u_mu_nC_{m,n}^\vee\frac{\pi}{a}(m_1-n_1)\cos\Bigl(\frac{|n_1+m_1|\pi x_1}{a}\Bigr)\sin\Bigl(\frac{|n_2-m_2|\pi x_2}{b}\Bigr)\notag\\
+&\sum_{\begin{subarray}{l}m,n\in\mathbb N_0^2\\m<n\\m_2>n_2\end{subarray}}u_mu_nC_{m,n}^\vee\frac{\pi}{a}(-m_1+n_1)\cos\Bigl(\frac{|n_1+m_1|\pi x_1}{a}\Bigr)\sin\Bigl(\frac{|n_2-m_2|\pi x_2}{b}\Bigr)\notag\\
+&\sum_{\begin{subarray}{l}m,n\in\mathbb N_0^2\\m<n\end{subarray}}u_mu_nC_{m,n}^\vee\frac{\pi}{a}(m_1+n_1)\cos\Bigl(\frac{|n_1-m_1|\pi x_1}{a}\Bigr)\sin\Bigl(\frac{|n_2+m_2|\pi x_2}{b}\Bigr)\notag\\
+&\sum_{\begin{subarray}{l}m,n\in\mathbb N_0^2\\m<n\\m_2<n_2\end{subarray}}u_mu_nC_{m,n}^\wedge\frac{\pi}{a}(-m_1-n_1)\cos\Bigl(\frac{|n_1-m_1|\pi x_1}{a}\Bigr)\sin\Bigl(\frac{|n_2-m_2|\pi x_2}{b}\Bigr)\notag\\
+&\sum_{\begin{subarray}{l}m,n\in\mathbb N_0^2\\m<n\\m_2>n_2\end{subarray}}u_mu_nC_{m,n}^\wedge\frac{\pi}{a}(m_1+n_1)\cos\Bigl(\frac{|n_1-m_1|\pi x_1}{a}\Bigr)\sin\Bigl(\frac{|n_2-m_2|\pi x_2}{b}\Bigr)\notag\\
+&\sum_{n\in\mathbb N_0^2} u_n^2\frac{\pi^2}{4ab}n_1n_2\frac{\pi}{a}(2n_1)\sin\Bigl(\frac{2n_2\pi x_2}{b}\Bigr)\label{unu2.3}.
\end{align}
Now we project $(u\cdot\nabla)u$ in $H$:\footnote{See section \ref{ScProj}.} First we put
\begin{eqnarray*}
n(--)m:=(n_1-m_1,n_2-m_2);
& &n(-+)m:=(n_1-m_1,n_2+m_2);\\
n(+-)m:=(n_1+m_1,n_2-m_2);
& &n(++)m:=(n_1+m_1,n_2+m_2);\\
& &(a,b)^+:=(|a|,\,|b|),\quad a,b\in\mathbb Z\quad\quad\footnotemark
\end{eqnarray*}\footnotetext{$\mathbb Z$ represents the set of integer numbers.}
and then the projection can be written as

\begin{align}
&P^\nabla[(u\cdot\nabla)u]\notag\\
=&\sum_{\begin{subarray}{l}m,n\in\mathbb N_0^2\\m<n\end{subarray}}\frac{u_mu_nC_{m,n}^\wedge}{\,\overline{(n(++)m)^+}\,}D_{++}W_{(n(++)m)^+}\notag\\
+&\sum_{\begin{subarray}{l}m,n\in\mathbb N_0^2\\m<n\\m_2<n_2\end{subarray}}\frac{u_mu_nC_{m,n}^\vee}{\,\overline{(n(+-)m)^+}\,}D_{+-}^1W_{(n(+-)m)^+}\notag\\
+&\sum_{\begin{subarray}{l}m,n\in\mathbb N_0^2\\m<n\\m_2>n_2\end{subarray}}\frac{u_mu_nC_{m,n}^\vee}{\,\overline{(n(+-)m)^+}\,}D_{+-}^2W_{(n(+-)m)^+}\notag\\
+&\sum_{\begin{subarray}{l}m,n\in\mathbb N_0^2\\m<n\end{subarray}}\frac{u_mu_nC_{m,n}^\vee}{\,\overline{(n(-+)m)^+}\,}D_{-+}W_{(n(-+)m)^+}\notag\\
+&\sum_{\begin{subarray}{l}m,n\in\mathbb N_0^2\\m<n\\m_2<n_2\end{subarray}}\frac{u_mu_nC_{m,n}^\wedge}{\,\overline{(n(--)m)^+}\,}D_{--}^1W_{(n(--)m)^+}\notag\\
+&\sum_{\begin{subarray}{l}m,n\in\mathbb N_0^2\\m<n\\m_2>n_2\end{subarray}}\frac{u_mu_nC_{m,n}^\wedge}{\,\overline{(n(--)m)^+}\,}D_{--}^2W_{(n(--)m)^+}\label{Punu};
\end{align}
where
\begin{align*}
D_{++}=\frac{(n(++)m)_2^+\pi}{b}\frac{\pi}{b}(m_2-n_2)-\frac{(n(++)m)_1^+\pi}{a}\frac{\pi}{a}(n_1-m_1)\\
D_{+-}^1=\frac{(n(+-)m)_2^+\pi}{b}\frac{\pi}{b}(m_2+n_2)-\frac{(n(+-)m)_1^+\pi}{a}\frac{\pi}{a}(m_1-n_1)\\
D_{+-}^2=\frac{(n(+-)m)_2^+\pi}{b}\frac{\pi}{b}(m_2+n_2)-\frac{(n(+-)m)_1^+\pi}{a}\frac{\pi}{a}(n_1-m_1)\\
D_{-+}=\frac{(n(-+)m)_2^+\pi}{b}\frac{\pi}{b}(m_2-n_2)-\frac{(n(-+)m)_1^+\pi}{a}\frac{\pi}{a}(m_1+n_1)\\
D_{--}^1=\frac{(n(--)m)_2^+\pi}{b}\frac{\pi}{b}(m_2+n_2)-\frac{(n(--)m)_1^+\pi}{a}\frac{\pi}{a}(-m_1-n_1)\\
D_{--}^2=\frac{(n(--)m)_2^+\pi}{b}\frac{\pi}{b}(m_2+n_2)-\frac{(n(--)m)_1^+\pi}{a}\frac{\pi}{a}(m_1+n_1).
\end{align*}
Equation \eqref{Punu} can be reduced to
\begin{align*}
&P^\nabla[(u\cdot\nabla)u]\\
=&\sum_{\begin{subarray}{l}m,n\in\mathbb N_0^2\\m<n\end{subarray}}\frac{u_mu_nC_{m,n}^\wedge}{\,\overline{(n(++)m)^+}\,}(\bar n-\bar m)W_{(n(++)m)^+}\\
+&\sum_{\begin{subarray}{l}m,n\in\mathbb N_0^2\\m<n\\m_2<n_2\end{subarray}}\frac{u_mu_nC_{m,n}^\vee}{\,\overline{(n(+-)m)^+}\,}(\bar m-\bar n)W_{(n(+-)m)^+}\\
+&\sum_{\begin{subarray}{l}m,n\in\mathbb N_0^2\\m<n\\m_2>n_2\end{subarray}}\frac{u_mu_nC_{m,n}^\vee}{\,\overline{(n(+-)m)^+}\,}(\bar n-\bar m)W_{(n(+-)m)^+}\\
+&\sum_{\begin{subarray}{l}m,n\in\mathbb N_0^2\\m<n\end{subarray}}\frac{u_mu_nC_{m,n}^\vee}{\,\overline{(n(-+)m)^+}\,}(\bar n-\bar m)W_{(n(-+)m)^+}\\
+&\sum_{\begin{subarray}{l}m,n\in\mathbb N_0^2\\m<n\\m_2<n_2\end{subarray}}\frac{u_mu_nC_{m,n}^\wedge}{\,\overline{(n(--)m)^+}\,}(\bar m-\bar n)W_{(n(--)m)^+}\\
+&\sum_{\begin{subarray}{l}m,n\in\mathbb N_0^2\\m<n\\m_2>n_2\end{subarray}}\frac{u_mu_nC_{m,n}^\wedge}{\,\overline{(n(--)m)^+}\,}(\bar n-\bar m)W_{(n(--)m)^+}.
\end{align*}
Hence
\begin{align*}
&P^\nabla[(u\cdot\nabla)u]\\
=&\sum_{\begin{subarray}{l}m,n\in\mathbb N_0^2\\m<n\end{subarray}}\frac{u_mu_nC_{m,n}^\wedge}{\,\overline{(n(++)m)^+}\,}(\bar n-\bar m)W_{(n(++)m)^+}\\
+&\sum_{\begin{subarray}{l}m,n\in\mathbb N_0^2\\m<n\end{subarray}}-\frac{u_mu_nC_{m,n}^\vee}{\,\overline{(n(+-)m)^+}\,}(\bar n-\bar m)\mathrm{sign}(n_2-m_2)W_{(n(+-)m)^+}\\
+&\sum_{\begin{subarray}{l}m,n\in\mathbb N_0^2\\m<n\end{subarray}}\frac{u_mu_nC_{m,n}^\vee}{\,\overline{(n(-+)m)^+}\,}(\bar n-\bar m)W_{(n(-+)m)^+}\\
+&\sum_{\begin{subarray}{l}m,n\in\mathbb N_0^2\\m<n\end{subarray}}-\frac{u_mu_nC_{m,n}^\wedge}{\,\overline{(n(--)m)^+}\,}(\bar n-\bar m)\mathrm{sign}(n_2-m_2)W_{(n(--)m)^+}
\end{align*}
or, what is the same
\begin{align}
&P^\nabla[(u\cdot\nabla)u]\notag\\
=&\sum_{\begin{subarray}{l}m,n\in\mathbb N_0^2\\m<n\end{subarray}}\frac{u_mu_nC_{m,n}^\wedge}{\,\overline{(n(++)m)^+}\,}(\bar n-\bar m)W_{(n(++)m)^+}\notag\\
+&\sum_{\begin{subarray}{l}m,n\in\mathbb N_0^2\\m<n\end{subarray}}-\frac{u_mu_nC_{m,n}^\vee}{\,\overline{(n(+-)m)^+}\,}(\bar n-\bar m)\mathrm{sign}(n_2-m_2)W_{(n(+-)m)^+}\notag\\
+&\sum_{\begin{subarray}{l}m,n\in\mathbb N_0^2\\m<n\end{subarray}}\frac{u_mu_nC_{m,n}^\vee}{\,\overline{(n(-+)m)^+}\,}(\bar n-\bar m)\mathrm{sign}(n_1-m_1)W_{(n(-+)m)^+}\notag\\
+&\sum_{\begin{subarray}{l}m,n\in\mathbb N_0^2\\m<n\end{subarray}}-\frac{u_mu_nC_{m,n}^\wedge}{\,\overline{(n(--)m)^+}\,}(\bar n-\bar m)\mathrm{sign}(n_1-m_1)\mathrm{sign}(n_2-m_2)W_{(n(--)m)^+}\label{PHunu}
\end{align}
\section{The Infinite ODE System.}
We are now able to write an infinite system of ODEs related with the N-S Equation \eqref{nse1}
$$
u_t+(u\cdot\nabla)u+\nabla p=\nu\Delta u+F+v.
$$
As we can seen in the strong formulation of the N-S Problem (Problem \ref{s.p.leray}), at each time $\tau$ the last equality is  an equality between elements of $V^\prime$. However if we project this equality onto the solenoidal space $H$ we obtain
$$
P^\nabla\Bigl(u_t+(u\cdot\nabla)u+\nabla p\Bigr)=P^\nabla\Bigl(\nu\Delta u+F+v\Bigr)
$$
and, by \eqref{uDA.AuH} and the fact that $\tau\in D(A)$ a.e., we obtain
\begin{equation}\label{PrNSE}
u_t(\tau)=-P^\nabla Bu(\tau)+\nu\Delta u(\tau)+F+v(\tau).
\end{equation}
We suppose $F$ is solenoidal, otherwise we just take its solenoidal part. 

From
$$
u=\sum_{k\in\mathbb N_0^2}u_kW_k
$$
we derive
\begin{align*}
u_t&=\sum_{\begin{subarray}{l}m,n\in\mathbb N_0^2\\m<n\end{subarray}}(u_k)_tW_k\\
\Delta u&=\sum_{\begin{subarray}{l}m,n\in\mathbb N_0^2\\m<n\end{subarray}}\bar ku_k
\end{align*}
Therefore we may write \eqref{PrNSE} in the form of the infinite system of ODEs:
\begin{align}
\dot u_k&=\sum_{\begin{subarray}{l}m,n\in\mathbb N_0^2\\m<n\\(n(++)m)^+=k\end{subarray}}-\frac{u_mu_nC_{m,n}^\wedge}{\bar k}(\bar n-\bar m)\notag\\
&+\sum_{\begin{subarray}{l}m,n\in\mathbb N_0^2\\m<n\\(n(--)m)^+=k\end{subarray}}\frac{u_mu_nC_{m,n}^\wedge}{\bar k}(\bar n-\bar m)\mathrm{sign}(n_1-m_1)\mathrm{sign}(n_2-m_2)\notag\\
&+\sum_{\begin{subarray}{l}m,n\in\mathbb N_0^2\\m<n\\(n(-+)m)^+=k\end{subarray}}-\frac{u_mu_nC_{m,n}^\vee}{\bar k}(\bar n-\bar m)\mathrm{sign}{n_1-m_1}\notag\\
&+\sum_{\begin{subarray}{l}m,n\in\mathbb N_0^2\\m<n\\(n(+-)m)^+=k\end{subarray}}\frac{u_mu_nC_{m,n}^\vee}{\bar k}(\bar n-\bar m)\mathrm{sign}(n_2-m_2)\notag\\
&+\nu\bar ku_k+F_k+v_k.\label{infsys}
\end{align}
\subsection{Galerkin Approximations and G-Saturating Sets.}
Trying to make the writting simpler we introduce some notation referring to the coeficients appearing in system \eqref{infsys}:
\begin{Definition}
\begin{align*}
C_{m,n}^{++}&:=-\frac{\pi^2}{4ab}\frac{m\wedge n}{\,\overline{(m(++)n)^+}\,}(\bar n-\bar m)\\
C_{m,n}^{--}&:=\frac{\pi^2}{4ab}\frac{m\wedge n}{\,\overline{(m(--)n)^+}\,}(\bar n-\bar m)\mathrm{sign}(n_1-m_1)\mathrm{sign}(n_2-m_2)\\
C_{m,n}^{-+}&:=-\frac{\pi^2}{4ab}\frac{m\vee n}{\,\overline{(m(-+)n)^+}\,}(\bar n-\bar m)\mathrm{sign}(n_1-m_1)\\
C_{m,n}^{+-}&:=\frac{\pi^2}{4ab}\frac{m\vee n}{\,\overline{(m(+-)n)^+}\,}(\bar n-\bar m)\mathrm{sign}(n_2-m_2).
\end{align*}
\end{Definition}
Let $g$ be a set of (coordinates of) $r$ linearly independent vectors in $H$. Write also $g$ for the matrix whose columns are those $r$ vectors.\\
Rewriting system \eqref{infsys} with the directions in $g$ as the controlled ones, we obtain
$$
\dot u=f(u)+gv\qquad v\in\mathbb R^r
$$
or, equivalently
$$
\dot u_k=(f(u))_k+(gv)_k\qquad v\in\mathbb R^r,\;k\in\mathbb N_0^2
$$
where
\begin{eqnarray*}
(f(u))_k&:=&\sum_{\begin{subarray}{l}m,n\in\mathbb N_0^2\\m<n\\(n(++)m)^+=k\end{subarray}}u_mu_nC_{m,n}^{++}\notag\\
&+&\sum_{\begin{subarray}{l}m,n\in\mathbb N_0^2\\m<n\\(n(--)m)^+=k\end{subarray}}u_mu_nC_{m,n}^{--}\notag\\
&+&\sum_{\begin{subarray}{l}m,n\in\mathbb N_0^2\\m<n\\(n(-+)m)^+=k\end{subarray}}u_mu_nC_{m,n}^{-+}\notag\\
&+&\sum_{\begin{subarray}{l}m,n\in\mathbb N_0^2\\m<n\\(n(+-)m)^+=k\end{subarray}}u_mu_nC_{m,n}^{+-}\notag\\
&+&\nu\bar ku_k+F_k\\
&=& B_k(u,\,u)+\nu\bar k u_k+F_k,\notag\\
\end{eqnarray*}
defining $B_k(u,\,u):=(f(u))_k-(\nu\bar k u_k+F_k)$. Let us apply a FCE Procedure to this infinite system. Applying factorization to this system we obtain,\footnote{See section \ref{SecFCE}.} the factorized system
\begin{eqnarray}
(f_{1X}(u))_k&=&(f_X(u))_k+gv^1=(f(u))_k+B_k(u,\,gX)+B_k(gX,\,u)\notag\\
&+&\nu\bar k (gX)_k+B_k((gX),\,(gX))+gv^1\label{FacSysGG}
\end{eqnarray}
Now we put
\begin{eqnarray}
(\mathcal V_0(gX))_k&:=&(f(u))_k+gv^1;\notag\\
(\mathcal V_1(gX))_k&:=&B_k(u,\,gX)+B_k(gX,\,u)+\nu\bar k (gX)_k;\label{V012}\\
(\mathcal V_2(gX))_k&:=&B_k((gX),\,(gX));\notag
\end{eqnarray}
and we note that $\mathcal V_0,\,\mathcal V_1,\,\,\text{and}\,\,\mathcal V_2$ are respectively, independent, linear and bilinear on $gX$ vector fields.\par  
Now, ginen $X\in\mathbb R^r$ we have
$$
\mathcal V_0(g(-X))=\mathcal V_0(gX);\quad \mathcal V_1(g(-X))=-\mathcal V_1(gX)\;\text{and}\;\mathcal V_2(g(-X))=\mathcal V_2(gX).
$$
So,
$$
f(u)+B(gX,\,gX)=\frac{1}{2}\Big(f_X(u)+f_{-X}(u)\Big)\in Conv\{f_X(u)\mid\,X\in\mathbb R^r\}.
$$
The set
$$
\{B(gX,\,gX)\mid\,X\in\mathbb R^r\}=\{B(Y,\,Y)\mid\,Y\in span(g)\}
$$
is a cone.\par
Next, from
$$
Conv(\{B(gX,\,gX)\mid\,X\in\mathbb R^r\}+G=\{B(Y,\,Y)\mid\,Y\in span(g)\})+G,
$$
where $G:=span(g)$, we extract the subspace
$$
G^1:=\Big(G+Conv\{B(Y,\,Y)\mid\,Y\in G\}\Big)\cap\Big(G-Conv\{B(Y,\,Y)\mid\,Y\in G\}\Big).
$$
Now we present two definitions:
\begin{Definition}
A finite set of vectors $g\subseteq D(A)\quad$\footnote{Recall that $D(A)\subset H$ has been defined in \eqref{spaces}.} is said {\bf Saturating} for system \eqref{infsys} if the sequence $(G^j)_{n\in\mathbb N}$ of subspaces of $H$ defined recursively by
\begin{align*}
G^0&:=G=span(g);\\
G^{j+1}&:=\Big(G^j+Conv\{B(Y,\,Y)\mid\,Y\in G^j\}\Big)\notag\\
&\qquad\qquad\cap\Big(G^j-Conv\{B(Y,\,Y)\mid\,Y\in G^j\}\Big)\cap D(A);\quad\footnotemark
\end{align*}
\footnotetext{We intersect with $D(A)$ because, it is helpful for the study of controllability issues of the infinite system.}
satisfies
$$
\overline{\bigcup_{i\in\mathbb N}G^i}=H.
$$\par
We say that $g$ is {\bf G-Saturating} for system \eqref{infsys} if $span(g)=span\{W_k\mid\,k\in\mathcal K^1\}$ for some finite set of modes $\mathcal K^1\in\mathcal{FP}(\mathbb N_0^2)$, and there exists a sequence of subspaces $H^j$ such that $H^0:=span(g)$ and
\begin{align*}
&H^{j+1}\subseteq\Big(H^j+Conv\{B(Y,\,Y)\mid\,Y\in H^j\}\Big)\cap\Big(H^j-Conv\{B(Y,\,Y)\mid\,Y\in H^j\}\Big);\\
&\overline{\bigcup_{i\in\mathbb N}H^i}=H;
\end{align*}
and, besides there exists a finite subset
$\mathcal K^{j+2}\in\mathcal{FP}(\mathbb N_0^2)$ such that $H^{j+1}=span\{W_k\mid\,k\in\mathcal K^{j+2}\}$.
\end{Definition}
\begin{Remark}
In the previuos definition the condition
\begin{align*}
(A)\qquad&\overline{\bigcup_{i\in\mathbb N}D^i}=H\\
\text{is equivalent to}&\\
(B)\qquad&\forall x\in H\,[j\to\infty\;\text{only if}\;|x-\Pi_jx|\to 0]
\end{align*}
where $D^i$ stays for either $G^i$ or $H^i$. Indeed if $(B)$ is satisfied and $x\in H$, we have $\Pi_ix\in D^i$ and so, $(\Pi_ix)_{i\in\mathbb N}$ is a sequence in $\bigcup_{i\in\mathbb N}D^i$ converging to $x$. Thus $\overline{\bigcup_{i\in\mathbb N}D^i}\supseteq H$.\par
Conversely if $(A)$ is satisfied and $x\in H$, we can find $j_0\in\mathbb N$ such that $|x^0-G^{j_0}|<1$ and, for each $n\in\mathbb N_0$ we can find $j_n\in\mathbb N$ and $x^n\in D^{j_n}$ such that $j_n>j_{n-1}$ and $|x-x^n|<\frac{1}{n+1}$, then also $|x-\Pi_{j_n}x|<\frac{1}{n+1}$.\;\footnote {Because $|x-\Pi_{j_n}x|=min\{|x-y|\,\mid\,y\in G^{j_n}\}$.}\; So $(|x-\Pi_{j_n}x|)_{n\in\mathbb N}$ (that is a subsequence of $(|x-\Pi_j x|)_{j\in\mathbb N}$) converges to zero. Since the sequence $(D^i)_{i\in\mathbb N}$ is increasing (in the inclusion sense), the sequence $(|x-\Pi_j x|)_{j\in\mathbb N}$ is decreasing so, it converges to zero (as does its subsequence $(|x-\Pi_{j_n}x|)_{n\in\mathbb N}$). 
\end{Remark}

\begin{Definition}
A {\bf Galerkin approximation} of system \eqref{infsys} is the same system with the additional condition
$$k,\,n,\,m\,\in\mathcal G\in\mathcal{FP}(\mathbb N_0^2).$$
\end{Definition}
From now we shall look for the existence of a G-Saturating set $\mathcal K^1$ for system \eqref{infsys}. The existence of such a set means that for any given finite set of modes $\mathcal O\subseteq\mathbb N_0^2$, we want to observe, there exists $p\in\mathbb N$ such that $\mathcal O\subseteq H^p$.\par
For $j<p$ the construction of $H^{j+1}$ depends only on vectors of $H^j\subseteq H^p$. Therefore from the $\mathcal K^{p+1}$-Galerkin approximation with $\mathcal K^j\quad(1\leq j\leq p)$ as set of controlled modes, i.e., with $H^{j-1}$ as set of controlled directions, we can arrive by a FCE Procedure to the same Galerkin approximation with $\mathcal K^{j+1}$ as set of controlled modes.\par
Iterating FCE procedures, after $p$ steps (starting with $\mathcal K^1$ as set of controlled modes), we arrive to the $\mathcal K^{p+1}$-Galerkin approximation with $\mathcal K^{p+1}$ as set of controlled modes and, the approximate controllability at time $t$ of such system is an immediate consequence of Corollary \ref{clAfxeV.Rn}. The approximate controllability at time $t$ of the $\mathcal K^{p+1}$-Galerkin approximation with $\mathcal K^1$ as set of controlled modes follows from the fact that a FCE Procedure does not change closure of attainable set at time $t$.\\
Hence the controlling of the modes in $\mathcal K^1$ in the  $\mathcal K^{p+1}$-Galerkin approximation we ``observe'' approximate controllability in the state space $H^p=span(\mathcal K^{p+1})$ and so also in the ``$\mathcal O$-space''.

\section{Looking for Saturating Sets.}\label{looksat}
We still working in the rectangle $[0,\,a]\times[0,\,b]$. We set two elements $n,\,m\in\mathbb N_0^2$ such that $m<n$. We have that 
$$
\bar n-\bar m=0\;\iff\; \frac{m_1^2}{a^2}+\frac{m_2^2}{b^2}=\frac{n_1^2}{a^2}+\frac{n_2^2}{b^2}
$$
and, in this equation $m_1=n_1\iff m_2=n_2$. So since $m<n$ we have that $\bar n-\bar m=0$ implies $m_1\ne n_1\And m_2\ne n_2$. Then 
\begin{equation}\label{van.n-m.1}
\bar n-\bar m=0\;\iff\; \frac{a^2}{b^2}=-\frac{m_1^2-n_1^2}{m_2^2-n_2^2}=\frac{(n_1-m_1)(n_1+m_1)}{(m_2-n_2)(m_2+n_2)}.
\end{equation}
Now since $m_1<n_1$ which, from \eqref{van.n-m.1}, implies $m_2>n_2$ we obtain
\begin{equation}\label{van.n-m}
\bar n-\bar m=0\And m<n\;\Rightarrow\; n_1>m_1\And n_2<m_2.
\end{equation}

By \eqref{van.n-m} we obtain
\begin{Corollary}\label{VanishC}
Under the condition $m<n\wedge (m_1=n_1\vee n_2\geq m_2)$ we have
$$
C_{m,n}^{++}=0\iff m\wedge n=0\iff C_{m,n}^{--}=0
$$
and
$$
C_{m,n}^{-+}\ne 0\ne C_{m,n}^{+-}.
$$
\end{Corollary}
\begin{Proposition}\label{k1sat}
The set $\mathcal K^1:=\{(n_1,\,n_2)\in\mathbb N_0^2\mid n_1,\,n_2\leq 3\}\setminus\{(3,3)\}$ is G-Saturating for system \eqref{infsys}.\footnote{We take $(3,3)$ off only because a writting reason. We shall apply a induction procedure in the proof and without these ``corners'' it becomes simpler.}
\end{Proposition}
\begin{proof}
For any $j\in\mathbb N_0$ put
$$
\mathcal K^j:=\{(n_1,\,n_2)\in\mathbb N_0^2\mid n_1,\,n_2\leq j+2\}\setminus\{(j+2,j+2)\}.
$$
We shall prove that the directions in $H^{j-1}:=span(\mathcal K^j)$ can be obtained by a suitable FCE Procedure from the directions in $H^{j-2}\quad(j\geq 2)$.\par 
For every pair $(n,\,m)\in\mathcal K^1,\,n\ne m$ and for every $\lambda\in\mathbb R^+$ we define the vectors $v_{n,m}^\lambda,\;w_{m,n}^\lambda\in\mathbb R^r$ defined by
\begin{eqnarray*}
& &
\begin{cases}
(v_{m,n}^\lambda)_n=\lambda;\\
(v_{m,n}^\lambda)_m=1;\\
(v_{m,n}^\lambda)_k=0,\quad k\in\mathcal K^1\setminus\{n,\,m\};
\end{cases}\\
\text{and}& &\\
& &
\begin{cases}
(w_{m,n}^\lambda)_n=\lambda;\\
(w_{m,n}^\lambda)_m=-1;\\
(w_{m,n}^\lambda)_k=0,\quad k\in\mathcal K^1\setminus\{n,\,m\}.
\end{cases}
\end{eqnarray*}
The vector fields
\begin{eqnarray*}
& &\frac{f_{v_{m,n}^\lambda}(u)+f_{-v_{m,n}^\lambda}(u)}{2}\,=\,f(u)+\lambda\gamma_{m,n}\\
\text{and,}& &\\
& &\frac{f_{w_{m,n}^\lambda}(u)+f_{-w_{m,n}^\lambda}(u)}{2}\,=\,f(u)-\lambda\gamma_{m,n},
\end{eqnarray*}
where
\begin{eqnarray*}
\gamma_{m,n}&=&C_{m,n}^{--}\partial_{(n(--)m)^+}+C_{m,n}^{-+}\partial_{(n(-+)m)^+}\\
&+&C_{m,n}^{+-}\partial_{(n(+-)m)^+}+C_{m,n}^{++}\partial_{(n(++)m)^+},
\end{eqnarray*}
belong to $Conv\{f_{\pm v_{m,n}^\lambda},\,f_{\pm w_{m,n}^\lambda}\mid\,\lambda\in\mathbb R^+,\,m,\,n\in\mathcal K^1,\,n\ne m\}$.\par
Now we shall extract from $\{f_{X(t)}\mid\,X(t)\in\mathbb R^r\}$ a family as follows:\\
First we define for each pair $(m,\,n)\in\mathcal K^1,\,n\ne m$ the vector $\delta_{m,n}$ in $\mathbb R^N$ by:
\begin{eqnarray}
\delta_{m,n}&=&C_{m,n}^{--}e_{(n(--)m)^+}+C_{m,n}^{-+}e_{(n(-+)m)^+}\notag\\
&+&C_{m,n}^{+-}e_{(n(+-)m)^+}+C_{m,n}^{++}e_{(n(++)m)^+}\label{def.deltamn}
\end{eqnarray}
Now we select the subfamily $F_{S_1}:=\{\delta_{m,n}\mid\,(m,\,n)\in S_1\subseteq\mathcal (K^1)^2\}$ of this vectors, where
\begin{align}
S_1=\{&((1,2),(2,1));((1,1),(2,3));((1,2),(2,2));\notag\\
&((1,1),(3,2));((2,1),(2,2));((1,1),(1,3));((1,1),(3,1))\}.\label{S1.extr}
\end{align}
The vectors of this family are precisely
\begin{align*}
\delta_{(1,2),(2,1)}&=\frac{9\pi^2(b^2-a^2)}{4ab(a^2+b^2)}e_{(1,1)}+\frac{15\pi^2(a^2-b^2)}{4ab(9a^2+b^2)}e_{(1,3)}\\
&+\frac{15\pi^2(a^2-b^2)}{4ab(a^2+9b^2)}e_{(3,1)}+\frac{\pi^2(b^2-a^2)}{4ab(a^2+b^2)}e_{(3,3)}\\
\delta_{(1,1),(2,3)}&=\frac{\pi^2(3b^2+8a^2)}{4ab(4a^2+b^2)}e_{(1,2)}-\frac{5\pi^2(8a^2+3b^2)}{4ab(16a^2+b^2)}e_{(1,4)}\\
&+\frac{5\pi^2(8a^2+3b^2)}{4ab(4a^2+9b^2)}e_{(3,2)}-\frac{\pi^2(3b^2+8a^2)}{4ab(16a^2+9b^2)}e_{(3,4)}\\
\delta_{(1,2),(2,2)}&=-\frac{9b\pi^2}{2a(16a^2+b^2)}e_{(1,4)}+\frac{3b\pi^2}{2a(16a^2+9b^2)}e_{(3,4)}\\
\delta_{(1,1),(3,2)}&=-\frac{\pi^2(8b^2+3a^2)}{4ab(a^2+4b^2)}e_{(2,1)}-\frac{5\pi^2(3a^2+8b^2)}{4ab(9a^2+4b^2)}e_{(2,3)}\\
&+\frac{5\pi^2(3a^2+8b^2)}{4ab(a^2+16b^2)}e_{(4,1)}+\frac{\pi^2(8b^2+3a^2)}{4ab(9a^2+16b^2)}e_{(4,3)}\\
\delta_{(2,1),(2,2)}&=-\frac{9a\pi^2}{2b(16b^2+a^2)}e_{(4,1)}-\frac{3a\pi^2}{2b(16b^2+9a^2)}e_{(4,3)}\\
\delta_{(1,1),(1,3)}&=\frac{2a\pi^2}{b(b^2+a^2)}e_{(2,2)}-\frac{a\pi^2}{b(b^2+4a^2)}e_{(2,4)}\\
\delta_{(1,1),(3,1)}&=-\frac{2b\pi^2}{a(b^2+a^2)}e_{(2,2)}+\frac{b\pi^2}{a(a^2+4b^2)}e_{(4,2)},
\end{align*}
Projecting the vectors in this subfamily on the space $span\{e_k\mid\,k\in\mathcal K^2\setminus\mathcal K^1\}$ we obtain
\begin{align}
\Pi_1\delta_{(1,2),(2,1)}&=\frac{\pi^2(b^2-a^2)}{4ab(a^2+b^2)}e_{(3,3)}\label{notSq}\\
\Pi_1\delta_{(1,1),(2,3)}&=-\frac{5\pi^2(8a^2+3b^2)}{4ab(16a^2+b^2)}e_{(1,4)}-\frac{\pi^2(3b^2+8a^2)}{4ab(16a^2+9b^2)}e_{(3,4)}\notag\\
\Pi_1\delta_{(1,2),(2,2)}&=-\frac{9b\pi^2}{2a(16a^2+b^2)}e_{(1,4)}+\frac{3b\pi^2}{2a(16a^2+9b^2)}e_{(3,4)}\notag\\
\Pi_1\delta_{(1,1),(3,2)}&=\frac{5\pi^2(3a^2+8b^2)}{4ab(a^2+16b^2)}e_{(4,1)}+\frac{\pi^2(8b^2+3a^2)}{4ab(9a^2+16b^2)}e_{(4,3)}\notag\\
\Pi_1\delta_{(2,1),(2,2)}&=-\frac{9a\pi^2}{2b(16b^2+a^2)}e_{(4,1)}-\frac{3a\pi^2}{2b(16b^2+9a^2)}e_{(4,3)}\notag\\
\Pi_1\delta_{(1,1),(1,3)}&=-\frac{a\pi^2}{b(b^2+4a^2)}e_{(2,4)}\notag\\
\Pi_1\delta_{(1,1),(3,1)}&=\frac{b\pi^2}{a(a^2+4b^2)}e_{(4,2)}.\notag
\end{align}
Now we consider the case $a\ne b$, i. e., if our rectangle $R$ is not a square (the case of the square will be considered below in Remark \ref{CaseSq}). If $a\ne b$ these projections are linearly independent so, by Lemma \ref{ind.pind} the $15$ vectors of the family $\{e_k\mid\,k\in\mathcal K^1\}\cup F_{S_1}$ are linearly independent in $span(\mathcal K^2)=H^1$ and then they span $span\{e_k\mid\,k\in\mathcal K^2\}\subset H$.\\
Now we extract the linear space
\begin{multline*}
\tilde H^1:=\Big(H^0+Conv\{\lambda\delta_{m,n}\mid\,(m,\,n)\in S_1,\,\lambda\in\mathbb R\}\Big)\\
\cap \Big(H^0-Conv\{\lambda\delta_{m,n}\mid\,(m,\,n)\in S_1,\,\lambda\in\mathbb R\}\Big)
\end{multline*}
from the space
$$
\Big(H^0+Conv\{B(Y,\,Y)\mid\,Y\in H^0\}\Big)\cap \Big(H^0-Conv\{B(Y,\,Y)\mid\,Y\in H^0\}\Big).
$$
Since $Conv\{\lambda\delta_{m,n}\mid\,(m,\,n)\in S_1,\,\lambda\in\mathbb R\}\Big)$ coincides with $span\{\lambda\delta_{m,n}\mid\,(m,\,n)\in S_1,\,\lambda\in\mathbb R\}\Big)$, we have that
$$
\tilde H^1=span(\mathcal K^1\cup\{\delta_{m,n}\mid\,(m,\,n)\in S_1\})=span(\mathcal K^2)=H^1.
$$

The following proposition completes the proof.
\end{proof}
\begin{Proposition}\label{Kj.to.Kj+1}
From the directions in $H^j=span(\mathcal K^{j+1}),\quad (j\geq 1)$\footnote{The case $j=0$ was already seen in the proof of proposition \ref{k1sat}.} we can obtain the directions in $H^{j+1}=span(\mathcal K^{j+2})$ by a FCE Procedure.
\end{Proposition}
\begin{proof}
We consider two cases ``$j$ even'' and ``$j$ odd''.\par
$\bullet$ $j$ even: In this case
$$
\mathcal K^{j+1}:=\{(n_1,\,n_2)\in\mathbb N_0^2\mid n_1,\,n_2\leq j+3\}\setminus\{(j+3,j+3)\}.
$$
can be written as
$$
\mathcal K^{j+1}:=\{(n_1,\,n_2)\in\mathbb N_0^2\mid n_1,\,n_2\leq 2p+1\}\setminus\{(2p+1,2p+1)\},
$$
setting $p=\frac{j+2}{2}$. Then $p\geq 2$.\\
As we did before in the case ``$j=0\to j=1$'', we extract a subfamily $F_{S_{j+1}}:=\{\delta_{m,n}\mid\,(m,\,n)\in S_{j+1}\subseteq\mathcal (K^{j+1})^2\}$ of $\{B(Y,\,Y)\mid\,Y\in H^j\}$ where now the ``selection'' is
\begin{align*}
S_{j+1}&=\{((1,2),(2p,2p-1))\}\\
&\cup\{((1,1),(2z,2p+1))\mid\,z=1,\,\dots,\,p\}\cup\{((1,p+1)(2,p+1))\}\\
&\cup\{((1,1),(2p+1,2z))\mid\,z=1,\,\dots,\,p\}\cup\{((p+1,1),(p+1,2))\}\\
&\cup\{((s,1),(s,2p+1))\mid\,s=1,\,\dots,\,p\}\\
&\cup\{((1,s),(2p+1,s))\mid\,s=1,\,\dots,\,p\}.
\end{align*}
If we write explicitely the vectors of $F_{S_{j+1}}$ we obtain quite long expressions, for example we have that $C_{(1,2),(2p,2p-1)}^{--}$ equals
$$
-\frac{\bigl(a^2(-3+2p)+b^2(-1+2p)\bigr)(\pi+2p\pi)^2\mathrm{sign}(-3+2p)\mathrm{sign}(-1+2p)}{4ab(b^2|1-2p|^2+a^2|3-2p|^2)};
$$
so, here we will not write those vectores explicitely. Anyway, those vectors are
\begin{align*}
\delta_{(1,2),(2p,2p-1)}=&C_{(1,2),(2p,2p-1)}^{--}e_{(2p-1,2p-3)}+C_{(1,2),(2p,2p-1)}^{-+}e_{(2p-1,2p+1)}\\
+&C_{(1,2),(2p,2p-1)}^{+-}e_{(2p+1,2p-3)}+C_{(1,2),(2p,2p-1)}^{++}e_{(2p+1,2p+1)};\\
\delta_{(1,1),(2z,2p+1)}=&C_{(1,1),(2z,2p+1)}^{--}e_{(2z-1,2p)}+C_{(1,1),(2z,2p+1)}^{-+}e_{(2z-1,2(p+1))}\\
+&C_{(1,1),(2z,2p+1)}^{+-}e_{(2z+1,2p)}+C_{(1,1),(2z,2p+1)}^{++}e_{(2z+1,2(p+1))}\\&\qquad\qquad\qquad\qquad z=1,\,\dots,\,p;\\
\delta_{(1,p+1),(2,p+1)}=&C_{(1,p+1),(2,p+1)}^{-+}e_{(1,2(p+1))}+C_{(1,p+1),(2,p+1)}^{++}e_{(3,2(p+1))};\\
\delta_{(1,1),(2p+1,2z)}=&C_{(1,1),(2p+1,2z)}^{--}e_{(2p,2z-1)}+C_{(1,1),(2p+1,2z)}^{-+}e_{(2p,2z+1)}\\
+&C_{(1,1),(2p+1,2z)}^{+-}e_{(2(p+1),2z-1)}+C_{(1,1),(2p+1,2z)}^{++}e_{(2(p+1),2z+1)}\\&\qquad\qquad\qquad\qquad z=1,\,\dots,\,p;\\
\delta_{(p+1,1),(p+1,2)}=&C_{(p+1,1),(p+1,2)}^{+-}e_{(2(p+1),1)}+C_{(p+1,1),(p+1,2)}^{++}e_{(2(p+1),3)};\\
\delta_{(s,1),(s,2p+1)}=&C_{(s,1),(s,2p+1)}^{+-}e_{(2s,2p)}+C_{(s,1),(s,2p+1)}^{++}e_{(2s,2(p+1))}\\&\qquad\qquad\qquad\qquad s=1,\,\dots,\,p;\\
\delta_{(1,s),(2p+1,s)}=&C_{(1,s),(2p+1,s)}^{-+}e_{(2p,2s)}+C_{(1,s),(2p+1,s)}^{++}e_{(2(p+1),2s)}\\&\qquad\qquad\qquad\qquad s=1,\,\dots,\,p;
\end{align*}
And, projecting them onto the space $span\{e_k\mid\,k\in\mathcal K^{j+2}\setminus\mathcal K^{j+1}\}$ we arrive to the family $\Pi_{j+1} F_{S_{j+1}}$ whose elements are
\begin{align}
\Pi_{j+1}\delta_{(1,2),(2p,2p-1)}=&C_{(1,2),(2p,2p-1)}^{++}e_{(2p+1,2p+1)};\notag\\
\Pi_{j+1}\delta_{(1,1),(2z,2p+1)}=&C_{(1,1),(2z,2p+1)}^{-+}e_{(2z-1,2(p+1))}+C_{(1,1),(2z,2p+1)}^{++}e_{(2z+1,2(p+1))}\notag\\&\qquad\qquad\qquad\qquad z=1,\,\dots,\,p;\label{v13.oE.z=1}\\
\Pi_{j+1}\delta_{(1,p+1),(2,p+1)}=&C_{(1,p+1),(2,p+1)}^{-+}e_{(1,2(p+1))}+C_{(1,p+1),(2,p+1)}^{++}e_{(3,2(p+1))};\\
\Pi_{j+1}\delta_{(1,1),(2p+1,2z)}=&C_{(1,1),(2p+1,2z)}^{+-}e_{(2(p+1),2z-1)}+C_{(1,1),(2p+1,2z)}^{++}e_{(2(p+1),2z+1)}\notag\\&\qquad\qquad\qquad\qquad z=1,\,\dots,\,p;\label{v13.Eo.z=1}\\
\Pi_{j+1}\delta_{(p+1,1),(p+1,2)}=&C_{(p+1,1),(p+1,2)}^{+-}e_{(2(p+1),1)}+C_{(p+1,1),(p+1,2)}^{++}e_{(2(p+1),3)};\\
\Pi_{j+1}\delta_{(s,1),(s,2p+1)}=&C_{(s,1),(s,2p+1)}^{++}e_{(2s,2(p+1))},\qquad s=1,\,\dots,\,p;\\
\Pi_{j+1}\delta_{(1,s),(2p+1,s)}=&C_{(1,s),(2p+1,s)}^{++}e_{(2(p+1),2s)},\qquad s=1,\,\dots,\,p.
\end{align}
No one of the coeficients appearing in these expressions vanishes because all pairs $(m,n)$ satisfy $m<n\wedge (m_1=n_1\vee n_2\geq m_2)$ and because no one of the following expressions vanish
\begin{align*}
(1,2)\wedge(2p,2p-1)&=-1-2p;\\
(1,1)\wedge(2z,2p+1)&=1+2p-2z,\quad z=1,\,\dots,\,p;\\
(1,p+1)\wedge(2,p+1)&=-(1+p);\\
(1,1)\wedge(2p+1,2z)&=-1-2p+2z,\quad z=1,\,\dots,\,p;\\
(p+1,1)\wedge(p+1,2)&=1+p;\\
(s,1)\wedge(s,2p+1)&=2ps,\quad s=1,\,\dots,\,p;\\
(1,s)\wedge(2p+1,s)&=-2ps,\quad s=1,\,\dots,\,p.
\end{align*}
Hence we can see that these vectors are linearly independent. Indeed it suffices to prove that:
\begin{itemize}
\item The vectors $\Pi_{j+1}\delta_{(1,1),(2,2p+1)},\quad (z=1\,\text{in \eqref{v13.oE.z=1}})$ and $\Pi_{j+1}\delta_{(1,p+1),(2,p+1)}$ are linearly independent; and\\
\item The vectors $\Pi_{j+1}\delta_{(1,1),(2p+1,2)},\quad (z=1\,\text{in \eqref{v13.Eo.z=1}})$ and $\Pi_{j+1}\delta_{(p+1,1),(p+1,2)}$ are linearly independent;
\end{itemize}
But that comes from  
\begin{align*}
&\frac{C_{(1,1),(2,2p+1)}^{-+}}{C_{(1,p+1),(2,p+1)}^{-+}}=\frac{(3+2p)(3b^2+4a^2p(1+p))}{9b^2(1+p)}\\
\ne&-\frac{(-1+2p)(3b^2+4a^2p(1+p))}{3b^2(1+p)}=\frac{C_{(1,1),(2,2p+1)}^{++}}{C_{(1,p+1),(2,p+1)}^{++}};\\
\text{and}&\\
&\frac{C_{(1,1),(2p+1,2z)}^{+-}}{C_{(p+1,1),(p+1,2)}^{+-}}=\frac{(3+2p)(3a^2+4b^2p(1+p))}{9a^2(1+p)}\\
\ne&-\frac{(-1+2p)(3a^2+4b^2p(1+p))}{3a^2(1+p)}=\frac{C_{(1,1),(2,2p+1)}^{++}}{C_{(1,p+1),(2,p+1)}^{++}}.
\end{align*}
By Lemma \ref{ind.pind} the $(2(p+1))^2-1$ of $F_{S_{j+1}}\cup\{e_k\mid\,k\in\mathcal K^{j+1}\}$ are linearly independent and then they span
\begin{eqnarray*}
& &span\{e_k\mid\,1\leq k_1,k_2\leq2(p+1)\}\setminus\{(2(p+1),2(p+1))\}\\
&=&span\{e_k\mid\,k\in\mathcal K^{j+2}\}.
\end{eqnarray*}
Since with $\delta_{m,n}$ also $\lambda\delta_{m,n}$ belongs to $\{B(Y,\,Y)\mid\,Y\in H^j\}$ for every $\lambda\in\mathbb R$, we can select  the linear space $\tilde H^{j+1}$ defined by
\begin{multline*}
\tilde H^{j+1}:=\Big(H^j+Conv\{\lambda\delta_{m,n}\mid\,(m,\,n)\in S_{j+1},\,\lambda\in\mathbb R\}\Big)\\
\cap \Big(H^j-Conv\{\lambda\delta_{m,n}\mid\,(m,\,n)\in S_{j+1},\,\lambda\in\mathbb R\}\Big)
\end{multline*}
from the space
$$
\Big(H^j+Conv\{B(Y,\,Y)\mid\,Y\in H^j\}\Big)\cap \Big(H^j-Conv\{B(Y,\,Y)\mid\,Y\in H^j\}\Big).
$$
Since $Conv\{\lambda\delta_{m,n}\mid\,(m,\,n)\in S_{j+1},\,\lambda\in\mathbb R\}\Big)$ coincides with $span\{\lambda\delta_{m,n}\mid\,(m,\,n)\in S_{j+1},\,\lambda\in\mathbb R\}\Big)$, we have that
$$
\tilde H^{j+1}=span(\mathcal K^{j+1}\cup\{\delta_{m,n}\mid\,(m,\,n)\in S_{j+1}\})=span(\mathcal K^{j+2})=H^{j+1}.
$$
\ \\
Now we study the other case\\
$\bullet$ $j$ odd: In this case
$$
\mathcal K^{j+1}:=\{(n_1,\,n_2)\in\mathbb N_0^2\mid n_1,\,n_2\leq j+3\}\setminus\{(j+3,j+3)\}.
$$
can be written as
$$
\mathcal K^{j+1}:=\{(n_1,\,n_2)\in\mathbb N_0^2\mid n_1,\,n_2\leq 2p\}\setminus\{(2p,2p)\},
$$
setting $p=\frac{j+3}{2}$. Then $p\geq 2$.\\
We extract the subfamily $F_{S_{j+1}}:=\{\delta_{m,n}\mid\,(m,\,n)\in S_{j+1}\subseteq\mathcal (K^{j+1})^2\}$ of  $\{B(Y,\,Y)\mid\,Y\in H^j\}$ where now the ``selection'' is

\begin{align*}
S_{j+1}&=\{((1,2p-1),(2p-1,1))\}\\
&\cup\{((1,1),(2z-1,2p))\mid\,z=2,\,\dots,\,p\}\cup\{((1,p)(3,p+1))\}\\
&\cup\{((1,1),(2p,2z-1))\mid\,z=2,\,\dots,\,p\}\cup\{((p,1),(p+1,3))\}\\
&\cup\{((1,1),(2s,2p))\mid\,s=1,\,\dots,\,p-1\}\cup\{((1,p),(2,p+1))\}\\
&\cup\{((1,1),(2p,2s))\mid\,s=1,\,\dots,\,p-1\}\cup\{((p,1),(p+1,2))\}.
\end{align*}
Those vectors are
\begin{align*}
\delta_{(1,2p-1),(2p-1,1)}=&C_{(1,2p-1),(2p-1,1)}^{--}e_{(2(p-1),2(p-1))}\\
+&C_{(1,2p-1),(2p-1,1)}^{-+}e_{(2(p-1),2p)}\\
+&C_{(1,2p-1),(2p-1,1)}^{+-}e_{(2p,2(p-1))} +C_{(1,2p-1),(2p-1,1)}^{++}e_{(2p,2p)}\\
\delta_{(1,1),(2z-1,2p)}=&C_{(1,1),(2z-1,2p)}^{--}e_{(2(z-1),2p-1)} +C_{(1,1),(2z-1,2p)}^{-+}e_{(2(z-1),2p+1)}\\
+&C_{(1,1),(2z-1,2p)}^{+-}e_{(2z,2p-1)} +C_{(1,1),(2z-1,2p)}^{++}e_{(2z,2p+1)}\\ &\qquad\qquad\qquad z=2,\,\dots,\,p\\
\delta_{(1,p)(3,p+1)}=&C_{(1,p)(3,p+1)}^{--}e_{(2,1)} +C_{(1,p)(3,p+1)}^{-+}e_{(2,2p+1)}\\
+&C_{(1,p)(3,p+1)}^{+-}e_{(4,1)}+C_{(1,p)(3,p+1)}^{++}e_{(4,2p+1)}\\
\delta_{(1,1),(2p,2z-1)}=&C_{(1,1),(2p,2z-1)}^{--}e_{(2p-1,2(z-1))} +C_{(1,1),(2p,2z-1)}^{-+}e_{(2p-1,2z)}\\
+&C_{(1,1),(2p,2z-1)}^{+-}e_{(2p+1,2(z-1))} +C_{(1,1),(2p,2z-1)}^{++}e_{(2p+1,2z)}\\
&\qquad\qquad\qquad z=2,\,\dots,\,p
\end{align*}
\begin{align*}
\delta_{(p,1),(p+1,3)}=&C_{(p,1),(p+1,3)}^{--}e_{(1,2)} +C_{(p,1),(p+1,3)}^{-+}e_{(1,4)}\\
+&C_{(p,1),(p+1,3)}^{+-}e_{(2p+1,2)} +C_{(p,1),(p+1,3)}^{++}e_{(2p+1,4)}\\
\delta_{(1,1),(2s,2p)}=&C_{(1,1),(2s,2p)}^{--}e_{(2s-1,2p-1)} +C_{(1,1),(2s,2p)}^{-+}e_{(2s-1,2p+1)}\\
+&C_{(1,1),(2s,2p)}^{+-}e_{(2s+1,2p-1)} +C_{(1,1),(2s,2p)}^{++}e_{(2s+1,2p+1)}\\
&\qquad\qquad\qquad s=1,\,\dots,\,p-1\\
\delta_{(1,p),(2,p+1)}=&C_{(1,p),(2,p+1)}^{--}e_{(1,1)} +C_{(1,p),(2,p+1)}^{-+}e_{(1,2p+1)}\\
+&C_{(1,p),(2,p+1)}^{+-}e_{(3,1)} +C_{(1,p),(2,p+1)}^{++}e_{(3,2p+1)}\\
\delta_{(1,1),(2p,2s)}=&C_{(1,1),(2p,2s)}^{--}e_{(2p-1,2s-1)} +C_{(1,1),(2p,2s)}^{-+}e_{(2p-1,2s+1)}\\
+&C_{(1,1),(2p,2s)}^{+-}e_{(2p+1,2s-1)} +C_{(1,1),(2p,2s)}^{++}e_{(2p+1,2s+1)}\\
&\qquad\qquad\qquad s=1,\,\dots,\,p-1\\
\delta_{(p,1),(p+1,2)}=&C_{(p,1),(p+1,2)}^{--}e_{(1,1)} +C_{(p,1),(p+1,2)}^{-+}e_{(1,3)}
\\+&C_{(p,1),(p+1,2)}^{+-}e_{(2p+1,1)} +C_{(p,1),(p+1,2)}^{++}e_{(2p+1,3)};
\end{align*}
And, projecting them onto the space $span\{e_k\mid\,k\in\mathcal K^{j+2}\setminus\mathcal K^{j+1}\}$ we arrive to the family $\Pi_{j+1} F_{S_{j+1}}$ which elements are
\begin{align}
\Pi_{j+1}\delta_{(1,2p-1),(2p-1,1)}=&C_{(1,2p-1),(2p-1,1)}^{++}e_{(2p,2p)}\notag\\
\Pi_{j+1}\delta_{(1,1),(2z-1,2p)}=&C_{(1,1),(2z-1,2p)}^{-+}e_{(2(z-1),2p+1)}+C_{(1,1),(2z-1,2p)}^{++}e_{(2z,2p+1)}\notag\\
&\qquad\qquad\qquad z=2,\,\dots,\,p\label{v24.eO.z=2}\\
\Pi_{j+1}\delta_{(1,p)(3,p+1)}=&C_{(1,p)(3,p+1)}^{-+}e_{(2,2p+1)}+C_{(1,p)(3,p+1)}^{++}e_{(4,2p+1)}\notag\\
\Pi_{j+1}\delta_{(1,1),(2p,2z-1)}=&C_{(1,1),(2p,2z-1)}^{+-}e_{(2p+1,2(z-1))} +C_{(1,1),(2p,2z-1)}^{++}e_{(2p+1,2z)}\notag\\
&\qquad\qquad\qquad z=2,\,\dots,\,p\label{v24.Oe.z=2}\\
\Pi_{j+1}\delta_{(p,1),(p+1,3)}=&C_{(p,1),(p+1,3)}^{+-}e_{(2p+1,2)} +C_{(p,1),(p+1,3)}^{++}e_{(2p+1,4)}\notag\\
\Pi_{j+1}\delta_{(1,1),(2s,2p)}=&C_{(1,1),(2s,2p)}^{-+}e_{(2s-1,2p+1)}+C_{(1,1),(2s,2p)}^{++}e_{(2s+1,2p+1)}\notag\\
&\qquad\qquad\qquad s=1,\,\dots,\,p-1\label{v13.oO.s=1}\\
\Pi_{j+1}\delta_{(1,p),(2,p+1)}=&C_{(1,p),(2,p+1)}^{-+}e_{(1,2p+1)}+C_{(1,p),(2,p+1)}^{++}e_{(3,2p+1)}\notag\\
\Pi_{j+1}\delta_{(1,1),(2p,2s)}=&C_{(1,1),(2p,2s)}^{+-}e_{(2p+1,2s-1)} +C_{(1,1),(2p,2s)}^{++}e_{(2p+1,2s+1)}\notag\\
&\qquad\qquad\qquad s=1,\,\dots,\,p-1\label{v13.Oo.s=1}\\
\Pi_{j+1}\delta_{(p,1),(p+1,2)}=&C_{(p,1),(p+1,2)}^{+-}e_{(2p+1,1)} +C_{(p,1),(p+1,2)}^{++}e_{(2p+1,3)};
\end{align}
No one of the coeficients appearing in these expressions vanishes because all pairs $(m,n)$ satisfy $m<n\wedge (m_1=n_1\vee n_2\geq m_2)$ and because no one of the following expressions vanish
\begin{align*}
(1,2p-1)\wedge(2p-1,1)&=1-(2p-1)^2;\\
(1,1)\wedge(2z-1,2p)&=1+2(p-z),\quad z=2,\,\dots,\,p;\\
(1,p)\wedge(3,p+1)&=1-2p;\\
(1,1)\wedge(2p,2z-1)&=-1-2(p-z),\quad z=2,\,\dots,\,p;\\
(p,1)\wedge(p+1,3)&=2p-1;\\
(1,1)\wedge(2s,2p)&=2(p-s),\quad s=1,\,\dots,\,p-1;\\
(1,p)\wedge(2,p+1)&=1-p;\\
(1,1)\wedge(2p,2s)&=2(s-p),\quad s=1,\,\dots,\,p-1\\
(p,1)\wedge(p+1,2)&=p-1.
\end{align*}
Hence we can see that these vectors are linearly independent. To see this is enough to see that:
\begin{itemize}
\item The vectors $\Pi_{j+1}\delta_{(1,1),(3,2p)},\quad (z=2\,\text{in \eqref{v24.eO.z=2}})$ and $\Pi_{j+1}\delta_{(1,p)(3,p+1)}$ are linearly independent;
\item The vectors $\Pi_{j+1}\delta_{(1,1),(2p,3)},\quad (z=2\,\text{in \eqref{v24.Oe.z=2}})$ and $\Pi_{j+1}\delta_{(p,1),(p+1,3)}$ are linearly independent;
\item The vectors $\Pi_{j+1}\delta_{(1,1),(2,2p)},\quad (s=1\,\text{in \eqref{v13.oO.s=1}})$ and $\Pi_{j+1}\delta_{(1,p)(2,p+1)}$ are linearly independent;
\item The vectors $\Pi_{j+1}\delta_{(1,1),(2p,2)},\quad (s=1\,\text{in \eqref{v13.Oo.s=1}})$ and $\Pi_{j+1}\delta_{(p,1),(p+1,2)}$ are linearly independent;
\end{itemize}
But that comes from  
\begin{align*}
&\frac{C_{(1,1),(3,2p)}^{-+}}{C_{(1,p)(3,p+1)}^{-+}}=\frac{(3+2p)(8b^2+a^2(4p^2-1))}{(1+4p)(8b^2+a^2(2p+1))}\\
\ne&-\frac{(2p-3)(8b^2+a^2(4p^2-1))}{(2p-1)(8b^2+a^2(2p+1))}=\frac{C_{(1,1),(3,2p)}^{++}}{C_{(1,p)(3,p+1)}^{++}};\\
&\frac{C_{(1,1),(2p,3)}^{+-}}{C_{(p,1),(p+1,3)}^{+-}}=\frac{(3+2p)(8a^2+b^2(4p^2-1))}{(1+4p)(8a^2+b^2(2p+1))}\\
\ne&-\frac{(2p-3)(8a^2+b^2(4p^2-1))}{(2p-1)(8a^2+b^2(2p+1))}=\frac{C_{(1,1),(2p,3)}^{++}}{C_{(p,1)(p+1,3)}^{++}};
\end{align*}
\begin{align*}
&\frac{C_{(1,1),(2,2p)}^{-+}}{C_{(1,p)(2,p+1)}^{-+}}=\frac{2(1+p)(3b^2+a^2(4p^2-1))}{(3p+1)(3b^2+a^2(2p+1))}\\
\ne&-\frac{2(3b^2+a^2(4p^2-1))}{3b^2+a^2(2p+1)}=\frac{C_{(1,1),(2,2p)}^{++}}{C_{(1,p)(2,p+1)}^{++}};\\
\text{and}&\\
&\frac{C_{(1,1),(2p,2)}^{+-}}{C_{(p,1),(p+1,2)}^{+-}}=\frac{2(1+p)(3a^2+b^2(4p^2-1))}{(3p+1)(3a^2+b^2(2p+1))}\\
\ne&-\frac{2(3a^2+b^2(4p^2-1))}{3b^2+a^2(2p+1)}=\frac{C_{(1,1),(2p,2)}^{++}}{C_{(p,1)(p+1,2)}^{++}}.
\end{align*}
By Lemma \ref{ind.pind} the $(2p+1)^2-1$ of $F_{S_{j+1}}\cup\{e_k\mid\,k\in\mathcal K^{j+1}\}$ are linearly independent and then they span
$$
span\{e_k\mid\,1\leq k_1,k_2\leq 2p+1\}\setminus\{(2p+1,2p+1)\}=span\{e_k\mid\,k\in\mathcal K^{j+2}\}.
$$
Again, since with $\delta_{m,n}$ also $\lambda\delta_{m,n}$ belongs to $\{B(Y,\,Y)\mid\,Y\in H^j\}$ for every $\lambda\in\mathbb R$, we can select  the linear space $\tilde H^{j+1}$ defined by
\begin{multline*}
\tilde H^{j+1}:=\Big(H^j+Conv\{\lambda\delta_{m,n}\mid\,(m,\,n)\in S_{j+1},\,\lambda\in\mathbb R\}\Big)\\
\cap \Big(H^j-Conv\{\lambda\delta_{m,n}\mid\,(m,\,n)\in S_{j+1},\,\lambda\in\mathbb R\}\Big)
\end{multline*}
from the space
$$
\Big(H^j+Conv\{B(Y,\,Y)\mid\,Y\in H^j\}\Big)\cap \Big(H^j-Conv\{B(Y,\,Y)\mid\,Y\in H^j\}\Big).
$$
Since $Conv\{\lambda\delta_{m,n}\mid\,(m,\,n)\in S_{j+1},\,\lambda\in\mathbb R\}\Big)$ coincides with $span\{\lambda\delta_{m,n}\mid\,(m,\,n)\in S_{j+1},\,\lambda\in\mathbb R\}\Big)$, we have that
$$
\tilde H^{j+1}=span(\mathcal K^{j+1}\cup\{\delta_{m,n}\mid\,(m,\,n)\in S_{j+1}\})=span(\mathcal K^{j+2})=H^{j+1}.
$$
\end{proof}
\begin{Remark}
Starting with $span(\mathcal K^2)$, by Proposition \ref{Kj.to.Kj+1}, we can obtain the directions on $span(\mathcal K^j),\quad(j\geq2)$ iterating FCE procedures. Since that Proposition is valid for any rectangle (including the square) we can say that $\mathcal K^2$ is G-Saturating for all rectangles. 
\end{Remark}

\begin{Remark}\label{CaseSq}
We can show that $\mathcal K^1$ is saturating for the square too and so complete the proof of Proposition \ref{Kj.to.Kj+1}. The only step that is not clear in the square, is how to arrive to $span(\mathcal K^2)$ (or to a superspace of its) from $span(\mathcal K^1)$.  In the case of the square the extracted family in the end of proof of Proposition \ref{k1sat} (Case $j\, even$: $j=0$ in our present case) does not lead to a family of linearly independent vectors when projected in $span\{e_k\mid\,k\notin\mathcal K^1\}$ because $\Pi\delta_{(1,2)2,1)}=0$. Then we proceed as follows: First we select the family $F_{S_1}\setminus\{\delta_{(1,2),(2,1)}\}$. Projecting this family onto $span(\mathcal K^2\setminus(\{(3,3)\}\cup\mathcal K^1))$ we obtain a family of linearly independent vectors. So ``adding'' the vectors in $\{e_k\mid\,k\in\mathcal K^1\}$ to those $F_{S_1}$ we obtain a family of $4^2-1-1$ linearly independent vectors spanning the space $span\{e_k\mid\,k\in\mathcal K^2\setminus\{(3,3)\}\,\}$. So proceeding as before we reach the directions $span(K^2\setminus\{(3,3)\})$. Then we select the family
\begin{align*}
\delta_{(1,1),(2,4)}=&\frac{3\pi^2(5a^2+b^2)}{2ab(9a^2+b^2)}e_{(1,3)} -\frac{9\pi^2(5a^2+b^2)}{2ab(25a^2-b^2)}e_{(1,5)}\\
+&\frac{\pi^2(5a^2+b^2)}{2ab(a^2+b^2)}e_{(3,3)}-\frac{3\pi^2(5a^2+b^2)}{2ab(25a^2+9b^2)}e_{(3,5)}\\
\delta_{(1,2),(2,3)}=&-\frac{\pi^2(5a^2+3b^2)}{4ab(a^2+b^2)}e_{(1,1)}-\frac{7\pi^2(5a^2+3b^2)}{4ab(25a^2+b^2)}e_{(1,5)}\\
+&\frac{7\pi^2(5a^2+3b^2)}{4ab(a^2+9b^2)}e_{(3,1)}+\frac{\pi^2(5a^2+3b^2)}{4ab(25a^2+9b^2)}e_{(3,5)}\\
\delta_{(1,4),(2,1)}=&\frac{21\pi^2(b^2-5a^2)}{4ab(9a^2+b^2)}e_{(1,3)}+\frac{27\pi^2(5a^2-b^2)}{4ab(25a^2+b^2)}e_{(1,5)}\\
+&\frac{3\pi^2(5a^2-b^2)}{4ab(a^2+b^2)}e_{(3,3)}+\frac{21\pi^2(b^2-5a^2)}{4ab(25a^2+9b^2)}e_{(3,5)}.
\end{align*}
Projecting in $span\{e_k\mid\,k\in\mathcal K^3\setminus(\mathcal K^2\setminus\{(3,3)\})\,\}$ we obtain
\begin{align*}
\Pi\delta_{(1,1),(2,4)}=&-\frac{9\pi^2(5a^2+b^2)}{2ab(25a^2-b^2)}e_{(1,5)}\\
+&\frac{\pi^2(5a^2+b^2)}{2ab(a^2+b^2)}e_{(3,3)}-\frac{3\pi^2(5a^2+b^2)}{2ab(25a^2+9b^2)}e_{(3,5)}\\
\Pi\delta_{(1,2),(2,3)}=&-\frac{7\pi^2(5a^2+3b^2)}{4ab(25a^2+b^2)}e_{(1,5)}+\frac{\pi^2(5a^2+3b^2)}{4ab(25a^2+9b^2)}e_{(3,5)}\\
\Pi\delta_{(1,4),(2,1)}=&\frac{27\pi^2(5a^2-b^2)}{4ab(25a^2+b^2)}e_{(1,5)}\\+&\frac{3\pi^2(5a^2-b^2)}{4ab(a^2+b^2)}e_{(3,3)}+\frac{21\pi^2(b^2-5a^2)}{4ab(25a^2+9b^2)}e_{(3,5)}.
\end{align*}
Since we are working in the square $a=b$ and no one of the coeficients appearing in the last expressions vanish. We compute
\begin{align*}
\mathrm{Det}
&\begin{pmatrix}
C_{(1,1),(2,4)}^{-+}&& C_{(1,1),(2,4)}^{+-}&& C_{(1,1),(2,4)}^{++}\\
C_{(1,2),(2,3)}^{-+}&& 0 &&C_{(1,2),(2,3)}^{++}\\
C_{(1,4),(2,1)}^{-+}&& C_{(1,4),(2,1)}^{+-}&& C_{(1,4),(2,1)}^{++}
\end{pmatrix}\\
=&-\frac{15\pi^2(125a^6+75a^4b^2-5a^2b^4-32b^6)}
{16a^3b^3(a^2+b^2)(25a^2+b^2)(25a^2+9b^2)}\\
(\text{since}\,\,a=b)&\\
=&-\frac{2880a^6}{28288a^{12}}=-\frac{45}{442a^6}\ne 0.
\end{align*}
Hence the vectors are linearly independent. If we joint to the three $\delta_{m,n}$ vectors the vectors in $\{e_k\mid\,k\in\mathcal K^2\setminus\{(3,3)\}\,\}$ we obtain a family of $4^2+1$ linearly independent vectors spanning the space $\{e_k\mid\,k\in\mathcal K^2\cup\{(1,5),(3,3),(3,5)\}\,\}$. We arrive in this way to $span(\mathcal K^2\cup\{(1,5),(3,3),(3,5)\})\supset span(\mathcal K^2)$. Since $\mathcal K^2$ is G-Saturating for the square so is $\mathcal K^2\cup\{(1,5),(3,3),(3,5)\}$. Hence $\mathcal K^1$ is saturating for the square.
\end{Remark}
\section{Exact Controllability of Galerkin Approximations.}
Write the $\mathcal K^N$-Galerkin approximation of the NSE system \eqref{infsys}, with $\mathcal K^1$ as set of excited modes, in the concise form
\begin{equation}\label{galN}
N:
\begin{cases}
\dot u_k=\mathcal B_k(u)+\nu\Delta u_k +F_k+v_k&\quad k\in\mathcal K^1\\
\dot u_k=\mathcal B_k(u)+\nu\Delta u_k+F_k&\quad k\in\mathcal K^N\setminus\mathcal K^1\\
u\in\mathbb R^{\kappa_N}.&
\end{cases}
\end{equation}
In \cite{wei.mat} E. Weinam and J. Mattingly proved the Full Lie Rank Property for the 2D NSE with periodic conditions and for some class of few low modes controls. Similarly we prove that our equation also is ``full Lie rank'', i.e., Lie brackets at each point span the ambient space $\mathbb R^{\kappa_N}$.\par
Before we have proved that for all $N\in\mathbb N_0$ and all $t>0$ the system [\eqref{galN}.N] is time-$t$ approximately controllable:
$$
\forall u\in\mathbb R^{\kappa_N}\quad\overline{\mathcal A_u(\mathcal F_N)(t)}=\mathbb R^{\kappa_N}
$$ 
where $\mathcal F_N$ is the family of vector fields of system [\eqref{galN}.N], i.e.,
$$
\mathcal F_N=\{\mathcal B(\cdot)+\nu\Delta(\cdot)+F^{\kappa_N} +v\mid\,v\in\mathbb R^{\kappa_1}\}.
$$
Next we prove the (exact) controllability of system [\eqref{galN}.N]. For that we need to compute some Lie brackets.
\subsection{Lie Brackets. Full Lie Rank Property.}
From $\mathcal F_N$, we set the vector fields
$$
V^0:=\mathcal B+\nu\Delta +F^{\kappa_N},\quad X^i:=V^0+\frac{\partial}{\partial u_i}
$$
where $F^{\kappa_N}$ is the projection of $F$ onto $\mathbb R^\kappa_N$ and $\kappa_N=\#\mathcal K^N$ and, compute
$$
V^i:=[X^i,\,V^0]=[\frac{\partial}{\partial u_i},\,V^0]=\frac{\partial V^0}{\partial u}\frac{\partial}{\partial u_i}=\frac{\partial V^0}{\partial u_i},
$$
so
\begin{align*}
V^i_k:=&\Bigl(\sum_{\begin{subarray}{l}k=(n++i)^+\\i<n\end{subarray}}u_n C_{i,n}^{++}+\sum_{\begin{subarray}{l}k=(n+-i)^+\\i<n\end{subarray}}u_n C_{i,n}^{+-}\\
+&\sum_{\begin{subarray}{l}k=(n-+i)^+\\i<n\end{subarray}}u_n C_{i,n}^{-+}+\sum_{\begin{subarray}{l}k=(n- -i)^+\\i<n\end{subarray}}u_n C_{i,n}^{- -}\Bigr)+\delta^{[k,i]}\nu\bar k,
\end{align*}
where $\delta^{[k,i]}$ is the Kronecker delta function: $\delta^{[k,i]}=\begin{cases}1&\text{if}\; k=i\\ 0&\text{if}\; k\ne i\end{cases}$.
$$
V^{j,i}:=[X^j,\,V^i]=[X^j,\,[X^i,\,V^0]]=\frac{\partial V^i}{\partial u_j},
$$
so,
$$
V^{j,i}=\gamma_{i,j},\quad j>i.
$$
At each given point $\bar u\in\mathbb R^{\kappa_N}$ the family of brackets\footnote{Including the elements of $\mathcal F_N$ we consider brackets of ``length'' $0$.}
$$
\mathcal H^1:=\{V^0\pm\frac{1}{2}\frac{\partial}{\partial u_n},\,\alpha\gamma_{i,j}\mid\,n\in\mathbb R^{\kappa_1},\,(i,\,j)\in S_1\}
$$
span a superspace of $\mathbb R^{\kappa_2}$. Indeed, $\frac{\partial}{\partial u_n}=\bigl(V^0+\frac{1}{2}\frac{\partial}{\partial u_n}\bigr)-\bigl(V^0-\frac{1}{2}\frac{\partial}{\partial u_n}\bigr)$ .\footnote{Recall that $S_1$ has been chosen so that $\{e_n,\,\delta_{i,j}\mid\,n\in\mathcal K^1,\,(i,\,j)\in S_1\}$ span $\mathbb R^{\kappa_2}$.}\par
Now we have a technical difficulty, the vector fields have coeficients with not so nice expressions and, if we compute Lie brackets envolving them we will obtain even more complicated expressions. To avoid these expressions we prove by (finite) induction on $i\in\{1,\,\dots,\,\kappa_N\}$ that each  constant vector fied $\frac{\partial}{\partial u_n},\quad n\in\mathcal K^i$ is a linear combination of brackets:\\
\begin{itemize}
\item For $i=1$ we take the family $\mathcal H^1$.
\item Inductive Step: The induction hypothesis is:\\
``There is a family of brackets $\mathcal H^{p-1}=\{W^j\mid\, j=1,\,\dots,\,M_{p-1}\}$ such that every constant vector field $\frac{\partial}{\partial u_i},\quad i\in\mathcal K^{p-1}$ can be written as a linear combination of its elements:
$$
\frac{\partial}{\partial u_i}=\sum_{j=1}^{M_{p-1}}\alpha_i^j W^j,\quad \alpha_i^j\in\mathbb R;
$$
Then for all $i\in\mathcal K^{p-1}$
$$
V^i:=\Bigl[\frac{\partial}{\partial u_i},\,V^0\Bigr]\in span\{[W_j,\,V^0]\mid\,j=1,\,\dots,\,M_{p-1}\}
$$
and, for each $i,\,n\in\mathcal K^{p-1}$:
$$
V^{n,i}:=\Bigl[\frac{\partial}{\partial u_n},\,V^i\Bigr]=\gamma_{i,n}\in span\{[W^k,\,[W_j,\,V^0]]\mid\,k,\,j=1,\,\dots,\,M_{p-1}\}.
$$
Since the vectors in $\{\frac{\partial}{\partial u_k},\,\gamma_{i,n}\mid\,k\in\mathcal K^{p-1},\,(i,\,n)\in S_{p-1}\}$ span $\mathbb R^{\kappa_p}$ and can be written as a linear combination of brackets, then also the vector fields $\frac{\partial}{\partial u_i},\;i\in\mathcal K^{p}\setminus\mathcal K^{p-1}$ are linear combinations of brackets. The wanted family is $\mathcal H^{p-1}\cup\{[W^k,\,[W_j,\,V^0]]\mid\,k,\,j=1,\,\dots,\,M_{p-1}\}$.
\end{itemize}
Therefore, for all $N\in\mathbb N_0$, system [\eqref{galN}.N] is a full-rank bracket generating system. From that and from its approximate controllability\footnote{Approximate controllability at time $t$ implies, trivialy approximate controllability.} we conclude its controllability. Unfortunately for fixed time the bracket generating property is not suficient to conclude controllability from approximate controllability. To achieve controllability at time $t$ we shall need some lemmas which proofs can be found in \cite{jur}.
\subsection{Zero Orbits and Zero Ideal.}
\begin{Definition}
A {\bf zero-time orbit} $N_{0u}$ through $u$  of a family of vector fields $\mathcal F$ is the set
$$
N_{0u}:=\{u\circ e^{t_1V_1}\circ\,\dots\,\circ e^{t_pV_p}\mid\,p\in\mathbb N_0,\,V_i\in\mathcal F,\,t_i\in\mathbb R,\,\sum_{i=1}^p t_i=0\}
$$ 
\end{Definition}
\begin{Definition}
The {\bf derived algebra} of $\mathcal F$, denoted $\mathcal{D}_{er}(\mathcal F)$, is the set of all linear combinations of iterated brackets\footnote{Brackets of ``length'' $\geq 1$, considering the elements of $\mathcal F$ brackets of length $0$.}\\
The {\bf zero-time ideal}, denoted $\mathcal{I}(\mathcal F)$, is the span of elements in $\mathcal{D}_{er}(\mathcal F)$ and differences of the form $X-Y$ with $X$ and $Y$ in $\mathcal F$. 
\end{Definition}
\begin{Lemma}\label{LemmaA}
Let $\mathcal F$ be any family of analytic vector fields on an analytic manifold $M$. Let $N$ be an orbit of $\mathcal F$ and, $N_0$ be a zero orbit of $\mathcal F$ contained in $N$. Then we have the following:\\
\begin{itemize} 
\item Each connected component of $N_0$ is an orbit of $\mathcal{I}(\mathcal F)$;
\item For each $u\in N_0$, the tangent space of $N_0$ at $u$ is equal to the evaluation of $\mathcal{I}(\mathcal F)$ at $u$;
\item The dimension of $\mathcal{I}_u(\mathcal F)$ is constant as $u$ varies on $N$. It is equal either to $dim(Lie_u(\mathcal F))-1$ or to $dim(Lie_u(\mathcal F))$;
\item $dim(Lie_u(\mathcal F))=dim(\mathcal I_u(\mathcal F))$ iff $X(u)\in \mathcal I_u(\mathcal F)$ for some $X\in\mathcal F$.
\end{itemize}
\end{Lemma}
\begin{Lemma}\label{LemmaB}
Suppose that $\mathcal F$ is a family of vector fields on $M$ such that both $\mathcal F$ and its zero-time ideal $\mathcal{I}(\mathcal F)$ are Lie-determined (the evaluation of Lie brackets at each point span the tangent space to the orbit). In addiction, assume that $\mathcal F$ contains a complete vector field. Then\begin{itemize}
\item  $\mathcal A_u(\mathcal F)(t)$ is a connected subset of some zero orbit $N_{0z}$ through some element $z\in M$. 
\item $\mathcal A_u(\mathcal F)(t)$ has a nonempty interior in the manifold topology of the zero-orbit where it is contained. Moreover, the set of interior points is dense in  $\mathcal A_u(\mathcal F)(t)$.
\end{itemize}
\end{Lemma}
Coming back to our system [\eqref{galN}.N], by Lemma \ref{LemmaA} and due to the fact that $V^0(0)=0\in\mathcal{I}(\mathcal F)$, we have
$$
dim(Lie_0(\mathcal F_N))=dim(\mathcal I_0(\mathcal F_N))=\kappa_N;
$$
which means that the zero-time orbit $N_0$ through $0$ has dimension $\kappa_N$ and, since that dimension is constant in all points in the unique orbit $\mathbb R^{\kappa_N}$ of the system, we conclude that $N_0$ is a union of connected components of dimension $\kappa_N$. Since the dimension of that components is $\kappa_N$ their topology coincide with that of $\mathbb R^{\kappa_N}$ and, from the fact that the zero-time orbits form a partition of $\mathbb R^{\kappa_N}$ we conclude that $\mathbb R^{\kappa_N}$ is a union of connected open sets. Therefore there is only one zero-orbit, it is the whole state space $\mathbb R^{\kappa_N}$.\par
By Lemma \ref{LemmaB}, and by the fact that $V^0$ is a complete vector field which follows from the estimate $|u(s)|\leq|u(0)|+\frac{s}{\nu}\|F\|_{V^\prime}^2$ (see estimate \eqref{|u|finite} with $\tilde F=F$), the interior $int\mathcal A_u(\mathcal F_N)(t)$ of the attainable set from $u$ at time $t$ is dense in $\mathcal A_u(\mathcal F)(t)$, where the interior and density are relative to the topology of $\mathbb R^{\kappa_N}$ because that is the topology of the zero-orbit. Hence we arrive to the equality
$$
\overline{int\mathcal A_u(\mathcal F_N)(t)}=\overline{\mathcal A_u(\mathcal F_N)(t)}=\mathbb R^{\kappa_N} 
$$
for all $t>0$.\par
Now we can prove the controllability at time $t$ of  system [\eqref{galN}.N]: Let $u,\,z$ be two elements in $\mathbb R^{\kappa_N}$. Since the intersection of two open dense sets stills open and dense, we may take a point
$$
w\in int\mathcal A_u(\mathcal F_N)(t/2)\cap int\mathcal A_z(\mathcal -F_N)(t/2).
$$
[Note that the family $-\mathcal F_N:=\{-V\mid\,V\in\mathcal F\}$ satisfies the requirements of lemmas \ref{LemmaA} and \ref{LemmaB}.]\par
Then we can write
\begin{align*}
w&=u\circ e^{t_1V_1}\circ\,\cdots\,\circ e^{t_nV_n},\qquad V_i\in\mathcal F_N,\,t_i\geq 0,\,\sum_{i=1}^n t_i=\frac{t}{2};\\
w&=z\circ e^{-s_1W_1}\circ\,\cdots\,\circ e^{-s_mW_m},\qquad W_i\in\mathcal F_N,\,s_i\geq 0,\,\sum_{i=1}^m t_i=\frac{t}{2};
\end{align*}
So, $z$ is reachable from $u$ in time $t$:
$$
z=u\circ e^{t_1V_1}\circ\,\cdots\,\circ e^{t_nV_n}\circ e^{s_mW_m}\circ\,\cdots\,\circ e^{s_1W_1}.
$$

\chapter{Controllability in Observed Component.}\label{Ch:COP}

\section{Controlled N-S Problem. Existence, Uniqueness and Continuity.}
In chapter \ref{Ch:EUC} we have presented the classical, weak and strong formulations of \eqref{c.s.nse}-\eqref{c.s.inicond}. So for the controlled version with $F+v$ in the place of $\tilde F$ we have existence, unicity and continuity in the data $(u_0,\,F+v,\,\nu)$, because our control is an essencially bounded function, so that $F+v$ will belong where $\tilde F$ did: ($L^2(0,\,T,\,V^\prime)$ or $L^2(0,\,T,\,H)$). If we consider the initial data as $(u_0,\,\tilde F,\,v,\,\nu)$, by Theorems \ref{ex.w.sol}, \ref{contindata} and \ref{contindata2} we easily conclude that (considering again the external force depending on time):
\begin{Corollary}\label{ex.w.solv}
Given 
$$
(u_0,\,\tilde F,\,v,\,\nu)\in H\times L^2(0,\,T,\,V^\prime)\times L^\infty(0,\,T,\,V^\prime)\times]0,\,+\infty[,
$$
there is at least one weak solution $u\in C([0,\,T],\,H)$ for Problem \ref{w.p.leray} with $\tilde F+v$ in the place of $\tilde F$.
\end{Corollary}
and,
\begin{Corollary}\label{contindatav}
The maps
\begin{align*}
\mathbb S:\,H\times L^2(0,\,T,\,V^\prime)\times L^\infty(0,\,T,\,V^\prime)\times]0,\,+\infty[&\to C([0,\,T],\,H)\\
(u_0,\,\tilde F,\,v,\,\nu)&\mapsto u
\end{align*}
and
\begin{align*}
\mathbb S_2:\,H\times L^2(0,\,T,\,V^\prime)\times L^\infty(0,\,T,\,V^\prime)\times]0,\,+\infty[&\to L^2(0,\,T,\,V)\\
(u_0,\,\tilde F,\,v,\,\nu)&\mapsto u
\end{align*}
are continuous.
\end{Corollary}
Similarly, by Theorems \ref{ex.s.sol}, \ref{scontindata} and \ref{scontindata2} we can easily see that
\begin{Corollary}\label{ex.s.solv}
Given 
$$
(u_0,\,\tilde F,\,v,\,\nu)\in V\times L^2(0,\,T,\,H)\times L^\infty(0,\,T,\,H)\times]0,\,+\infty[,
$$
there is at least one weak solution $u\in C([0,\,T],\,V)$ for Problem \ref{s.p.leray} with $\tilde F+v$ in the place of $\tilde F$.
\end{Corollary}
and,
\begin{Corollary}\label{scontindatav}
The maps
\begin{align*}
\mathbb S_s:\,V\times L^2(0,\,T,\,V^\prime)\times L^\infty(0,\,T,\,H)\times]0,\,+\infty[&\to C([0,\,T],\,V)\\
(u_0,\,\tilde F,\,v,\,\nu)&\mapsto u
\end{align*}
and
\begin{align*}
\mathbb S_{2s}:\,V\times L^2(0,\,T,\,V^\prime)\times L^\infty(0,\,T,\,H)\times]0,\,+\infty[&\to L^2(0,\,T,\,D(A))\\
(u_0,\,\tilde F,\,v,\,\nu)&\mapsto u
\end{align*}
are continuous.
\end{Corollary}

\section{Change of Variables.}\label{S:ch.var}
If we make the change of variables
$$
u=y+\mathbb Iv
$$
where $\mathbb I$ is the primitive operator --- $[\mathbb Iv](t)=\int_0^t v(\tau)\,d\tau$,from
$$
u^\prime=-\nu Au-Bu+\tilde F+v
$$
we arrive to the equation
$$
y^\prime=-\nu A(y+\mathbb Iv)-B(y+\mathbb Iv)+\tilde F.
$$
Note that the function $v$ appears only implicitly in the last equation. Now we forget that $\mathbb Iv$ is a primitive of an essentially bounded function and replace it by $P$ in the equation. Since $v$ is a low modes forcing it takes value in a finite-dimensional space and, $\mathbb Iv$ being a primitive we have $\mathbb Iv\in C([0,\,T],\,D(A))$. But we take $P$ in the larger space $L^4(0,\,T,\,D(A))$.\\
\section{Weak Case.}
We want to study the following equivalent problems \ref{yleray} and \ref{yAB}.
\begin{Problem}\label{yleray}
Given
\begin{align}
&\quad \tilde F\in L^2(0,\,T,\,V^\prime),\quad P\in L^4(0,\,T,\,D(A))\label{yleray.FP}\\
&\And\notag\\
&\quad y_0\in H\label{yleray.y0}\\
&to\, find\notag\\
&\quad y\in L^2(0,\,T,\,V)\label{yleray.y}\\
&\text{satisfying (in the distribution sense)}\notag\\
&\quad \forall v\in V:\notag\\
&\quad \frac{d}{dt}(y,\,v)+\nu((y+P,\,v))+b(y+P,\,y+P,\,v)=<\tilde F,\,v>,\label{yleray.nse}\\
&\text{and}\notag\\
&\quad y(0)=y_0\label{yleray.inicond}.
\end{align}
\end{Problem}
\begin{Problem}\label{yAB}
Given
\begin{align}
&\quad \tilde F\in L^2(0,\,T,\,V^\prime),,\quad P\in L^4(0,\,T,\,D(A))\label{yAB.FP}\\
&\And\notag\\
&\quad y_0\in H,\label{yAB.y0}\\
&\text{to find}\notag\\
&\quad y\in L^2(0,\,T,\,V),\qquad y^\prime\in L^1(0,\,T,\,V^\prime)\label{yAB.yy'}\\
&\text{satisfying}\notag\\
&\quad y^\prime+\nu A(y+P)+B(y+P)=\tilde F\quad\text{on}\quad]0,\,T[,\label{yAB.nse}\\
&\text{and}\notag\\
&\quad y(0)=y_0\label{yAB.inicond}.
\end{align}
\end{Problem}
\begin{Remark}
The equivalence of these problems can be shown the way we proved the equivalence of problems \ref{w.p.leray} and \ref{w.p.AB} (see chapter \ref{Ch:EUC}). This equivalence follows from
$$
\begin{cases}
&\bullet\;\tilde F\in L^2(0,\,T,\,V^\prime);\\
&\bullet\;A(y+P)\in L^2(0,\,T,\,V^\prime).\\
&\qquad\qquad\text{Indeed $\|A(y+P)\|_{V^\prime}\leq\|y+P\|$ and $y+P\in L^2(0,\,T,\,V)$};\\
&\bullet\;B(y+P)\in L^1(0,\,T,\,V^\prime).\\
&\qquad\qquad\text{Indeed $\|B(y+P)\|_{V^\prime}\leq C\|y+P\|^2$ and $y+P\in L^2(0,\,T,\,V)$}.
\end{cases}
$$
So that $\tilde F-\nu A(y+P)-B(y+P)\in L^1(0,\,T,\,V^\prime)$.
\end{Remark}

\section{Existence.}
We have the Theorem
\begin{Theorem}\label{ex.soly}
Given $\tilde F,\;P$ and $y_0$ satisfying \eqref{yAB.FP} and
\eqref{yAB.y0}.
There is at least one function $y$ satisfying
\eqref{yAB.yy'}-\eqref{yAB.inicond}.
\end{Theorem}
The proof is analogous to that of Theorem
\ref{ex.w.sol}.
Basically it differs only in some estimates we
compute now: Following the proof of Theorem \ref{ex.w.sol}
presented in chapter \ref{Ch:EUC}, with the suitable adaptations to
problems \ref{yleray}-\ref{yAB}, we define an approximate
solution
$$
y^m=\sum_{\max\{i_1,\,i_2\}\leq m}y_i^m(t)W_i
$$ for each $m\in\mathbb N_0$ and arrive to the equation
\begin{multline}
<(y^m)^\prime,\,y^m>+\nu<A(y^m+P^m),\,y^m>\\+<B(y^m+P^m),\,y^m>=<\tilde F,\,y^m>.\label{pp.est}\quad\footnotemark
\end{multline}
\footnotetext{Where $P^m$ is the projection of $P$ onto $span\{W_i\mid\,\max\{i_1,\,i_2\}\leq m\}$.}
We note that
\begin{eqnarray*}
&\bullet&<A(y^m+P^m),\,y^m>=\|y^m\|^2+((P^m,\,y^m))\\
&\bullet&b(y^m+P^m,\,y^m+P^m,\,y^m)\\
& &=b(y^m+P^m,\,y^m+P^m,\,y^m+P^m)-b(y^m+P^m,\,y^m+P^m,\,P^m)\\
& &= -b(y^m+P^m,\,y^m+P^m,\,P^m)\\
& &=-b(y^m+P^m,\,y^m,\,P^m)-b(y^m+P^m,\,P^m,\,P^m)\\
& &=-b(y^m+P^m,\,y^m,\,P^m)=-b(y^m,\,y^m,\,P^m)-b(P^m,\,y^m,\,P^m).
\end{eqnarray*}
Hence from \eqref{pp.est} we have
\begin{multline*}
\frac{d}{dt}|y^m|^2+2\nu\|y^m\|^2=-2\nu((P^m,\,y^m))\\
-b(y^m,\,P^m,\,y^m)+b(P^m,\,y^m,\,P^m)+<\tilde F,\,y^m>;
\end{multline*}
hence
$$
\frac{d}{dt}|y^m|^2+2\nu\|y^m\|^2\leq 2\nu\|P\|\|y^m\|+C|y^m|\|y^m\|\|P\|+C\|P\|^2\|y^m\|+\|\tilde F\|_{V^\prime}\|y^m\|;
$$
thus
\begin{equation}\label{p.est}
\frac{d}{dt}|y^m|^2+\nu\|y^m\|^2\leq 4\nu\|P\|^2+\frac{C^2}{\nu}|y^m|^2\|P\|^2+\frac{C^2}{\nu}\|P\|^4+\frac{1}{\nu}\|\tilde F\|_{V^\prime}^2.
\end{equation}
From equation \eqref{p.est} and from Gronwall Inequality we can
derive the estimates \eqref{E1} and \eqref{E2} below:
\begin{align}
&|y^m(s)|^2\leq\exp\int_0^T\frac{C^2}{\nu}\|P(t)\|^2\,dt\biggl(|y_0|^2\notag\\
+& \int_0^T 4\Bigl(\nu\|P(t)\|^2+\frac{C^2}{\nu}\|P(t)\|^4+\frac{1}{\nu}\|\tilde F(t)\|^2\Bigr)\,dt\biggr)\label{E1}.
\end{align}
Since $P\in L^4(0,\,T,\,D(A))\subset L^4(0,\,T,\,V)\subset L^2(0,\,T,\,V)$ and $\tilde F\in L^2(0,\,T,\,V)$, \eqref{E1} shows that
\begin{equation}\label{ymbddLinfH}
y^m\; remains\; in\; a\; bounded\; set\; of\; L^\infty(0,\,T,\,H).
\end{equation}
From
\begin{multline}
|y^m(T)|^2-|y^m(0)|^2+\nu\int_0^T\|y^m(t)\|^2\,dt\\
\leq\int_0^T \Bigl(4\nu\|P(t)\|^2+\frac{C^2}{\nu}\|P(t)\|^4+\frac{1}{\nu}\|\tilde F(t)\|^2\Bigr)\,dt+\int_0^T\frac{C^2}{\nu}|y^m(t)|^2\|P(t)\|^2\,dt.\label{E2}
\end{multline}
By \eqref{ymbddLinfH} the last integral is finite. Then \eqref{E2} shows that
\begin{equation}\label{ymbddL2V}
y^m\; remains\; in\; a\; bounded\; set\; of\; L^2(0,\,T,\,V).
\end{equation}
The rest of the proof is completely analogous.
\begin{Remark}\label{PL2V}
As we can see in the estimates \eqref{E1} and \eqref{E2}, considering $P\in L^4(0,\,T,\,V)$ is sufficient to guarantee existence of weak solutions.
\end{Remark}
\section{Uniqueness.}
\begin{Theorem}\label{cont.yH}
The solution of Problems \ref{yleray}-\ref{yAB} given by Theorem \ref{ex.soly} is unique. Moreover it is a.e. equal to a continuous function from $[0,\,T]$ into $H$.
\end{Theorem}
The continuity follows from
\begin{eqnarray*}
&\bullet&A(y+P)\in L^2(0,\,T,\,V^\prime)\\
&\bullet&\tilde F\in L^2(0,\,T,\,V^\prime)\\
&\bullet&B(y+P)\in L^2(0,\,T,\,V^\prime).
\end{eqnarray*}
Indeed these expressions imply
$$
y\in L^2(0,\,T,\,V)\And y^\prime\in L^2(0,\,T,\,V^\prime).
$$
The first two expressions we already know to be true. The last
one follows from
\begin{eqnarray*}
& &\|B(y+P)\|\leq C|y|\|P\|+C\|P\|^2;\\
& &P\in L^4(0,\,T,\,D(A))\subset L^4(0,\,T,\,V);\quad y\in L^\infty(0,\,T,\,H).
\end{eqnarray*}
To conclude the uniqueness we consider two solutions $y$ and $z$
of problems \eqref{yleray}-\eqref{yAB}. The difference $w:=y-z$
satisfies
$$
w^\prime=-\nu Aw-B(y+P)+B(z+P)
$$
from which we derive
\begin{align*}
&\;\frac{d}{dt}|w|^2+2\nu\|w\|^2\leq 2C|w|\|w\|\|z+P\|\\
&\;\frac{d}{dt}|w|^2\leq\frac{C^2}{2\nu}|w|^2\|z+P\|^2\\
&\;|w(s)|^2\leq|w(0)|^2\exp\int_0^T\frac{C^2}{2\nu}\|z(t)+P(t)\|^2\,dt=0.
\end{align*}

\section{Continuity.}
\begin{Theorem}\label{ycontindata}
The map
\begin{align*}
\mathbb Y:\,H\times L^2(0,\,T,\,V^\prime)\times L^4(0,\,T,\,D(A))\times]0,\,+\infty[&\to C([0,\,T],\,H)\\
(y_0,\,\tilde F,\,P,\,\nu)&\mapsto y
\end{align*}
is continuous. Where $y$ is the unique solution of problems \eqref{yleray}-\eqref{yAB} corresponding to the data $(y_0,\,\tilde F,\,P,\,\nu)$.
\end{Theorem}
\begin{proof}
The proof is analogous to that of Theorem \ref{contindatav}. We fix a quadruple
$$
(y_0,\,\tilde F,\,P,\,\nu)\in H\times L^2(0,\,T,\,V^\prime)\times L^4(0,\,T,\,D(A))\times]0,\,+\infty[
$$
and $\varepsilon>0$. Then consider another quadruple
$$
(z_0,\,G,\,Q,\,\eta)\in H\times L^2(0,\,T,\,V^\prime)\times L^4(0,\,T,\,D(A))\times]0,\,+\infty[.
$$
Put
$$
y:=\mathbb Y(y_0,\,\tilde F,\,P,\,\nu)\And z:=\mathbb Y(z_0,\,G,\,Q,\,\eta)
$$ so that
$$
y^\prime + \nu A(y+P) + B(y+P) = \tilde F\And z^\prime + \eta A(z+Q) + B(z+Q) = G.
$$
Putting $w:=z-y$ we obtain
$$
w^\prime=G-\tilde F-\eta Aw -\eta A(Q-P)+(\nu-\eta)A(y+P)-B(z+Q)+B(y+P).
$$
Taking the scalar product with $w$ we obtain
\begin{align}
<w^\prime,\,w>&=<G-\tilde F,\,w>-\eta\|w\|^2-\eta((Q-P,\,w))\notag\\
&+(\nu-\eta)((y+P,\,w))+b(y+P,\,y+P,\,w)-b(z+Q,\,z+Q,\,w).\label{n.simp.b}
\end{align}
Now we estimate the last term of the last equality:
\begin{align*}
&b(y+P,\,y+P,\,w)-b(z+Q,\,z+Q,\,w)\\
&\qquad=b(y,\,y,\,w)+b(y,\,P,\,w)+b(P,\,y,\,w)+b(P,\,P,\,w)\\
&\qquad-b(z,\,z,\,w)-b(z,\,Q,\,w)-b(Q,\,z,\,w)-b(Q,\,Q,\,w)\\
&\qquad=b(y,\,y,\,w)-b(z,\,z,\,w)+b(y,\,P,\,w)-b(z,\,Q,\,w)\\
&\qquad+b(P,\,y,\,w)-b(Q,\,y,\,w)+b(P,\,P,\,w)-b(Q,\,Q,\,w).
\end{align*}
We put
\begin{align*}
Z_1&:=b(y,\,y,\,w)-b(z,\,z,\,w)\\
Z_2&:=b(y,\,P,\,w)-b(z,\,Q,\,w)\\
Z_3&:=b(P,\,y,\,w)-b(Q,\,z,\,w)\\
Z_4&:=b(P,\,P,\,w)-b(Q,\,Q,\,w).
\end{align*}
Then
\begin{equation}\label{b<Zs}
|b(y+P,\,y+P,\,w)-b(z+Q,\,z+Q,\,w)|\leq|Z_1|+|Z_2|+|Z_3|+|Z_4|
\end{equation}
Since
\begin{align}
|Z_1|&=|b(y,\,y,\,w)-b(z,\,z,\,w)|=|-b(w,\,y,\,w)|\notag\\
&\leq C|w|\|w\|\|y\|,\quad\text{by \eqref{E.01.1.01}}\label{B.Z_1}\\
|Z_2|&=|b(y,\,P,\,w)-b(z,\,Q,\,w)|=|b(y,\,P,\,w)-b(y,\,Q,\,w)-b(w,\,Q,\,w)|\notag\\
&=|b(y,\,P-Q,\,w)-b(w,\,Q,\,w)|\leq|b(y,\,w,\,P-Q)|+|b(w,\,Q,\,w)|\notag\\
&\leq C|y|\|w\||P-Q|^\frac{1}{2}|P-Q|_{[2]}^\frac{1}{2}+C|w|\|w\|\|Q\|\quad\text{by \eqref{E.0.1.02} and \eqref{E.01.1.01}},\notag\\
&\leq C_1|y|\|w\||P-Q|_{[2]}+C|w|\|w\|\|Q\|,\quad\text{by $|\cdot|\leq\|\cdot\|_2$}\label{B.Z_2}\\
|Z_3|&=|b(P,\,y,\,w)-b(Q,\,z,\,w)|=|b(P,\,y,\,w)-b(Q,\,y,\,w)-b(Q,\,w,\,w)|\notag\\
&=|-b(P-Q,\,w,\,y)|\leq C|P-Q|^\frac{1}{2}|P-Q|_{[2]}^\frac{1}{2}\|w\||y|\quad\text{by \eqref{E.02.1.0}}\notag\\
&\leq C_1|P-Q|_{[2]}\|w\||y|,\label{B.Z_3}\\
|Z_4|&=|b(P,\,P,\,w)-b(Q,\,Q,\,w)|=|b(P-Q,\,P,\,w)+b(Q,\,P,\,w)\\
&-b(Q,\,Q,\,w)|\leq|b(P-Q,\,P,\,w)|+|b(Q,\,P-Q,\,w)|\notag\\
&\leq C\|P-Q\|\|P\|\|w\|+C\|Q\|\|P-Q\|\|w\|\quad\text{by \eqref{E.1.1.1}}.\label{B.Z_4}
\end{align}
From \eqref{b<Zs}, \eqref{B.Z_1}-\eqref{B.Z_4} and \eqref{n.simp.b} we obtain
\begin{multline*}
\frac{d}{dt}|w|^2\leq 2\|G-\tilde F\|_{V^\prime}\|w\|-2\eta\|w\|^2+2\eta\|Q-P\|\|w\|+2|\nu-\eta|\|y+P\|\|w\|\\+2C|w|\|w\|\|y\|+2C_1|y|\|w\||P-Q|_{[2]}+2C|w|\|w\|\|Q\|\\+2C_1|P-Q|_{[2]}\|w\||y|+2C\|P-Q\|\|P\|\|w\|+2C\|Q\|\|P-Q\|\|w\|;
\end{multline*}
Hence
\begin{multline*}
\frac{d}{dt}|w|^2+\eta\|w\|^2\leq
\frac{8}{\eta}\|G-\tilde F\|_{V^\prime}^2+8\eta\|Q-P\|^2+
\frac{8}{\eta}|\nu-\eta|^2\|y+P\|^2\\
+\frac{8}{\eta}C^2|w|^2(\|y\|^2+\|Q\|^2)+\frac{8}{\eta}4C_1^2|y|^2|P-Q|_{[2]}^2
+\frac{8}{\eta}C^2\|P-Q\|^2(\|P\|^2+\|Q\|^2).
\end{multline*}
Now if $\eta$ satisfies
$$
|\nu-\eta|<\frac{\nu}{2}
$$
we have $\eta\in]\frac{\nu}{2},\,\frac{3\nu}{2}[$. Hence, from the last equation we obtain
\begin{multline}
\frac{d}{dt}|w|^2+\frac{\nu}{2}\|w\|^2\leq
\frac{16}{\nu}\|G-\tilde F\|_{V^\prime}^2+|Q-P|_{[2]}^2(12\nu C_2+\frac{16}{\nu}4C_1^2|y|^2)\\
+\frac{16}{\nu}|\nu-\eta|^2\|y+P\|^2+\frac{16}{\nu}C^2|w|^2(\|y\|^2+\|Q\|^2)
+\frac{16}{\nu}C^2\|P-Q\|^2(\|P\|^2+\|Q\|^2).\label{estYHV1}
\end{multline}
By the Gronwall Inequality:
\begin{multline}\label{est.w.Y}
|w(s)|^2\leq\exp\biggl(\int_0^T\frac{16}{\nu}C^2\bigl(\|y(t)\|^2+\|Q(t)\|^2
\bigr)\,dt\biggr)\Biggl(|w(0)|^2\\
+\int_0^T
\Bigl[\frac{16}{\nu}\|G(t)-\tilde F(t)\|_{V^\prime}^2+|Q(t)-P(t)|_{[2]}^2(12\nu C_2+\frac{16}{\nu}4C_1^2|y|^2)\\
+\frac{16}{\nu}|\nu-\eta|^2\|y(t)+P(t)\|^2+\frac{16}{\nu}C^2\|P(t)-Q(t)\|^2\bigl(\|P(t)\|^2+\|Q(t)\|^2\bigr)\Bigr]\,dt\Biggr).
\end{multline}
Now we consider two cases $P=0\And P\ne 0$.\par
If $P\ne 0$ and if  $Q$ satisfies
$$
\|P-Q\|_{L^4(0,T,D(A))}<\frac{\|P\|_{L^4(0,T,D(A))}}{2}
$$
and since $\|y\|_{L^\infty(0,\,T,\,H)}=:D<+\infty$ we have that
\begin{multline*}
|w(s)|^2\leq\exp\biggl[C_3\biggl(\int_0^T\|y(t)\|^2\,dt+\Bigl(\int_0^T\|P(t)\|^4\,dt\Bigr)^\frac{1}{2}\biggr)\,\biggr]\Biggl(|w(0)|^2\\
+C_4\biggr[\int_0^T \|G(t)-\tilde F(t)\|_{V^\prime}^2\,dt+\Bigl(\int_0^T|Q(t)-P(t)|_{[2]}^4\,dt\Bigr)^\frac{1}{2}\\
+|\nu-\eta|^2\int_0^T\|y(t)+P(t)\|^2\,dt+\Bigl(\int_0^T\|P-Q\|^4\,dt\Bigr)^\frac{1}{2}\Bigl(\int_0^T\|P\|^4\,dt\Bigr)^\frac{1}{2}\biggr]\Biggr).
\end{multline*}
Thus
\begin{multline*}
|w(t)|^2\leq C_5\biggl(|w(0)|^2+\|G-\tilde F\|_{L^2(0,\,T,\,V^\prime)}^2+\|Q-P\|_{L^4(0,\,T,\,D(A))}^2\\+|\nu-\eta|^2+\|Q-P\|_{L^4(0,\,T,\,D(A))}^2\biggr).
\end{multline*}
We rewrite the last equation as
$$
|w(t)|^2\leq C_6^2\biggl(|w(0)|^2+\|G-\tilde F\|_{L^2(0,\,T,\,V^\prime)}^2+\|Q-P\|_{L^4(0,\,T,\,D(A))}^2+|\nu-\eta|^2\biggr).
$$
Then (in the case $P\ne0$), if the quadruple $(z_0,\,G,\,Q,\,\eta)$ satisfies
\begin{align*}
&|\nu-\eta|<\min\{\frac{\nu}{2},\,\frac{\varepsilon}{2C_6}\};\\
&\|P-Q\|_{L^4(0,\,T,\,D(A))}<\min\Bigl\{\frac{\|P\|_{L^4(0,\,T,\,D(A))}}{2},\,\frac{\varepsilon}{2C_6}\Bigr\};\\
&|z_0-y_0|=|w(0)|<\frac{\varepsilon}{2C_6};\quad\|G-\tilde F\|_{L^2(0,\,T,\,V^\prime)}<\frac{\varepsilon}{2C_6};
\end{align*}
we have that $|w(s)|<\varepsilon$. Therefore we have the continuity of $\mathbb Y$ in all quadruples such that $P\ne0$. Note that $C_6$ depends only in the fixed quadruple $(y_0,\,\tilde F,\,P,\,\nu)$ where we are studying the continuity.\;\footnote{$C_6$ depends in $R$ and $T$ too, but $T$ and $R$ are fixed in the statement of the Theorem.}\par
In the case $P=0$ we can start from a bound for $Q$, say that $Q$ satisfies $\|Q\|_{L^4(0,\,T,\,D(A))}<C_7$ (so that the argument of the exponential is bounded by a constant depending only on $C_7$ and $(y_0,\,\tilde F,\,P,\,\nu)$, but not on $Q$). Then from \eqref{est.w.Y} we can conclude that
\begin{multline*}
|w(t)|^2\leq C_8\biggl(|w(0)|^2+\|G-\tilde F\|_{L^2(0,\,T,\,V^\prime)}^2\\+\|Q\|_{L^4(0,\,T,\,D(A))}^2+|\nu-\eta|^2+\|Q\|_{L^4(0,\,T,\,D(A))}^4\biggr).
\end{multline*}
If $\|Q\|_{L^4(0,\,T,\,D(A))}<1$ we can rewrite the last equation obtaining
$$
|w(t)|^2\leq C_9^2\biggl(|w(0)|^2+\|G-\tilde F\|_{L^2(0,\,T,\,V^\prime)}^2+\|Q\|_{L^4(0,\,T,\,D(A))}^2+|\nu-\eta|^2\biggr)
$$
So in the case $P=0$ we choose, for the quadruple $(z_0,\,G,\,Q,\,\eta)$, the bounds
\begin{align*}
|\nu-\eta|<\min\{\frac{\nu}{2},\,\frac{\varepsilon}{2C_9}\};\quad&\|Q\|_{L^4(0,\,T,\,D(A))}<\min\{1,\,\frac{\varepsilon}{2C_9}\};\\
|z_0-y_0|=|w(0)|<\frac{\varepsilon}{2C_9};\quad&\|G-\tilde F\|_{L^2(0,\,T,\,V^\prime)}<\frac{\varepsilon}{2C_9};
\end{align*}
and $\mathbb Y$ is continuous in this case too.
\end{proof}

\begin{Theorem}\label{ycontindata2}
The map
\begin{align*}
\mathbb Y_2:\,H\times L^2(0,\,T,\,V^\prime)\times
L^4(0,\,T,\,D(A))\times]0,\,+\infty[&
\to L^2([0,\,T],\,V)\\
(y_0,\,\tilde F,\,P,\,\nu)&\mapsto y
\end{align*}
is continuous. Where $y$ is the unique solution of problems
\eqref{yleray}-\eqref{yAB} corresponding to the data
$(y_0,\,\tilde F,\,P,\,\nu)$.
\end{Theorem}
\begin{proof}
We Fix $\varepsilon>0$, $(y_0,\,\tilde F,\,P,\,\nu)\in H\times
L^2(0,\,T,\,V^\prime)\times L^4(0,\,T,\,D(A))\times]0,\,+\infty[$
and put $y:=\mathbb Y_2(y_0,\,\tilde F,\,P,\,\nu)$.\\
Let $(z_0,\,G,\,Q,\,\eta)$ be another quadruple in the product $H\times
L^2(0,\,T,\,V^\prime)\times
L^4(0,\,T,\,D(A))\times]0,\,+\infty[$. Put $z:=\mathbb
Y_2(z_0,\,G,\,Q,\,\eta)$ and $w:=z-y$ and like in the proof of
Theorem \ref{ycontindata} we arrive to \eqref{estYHV1}:

 \begin{multline*}
\frac{d}{dt}|w|^2+\frac{\nu}{2}\|w\|^2\leq
\frac{16}{\nu}\|G-\tilde F\|_{V^\prime}^2+\frac{16}{\nu}\|Q-P\|_2^2(1+4C_1^2|y|^2)\\
+\frac{16}{\nu}|\nu-\eta|^2\|y+P\|^2+\frac{16}{\nu}C^2|w|^2(\|y\|^2+\|Q\|^2)
+\frac{16}{\nu}C^2\|P-Q\|^2(\|P\|^2+\|Q\|^2).
\end{multline*}
Integrating over $[0,\,T]$

\begin{multline*}
|w(T)|^2-|w(0)|^2+\int_0^T\frac{\nu}{2}\|w(t)\|^2\,dt\leq
\int_0^T\biggl[\frac{16}{\nu}\|G(t)-\tilde F(t)\|_{V^\prime}^2\\+
\frac{16}{\nu}\|Q(t)-P(t)\|_2^2(1+4C_1^2|y(t)|^2)
+\frac{16}{\nu}|\nu-\eta|^2\|y(t)+P(t)\|^2\\+\frac{16}{\nu}C^2|w(t)|^2(\|y(t)\|^2+\|Q(t)\|^2)
+\frac{16}{\nu}C^2\|P(t)-Q(t)\|^2(\|P(t)\|^2+\|Q(t)\|^2)\biggr]\,dt.
\end{multline*}
Looking at the last two terms we see that
\begin{align*}
\int_0^T \|y(t)\|^2+\|Q(t)\|^2\,dt&=\int_0^T \|y(t)\|^2\,dt+\int_0^T\|Q(t)\|^2\,dt\\
&\leq D_1+D_2\|Q\|_{L^4(0,\,T,\,D(A))}^2
\end{align*}
and
\begin{align*}
&\int_0^T\|P(t)-Q(t)\|^2\|P(t)\|^2\,dt+\int_0^T\|P(t)-Q(t)\|^2\|Q(t)\|^2\,dt\\
\leq&D_3\|P-Q\|_{L^4(0,\,T,\,D(A))}^2\biggl(\|P\|_{L^4(0,\,T,\,D(A))}^2+\|Q\|_{L^4(0,\,T,\,D(A))}^2\biggr)
\end{align*}
so, we arrive to
\begin{multline*}
\int_0^T\|w(t)\|^2\,dt\leq
C_1\|G-\tilde F\|_{L^2(0,\,T,\,V^\prime}^2+C_2\|Q-P\|_{L^4(0,\,T,\,D(A)}^2\\
+C_3|\nu-\eta|^2+C_4|w|_{C([0,\,T],\,H)}^2(1+\|Q\|_{L^4(0,\,T,\,D(A)}^2)\\
+C_5\|P-Q\|_{L^4(0,\,T,\,D(A)}^2(\|P\|_{L^4(0,\,T,\,D(A)}^2+\|Q\|_{L^4(0,\,T,\,D(A)}^2)
\end{multline*}
where the constants do not depend on the quadruple $(z_0,\,G,\,Q,\,\eta)$.\\
If we choose $Q$ close to $P$, say $\|P-Q\|_{L^4(0,\,T,\,D(A)}<1$
we have
\begin{multline*}
\int_0^T\|w(t)\|^2\,dt\leq
C_1\|G-\tilde F\|_{L^2(0,\,T,\,V^\prime}^2+C_2\|Q-P\|_{L^4(0,\,T,\,D(A)}^2\\
+C_3|\nu-\eta|^2+D_4|w|_{C([0,\,T],\,H)}^2\biggl(1+\|P\|_{L^4(0,\,T,\,D(A)}\biggr)^2\\
+D_5\|P-Q\|_{L^4(0,\,T,\,D(A)}^2\biggl(1+\|P\|_{L^4(0,\,T,\,D(A)}\biggr)^2.
\end{multline*}
or
\begin{multline*}
\int_0^T\|w(t)\|^2\,dt\leq
C_1\|G-\tilde F\|_{L^2(0,\,T,\,V^\prime)}^2\\+
D_6\|Q-P\|_{L^4(0,\,T,\,D(A)}^2
+C_3|\nu-\eta|^2+D_7|w|_{C([0,\,T],\,H)}^2
\end{multline*}
(where the constants do not depend on the quadruple $(z_0,\,G,\,Q,\,\eta)$).\\
By Theorem \ref{ycontindata} there is a $\delta>0$ such that if
both $|y_0-z_0|$, $\|G-\tilde F\|_{L^2(0,\,T,\,V^\prime)}$, $\|P-Q\|_{L^4(0,\,T,\,D(A)}$ and $|\eta-\nu|$
are less than $\delta$, then
$|w|_{C([0,\,T],\,H)}<\sqrt\frac{\varepsilon}{4D_7}$. Hence for
some $\delta_1$ smaller that $\delta$ we have that if both
$|y_0-z_0|$, $\|G-\tilde F\|_{L^2(0,\,T,\,V^\prime)}$, $\|P-Q\|_{L^4(0,\,T,\,D(A)}$ and $|\eta-\nu|$ are
less than $\delta_1$, then
$\Bigl(\int_0^T\|w(t)\|^2\,dt\Bigr)^\frac{1}{2}<\varepsilon$.\\
Therefore the map $\mathbb Y_2$ is continuous.
\end{proof}
\begin{Remark}
As said in Remark \ref{PL2V} considering $P\in L^2(0,\,T,\,V)$ is enough to guarantee the existence of weak solutions but, to have the continuity of $\mathbb Y$ or $\mathbb Y_2$ we have to consider $P$ varying in $L^2(0,\,T,\,D(A))$-norm. Variation in $L^2(0,\,T,\,V)$-norm seems to be not stronger enough.
\end{Remark}
\section{Strong Case.}
Consider the equivalent problems \ref{syleray} and \ref{syAB}:
\begin{Problem}\label{syleray}
Given
\begin{align}
&\quad \tilde F\in L^2(0,\,T,\,H),\quad P\in L^4(0,\,T,\,D(A))\label{syleray.FP}\\
&\And\notag\\
&\quad y_0\in V\label{syleray.y0}\\
&to\, find\notag\\
&\quad y\in L^2(0,\,T,\,D(A))\cap L^\infty(0,\,T,\,V)\label{syleray.y}\\
&\text{satisfying (in the distribution sense)}\notag\\
&\quad \forall v\in V:\notag\\
&\quad \frac{d}{dt}(y,\,v)+\nu((y+P,\,v))+b(y+P,\,y+P,\,v)=(\tilde F,\,v),\label{syleray.nse}\\
&\text{and}\notag\\
&\quad y(0)=y_0\label{syleray.inicond}.
\end{align}
\end{Problem}
\begin{Problem}\label{syAB}
Given
\begin{align}
&\quad \tilde F\in L^2(0,\,T,\,H),,\quad P\in L^4(0,\,T,\,D(A))\label{syAB.FP}\\
&\And\notag\\
&\quad y_0\in V,\label{syAB.y0}\\
&\text{to find}\notag\\
&\quad y\in L^2(0,\,T,\,D(A))\cap L^\infty(0,\,T,\,V),\qquad
y^\prime\in L^2(0,\,T,\,H)\label{syAB.yy'}\\
&\text{satisfying}\notag\\
&\quad y^\prime+\nu A(y+P)+B(y+P)=\tilde F\quad\text{on}\quad]0,\,T[,\label{syAB.nse}\\
&\text{and}\notag\\
&\quad y(0)=y_0\label{syAB.inicond}.
\end{align}
\end{Problem}

The equivalence of these problems follows from
$$
\begin{cases}
&\bullet\;\text{A solution of Problem \ref{syAB} is a solution of
Problem \ref{syleray}};\\
&\bullet\;\text{Given a solution of Problem \ref{syleray}}\\
&--\;\tilde F\in L^2(0,\,T,\,H);\\
&--\;A(y+P)\in L^2(0,\,T,\,H).\\
&\qquad\qquad\text{Indeed $\|A(y+P)\|^2=|y+P|_{[2]}^2$ and $y+P\in L^2(0,\,T,\,D(A))$};\\
&--\;B(y+P)\in L^2(0,\,T,\,H).\\
&\qquad\qquad\text{See Remark \ref{ByPL2H} below.}
\end{cases}
$$
So that $\tilde F-\nu A(y+P)-B(y+P)\in L^2(0,\,T,\,H)\subseteq
L^2(0,\,T,\,V^\prime)$. By \ref{syleray.nse} and $y\in
L^2(0,\,T,\,V)$ we have $y^\prime=\tilde F-\nu A(y+P)-B(y+P)$ a.e., and
$y\in C([0,\,T],\,V^\prime)$ ``a.e.''.
\begin{Remark}\label{ByPL2H}
By definition and trilinearity we have
\begin{align*}
B(y+P)(h)&=b(y+P,\,y+P,\,h)\\
&=b(y,\,y,\,h)+b(y,\,P,\,h)+b(P,\,y,\,h)+b(P,\,P,\,h)
\end{align*}
or
\begin{equation}\label{expByP}
B(y+P)=B(y,\,y)+B(y,\,P)+B(P,\,y)+B(P,\,P).
\end{equation}
From \eqref{E.01.12.0} we obtain\\
\begin{align*}
& |B(y,\,y)|\leq C|y|^\frac{1}{2}\|y\||y|_{[2]}^\frac{1}{2};\\
& |B(y,\,P)|\leq C|y|^\frac{1}{2}\|y\|^\frac{1}{2}\|P\|^\frac{1}{2}|P|_{[2]}^\frac{1}{2};\\
& |B(y,\,P)|\leq C|P|^\frac{1}{2}\|P\|^\frac{1}{2}\|y\|^\frac{1}{2}|y|_{[2]}^\frac{1}{2};\;\text{and}\\
& |B(P,\,P)|\leq C|P|^\frac{1}{2}\|P\||P|_{[2]}^\frac{1}{2}.
\end{align*}
Thus for some constant $C_1$
\begin{align*}
& |B(y,\,y)|^4\leq C_1\|y\|^6|y|_{[2]}^2\,\in L^1(0,\,T,\,\mathbb R);\\
& |B(y,\,P)|^4\leq C_1\|y\|^4|P|_{[2]}^4\,\in L^1(0,\,T,\,\mathbb R);\\
& |B(y,\,P)|^2\leq C_1|P|_{[2]}^2\|y\||y|_{[2]}\,\in L^1(0,\,T,\,\mathbb R);\;\text{and}\\
& |B(P,\,P)|^2\leq C_1|P|_{[2]}^4\,\in L^1(0,\,T,\,\mathbb R).
\end{align*}
So, $B(y,\,y)$ and $B(y,\,P)$ are elements of $L^4(0,\,T,\,H)$ and,  $B(P,\,y)$ and $B(P,\,P)$ are elements of $L^2(0,\,T,\,H)$ thus, by \eqref{expByP}, $B(y+P)$ is a sum of elements of $L^2(0,\,T,\,H)$. 
\end{Remark}
\subsection{Existence.}
\begin{Theorem}\label{ex.solsy}
Given $\tilde F,\;P$ and $u_0$ satisfying \eqref{syAB.FP} and
\eqref{syAB.y0}. There is at least one function $y$ satisfying
\eqref{syAB.yy'}-\eqref{syAB.inicond}.
\end{Theorem}
\begin{proof}
As we have done for Theorem \ref{ex.soly} starting from
approximate solutions we arrive to the conclusions
\eqref{ymbddLinfH} and \eqref{ymbddL2V}.\\
An approximate solution $y^m:=\sum_{\max\{i_1,\,i_2\}\leq m}
y_i^m(t)W_i$ satisfy
$$
((y^m)^\prime,\,v)+\nu(A(y^m+P^m),\,v)+(B(y^m+P^m),\,v)=(\tilde F,\,v),\qquad
v\in V,
$$
setting $y_i^m\bar\imath W_i$ for $v$ and summing up we arrive to
\begin{equation}\label{estsyVDA1}
((y^m)^\prime,\,Ay^m)+\nu(A(y^m+P^m),\,Ay^m)+(B(y^m+P^m),\,Ay^m)=(\tilde F,\,Ay^m).
\end{equation}
We have that
\begin{align*}
((y^m)^\prime,\,Ay^m)&=(((y^m)^\prime,\,y^m))=\frac{1}{2}\frac{d}{dt}\|y^m\|^2;\\
(A(y^m+P^m),\,Ay^m)&=|y^m|_{[2]}^2+((P^m,\,y^m));\\
\end{align*}
and for the term with $B$ we have
\begin{align*}
&(B(y^m+P^m),\,Ay^m)\\
=&b(y^m+P^m,\,y^m+P^m,\,Ay^m)\\
=&b(y^m,\,P^m,\,Ay^m)+b(y^m,\,y^m,\,Ay^m)+b(P^m,\,P^m,\,Ay^m)+b(P^m,\,y^m,\,Ay^m).
\end{align*}
By \eqref{E.02.1.0} we have
\begin{align*}
|b(y^m,\,P^m,\,Ay^m)|&\leq
C|y^m|^\frac{1}{2}|y^m|_{[2]}^\frac{1}{2}\|P\||y^m|_{[2]}
=C|y^m|^\frac{1}{2}\|P\||y^m|_{[2]}^\frac{3}{2}\\
|b(y^m,\,y^m,\,Ay^m)|&\leq
C|y^m|^\frac{1}{2}|y^m|_{[2]}^\frac{1}{2}\|y^m\||y^m|_{[2]}
=C|y^m|^\frac{1}{2}\|y^m\||y^m|_{[2]}^\frac{3}{2}\\
|b(P^m,\,P^m,\,Ay^m)|&\leq
C|P|^\frac{1}{2}|P|_{[2]}^\frac{1}{2}\|P\||y^m|_{[2]}\\
|b(P^m,\,y^m,\,Ay^m)|&\leq
C|P|^\frac{1}{2}|P|_{[2]}^\frac{1}{2}\|y^m\||y^m|_{[2]}.
\end{align*}
From \eqref{estsyVDA1} we obtain
\begin{align*}
\frac{d}{dt}\|y^m\|^2+2\nu|y^m|_{[2]}^2&\leq 2\nu|P|_{[2]}|y^m|_{[2]}+2|\tilde F||y^m|_{[2]};\\
&+2C|y^m|^\frac{1}{2}\|P\||y^m|_{[2]}^\frac{3}{2}
+2C|y^m|^\frac{1}{2}\|y^m\||y^m|_{[2]}^\frac{3}{2}\\
&+2C|P|^\frac{1}{2}|P|_{[2]}^\frac{1}{2}\|P\||y^m|_{[2]}
+2C|P|^\frac{1}{2}|P|_{[2]}^\frac{1}{2}\|y^m\||y^m|_{[2]};
\end{align*}
by Young inequalities
\begin{align}
\frac{d}{dt}\|y^m\|^2+\nu|y^m|_{[2]}^2&\leq C_1|P|_{[2]}^2+C_1|\tilde F|^2;\notag\\
&+C_2|y^m|^2\|P\|^4+C_2|y^m|^2\|y^m\|^4\notag\\
&+C_3|P||P|_{[2]}\|P\|^2+C_3|P||P|_{[2]}\|y^m\|^2.\label{estsyVDA2}
\end{align}
By Gronwall Inequality
\begin{align*}
\|y^m(t)\|^2&\leq \exp\Bigl(\int_0^T
C_2|y^m(s)|^2\|y^m(s)\|^2+ C_4|P(s)|_{[2]}^2\,ds\Bigr)\biggl(\|y(0)\|^2\\
&+\int_0^T
C_1|P(s)|_{[2]}^2+C_1|\tilde F(s)|^2+C_2|y^m(s)|^2\|P(s)\|^4+C_5|P(s)|_{[2]}^4\,ds\biggr).
\end{align*}
By \eqref{ymbddLinfH}, \eqref{ymbddL2V}, $P\in L^4(0,\,T,\,D(A))$
and $\tilde F\in L^2(0,\,T,\,H)$ we have
$$
\|y^m(t)\|^2\leq K_1\qquad\text{for some constant $K_1$}.
$$
Hence
\begin{equation}\label{ymbddLinfV}
y^m\; remains\; in\; a\; bounded\; set\; of\; L^\infty(0,\,T,\,V).
\end{equation}
Integrating \eqref{estsyVDA2} over ${[0,\,T]}$ we obtain
\begin{multline*}
\|y^m(T)\|^2-\|y^m(0)\|^2+\nu\int_0^T|y^m(t)|_{[2]}^2\,dt\leq
\int_0^TC_1|P(t)|_{[2]}^2+C_1|\tilde F(t)|^2\,dt\\
+\int_0^T C_2|y^m|^2\|P\|^4+C_2|y^m|^2\|y^m\|^4+C_5|P|_{[2]}^4+C_6|P|_{[2]}^2\|y^m\|^2\,dt.
\end{multline*}
Again, by \eqref{ymbddLinfH}, \eqref{ymbddL2V}, \eqref{ymbddLinfV},
$P\in L^4(0,\,T,\,D(A))$ and $\tilde F\in L^2(0,\,T,\,H)$ we arrive to
$$
\int_0^T|y^m(t)|_{[2]}^2<K_2\qquad\text{for some constant $K_2$}.
$$
Hence
\begin{equation}\label{ymbddL2DA}
y^m\; remains\; in\; a\; bounded\; set\; of\; L^2(0,\,T,\,D(A)).
\end{equation}
The rest of the proof is analogous to that of Corollary
\ref{ex.w.solv}.
\end{proof}
\begin{Remark}
Contrary to the case of weak solutios (see Remark \ref{PL2V}), we can not guarantee from the previous sketch that considering $P\in L^2(0,\,T,\,V)$ is enough to guarantee the existence of a strong solution.
\end{Remark}
\subsection{Uniqueness.}
\begin{Theorem}\label{cont.yV}
The solution of Problems \ref{syleray}-\ref{syAB} given by Theorem
\ref{ex.solsy} is unique. Moreover it is a.e. equal to a
continuous function from $[0,\,T]$ into $V$.
\end{Theorem}
\begin{proof}
The uniqueness follows from Theorem \ref{cont.yH} and from the
fact that a solution of problems \ref{syleray}-\ref{syAB} is a
solution of problems \ref{yleray}-\ref{yAB}.\\
By $y^\prime\in L^2(0,\,T,\,H)$ we have $y^\prime\in
L^2(0,\,T,\,D(A)^\prime)$. Since $y\in L^2(0,\,T,\,D(A))$ and the
inclusions
$$
D(A)\subseteq V\subseteq D(A)^\prime
$$
are dense and continuous we have that $y\in C([0,\,T],\,V)$. [See Lemma 1.2
in \cite{temams} section III.1].
\end{proof}
\subsection{Continuity.}
\begin{Theorem}\label{sycontindata}
The map
\begin{align*}
\mathbb Y_s:\,V\times L^2(0,\,T,\,H)\times L^4(0,\,T,\,D(A))\times]0,\,+\infty[
&\to C([0,\,T],\,V)\\
(y_0,\,\tilde F,\,P,\,\nu)&\mapsto y
\end{align*}
is continuous. Where $y$ is the unique solution of problems
\eqref{syleray}-\eqref{syAB} corresponding to the data
$(y_0,\,\tilde F,\,P,\,\nu)$.
\end{Theorem}
\begin{proof}
The proof is analogous to that of Theorem \ref{ycontindata}. We
fix a quadruple
$$
(y_0,\,\tilde F,\,P,\,\nu)\in V\times L^2(0,\,T,\,H)\times
L^4(0,\,T,\,D(A))\times]0,\,+\infty[
$$
and $\varepsilon>0$. Then consider another quadruple
$$
(z_0,\,G,\,Q,\,\eta)\in V\times L^2(0,\,T,\,H)\times
L^4(0,\,T,\,D(A))\times]0,\,+\infty[.
$$
Put
$$
y:=\mathbb Y_s(y_0,\,\tilde F,\,P,\,\nu)\And z:=\mathbb
Y_s(z_0,\,G,\,Q,\,\eta)
$$
Putting $w:=z-y$ we obtain
$$
w^\prime=G-\tilde F-\eta Aw -\eta A(Q-P)+(\nu-\eta)A(y+P)-B(z+Q)+B(y+P).
$$
Taking the scalar product with $Aw$ we obtain
\begin{align}
(w^\prime,\,Aw)&=<G-\tilde F,\,Aw>-\eta|w|_{[2]}^2-\eta(Q-P,\,w)_{[2]}\notag\\
&+(\nu-\eta)(y+P,\,w)_{[2]}+b(y+P,\,y+P,\,w)-b(z+Q,\,z+Q,\,w).\label{sn.simp.b}
\end{align}
Now we estimate the last difference:
\begin{align*}
&b(y+P,\,y+P,\,w)-b(z+Q,\,z+Q,\,w)\\
&\qquad=b(y,\,y,\,Aw)+b(y,\,P,\,Aw)+b(P,\,y,\,Aw)+b(P,\,P,\,Aw)\\
&\qquad-b(z,\,z,\,Aw)-b(z,\,Q,\,Aw)-b(Q,\,z,\,Aw)-b(Q,\,Q,\,Aw).
\end{align*}
We put
\begin{align*}
Z_1&:=b(y,\,y,\,Aw)-b(z,\,z,\,Aw)\\
&=b(y,\,y,\,Aw)-\Bigl(b(w,\,z,\,Aw)+b(y,\,z,\,Aw)\Bigr)\\
&=b(y,\,y,\,Aw)-\Bigl(b(w,\,w,\,Aw)+b(w,\,y,\,Aw)+b(y,\,w,\,Aw)+b(y,\,y,\,Aw)\Bigr)\\
&=-b(w,\,w,\,Aw)-b(w,\,y,\,Aw)-b(y,\,w,\,Aw);\\
Z_2&:=b(y,\,P,\,Aw)-b(z,\,Q,\,Aw)\\
&=b(y,\,P,\,Aw)-\Bigl(b(w,\,Q,\,Aw)+b(y,\,Q,\,Aw)\Bigr)\\
&=b(y,\,P,\,Aw)-\Bigl(b(w,\,Q,\,Aw)+b(y,\,Q-P,\,Aw)+b(y,\,P,\,Aw)\Bigr)\\
&=-b(w,\,Q,\,Aw)-b(y,\,Q-P,\,Aw);\\
Z_3&:=b(P,\,y,\,Aw)-b(Q,\,z,\,Aw)\\
&=b(P,\,y,\,Aw)-\Bigl(b(Q,\,w,\,Aw)+b(Q,\,y,\,Aw)\Bigr)\\
&=b(P,\,y,\,Aw)-\Bigl(b(Q,\,w,\,Aw)+b(Q-P,\,y,\,Aw)+b(P,\,y,\,Aw)\Bigr)\\
&=-b(Q,\,w,\,Aw)-b(Q-P,\,y,\,Aw);\\
Z_4&:=b(P,\,P,\,Aw)-b(Q,\,Q,\,Aw)\\
&=b(P-Q,\,P,\,Aw)+\Bigl(b(Q,\,P,\,Aw)-b(Q,\,Q,\,Aw)\Bigr)\\
&=b(P-Q,\,P,\,Aw)+b(Q,\,P-Q,\,Aw);
\end{align*}
We have
\begin{equation}\label{sb<Zs}
|b(y+P,\,y+P,\,w)-b(z+Q,\,z+Q,\,w)|\leq|Z_1|+|Z_2|+|Z_3|+|Z_4|
\end{equation}
and, by \eqref{E.02.1.0}
\begin{align}
|Z_1|&\leq C|w|^\frac{1}{2}|w|_{[2]}^\frac{1}{2}\|w\||w|_{[2]}
+C|w|^\frac{1}{2}|w|_{[2]}^\frac{1}{2}\|y\||w|_{[2]}
+C|y|^\frac{1}{2}|y|_{[2]}^\frac{1}{2}\|w\||w|_{[2]},\label{sB.Z_1}\\
|Z_2|&\leq C|w|^\frac{1}{2}|w|_{[2]}^\frac{1}{2}\|Q\||w|_{[2]}
+C|y|^\frac{1}{2}|y|_{[2]}^\frac{1}{2}\|P-Q\||w|_{[2]},\label{sB.Z_2}\\
|Z_3|&\leq C|Q|^\frac{1}{2}|Q|_{[2]}^\frac{1}{2}\|w\||w|_{[2]}
+C|P-Q|^\frac{1}{2}|P-Q|_{[2]}^\frac{1}{2}\|y\||w|_{[2]},\label{sB.Z_3}\\
|Z_4|&\leq C|P-Q|^\frac{1}{2}|P-Q|_{[2]}^\frac{1}{2}\|P\||w|_{[2]}
+C|Q|^\frac{1}{2}|Q|_{[2]}^\frac{1}{2}\|P-Q\||w|_{[2]}.\label{sB.Z_4}
\end{align}

From \eqref{sb<Zs} and \eqref{sn.simp.b} we obtain
\begin{multline*}
\frac{d}{dt}\|w\|^2\leq
2|G-\tilde F||w|_{[2]}-2\eta|w|_{[2]}^2+2\eta|Q-P|_{[2]}|w|_{[2]}\\
+2|\nu-\eta||y+P|_{[2]}|w|_{[2]}+2\sum_{i=1}^4|Z_i|.
\end{multline*}
If $\eta$ satisfies $|\eta-\nu|<\frac{\nu}{2}$, i.e.,
$\frac{\nu}{2}<\eta<\frac{3\nu}{2}$ we obtain
\begin{multline*}
\frac{d}{dt}\|w\|^2+\nu|w|_{[2]}^2\leq
2|G-\tilde F||w|_{[2]}+3\nu|Q-P|_{[2]}|w|_{[2]}\\
+2|\nu-\eta||y+P|_{[2]}|w|_{[2]}+2\sum_{i=1}^4|Z_i|
\end{multline*}
so, applying Young inequalities, for appropriate
constants (independent of the quadruple $(z_0,\,G,\,Q,\,\eta)$) we have
\begin{multline}\label{syestVDA}
\frac{d}{dt}\|w\|^2+\frac{\nu}{2}|w|_{[2]}^2\leq
C_1|G-\tilde F|^2+C_2|Q-P|_{[2]}^2+C_3|\nu-\eta|^2|y+P|_{[2]}^2\\
+C_4\biggl(|w|^2\|w\|^4+|w|^2\|y\|^4+|y||y|_{[2]}\|w\|^2
+|w|^2\|Q\|^4+|y||y|_{[2]}\|P-Q\|^2\\
+|Q||Q|_{[2]}\|w\|^2+|P-Q||P-Q|_{[2]}\|y\|^2
+|P-Q||P-Q|_{[2]}\|P\|^2+|Q||Q|_{[2]}\|P-Q\|^2\biggr).
\end{multline}

By the Gronwall Inequality:
\begin{align}
&\|w(s)\|^2\notag\\
\leq&\exp\biggl(\int_0^T
C_4\bigl(|w(t)|^2\|w(t)\|^2+|y(t)||y(t)|_{[2]}+|Q(t)||Q(t)|_{[2]}
\bigr)\,dt\biggr)\Biggl(\|w(0)\|^2\notag\\
+&C_5\int_0^T
\Biggl[|Q(t)-P(t)|_{[2]}^2\biggl(1+|y(t)||y(t)|_{[2]}+\|y(t)\|^2
\\+&\|P(t)\|^2+|Q(t)||Q(t)|_{[2]}\biggr)+|G(t)-\tilde F(t)|^2+|\nu-\eta|^2|y(t)+P(t)|_{[2]}^2\notag\\
+&|w(t)|^2\|y(t)\|^4+|w(t)|^2\|Q(t)\|^4\Biggr]\,dt\Biggr).\label{sest.w.Y}
\end{align}
Given a constant $E>0$, by Theorems \ref{ycontindata},
\ref{ycontindata2} and by the existence of a constant $C$ such that\footnote{The existence of such a constant comes from the continuity of the inclusions $V\subseteq H\subseteq V^\prime$.}
\begin{multline*}
|a|\leq C\|a\|\And \|f\|_{L^2(0,\,T,\,V^\prime)}\leq C\|f\|_{L^2(0,\,T,\,H)}\\
\text{for all}\; a\in V,\,f\in L^2(0,\,T,\,H).
\end{multline*}
 there is $\delta>0$ such that if both
$\|y_0-z_0\|$, $\|G-\tilde F\|_{L^2(0,\,T,\,H)}$, $\|P-Q\|_{L^4(0,\,T,\,D(A))}$ and $|\eta-\nu|$ are less than $\delta$ we
have
$$
\|w\|_{C([0,\,T],\,H)}<E \qquad
\text{and}\qquad\int_0^T\|w(t)\|^2\,dt<E.
$$
Therefore, if we choose in addition $\delta<\frac{\nu}{2}$, we have

\begin{multline*}
\|w(s)\|^2\leq D_1\Biggl[\|w(0)\|^2+\|G-\tilde F\|_{L^2(0,\,T,\,H)}^2\\+\int_0^T|Q(t)-P(t)|_{[2]}^2\Bigl(1+|y(t)|_{[2]}
+|P(t)|_{[2]}^2+|Q(t)|_{[2]}^2\Bigr)\,dt\\
+|\nu-\eta|^2+\|w\|_{C([0,\,T],\,H)}^2(1+\|Q\|_{L^4(0,\,T,\,D(A))}^4)\Biggr]
\end{multline*}
and, because
\begin{align*}
&\int_0^T|Q(t)-P(t)|_{[2]}^2\Bigl(1+|y(t)|_{[2]} +|P(t)|_{[2]}^2+|Q(t)|_{[2]}^2\Bigr)\,dt\\
\leq&\|Q-P\|_{L^4(0,\,T,\,D(A))}^2\Biggl(\|1\|_{L^2(0,\,T,\,\mathbb R)}+
\|y\|_{L^2(0,\,T,\,D(A))}\\
+&\|P\|_{L^4(0,\,T,\,D(A))}^2+\|Q\|_{L^4(0,\,T,\,D(A))}^2\Biggr)\\
\leq&K_1\|Q-P\|_{L^4(0,\,T,\,D(A))}^2\Biggl(1+\Bigl(\delta+\|P\|_{L^4(0,\,T,\,D(A))}\Bigr)^2\Biggr)\\
\leq&K_2\|Q-P\|_{L^4(0,\,T,\,D(A))}^2.
\end{align*}
we have
\begin{multline}
\|w(s)\|^2\\
\leq K_3\Bigl[\|w(0)\|^2+\|G-\tilde F\|_{L^2(0,\,T,\,H)}^2+\|Q-P\|_{L^4(0,\,T,\,D(A))}^2
+|\nu-\eta|^2+\|w\|_{C([0,\,T],\,H)}^2\Bigr]
\end{multline}
(with $K_3$ independent of the quadruple $(z_0,\,G,\,Q,\,\eta)$).\\
Then for some $\delta_1$ smaller than $\delta$, and using Theorem \ref{ycontindata}, we have
$\|w(s)\|<\varepsilon$ if both $\|y_0-z_0\|$, $\|\tilde F-G\|_{L^2(0,\,T,\,H)}$, $\|P-Q\|_{L^4(0,\,T,\,D(A))}$ and
$|\eta-\nu|$ are less than $\delta_1$. Thus the map $\mathbb Y_s$
is continuous.
\end{proof}
\begin{Theorem}\label{sycontindata2}
The map
\begin{align*}
\mathbb Y_{2s}:\,V\times L^2(0,\,T,\,H)\times
L^4(0,\,T,\,D(A))\times]0,\,+\infty[&
\to L^2([0,\,T],\,D(A))\\
(y_0,\,\tilde F,\,P,\,\nu)&\mapsto y
\end{align*}
is continuous. Where $y$ is the unique solution of problems
\eqref{syleray}-\eqref{syAB} corresponding to the data
$(y_0,\,\tilde F,\,P,\,\nu)$.
\end{Theorem}
\begin{proof}
We fix $(y_0,\,\tilde F,\,P,\,\nu)\in V\times
L^2(0,\,T,\,H)\times L^4(0,\,T,\,D(A))\times]0,\,+\infty[$, $\varepsilon>0$ and, put $y:=\mathbb Y_{2s}(y_0,\,\tilde F,\,P,\,\nu)$.\\
Let $(z_0,\,G,\,Q,\,\eta)$ be another quadruple in the product space $V\times
L^2(0,\,T,\,H)\times L^4(0,\,T,\,D(A))\times]0,\,+\infty[$. Put
$z:=\mathbb Y_{2s}(z_0,\,G,\,Q,\,\eta)$ and $w:=z-y$ and like in
the proof of Theorem \ref{sycontindata} we arrive to
\eqref{syestVDA}:
\begin{multline*}
\frac{d}{dt}\|w\|^2+\frac{\nu}{2}|w|_{[2]}^2\leq
C_1|G-\tilde F|^2+C_2|Q-P|_{[2]}^2+C_3|\nu-\eta|^2|y+P|_{[2]}^2\\
+C_4\Biggl(|w|^2\|w\|^4+|w|^2\|y\|^4+|y||y|_{[2]}\|w\|^2
+|w|^2\|Q\|^4\\+|y||y|_{[2]}\|P-Q\|^2+|Q||Q|_{[2]}\|w\|^2+|P-Q||P-Q|_{[2]}\|y\|^2\\
+|P-Q||P-Q|_{[2]}\|P\|^2+|Q||Q|_{[2]}\|P-Q\|^2\Biggr).
\end{multline*}
or, 
\begin{multline*}
\frac{d}{dt}\|w\|^2+\frac{\nu}{2}|w|_{[2]}^2\leq
C_1|G-\tilde F|^2+C_2|Q-P|_{[2]}^2+C_3|\nu-\eta|^2|y+P|_{[2]}^2\\
+C_5\Biggl(\|w\|^6+\|w\|^2+|y|_{[2]}^2\|w\|^2
+\|w\|^2\|Q\|^4+|Q|_{[2]}^2\|w\|^2\\
+|P-Q|_{[2]}^2\Bigl(1+|y|_{[2]}+|P|_{[2]}^2+|Q|_{[2]}^2\Bigr)
\Biggr)
\end{multline*}
(where the constants are independent of the quadruple $(z_0,\,G,\,Q,\,\eta)$).\\
Integrating over $[0,\,T]$
\begin{align*}
\int_0^T|w(t)|_{[2]}^2\,dt&\leq
D\Biggl[\|w\|_{C([0,\,T],\,V)}^2+\|\tilde F-G\|_{L^2(0,\,T,\,H)}^2+\|Q-P\|_{L^4(0,\,T,\,D(A))}^2\\
&+|\nu-\eta|^2+\|w\|_{C([0,\,T],\,V)}^2\Bigl(\|w\|_{C([0,\,T],\,V)}^4+1+\|y\|_{L^2(0,\,T,\,D(A))}^2\\
&+\|Q\|_{L^4(0,\,T,\,D(A))}^4+\|Q\|_{L^4(0,\,T,\,D(A))}^2\Bigr)\\
&+\|P-Q\|_{L^4(0,\,T,\,D(A))}^2\Bigl(1+\|y\|_{L^2(0,\,T,\,D(A)}\\
&+\|P\|_{L^4(0,\,T,\,D(A))}^2+\|Q\|_{L^4(0,\,T,\,D(A))}^2\Bigr)\Biggr]
\end{align*}
By Theorem \ref{sycontindata}, given a constant $E>0$, there is a
$\delta>0$ such that, if both $\|y_0-z_0\|$,
$\|P-Q\|_{L^4(0,\,T,\,D(A))}$, $\|\tilde F-G\|_{L^2(0,\,T,\,H)}$, and $|\eta-\nu|$ are less than $\delta$ we have
$$
\|w\|_{C([0,\,T],\,V)}<E.
$$
Thus, 
\begin{align*}
\int_0^T|w(t)|_{[2]}^2\,dt&\leq D_1\Biggl[\|\tilde F-G\|_{L^2(0,\,T,\,H)}^2+\|Q-P\|_{L^4(0,\,T,\,D(A))}^2
+|\nu-\eta|^2\\
&+\|w\|_{C([0,\,T],\,V)}^2\Bigl(\|w\|_{C([0,\,T],\,V)}^4+1+\|y\|_{L^2(0,\,T,\,D(A))}^2\\
&+\Bigl(\delta+\|P\|_{L^4(0,\,T,\,D(A))}\Bigr)^4+\Bigl(\delta+\|P\|_{L^4(0,\,T,\,D(A))}\Bigr)^2\\
&+\|P-Q\|_{L^4(0,\,T,\,D(A))}^2\Bigl(1+\|y\|_{L^2(0,\,T,\,D(A)}\\&+\|P\|_{L^4(0,\,T,\,D(A))}^2+\bigl(\delta+\|P\|_{L^4(0,\,T,\,D(A))}\bigr)^2\Bigr)\Biggr]
\end{align*}
or
\begin{multline*}
\int_0^T|w(t)|_{[2]}^2\,dt\leq D_2\Biggl[\|\tilde F-G\|_{L^2(0,\,T,\,H)}^2+\|Q-P\|_{L^4(0,\,T,\,D(A))}^2\\+
|\nu-\eta|^2+\|w\|_{C([0,\,T],\,V)}^2\Biggr]
\end{multline*}
(with $D_2$ independent of the quadruple $(z_0,\,G,\,Q,\,\eta)$).\\
Hence for some $\delta_1$ smaller than $\delta$ we have
$\int_0^T|w|_{[2]}^2<\varepsilon^2$ if both $\|y_0-z_0\|$, $\|\tilde F-G\|_{L^2(0,\,T,\,H)}$, $\|P-Q\|_{L^4(0,\,T,\,D(A))}$
and $|\eta-\nu|$ are less than $\delta_1$. Thus the map $\mathbb
Y_{2s}$ is continuous.
\end{proof}
\section{Continuity in Relaxation Metric.}
We begin with a definition:
\begin{Definition}
The {\bf relaxation metric} in $L^1([0,\,T],\,\mathbb R^d)$ is defined by the norm
\begin{equation}\label{rxm}
\|g\|_{rx}:=\max_{t_1,\,t_2\in [0,\,T]}\biggl\|\int_{t_1}^{t_2} g(\tau)\,d\tau\biggr\|_{\mathbb R^d},
\end{equation}
where $\|\cdot\|_{\mathbb R^d}$ is the $l_1$ norm in $\mathbb R^d$.
\end{Definition}
If nothing in contrary is stated we consider the spaces $\mathbb R^d\quad (d\in\mathbb N_0$) endowed with $l_1$-norm --- $\|x\|_{\mathbb R^d}=\|x\|_{l_1}:=\sum_{i=1}^d|x_i|$.
\begin{Remark}
It is easy to check that \eqref{rxm} is a semi-norm and, since functions in $L^1([0,\,T],\,\mathbb R^d)$ coinciding on a set of measure $T$ are identified we can conclude that \eqref{rxm} is a norm. 
\end{Remark}

Consider, also, the {\bf w-relaxation metric} on $L^1([0,\,T],\,\mathbb R^d)$ defined by the norm
\begin{equation}\label{wrxm}
\|g\|_{wrx}:=\max_{t\in [0,\,T]}\biggl\|\int_0^t g(\tau)\,d\tau\biggr\|_{\mathbb R^d}.
\end{equation}
\begin{Remark}\label{wrx.rx}
The identity map
$$
\Bigl(L^1([0,\,T],\,\mathbb R^d),\,\|\cdot\|_{rx}\Bigr)\mapsto\Bigl(L^1([0,\,T],\,\mathbb R^d),\,\|\cdot\|_{wrx}\Bigr)
$$
is continuous.
\end{Remark}

Now we fix a finite subset $\mathbb F\subset\mathbb N_0^2$, put $\mathcal F:=\#\mathbb F$ and define the space $S_\mathbb F:=span\{W_k\mid\,k\in\mathbb F\}$. From the fact that the dimension of $S_\mathbb F$ is finite, we conclude the equivalence of all norms defined in it. In particular the norms  $\|\cdot\|_{\mathbb R^{\mathcal F}}$ and $|\cdot|_{[2]}$ are equivalent on  $S_\mathbb F$. From this equivalence we have the following:
\begin{Proposition}\label{Iwrxc.cont}
The map
\begin{align*}
\mathbb I:L_{wrx}^\infty([0,\,T],\,\mathbb R^\mathcal F)&\to C([0,\,T],\,S_\mathbb F)\\
(v_k(t))_{k\in\mathbb F}&\mapsto \sum_{k\in\mathbb F}\mathbb Iv_k(t)W_k
\end{align*}
is continuous. Where $\mathbb Iv_k(t):=\int_0^t v_k(\tau)\,d\tau$ and the subscript ``wrx'' means that we are considering relaxation metric on the set $L^\infty([0,T],\,\mathbb R^{\mathcal F})$.
\end{Proposition}
\begin{proof}
The continuity of the map
\begin{align*}
\mathbb I:L_{wrx}^\infty([0,\,T],\,\mathbb R^\mathcal F)&\to C([0,\,T],\,\mathbb R^\mathcal F)\\
(v_k(t))_{k\in\mathbb F}&\mapsto (\mathbb Iv_k(t))_{k\in\mathbb F}
\end{align*}
is trivial.
\end{proof}
Recall that by definition, the map $\mathbb S$ of Corollary \ref{contindatav} gives us the weak solution, belonging to $C([0,\,T],\,H)$, of the NSE for an initial data in $\Pi:=H\times L^2(0,\,T,\,V^\prime)\times L^\infty([0,\,T],S_\mathbb F)\times\mathbb R^+$. Changing the topology on the third factor of the previous product to the w-relaxation one we arrive to the space $L_{wrx}^\infty([0,\,T],\,\mathbb R^\mathcal F)$ and we define the function $\mathbb S_{wrx}$ as the function defined in the product $\Pi_{wrx}:=H\times L^2(0,\,T,\,V^\prime)\times L_{wrx}^\infty([0,\,T],\mathbb R^{\mathcal F})\times\mathbb R^+$  and taking the same values as $\mathbb S$.
\begin{Proposition}
The map $\mathbb S_{wrx}$ is continuous.
\end{Proposition}
\begin{proof}
As in the beginning of section \ref{S:ch.var} we put $\mathbb Y_{wrx}:=\mathbb S_{wrx}-\mathbb I_\circ$. But here we want an equality between functions defined on $\Pi_{wrx}$ so, we just put $\mathbb I_\circ(u_0,\,\tilde F,\,v,\,\nu):=\mathbb Iv.$\\
By Proposition \ref{Iwrxc.cont} and Theorem \ref{ycontindata} the map 
\begin{align*}
\mathbb Y_{wrx}:\Pi_{wrx}&\to C([0,\,T],\,H)\\
(u_0,\,\tilde F,\,v,\,\nu)&\mapsto \mathbb Y(u_0,\,\tilde F,\,\mathbb Iv,\,\nu)=\mathbb Y\circ \mathbb I^\circ(u_0,\,\tilde F,\,v,\,\nu)
\end{align*}
is continuous. Where $\mathbb I^\circ(u_0,\,\tilde F,\,v,\,\nu):=(u_0,\,\tilde F,\,\mathbb Iv,\,\nu)$.\\By the equality $\mathbb S_{wrx}=\mathbb Y_{wrx}+\mathbb I_\circ$ we conclude the continuity of $\mathbb S_{wrx}$.  
\end{proof}
Analogously, using Proposition \ref{Iwrxc.cont} and theorems \ref{ycontindata2}, \ref{sycontindata} and \ref{sycontindata2}, we can prove the continuity on relaxation metric of the maps $\mathbb S_{2}$, $\mathbb S_{s}$ and $\mathbb S_{2s}$ arriving to the Proposition
\begin{Proposition}
The maps $\mathbb S_{wrx}$, $\mathbb S_{2wrx}$, $\mathbb S_{swrx}$ $\mathbb S_{2swrx}$ are all continuous.
\end{Proposition}
By Remark \ref{wrx.rx} we obtain
\begin{Corollary}\label{NSsol.cont.rx}
The maps $\mathbb S_{rx}$, $\mathbb S_{2rx}$, $\mathbb S_{srx}$ $\mathbb S_{2srx}$ are all continuous.\quad\footnotemark
\end{Corollary}
\footnotetext{These ``rx''-maps are defined similarly as the ``wrx'' ones, just considering the ``rx''-topology in the factor of essentially bounded functions.}
\begin{Remark}
If, instead of a subspace spanned by a finite number of eigenfunctions, we consider any (with either finite or infinite dimension) subspace $\mathbb F\subseteq D(A)$ and, instead of defining the relaxation metric using the $l_1$-norm we define it on $L^1(0,\,T,\,\mathbb F)$ using the norm induced by $D(A)$:
$$
\|g\|_{rx[2]}:=\max_{t_1,\,t_2\in [0,\,T]}\biggl\|\int_{t_1}^{t_2} g(\tau)\,d\tau\biggr\|_{[2]},
$$
and, its corresponding weak form by
$$
\|g\|_{wrx[2]}:=\max_{t\in [0,\,T]}\biggl\|\int_0^t g(\tau)\,d\tau\biggr\|_{[2]}
$$
we still to have the continuity of the corresponding maps  $\mathbb S_{rx[2]}$, $\mathbb S_{2rx[2]}$, $\mathbb S_{srx[2]}$ $\mathbb S_{2srx[2]}$, i.e., the solutions of the NSE vary continuously when the control varies continuously in $rx[2]$-metric.\\
The reason why we have considered a finite space is that we are interested in the case where the control in a space spanned by a finite number of eigenfunctions and, the reason to consider $l_1$-metric in the definition of relaxation metric is that it will be convenient later.
\end{Remark}
\section{A Remark on Finite System versus Infinite System.}
From now we consider only strong solutions and, the external force $F$ and the coefficient of viscosity $\nu$ are fixed.
\begin{Definition}
Given $T>0$, we say that the N-S system is {\bf time-$T$  approximately controllable in observed component} if for any $\tilde\phi\in V$ and any finite subset of modes $\mathcal O\subseteq \mathbb N_0^2$, the projection of the closure of the attainable set at time $T$ from $\tilde\phi\in V$ onto $span\{W_k\mid k\in\mathcal O\}$ is surjective. It is {\bf  time-$T$ controllable in observed component} if the projection of the attainable set at time $T$ from $\tilde\phi\in V$ on $span\{W_k\mid k\in\mathcal O\}$ is surjective.
\end{Definition}
In this section we ask ourselves if can we prove approximate controllability in observed component similarly as we have proved approximated controllability of Galerkin approximations.\par
As we have seen in \eqref{infsys} from the N-S equation
$$
u_t+\nu Au+P^\nabla Bu=F+v\quad\footnotemark
$$
\footnotetext{Recall that $P^\nabla$ is the projection on the space of divergence free functions.}
we can derive the infinite ODE system
\begin{align}
\dot u_k&=\sum_{\begin{subarray}{l}m,n\in\mathbb N_0^2\\m<n\\(n(++)m)^+=k\end{subarray}}-\frac{u_mu_nC_{m,n}^\wedge}{\bar k}(\bar n-\bar m)\notag\\
&+\sum_{\begin{subarray}{l}m,n\in\mathbb N_0^2\\m<n\\(n(--)m)^+=k\end{subarray}}\frac{u_mu_nC_{m,n}^\wedge}{\bar k}(\bar n-\bar m)\mathrm{sign}(n_1-m_1)\mathrm{sign}(n_2-m_2)\notag\\
&+\sum_{\begin{subarray}{l}m,n\in\mathbb N_0^2\\m<n\\(n(-+)m)^+=k\end{subarray}}-\frac{u_mu_nC_{m,n}^\vee}{\bar k}(\bar n-\bar m)\mathrm{sign}{n_1-m_1}\notag\\
&+\sum_{\begin{subarray}{l}m,n\in\mathbb N_0^2\\m<n\\(n(+-)m)^+=k\end{subarray}}\frac{u_mu_nC_{m,n}^\vee}{\bar k}(\bar n-\bar m)\mathrm{sign}(n_2-m_2)\notag\\
&+\nu\bar ku_k+F_k+v_k.\label{infsysu1}
\end{align}
Now consider the subset $\mathcal K^1=\{(n_1,\,n_2)\in\mathbb N_0^2\mid\,n_1,\,n_2\leq 3\}\setminus\{(3,\,3)\}$ as the set of forced modes. Put $\kappa_1:=\#\mathcal K^1$.\\
If we put $y:=u-\mathbb Iv$ we can write $u$ in the form
\begin{equation}\label{u=yIv}
u=y+\mathbb Iv
\end{equation}
and $y$ satisfies
$$
y_t=-\nu A(y+\mathbb Iv)-P^\nabla B(y+\mathbb Iv)+F.
$$
Writing system \eqref{infsysu1} (with $\mathcal K^1$ as set of controlled modes) in the form
\begin{equation}\label{infsysu}
\begin{cases}
\dot u_k=\mathcal B_k(u)+\nu\bar ku_k+F_k+v_k\quad&k\in\mathcal K^1\\
\dot u_k=\mathcal B_k(u)+\nu\bar ku_k+F_k\quad&k\notin\mathcal K^1
\end{cases}
\end{equation}
in the new variable ``$y$'' the system becomes
\begin{equation}\label{infsysy}
\dot y_k=\mathcal B_k(y+\mathbb Iv)+\nu \bar k(y_k+\mathbb Iv_k)+F_k\qquad v_k=0\; \text{if}\;k\notin\mathcal K^1.
\end{equation}

For any time $T>0$ the closure of the attainable set at time $T$ of system \eqref{infsysu}--- $\overline{\mathcal A_{u_0}(u)(T)}$ --- and the closure of the attainable set at time $T$ of system \eqref{infsysy} ---$\overline{\mathcal A_{u_0}(y)(T)}$ ---\footnote{Note that, by \eqref{u=yIv}, at initial time we have $u_0=y_0$.} are related by
$$
\overline{\mathcal A_{u_0}(u)(T)}=\overline{\mathcal A_{u_0}(y)(T)\times\mathbb R^{\kappa_1}}=\overline{\mathcal A_{u_0}(y)(T)}\times\mathbb R^{\kappa_1}.
$$
Indeed the inclusion $\overline{\mathcal A_{u_0}(u)(T)}\subseteq\overline{\mathcal A_{u_0}(y)(T)\times\mathbb R^{\kappa_1}}$ follows from \eqref{u=yIv} and the reverse one follows from the density of the map
$$
L^\infty([0,\,T],\,\mathbb R^{\kappa_1})\ni v\mapsto(\mathbb Iv,\,\mathbb Iv(T))\in L^4(0,\,T,\,\mathbb R^{\kappa_1})\times\mathbb R^{\kappa_1}
$$
and from the continuity of the map
$$
L^4(0,\,T,\,\mathbb R^{\kappa_1})\ni P\mapsto\mathbb Y_s(u_0,\,F,\,P,\,\nu)(T)\in V
$$
for fixed $u_0$ that is a consequence of Theorem \ref{sycontindata}. Similarly if we rewrite the NS Equation as
$$
u_t+\nu Au+P^\nabla Bu=F+v^1+v^2
$$
we obtain
$$
\overline{\mathcal A_{u_0}(u)(T)}=\overline{\mathcal A_{u_0}(y^1)(T)}\times\mathbb R^{\kappa_1}.
$$
where $\mathcal A_{u_0}(y^1)(T)$ stays for the attainable set at time $T$ from $u_0$ of the system
$$
y_t=-\nu A(y+\mathbb Iv^2)-P^\nabla B(y+\mathbb Iv^2)+F+v^1.
$$
So Factorization works like in the finite dimensional case.\\
In the step of convexification a problem arises: In section \ref{looksat} we obtained the new vector fields $\gamma_{m,n}$ from the vector fields $v_{m,n}^\lambda$ and $w_{m,n}^\lambda$ and from the expressions
\begin{multline*}
\frac{f_{v_{m,n}^\lambda}(u)+f_{-v_{m,n}^\lambda}(u)}{2}=f(u)+\lambda\gamma_{m,n}\\
\text{and}\\
\frac{f_{w_{m,n}^\lambda}(u)+f_{-w_{m,n}^\lambda}(u)}{2}=f(u)-\lambda\gamma_{m,n},
\end{multline*}
i.e., we have extracted a vector field from the convexification of the vector fields $f_{\pm v_{m,n}^\lambda}$ and $f_{\pm w_{m,n}^\lambda}$ and, in the finite dimensional case we know that we can convex without  changing the closure of attainable set at time $t$. From Factorization we know that we can follow the vector fields $f_{\pm v_{m,n}^\lambda}$ and $f_{\pm w_{m,n}^\lambda}$ without changing closure of attainable set but, in the infinite dimensional case we do not know if we can convex these new vector fields without changing the closure of attainable set.\par
If using the vector fields $f(u)+\lambda\gamma_{m,n}$ and $f(u)-\lambda\gamma_{m,n}$ we would not change closure of attainable set then we could add these new directions to the old ones --- $span(\mathcal K^1)$ --- without changing the closure of attainable set.\par
 It is known that any control $v_0\in span\{\delta_{m,n},\,e_k\mid\,(m,n)\in S_1,\,k\in\mathcal K^1\}$\;\footnote{See \eqref{S1.extr} for the definition of $S_1$.} can be approximated in relaxation metric by controls taking values on  $\{\alpha\delta_{m,n},\,\alpha e_k\mid\,(m,n)\in S_1,\,k\in\mathcal K^1\}$ where $\alpha\in\mathbb R$ is a positive constant depending on $v_0$ (as we will see in Lemma \ref{app.lemma} below). By the continuity of $\mathbb S_{srx}$ we do not change the closure of attainable set using controls in $span\{\delta_{m,n},\,e_k\mid\,(m,n)\in S_1,\,k\in\mathcal K^1\}$. We have just applied a step of Convexification. Note that this last procedure of Convexification is quite different from that after Factorization Procedure that is more complicated: For example setting $m:=(1,\,1),\,n:=(2,\,1)$ and considering the vector field $f_{v_{m,n}^\lambda}(u)=f(u)+\lambda\gamma_{m,n}+\mathcal V_1(u,\,v_{m,n}^\lambda)$ (see \eqref{V012}), we see that the candidate to new control --- $\lambda\gamma_{m,n}+\mathcal V_1(u,\,v_{m,n}^\lambda)$, contrary to what happens in the last step of Convexification depends on $u$ and does not take values on a compact subset, indeed if we put $u^r:=rW_m\quad r\in\mathbb R$ we obtain
\begin{align*}
&\lambda\gamma_{m,n}+\mathcal V_1(u^r,\,v_{m,n}^\lambda)\\
=&\lambda\gamma_{m,n}-\lambda r \bigl[P^\nabla B(W_n,\,W_m)+P^\nabla B(W_m,\,W_n)\bigr]+r\nu\bar m.
\end{align*}
Since $\lambda\gamma_{m,n}-\lambda r \bigl[P^\nabla B(W_m,\,W_n)+P^\nabla B(W_m,\,W_n)\bigr]$ takes values  on
$$
span\{W_{(1,\,2)},\,W_{(3,\,2)}\}
$$
we obtain $r\nu\overline{(1,\,1)}$ for the projection of $\lambda\gamma_{m,n}+\mathcal V_1(u,\,v_{m,n}^\lambda)$  onto $span\{W_{(1,\,1)}\}$. So, the projection goes to $\infty$ as $r$ does.\\
Without compactness we are not able to use the Approximation Lemma \ref{app.lemma} below.\par
In the case of Galerkin approximations from approximate controllability we could derive (exact) controllability using the bracket generating property of the system but, in the infinite case we do not have a property like that so, we would not be able to conclude controllability in observed component immediately from approximate controllability (in the case we could somehow prove approximate controllability).\\
In the next Section, using some more tools, we prove the so called {\em solid controllability in observed component} for system \eqref{infsysu} which implies controllability in observed component. 
\section{Solid Controllability in Observed Component}
\begin{Definition}
Let $\phi^0:\,M^1\to M^2$ be a continuous map between two finite dimensional $C^0$-manifolds, $\Omega\subset M^1$ be an open subset with compact closure and, $S\subseteq M^2$ be any subset. We say that $\phi^0(\Omega)$ {\bf covers} $S$ {\bf solidly}, if for some $C^0$-neighborhood $\mathcal N$ of $\phi^0\mid_{\overline{\Omega}}$ there holds: $S\subseteq\phi(\Omega)$. 
\end{Definition}
Let $\mathcal O\subset\mathbb N_0^2$ be the finite set of modes we want to observe and, $\Pi_\mathcal O$ be the projection map from $V$ onto $span\{W_k\mid\,k\in\mathcal O\}$. Define, for each $T>0$ and each finite subset $\mathbb F\subset\mathbb N_0^2$, the ``end point'' map
\begin{align*}
\mathbb E_T:\,V\times L^\infty([0,\,T],\,\mathbb R^{\#\mathbb F})&\to \mathcal O\\
(u_0,\,v)&\mapsto \Pi_\mathcal O\circ\mathbb S_s(u_0,\,F,\,v,\,\nu)(T).
\end{align*}
For any $N\in\mathbb N_0$ define, also, the system
\begin{equation}\label{sysN}
N:
\begin{cases}
\dot u_k=\mathcal B_k(u)+\nu\bar ku_k+F_k+v_k;&\quad k\in\mathcal K^N\\
\dot u_k=\mathcal B_k(u)+\nu\bar ku_k+F_k;&\quad k\notin\mathcal K^N.
\end{cases}
\end{equation}
that is the same as system \eqref{infsysu} with $\mathcal K^N$ as the finite set of controlled modes.\par
\begin{Definition}
 We shall say that system $[\eqref{sysN}.N]$ is {\bf time-$T$ solidly controllable in observed component} if for any $u_0\in V$ and $R>0$ there exists a family
$$
\mathcal V_{u_0,R}:=\{v_b\in L^\infty([0,\,T],\mathbb R^{\kappa_N})\mid\,b\in B_{u_0,R}\}
$$
such that $\mathbb E_T(u_0,B_{u_0,R}):=\mathbb E_T(u_0,\mathcal V_{u_0,R})$ covers $\overline{\mathcal O}_R(u_0^{\#\mathcal O})$ solidly. Where, by $y^{\#\mathcal O}$ we mean the projection of $y$ onto $\mathbb R^{\#\mathcal O}=\mathcal O$; $B_{u_0,R}$ is an open relatively compact subset of a $C^0$-manifold and; $\overline{{\mathcal O}}_R(y)$ is the closed ball
$$
\{x\in {\mathcal O}\mid\,\|x-y\|_{l_1}\leq R\}:=\{x\in\mathbb R^{\#{\mathcal O}}\mid\,\|x-y\|_{l_1}\leq R\}.
$$
Below we will need also open balls we define them by
$$
{\mathcal O}_R(y)=
\{x\in {\mathcal O}\mid\,\|x-y\|_{l_1}< R\}:=\{x\in\mathbb R^{\#{\mathcal O}}\mid\,\|x-y\|_{l_1}< R\}.
$$
\end{Definition}
\begin{Proposition}\label{solcop1}
System $[\eqref{sysN}.1]$ is time-T solidly controllable in observed component.
\end{Proposition}
\begin{Remark}
Proposition \ref{solcop1} implies controllability in observed component and, it follows from Proposition \ref{solcopstN} (with $N=1$) below. Indeed given $R>0$ and $u_0\in V$, if $T\leq T^0$ it is included in the statement of Proposition \ref{solcopstN} (with $N=1$), otherwise if $T>T^0$ we apply any control $\bar v\in L^\infty([0,\,T],\mathbb R^{\kappa_1})$ (for example $\bar v=0$ --- no control) up to time $T-T^0$ arriving to some point $y\in V$. Put $\bar R:=R+\|y^{\kappa_1}-u_0^{\kappa_1}\|$. Then apply first part of Proposition \ref{solcopstN} (with $N=1$ and $T=T^0$) to the pair $(y,\bar R)\in V\times ]0,\,+\infty[$. The family $\mathcal V_{y,\bar R}\circ\bar v$ will do.
\end{Remark}
\begin{Proposition}\label{solcopstN}\ \\
\begin{enumerate}
\item For some $T^0>0$, every $0<T\leq T^0$ and every $N\in\mathbb N_0$ the system $[\eqref{sysN}.N]$ is time-T solid controllable in observed component;\\
\item For each pair $(u_0,\,R)\in V\times [0,\,+\infty[$ the family
$$\mathcal V_{u_0,R}:=\{v_b\mid\,b\in B_{u_0,R}\}$$
can be chosen satisfying:\\
\begin{itemize}
\item The map $b\mapsto v_b$ is $(B,\,L^2(0,\,T,\,\mathbb R^{\kappa_N}))$-continuous and;\\
\item The controls $v_b(t)$ are uniformly (w.r.t. $b$ and $t$) $l_1$-bounded:
$$
\|v_b(t)\|_{l_1}\leq A=A(T,R,u_0).
$$
\end{itemize}
\end{enumerate}
\end{Proposition}
Since $\mathcal K^1$ is saturating we have $\mathcal O\subseteq \mathcal K^M$ from some $M\in\mathbb N_0$. Fix $M$ with this property. We shall prove Proposition \ref{solcopstN} in two steps. Prove it in the case $N\geq M$ and prove the ``back-induction'' step `` it holds for $N$ implies it holds for $N-1$'' $(N=2,\,\dots,\,M)$. These steps are the following subsections \ref{S:1step} and \ref{S:2step}.
\subsection{First Step. Proposition \ref{solcopstN}: $N$ Big}\label{S:1step}
In this subsection we shall prove that the statement of Proposition \ref{solcopstN} holds for $N\geq M$.\par
Decompose $u\in V$ as $u=u^{\kappa_N}+U^{\kappa_N}$ where $u^{\kappa_N}:=P^{\kappa_N}u$, i.e., $u^{\kappa_N}$ is the projection of $u$ onto $\mathbb R^{\kappa_N}$. So $U^{\kappa_N}=P^{-\kappa_N}u$ is the projection of $u$ onto $(\mathbb R^{\kappa_N})_V^\bot$.\\
Now write system $[\eqref{sysN}.N]$ as
\begin{align*}
u_t^{\kappa_N}&=P^{\kappa_N}(-\nu Au-P^\nabla Bu+F+v)\\
U_t^{\kappa_N}&=P^{-\kappa_N}(-\nu Au-P^\nabla Bu+F+v),
\end{align*}
i.e.,
\begin{align*}
u_t^{\kappa_N}&=-\nu Au^{\kappa_N}-P^{\kappa_N}P^\nabla Bu+P^{\kappa_N}F+v=:g(u)+v\\
U_t^{\kappa_N}&=-\nu AU^{\kappa_N}-P^{-\kappa_N}P^\nabla Bu+P^{-\kappa_N}F=:G(u).
\end{align*}
Concisely we have
\begin{equation}\label{NSgGN}
u_t^{\kappa_N}=g(u)+v,\qquad U_t^{\kappa_N}=G(u)\qquad t\in[0,\,T].
\end{equation}
Let $(u_0,\,R)$ be an element of $V\times]0,\,+\infty[$. Fix $\gamma>1$. Let $T>0$ be a positive real number.\\
For each $p\in\gamma\mathcal O_R(0)\subseteq\mathbb R^{\kappa_N}$ define, on $[0,\,T]$ the constant control
$$
v_p(t):=T^{-1}p.
$$
Since $\gamma>1$ we have that $\overline{\mathcal O}_R(0)\subseteq\gamma\mathcal O_R(0)$ and, $\int_0^T v_p(t)\,dt=\int_0^T T^{-1}p\,dt=p$.\\
The family $\mathcal V_p:=\{v_p\mid\,p\in\gamma\mathcal O_R(0)\}$ is parametrized continuously in $L^q$-norm ($q>0$), i.e., the map $p\mapsto v_p$ is $(l_1,\,L^q(0,\,T,\,\mathbb R^{\#\mathcal O}))$-continuous. Indeed for $q>0$:
$$
\|T^{-1}p_2-T^{-1}p_1\|_{L^q}^q=\int_0^T T^{-q}\|p_2-p_1\|_{l_1}^q=T^{1-q}\|p_2-p_1\|_{l_1}^q
$$ 
Hence
$$
\|T^{-1}p_2-T^{-1}p_1\|_{L^q}=T^{\frac{1-q}{q}}\|p_2-p_1\|_{l_1}.
$$
For $q=\infty$: $\|T^{-1}p_2-T^{-1}p_1\|_{L^\infty}=T^{-1}\|p_2-p_1\|_{l_1}$.\par
We can also see that $\|v_p(t)\|_{l_1}< T^{-1}\gamma R$. To prove that the family $\mathcal V_p$ is the one we are looking for it remains to check that $\mathbb E_T(u_0,\,\mathcal V_p)$ covers $\overline{\mathcal O}_R(u_0^{\#\mathcal O})$ solidly.\\
By the $(l_1,\,L^2)$-continuity of $p\mapsto v_p$ and $\bigl(L^2,\,C([0,\,T],\,V)\bigr)$-continuity of $v\mapsto \mathbb S_s(u_0,\,F,\,v,\,\nu)$ we conclude the $(l_1,\,l_1)$-continuity of
$$
p\mapsto \Pi_{\mathcal O}\circ\mathbb S_s(u_0,\,F,\,v_p,\,\nu)(T)=\mathbb E_T(u_0,\,v_p).
$$
Rescaling time: $t=T\xi,\quad\xi\in[0,\,1]$. From \eqref{NSgGN} we obtain the system
\begin{equation}\label{NSgGNxi} 
u_\xi^{\kappa_N}=T(g(u)+v)\qquad U_\xi^{\kappa_N}=TG(u)\qquad u(0)=u_0,\,\xi\in[0,\,1].
\end{equation}
which solutions, for $v=v_p$, will be compared with those of the following system:
\begin{equation}\label{NSgGNxi.0} 
y_\xi^{\kappa_N}=p\qquad Y_\xi^{\kappa_N}=0;\qquad y(0)=y_0,\;\xi\in[0,\,1].
\end{equation}
 Put $z=u-y$. Then $z$ satisfies
\begin{equation*} 
z_\xi^{\kappa_N}=T(g(u)+v_p)-p\qquad Z_\xi^{\kappa_N}=TG(u)\qquad\xi\in[0,\,1],
\end{equation*}
i.e.,
\begin{align*}
\frac{1}{T}z_\xi^{\kappa_N}&=-\nu Au^{\kappa_N}-P^{\kappa_N}P^\nabla Bu+P^{\kappa_N}F\\
\frac{1}{T}Z_\xi^{\kappa_N}&=-\nu AU^{\kappa_N}-P^{-\kappa_N}P^\nabla Bu+P^{-\kappa_N}F.
\end{align*}
which is equivalent to
$$
\frac{1}{T}z_\xi=-\nu Au-P^\nabla Bu+F.
$$
Multiplying by $z$ we arrive to
\begin{align*}
\frac{1}{2T}\frac{d}{d\xi}|z|^2
&\leq-\nu((u,\,z))+|(P^\nabla Bu,\,z)|+|F||z|\\
&\leq-\nu\|z\|^2-\nu((y,\,z))+|(Bu,\,z)|+C_1|F|\|z\|
\end{align*}
and, since
\begin{align*}
b(u,\,u,\,z)&=b(y+z,\,y+z,\,z)=b(y+z,\,y,\,z)\\
&=b(y,\,y,\,z)+b(z,\,y,\,z)
\end{align*}
we arrive to
\begin{align*}
\frac{1}{2T}\frac{d}{d\xi}|z|^2
&\leq-\nu\|z\|^2-\nu((y,\,z))+C\|y\|^2\|z\|+C|z|\|z\|\|y\|+C_1|F|\|z\|
\end{align*}
from which we obtain
$$
\frac{1}{T}\frac{d}{d\xi}|z|^2+\nu\|z\|^2\leq C_2\|y\|^2+C_2\|y\|^4+C_2|z|^2\|y\|^2+C_2|F|^2.
$$
By Gronwall inequality:
\begin{align*}
|z(s)|^2&\leq\exp\Bigl\{TC_2\int_0^1\|y(\xi)\|^2\,d\xi\Bigr\}\biggl(|z(0)|^2+TC_2\int_0^1\|y(\xi)\|^2+\|y(\xi)\|^4+|F|^2\,d\xi\biggr)\\
&\leq\exp(T)D_1(|z(0)|^2+TD_2)
\end{align*}
where $D_1$ and $D_2$ depend only on $\gamma,\,R$ and $\|y_0\|$. Indeed
$y(\xi)=y_0+p\xi$ and $\|y_0+p\xi\|\leq\|y_0\|+C\|p\xi\|_{l_1}$ and, $\|p\xi\|_{l_1}<\gamma R$. In particular we have (for fixed $u_0$, $\gamma$ and $R$):
\begin{Corollary}\label{rt.comp.tr}\ \\
\begin{enumerate}
\item If $y_0=u_0$ then $|u-y|\leq [T\exp(T)]^\frac{1}{2}K$\label{rt.comp.tr.ini}\\
\item For ``bounded'' $T$ and $y_0$, say $T\leq T_1$ and $\|y_0-u_0\|<\beta$ we have
$$
|u-y|\leq K\Bigl(|u_0-y_0|^2+1\Bigr)^\frac{1}{2}
$$
with $K$ independent of $T$ ($K$ depends only on $\gamma,\,R$ and $\|y_0\|$).\label{rt.comp.tr.st}
\end{enumerate}
\end{Corollary}
Note that
$$
\mathbb E_T(u_0,\,v_p)=\Pi_{\mathcal O}\circ\Phi^T(u_0,\,p)(1)
$$
where $\Phi^T(u_0,\,p)$ is the solution of system \eqref{NSgGNxi}.\\
Represent by $\Phi^0(u_0,\,p)$ the solution of system \eqref{NSgGNxi.0}. Then $\Pi_{\mathcal O}\circ\Phi^0(u_0,\,p)(1)=u_0^{\#\mathcal O}+p$\par
Now we compute
\begin{align*}
&\|\Pi_{\mathcal O}\circ\Phi^T(u_0,\,p)(1)-\Pi_{\mathcal O}\circ\Phi^0(u_0,\,p)(1)\|_{l_\infty}\\
\leq&\|\Pi_{\mathcal O}\circ\Phi^T(u_0,\,p)(1)-\Pi_{\mathcal O}\circ\Phi^0(u_0,\,p)(1)\|_{l_1}\\
\leq& C\|\Phi^T(u_0,\,p)(1)-\Phi^0(u_0,\,p)(1)\|_H\leq C\|\Phi^T(u_0,\,p)-\Phi^0(u_0,\,p)\|_{C([0,\,1],H)}.
\end{align*}
By item \ref{rt.comp.tr.ini} of of corollary \ref{rt.comp.tr}:
\begin{equation*}
|\Pi_{\mathcal O}\circ\Phi^T(u_0,\,p)(1)-\Pi_{\mathcal O}\circ\Phi^0(u_0,\,p)(1)\|_{l_\infty}\leq CK[T\exp(T)]^\frac{1}{2}
\end{equation*}
where $K$ is independent of $T$ and $p$. Thus
\begin{equation*}
\|\Pi_{\mathcal O}\circ\Phi^T(u_0,\,\cdot)(1)-\Pi_{\mathcal O}\circ\Phi^0(u_0,\,\cdot)(1)\|_{C(\gamma\mathcal O_R(0),\mathbb R_\infty^{\kappa_N})}\leq CK[T\exp(T)]^\frac{1}{2}\quad\footnotemark
\end{equation*}
\footnotetext{Where the subscript ``$\infty$'' means that we are considering the ``$\max$''-norm --- $l_\infty$.}
or, defining in $\gamma\mathcal O_R(0)$:
$$
\mathbb G_T(p):=\mathbb E_T(u_0,p);\quad\mathbb G_0(p):=\Pi_{\mathcal O}\circ\Phi^0(u_0,\,p)(1);
$$
\begin{equation}
\|\mathbb G_T-\mathbb G_0\|_{C(\gamma\mathcal O_R(0),\mathbb R_\infty^{\kappa_N})}\leq CK[T\exp(T)]^\frac{1}{2}.
\end{equation}
Put $T^0=$ the unique solution of $T^0\exp(T^0)=\Bigl(\frac{(\gamma-1)R}{2\#\mathcal OCK}\Bigr)^2$. Then
\begin{equation}\label{GtG0<R,gamma}
\forall T\in ]0,\,T^0]\;\Bigl[\|\mathbb G_T-\mathbb G_0\|_{C(\gamma\mathcal O_R(0),\mathbb R_\infty^{\#\mathcal O})}\leq \frac{R(\gamma-1)}{2\#\mathcal O}\Bigr].
\end{equation}
Note that the map $\mathbb G_0$ is just the restriction of a translation in $\mathbb R^{\#\mathcal O}$ --- $\mathbb G_0(p)=u_0^{\#\mathcal O}+p$ --- restricted to $\gamma\mathcal O_R(0)\equiv\mathcal O_{\gamma R}(0)$. By the Degree Theory\footnote{See, for example, \cite{deg.th}.} we have that
\begin{multline}\label{deg1}
\Bigl[p\notin\partial(\mathbb G_0\mathcal O_{\gamma R}(0))\equiv \partial\mathcal O_{\gamma R}(u_0^{\#\mathcal O})\Bigr]\\
\Rightarrow\Bigl[deg(\mathbb G_0,\,\mathcal O_{\gamma R}(0),\,p)=deg(I_{\gamma R},\,\mathcal O_{\gamma R}(0),\,p-u_0^{\#\mathcal O})=1\Bigr].
\end{multline}
Where $I_{\gamma R}$ is the identity function on $\mathcal O_{\gamma R}(0)$.\par Yet by the Degree Theory we know that for every $p\notin\partial\mathcal O_{\gamma R}(u_0^{\#\mathcal O})$ and every continuous function $\psi:\mathcal O_{\gamma R}(0)\to\mathbb R^{\#\mathcal O}$ such that
\begin{multline}\label{deg2}
\Bigl[\|\psi-\mathbb G_0\|_{C(\mathcal O_{\gamma R}(0),\,\mathbb R_\infty^{\#\mathcal O})}<\|p-\partial\mathcal O_{\gamma R}(u_0^{\#\mathcal O})\|_{l_\infty}\Bigr]\\
\Rightarrow
\Bigl[deg(\psi,\,\mathcal O_{\gamma R}(0),\,p)=deg(\mathbb G_0,\,\mathcal O_{\gamma R}(0),\,p)\Bigr]\footnotemark
\end{multline}
\footnotetext{Where $\|p,\,A\|$ means the distance from the element $p$ to the set $A$, i.e., $\|p,\,A\|:=\inf\{\|p-a\|\,\mid\,a\in A\}$.}
We claim that for $T\leq T^0$ we have that $\mathbb G_T(\mathcal O_{\gamma R}(0))$ covers $\overline{\mathcal O}_R(u_0^{\#\mathcal O})$ solidly. Indeed, given $T\leq T^0$ and a continuous function $\phi:\mathcal O_{\gamma R}(0)\to\mathbb R^{\#\mathcal O}$ such that $\|\phi-\mathbb G_T\|_{C(\mathcal O_{\gamma R}(0),\,\mathbb R_\infty^{\#\mathcal O})}<\frac{(\gamma-1)R}{2\#\mathcal O}$ we have, using \eqref{GtG0<R,gamma}, 
\begin{align*}
\|\phi-\mathbb G_0\|_{C(\mathcal O_{\gamma R}(0),\,\mathbb R_\infty^{\#\mathcal O})}&\leq\|\phi-\mathbb G_T\|_{C(\mathcal O_{\gamma R}(0),\,\mathbb R_\infty^{\#\mathcal O})}+\|\mathbb G_T-\mathbb G_0\|_{C(\mathcal O_{\gamma R}(0),\,\mathbb R_\infty^{\#\mathcal O})}\\
&<\frac{(\gamma-1)R}{2\#\mathcal O}+\frac{(\gamma-1)R}{2\#\mathcal O}=\frac{(\gamma-1)R}{\#\mathcal O}.
\end{align*}
For $w\in\overline{\mathcal O}_R(u_0^{\#\mathcal O})$ we have:
$$
\|w-\partial\mathcal O_{\gamma R}(u_0^{\#\mathcal O})\|_{l_\infty}\geq\frac{1}{\#\mathcal O}\|w-\partial\mathcal O_{\gamma R}(u_0^{\#\mathcal O})\|_{l_1}\geq \frac{(\gamma-1)R}{\#\mathcal O}
$$
so,
$$
\|\phi-\mathbb G_0\|_{C(\mathcal O_{\gamma R}(0),\,\mathbb R_\infty^{\#\mathcal O})}\|<\|w-\partial\mathcal O_{\gamma R}(u_0^{\#\mathcal O})\|_{l_\infty}.
$$
By \eqref{deg1} and \eqref{deg2} we conclude that
$$
deg(\phi,\,\mathcal O_{\gamma R}(0),\,w)=1.
$$
which means, in particular, that the equation $\phi(y)=w$ has a solution on $\mathcal O_{\gamma R}(0)$, i.e., $\phi(\mathcal O_{\gamma R}(0))$ covers $\overline{\mathcal O}_R(u_0^{\#\mathcal O})$. 
\subsection{Second Step. Proposition \ref{solcopstN}: ``Back-Induction''.}\label{S:2step}
In this subsection we ``imitate'' a driving using controls on $\mathbb R^{\kappa_N}$ by a driving using controls on $\mathbb R^{\kappa_{N-1}}$, $N=2,\,\dots,\,M$, $M$ is fixed and satisfies $\mathcal O\subseteq\mathcal K^M$. Both drivings leading to the same projection onto $\mathbb R^{\kappa_{N-1}}$ at final time but, possibly going by paths with projections ``far from each other'' in the middle. The projection onto the orthogonal space $(\mathbb R^{\kappa_N})_V^\bot$ of the paths will be $H$-close to each other so, at time $T$ the two drivings lead to points close in $H$-metric. Hence the end points of the projection onto the finite dimensional observed space $\mathcal O$ are close. Solid controllability will follow from this closeness and (again) from a Degree Theory argument.\\
Such imitation is then the key for the prove that if the system [\eqref{sysN}.N] is solid controllable in observed component then so is system [\eqref{sysN}.N-1].\\
After we prove this ``$N\to N-1$'' step it will be clear, from the fact that [\eqref{sysN}.M] is solid controllable in observed component (see subsection \ref{S:1step}), that system \eqref{infsysu} is solid controllable in observed component, we just note that the systems [\eqref{sysN}.1] and \eqref{infsysu} are the same system.\par
To prove the ``back-induction'' step ``$N\to N-1$'' we shall need some lemmas:\\
For the next Lemma we may consider again the case where the external force depend on time:
\begin{Lemma}\label{NS:q.to.v}
Given:\\
\begin{itemize}
\item A finite subset $J\subset\mathbb N_0^2;\,\mathbb J:=span\{W_k\mid\,k\in J\}$\\
\item A function $q\in W^{1,\infty}([t_i,\,t_f],\,\mathbb J)$, such that $q(t_i)=q_i$\\
\item An element $Q_i\in \mathbb J_V^\bot$, the orthogonal space to $\mathbb J$ in $V$.
\end{itemize}
Then there exists a control $v^J(q,Q_i)\in L^\infty([t_i,\,t_f],\,\mathbb J)$ depending on $q$ and $Q_i$ such that the projection onto $\mathbb J$ of the solution of the NSE
$$
u_t=-\nu Au-P^\nabla Bu+\tilde F+v^J(q,\,Q_i),\qquad u(t_i)=q_i+Q_i
$$ 
equals $q$ on $[t_i,\,t_f]$.\\
Moreover the map $v^J:(q,\,Q_i)\mapsto v^J(q,\,Q_i)$ is $(W^{1,2}\times \mathbb J_V^\bot,\,L^2(t_i,\,t_f,\,\mathbb J)$-continuous.
\end{Lemma}
\begin{proof}
Let $q,\,q_i,\,Q_i$ and $\mathbb J$ be like in the statement of the Lemma. Consider the (non controlled) NSE: $u_t=-\nu Au-P^\nabla Bu+\tilde F$ with initial condition $u(t_i)=q_i+Q_i=:u_i$ and split it into
\begin{align*}
P^J(u_t)&=P^J(-\nu Au-P^\nabla Bu+\tilde F)\\
P^{-J}(u_t)&=P^{-J}(-\nu Au-P^\nabla Bu+\tilde F)\\
u(t_i)&=u_i
\end{align*}
with the same initial condition, where $P^J$ (resp. $P^{-J}$) is the projection $H\to\mathbb J$ (resp. $H\to \mathbb J_H^\bot$). If we put $u^J:=P^Ju$ and $U^J:=P^{-J}u$ we arrive to the systems
\begin{align}
&\begin{cases}
u_t^J&=-\nu Au^J-P^JP^\nabla Bu+P^J\tilde F\\
u^J(t_i)&=q_i
\end{cases}\\
&\begin{cases}
U_t^J&=-\nu AU^J-P^{-J}P^\nabla Bu+P^{-J}\tilde F\\
U^J(t_i)&=Q_i.
\end{cases}\label{NSqtoQ1}
\end{align}
In the system \eqref{NSqtoQ1}, for each $k\in\mathbb J$ replace $u_k(t)$ by $q_k(t)$ arriving, in this way, to the (closed) system
\begin{equation}\label{NSqtoQ}
\begin{cases}
U_t^J=-\nu AU^J-P^{-J}P^\nabla B(U^J+q)+P^{-J}\tilde F\\
U^J(t_i)=Q_i.
\end{cases}
\end{equation}
We can prove existence and uniqueness of a strong solution for this system as we prove existence and uniqueness for the ``full'' equation. We indicate only how to find some estimates: Starting from approximate solutions
\begin{align}
U^{J,L}:=&\sum_{\begin{subarray}{l}k\in\mathcal K^L\setminus\mathbb J\\\mathbb J\subseteq\mathcal K^L\end{subarray}}U_k^{J,L}W_k\notag\\
(U^{J,L}_t,\,W_k)
=&-\nu(AU^{J,L},\,W_k)-(P^{-J}P^\nabla B(U^{J,L}+q),\,W_k)+(P^{-J}\tilde F,\,W_k)\notag\\
=&-\nu((U^{J,L},\,W_k))-(B(U^{J,L}+q),\,W_k)+(\tilde F,\,W_k)\label{UL.Wk}\\
&\qquad\forall k\in\mathcal K^L\setminus\mathbb J\\
U^{J,L}(t_i)&=Q_i^{\kappa_L-J}=\,\text{projection of $Q_i$ onto $\mathbb R^{\kappa_L-J}$}\notag,
\end{align}
 from which we obtain the ODE
\begin{align*}
-\bar k\frac{ab}{4}\dot U^{J,L}_k&=-\nu\bar k^2\frac{ab}{4} U^{J,L}_k-\sum_{m,n\in\mathcal K^L}u_m^{J,L}u_n^{J,L}b(W_m,\,W_n,\,W_k)-\bar k\frac{ab}{4}\tilde F_k\\
U_i^{J,L}&=Q_i^{\kappa_L-J}
\end{align*}
that has a maximal solution defined on $[t_i,t_{max}[$. Now we compute some estimates that, in particular, imply $t_{max}=t_f$:\\
Multiplying, for each $k$, the equation \eqref{UL.Wk} by $U_k^{J,L}$ and summing up we obtain:
$$
\frac{1}{2}\frac{d}{dt}|U^{J,L}|^2\leq -\nu\|U^{J,L}\|^2+|(B(U^{J,L}+q),U^{J,L})|+|\tilde F||U^{J,L}|
$$ 
and, after simplify the term with $B$, using some estimates from section \ref{S:Est.b} and Young inequalities, we obtain
$$
\frac{d}{dt}|U^{J,L}|^2+\nu\|U^{J,L}\|^2\leq C|U^{J,L}|^2\|q\|^2+C\|q\|^4+C|\tilde F|^2$$
from which we can conclude that for every $s\in[t_i,\,t_f]$:
\begin{align*}
|U^{J,L}(s)|^2\leq& \exp\int_{t_i}^{t_f}C\|q(t)\|^2\,dt\biggl(|Q_i|^2+C\int_{t_i}^{t_f}\|q(t)\|^4+|\tilde F(t)|^2\,dt\biggr)\\
\leq& D_1(\|Q_i\|^2+1)
\end{align*}
where $D_1$ can be taken depending only in $\|q\|_{C([t_i,\,t_f],\,\mathbb J)}$\;\footnote{$\tilde F$ is fixed.}. So for a constant $C_1$ depending only in $\|q\|_{C([t_i,\,t_f],\,\mathbb J)}$ and $\|Q_i\|$ we have
\begin{equation}\label{qtoQ.bdH}
\|U^{J,L}\|_{L^\infty(t_i,\,t_f,\,\mathbb J_H^\bot)}\leq C_1.
\end{equation}
We also have
\begin{align*}
\int_{t_i}^{t_f}\|U^{J,L}(t)\|^2\,dt&\leq D_3\|U^{J,L}\|_{L^\infty(t_i,\,t_f,\,\mathbb J_H^\bot)}^2\Bigl(1+\int_{t_i}^{t_f}\|q(t)\|^2\,dt\Bigr)\\
&+D_3\int_{t_i}^{t_f}\|q(t)\|^4+|\tilde F(t)|^2\,dt.
\end{align*}
Hence for a constant $C_2$ depending only in $\|q\|_{C([t_i,\,t_f],\,\mathbb J)}$ and $\|Q_i\|$ we have
\begin{equation}\label{qtoQ.bd2V}
\|U^{J,L}\|_{L^2(t_i,\,t_f,\,\mathbb J_V^\bot)}\leq C_2.
\end{equation}
From \eqref{qtoQ.bdH} and \eqref{qtoQ.bd2V} we have that
\begin{equation}\label{UJLbddLinfHL2V}
\begin{cases}
(U^{J,L})_L\,\text{\emph remains in a bounded subset of}\,L^\infty([t_i,\,t_f],\,H)\\
(U^{J,L})_L\,\text{\emph remains in a bounded subset of}\,L^2(t_i,\,t_f,\,V).
\end{cases}
\end{equation}
Analogously if we multiply, for each $k$, the equation \eqref{UL.Wk} by $-\bar kU_k^{J,L}$ and summing up we obtain:
$$
\frac{1}{2}\frac{d}{dt}\|U^{J,L}\|^2\leq -\nu|U^{J,L}|_{[2]}^2+|(B(U^{J,L}+q),AU^{J,L})|+|\tilde F||U^{J,L}|_{[2]}
$$
and,
$$
\frac{d}{dt}\|U^{J,L}\|^2+\nu|U^{J,L}|_{[2]}^2\leq C|U^{J,L}|^2\|U^{J,L}\|^4+C|q|_{[2]}^4+C|q|_{[2]}^2\|U^{J,L}\|^2+C|\tilde F|^2
$$
from which, using \eqref{qtoQ.bdH} and \eqref{qtoQ.bd2V}, we conclude that
for some constants $C_3$ and $C_4$ depending only in $\|q\|_{C([t_i,\,t_f],\,\mathbb J)}$ and $\|Q_i\|$:
\begin{align}
&\|U^{J,L}\|_{L^\infty(t_i,\,t_f,\,\mathbb J_V^\bot)}\leq C_3;\label{qtoQ.bdV}\\
&\|U^{J,L}\|_{L^2(t_i,\,t_f,\,\mathbb J_{D(A)}^\bot)}\leq C_4.\label{qtoQ.bd2D(A)}
\end{align}
\begin{equation}\label{UJLbddLinfVL2H}
\begin{cases}
(U^{J,L})_L\,\text{\emph remains in a bounded subset of}\,L^\infty([t_i,\,t_f],\,V)\\
(U^{J,L})_L\,\text{\emph remains in a bounded subset of}\,L^2(t_i,\,t_f,\,D(A))
\end{cases}
\end{equation}
 The existence of a strong solution follows some classical compactness theorems.\par
For the uniqueness we consider the difference $w$ of two solutions $V^J=V^J(q,\,Q_i)$ and $U^J=(q,\,Q_i)$ --- $w:=V^J-U^J$. Then from
$$
w_t=-\nu Aw-P^{-J}P^\nabla B(V^J+q)+P^{-J}P^\nabla B(U^J+q),
$$
multiplying by $w$, using some estimates from section \ref{S:Est.b} and appropriate inequalities, we obtain
$$
\frac{d}{dt}|w|^2+\nu\|w\|^2\leq C|w|^2\|U^J+q\|^2.
$$
So, $|w(s)|^2\leq|w(t_i)|^2\exp\int_{t_i}^{t_f}C\|U^J(t)+q(t)\|\,dt=0$. Then a weak solution is unique and, so is the strong one.\par
We have just proved that the map $(q,\,Q_i)\mapsto U^J(q,\,Q_i)$ is well defined. We claim that it is $(L^4(t_i,\,t_f,\,\mathbb J)\times (\mathbb J_V^\bot),\,X)$-continuous, where $X$ is either $L^\infty([t_i,\,t_f],\,\mathbb J_V^\bot)$ or $L^2(t_i,\,t_f,\,\mathbb J_{D(A)}^\bot)$. To prove these continuities we proceed as usually: Fix a pair $(q,\,Q_i)\in (W^{1,\infty}\times \mathbb J_V^\bot)$ and consider another one $(p,\,P_i)$ in the same product space. Define $w:=U^J(q,\,Q_i)-U^J(p,\,P_i)$. Then we obtain the equation for $w$:
$$
\dot w=-\nu Aw-P^{-J}P^\nabla B(U^J(q,\,Q_i)+q)+P^{-J}P^\nabla B(U^J(p,\,P_i)+p).
$$
To simplify the writing we put $Q:=U^J(q,\,Q_i)$ and $P:=U^J(p,\,P_i)$ so,
\begin{equation}\label{for.est.qtoQ}
\dot w=-\nu Aw-P^{-J}P^\nabla B(Q+q)+P^{-J}P^\nabla B(P+p)
\end{equation}
and, multiplying by $w$:
\begin{align*}
\frac{1}{2}\frac{d}{dt}|w|^2\leq&-\nu\|w\|^2+|b(P,\,P,\,w)-b(Q,\,Q,\,w)|+|b(P,\,p,\,w)-b(Q,\,q,\,w)|\\
&+|b(p,\,P,\,w)-b(q,\,Q,\,w)|+|b(p,\,p,\,w)-b(q,\,q,\,w)|\\
\leq &-\nu\|w\|^2+|b(w,\,w,\,Q)|+\Bigl\{|b(-w,\,p,\,w)+b(Q,\,p-q,\,w)|\Bigr\}\\
&+|b(p-q,\,Q,\,w)|+\Bigl\{|b(p-q,\,p,\,w)+b(q,\,p-q,\,w)|\Bigr\}
\end{align*}
from which we obtain
\begin{align}
\frac{d}{dt}|w|^2+\nu\|w\|^2&\leq C|w|^2\|Q\|^2+C\|Q\|^2\|p-q\|^2+C|w|^2\|p\|^2\notag\\
&+C(\|p\|^2+\|q\|^2)\|p-q\|^2.\label{est.cont.qtoQ}
\end{align}
Then by Gronwall Inequality
\begin{align*}
\|w\|_{C([t_i,\,t_f],\,\mathbb J_H^\bot)}^2&\leq \exp\Bigl[C\int_{t_i}^{t_f}\|Q(t)\|^2+\|p(t)\|^2\,dt\Bigr]\biggl(|w(t_i)|^2\\
&+C\int_{t_i}^{t_f}\|p-q\|^2(\|Q(t)\|^2+\|p(t)\|^2+\|q(t)\|^2)\,dt\biggr).
\end{align*}
For $\|p-q\|_{L^4}\leq 1$ we obtain\footnote{Of course we could ask for $\|p-q\|_{L^4}\leq D$ for any $D>0$. What we need is a first bound for $\|p\|_{L^4}$.}
\begin{equation}\label{cont.qtoQ.CH}
\|w\|_{C([t_i,\,t_f],\,\mathbb J_H^\bot)}^2\leq C_1\|Q_i-P_i\|^2+C_1\|p-q\|_{L^4}^2
\end{equation}
and conclude that the map $U^J$ is $(L^4\times \mathbb J_V^\bot,\,C([t_i,\,t_f],\,\mathbb J_H^\bot))$-continuous.\par
From \eqref{est.cont.qtoQ} we obtain
\begin{align*}
\int_{t_i}^{t_f}\|w(t)\|^2\,dt&\leq C_0\|w\|^2_{C([t_i,\,t_f],\,\mathbb J_H^\bot)}\Big(1+\int_{t_i}^{t_f}\|Q(t)\|^2+\|p(t)\|^2\,dt\Bigr)\\&+C_0 \|p-q\|_{L^4}^2\Bigl(\|Q\|_{L^4}^2+\|p\|_{L^4}^2+\|q\|_{L^4}^2\Bigr)\\
&\leq C_2\|w\|_{C([0,\,T],\,\mathbb J_H^\bot)}^2+C_2\|p-q\|_{L^4}^2,
\end{align*}
for $\|p-q\|_{L^4}\leq 1$. Hence by \eqref{cont.qtoQ.CH} we arrive to
\begin{equation}\label{cont.qtoQ.L2V}
\Bigl[\|p-q\|_{L^4}\leq 1\Bigr]\Rightarrow\Bigl[|w|_{L^2(t_i,\,t_f,\,\mathbb J_V^\bot)}^2\leq C_3\|Q_i-P_i\|^2+C_3\|p-q\|_{L^4}^2\Bigr].
\end{equation}
and conclude that the map $U^J$ is $(L^4\times \mathbb J_V^\bot,\,L^2(t_i,\,t_f,\,\mathbb J_V^\bot))$-continuous.\par
Now from equation \eqref{for.est.qtoQ}, multiplying it by $Aw$, we obtain
\begin{multline}
\frac{1}{2}\frac{d}{dt}\|w\|^2+\nu|w|_{[2]}^2\leq |(B(P+p),\,Aw)-(B(Q+q),\,Aw)|\\
\leq|b(P,\,P,\,Aw)-b(Q,\,Q,\,Aw)|+|b(P,\,p,\,Aw)-b(Q,\,q,\,Aw)|\\
+|b(p,\,P,\,Aw)-b(q,\,Q,\,Aw)|+|b(p,\,p,\,Aw)-b(q,\,q,\,Aw)|.\label{nsimp.bs}
\end{multline}
Note that
\begin{align*}
&b(P,\,P,\,Aw)-b(Q,\,Q,\,Aw)&&=b(w,\,w,\,Aw)-b(w,\,Q,\,Aw)-b(Q,\,w,\,Aw)\\
&b(P,\,p,\,Aw)-b(Q,\,q,\,Aw)&&=-b(w,\,p,\,Aw)+b(Q,\,p-q,\,Aw)\\
&b(p,\,P,\,Aw)-b(q,\,Q,\,Aw)&&=-b(p,\,w,\,Aw)+b(p-q,\,Q,\,Aw)\\
&b(p,\,p,\,Aw)-b(q,\,q,\,Aw)&&=b(p-q,\,p,\,Aw)+b(q,\,p-q,\,Aw).
\end{align*}
Hence from \eqref{nsimp.bs}, from the estimates of section \ref{S:Est.b} and from the continuity of the inclusions $D(A)\mapsto V\mapsto H$ we arrive to
\begin{multline*}
\frac{1}{2}\frac{d}{dt}\|w\|^2+\nu|w|_{[2]}^2
\leq C\Bigl(|w|^\frac{1}{2}\|w\||w|_{[2]}^\frac{3}{2}+\|w\|^\frac{1}{2}\|Q\||w|_{[2]}^\frac{3}{2}\Bigr)\\+C\Bigl(\|w\||p|_{[2]}|w|_{[2]}+\|Q\||p-q|_{[2]}|w|_{[2]}\Bigr)
+C\Bigl(|p|_{[2]}\|w\||w|_{[2]}+|p-q|_{[2]}\|Q\||w|_{[2]}\Bigr)\\+C\Bigl(|p-q|_{[2]}\|p\||w|_{[2]}+\|q\||p-q|_{[2]}|w|_{[2]}\Bigr).
\end{multline*}
By appropriate Young inequalities:
\begin{align}
&\frac{d}{dt}\|w\|^2+\nu|w|_{[2]}^2\notag\\
\leq& D\Bigl(|w|^2\|w\|^4+\|w\|^2\|Q\|^4\Bigr)+D\Bigl(\|w\|^2|p|_{[2]}^2+\|Q\|^2|p-q|_{[2]}^2\Bigr)\notag\\
+&D\Bigl(|p|_{[2]}^2\|w\|^2+|p-q|_{[2]}^2\|Q\|^2\Bigr)+D\Bigl(|p-q|_{[2]}^2\|p\|^2+\|q\|^2|p-q|_{[2]}^2\Bigr).\label{est.cont.qtoQ.2}
\end{align}
Thus
\begin{align*}
\|w\|_{C([t_i,\,t_f],\,\mathbb J_V^\bot)}^2&\leq \exp\Bigl[D\int_{t_i}^{t_f}|w(t)|^2\|w(t)\|^2+\|Q(t)\|^4+|p(t)|_{[2]}^2\,dt\Bigr]\biggl(\|w(t_i)\|^2\\
&+D\int_{t_i}^{t_f}|p(t)-q(t)|_{[2]}^2(\|Q(t)\|^2+\|p(t)\|^2+\|q(t)\|^2)\,dt\biggr).
\end{align*}
Then using \eqref{cont.qtoQ.CH} and \eqref{cont.qtoQ.L2V} 
\begin{equation}\label{cont.qtoQ.CV}
\begin{cases}
&\|p-q\|_{L^4}\leq 1\\
&\|P_i-Q_i\|_{\mathbb J_V^\bot}\leq 1
\end{cases}\Rightarrow
\Bigl[\|w\|_{C([t_i,\,t_f],\,\mathbb J_V^\bot)}^2\leq D_1\|Q_i-P_i\|^2+D_1\|p-q\|_{L^4}^2\Bigr].
\end{equation}
and conclude that the map $U^J$ is $(L^4\times (\mathbb J_V^\bot),\,C([t_i,\,t_f],\,\mathbb J_V^\bot))$-continuous.\par
From \eqref{est.cont.qtoQ.2} we can also obtain
\begin{align*}
&\quad\int_{t_i}^{t_f}|w(t)|_{[2]}^2\,dt\\
&\leq D_2\|w\|^2_{C([0,\,T],\,\mathbb J_V^\bot)}\Bigl(1+\int_{t_i}^{t_f}|w(t)|^2\|w(t)\|^2+\|Q(t)\|^4+|p(t)|_{[2]}^2\,dt\Bigr)\\
&+D_2\|p-q\|_{L^4}^2\Big(\|Q\|_{L^4}^2+\|p\|_{L^4}^2+\|q\|_{L^4}^2).
\end{align*}
Hence by \eqref{cont.qtoQ.CH} and \eqref{cont.qtoQ.L2V} we have
\begin{multline*}
\begin{cases}
&\|p-q\|_{L^4}\leq 1\\
&\|P_i-Q_i\|_{\mathbb J_V^\bot}\leq 1
\end{cases}\\
\Rightarrow
\Bigl[\|w\|_{L^2(t_i,\,t_f,\,\mathbb J_{D(A)}^\bot)}^2\leq D_3\|w\|_{C([t_i,\,t_f],\,\mathbb J_V^\bot)}^2+D_3\|p-q\|_{L^4}^2\Bigr]
\end{multline*}
and, by \eqref{cont.qtoQ.CV} we obtain
\begin{equation}\label{cont.qtoQ.L2D(A)}
\begin{cases}
&\|p-q\|_{L^4}\leq 1\\
&\|P_i-Q_i\|_{\mathbb J_V^\bot}\leq 1
\end{cases}\Rightarrow
\Bigl[\|w\|_{L^2(t_i,\,t_f,\,\mathbb J_{D(A)}^\bot)}^2\leq D_4\|Q_i-P_i\|^2+D_4\|p-q\|_{L^4}^2\Bigr]
\end{equation}
and conclude that the map $U^J$ is $(L^4\times (\mathbb J_V^\bot),\,L^2(t_i,\,t_f,\,D(A)))$-continuous.\par
Now we define another map on the product $W^{1,\infty}([t_i,\,t_f],\,\mathbb J)\times \mathbb J_V^\bot$ taking values on H:
$$
\Gamma^J:(q,\,Q_i)\mapsto -\nu A(U^J(q,\,Q_i)+q)-P^\nabla B(U^J(q,\,Q_i)+q)+\tilde F
$$
Fix $(q,\,Q_i)\in W^{1,2}([t_i,\,t_f],\,\mathbb J)\times \mathbb J_V^\bot$ and, consider another pair $(p,\,P_i)$ in the same space. Again, as we have done before put $Q=U^J(q,\,Q_i)$ and $P=U^J(p,\,P_i)$. We compute the norm of the difference
\begin{align*}
&|\Gamma^J(q,\,Q_i)-\Gamma^J(p,\,P_i)|\\
=&|\nu A(P+p)+P^\nabla B(P+p)-\nu A(Q+q)-P^\nabla B(Q+q)|\\
\leq& |\nu A(P+p-Q-q)|+|P^\nabla[B(P+p)-B(Q+q)]|\\
\leq& |\nu A(P+p-Q-q)|+|B(P+p)-B(Q+q)|\\
\leq& \nu|P-Q|_{[2]}+\nu|p-q|_{[2]}+|BP-BQ|\\
+&|B(P,\,p)-B(Q,\,q)|+|B(p,\,P)-B(q,\,Q)|+|Bp-Bq|\\
\leq& \nu|P-Q|_{[2]}+\nu|p-q|_{[2]}+C|P-Q|_{[2]}(\|P\|+\|Q\|)\\
+&C\Bigl(|P-Q|_{[2]}\|q\|+\|P\||p-q|_{[2]}\Bigr)\\
+&C\Bigl(\|P\||p-q|_{[2]}+|P-Q|_{[2]}\|q\|\Bigr)\\
+&C|p-q|_{[2]}(\|p\|+\|q\|).
\end{align*}
Therefore
\begin{align*}
&|\Gamma^J(q,\,Q_i)-\Gamma^J(p,\,P_i)|\\
\leq& \nu|P-Q|_{[2]}+\nu|p-q|_{[2]}+C_1|P-Q|_{[2]}\Bigl(\|P\|+\|Q\|+\|q\|\Bigr)\\
+&C_1|p-q|_{[2]}\Bigl(\|P\|+\|p\|+\|q\|\Bigr)\\
\leq& C_0|P-Q|_{[2]}\Bigl(1+\|P\|+\|Q\|+\|q\|\Bigr)+C_0|p-q|_{[2]}\Bigl(1+\|P\|+\|p\|+\|q\|\Bigr).
\end{align*}
For $\|p-q\|_{L^4}\leq 1$, and $\|P_i-Q_i\|\leq 1$, using \eqref{cont.qtoQ.CV} we obtain $\|P(t)-Q(t)\|^2\leq 2D_1$ and then $\|P(t)\|\leq \sqrt{2D_1}+\|Q(t)\|$. So we can arrive to
\begin{equation}\label{for.gamma.conts}
|\Gamma^J(q,\,Q_i)-\Gamma^J(p,\,P_i)|\leq C_2|P-Q|_{[2]}+C_2|p-q|_{[2]}\Bigl(1+\|p\|\Bigr)
\end{equation}
Note that we have used the fact that $q\in W^{1,\infty}$ and so, $\|q\|$ is bounded. But we can not ``replace'' $\|p\|$ by a constant in the last member because, unlike as $q$, $p$ is not fixed and we are considering that $p$ varying in $L^4$ topology.
\begin{align*}
&\int_{t_i}^{t_f}|\Gamma^J(q,\,Q_i)-\Gamma^J(p,\,P_i)|^2\,dt\\
\leq& 2C_2^2\int_{t_i}^{t_f}|P-Q|_{[2]}^2\,dt+2C_2^2\int_{t_i}^{t_f}|p-q|_{[2]}^2\Bigl(1+\|p\|\Bigr)^2\,dt\\
\leq& 2C_2^2\|P-Q\|_{L^2(t_i,\,t_f,\,D(A))}^2+2C_2^2\int_{t_i}^{t_f}|p-q|_{[2]}^2 2\Bigl(1+\|p\|^2\Bigr)\,dt\\
\leq& C_4\|P-Q\|_{L^2(t_i,\,t_f,\,D(A))}^2+C_4\|p-q\|_{L^4}^2(1+\|p\|_{L^4}^2).
\end{align*}
Therefore, using \eqref{cont.qtoQ.L2D(A)}:
\begin{multline}\label{cont.qtoGamma.L2H}
\begin{cases}
&\|p-q\|_{L^4}\leq 1\\
&\|P_i-Q_i\|_{\mathbb J_V^\bot}\leq 1
\end{cases}\\
\Rightarrow
\|\Gamma^J(q,\,Q_i)-\Gamma^J(p,\,P_i)\|_{L^2(t_i,\,t_f,\,H)}^2\,dt\leq
C_5\|P_i-Q_i\|^2+C_5\|p-q\|_{L^4}^2
\end{multline}

and we conclude the $(L^4\times (\mathbb J_V^\bot),\,L^2(t_i,\,t_f,\,H))$-continuity of $\Gamma^J$.\par
Now we can indicate which is the control $v^J$ appearing in the statement of the proposition: In fact
$$
v^J(q,\,Q_i):=\dot q - P^J\Gamma^J(q,\,Q_i)
$$
satisfies the statement. Indeed its $(W^{1,2}([t_i,\,t_f],\,\mathbb J)\times (\mathbb J_V^\bot),\,L^2(t_i,\,t_f,\,\mathbb J)$-continuity follows from the  $(L^4\times (\mathbb J_V^\bot),\,L^2(t_i,\,t_f,\,H)$-continuity of $\Gamma^J$ and from the\\ $(W^{1,2}([t_i,\,t_f],\,\mathbb J),\,L^2(t_i,\,t_f,\,\mathbb J)$-continuity of $q\mapsto\dot q$.\par
To prove that the projection of the solution of the system
\begin{equation}\label{NS.vJ}
u_t=-\nu Au-Bu+\tilde F+v^J(q,\,Q_i),\qquad u_i=q(t_i)+Q_i
\end{equation}
coincides with $q$ we differentiate $q+U^J(q,\,Q_i)$ obtaining
\begin{align*}
[q+U^J(q,\,Q_i)]_t&=\dot q-\nu A(U^J(q,\,Q_i))-P^{-J}P^\nabla B(U^J(q,\,Q_i)+q)+P^{-J}\tilde F\\
&=-\nu A(q+U^J(q,\,Q_i))-P^\nabla B(q+U^J(q,\,Q_i))+\tilde F\\
&+\dot q -(-\nu Aq-P^JP^\nabla B(q+U^J(q,\,Q_i))+P^J\tilde F)\\
&=-\nu A(q+U^J(q,\,Q_i))-P^\nabla B(q+U^J(q,\,Q_i))+\tilde F+v^J(q,\,Q_i)
\end{align*}
showing that $q+U^J(q,\,Q_i)$ is the (unique) solution of \eqref{NS.vJ}.\par
To finish the prove remains to verify that $v^J\in L^\infty([t_i,\,t_f],\,\mathbb J)$. Since $q\in W^{1,\infty}([t_i,\,t_f],\,\mathbb J)$ we have $\dot q\in L^\infty([t_i,\,t_f],\,\mathbb J)$ and by
\begin{align}
\|\Gamma^J(q,\,Q_i)\|_{V^\prime}&\leq\nu\|Q+q\|+\|Q+q\|^2+\|\tilde F\|_{V^\prime}\notag\\
&\leq C\Bigl(\|Q+q\|+\|Q+q\|^2+|\tilde F|\Bigr)\leq C_1.\label{bd.PGamma.Vprime}
\end{align}
Hence $\Gamma^J(q,\,Q_i)\in L^\infty([t_i,\,t_f],\,V^\prime)$ and then, $P^J\Gamma^J(q,\,Q_i)\in L^\infty([t_i,\,t_f],\,P^JV^\prime)$, i.e., $P^J\Gamma^J(q,\,Q_i)\in L^\infty([t_i,\,t_f],\,\mathbb J)$.\footnote{By the equivalence of $V^\prime$-norm and $l_1$-norm in $\mathbb J$. One can check that $V^\prime$ coincides with the domain $D(A^{-\frac{1}{2}})$ of the operator $A^{-\frac{1}{2}}$ and its Fourier characterization is $V^\prime=\{\sum_{k\in\mathbb N_0^2}u_kW_k\mid\,\sum_{k\in\mathbb N_0^2}u_k^2<+\infty\}$. For $v\in V^\prime,\,\|v\|_{V^\prime}=\frac{ab}{4}\sum_{k\in\mathbb N_0^2}u_k^2$.}\\
Moreover we can see that
\begin{align*}
\|v^J(q,\,Q_i)\|_{L^\infty([t_i,\,t_f],\,\mathbb J)}\leq& \|\dot q\|_{L^\infty([t_i,\,t_f],\,\mathbb J)}+\|\Gamma^J(q,\,Q_i)\|_{L^\infty([t_i,\,t_f],\,\mathbb J)}\\
\leq& \|\dot q\|_{L^\infty([t_i,\,t_f],\,\mathbb J)}+D_1\|\Gamma^J(q,\,Q_i)\|_{L^\infty([t_i,\,t_f],\,V^\prime)}
\end{align*}
and, by \eqref{bd.PGamma.Vprime}, we obtain
\begin{equation}\label{bd.vJLinf}
\|v^J(q,\,Q_i)\|_{L^\infty([t_i,\,t_f],\,\mathbb J)}\leq D_2
\end{equation}
where $D_2$ depends only on the norms $\|\dot q\|_{L^\infty([t_i,\,t_f],\,\mathbb J)}$, $\|q\|_{L^\infty([t_i,\,t_f],\,\mathbb J)}$ and\\ $\|Q\|_{L^\infty([t_i,\,t_f],\,\mathbb J_V^\bot)}$. Then using \eqref{qtoQ.bdV} the constant $D_2$ can be chosen depending only on $\|q\|_{W^{1,\infty}([t_i,\,t_f],\,\mathbb J)}$ and $\|Q_i\|$.
\end{proof}
\begin{Definition}
We call $\delta$-{\bf metric} the function defined on the product space $\bigl(L^\infty([0,\,T],\mathbb R^d)\bigr)^2$ by
$$
\delta\bigl(u,\,v\bigr):=meas\{t\in[0,\,T]\mid\,u(t)\ne v(t) \}.
$$
\end{Definition}
\begin{Remark}
The $\delta$-{\bf metric}  is the restriction, to the space of ordinary controls, of the strong metric defined on the space of relaxed controls (see \cite{gam}).  
\end{Remark}
In \cite{gam} (Chapter 3) we can find the so called {\em Approximation Lemma} that says that a strongly continuous family of relaxed controls can be weakly approximated with arbitrary accuracy by a strongly continuous family of ordinary controls. At the end of that chapter we can find a ``Remark on the Terminology'' leading us to the following Lemma:    
\begin{Lemma}[Approximation Lemma; \cite{gam}]\label{app.lemma}
Let $A\subseteq \mathbb R^d$ be the convexification of a finite set of points:
$$
A:=Conv\{p_1,\,p_2,\,\dots,\,p_r\}
$$
and, $\mathcal V:=\{v(t,\,b)\in L^\infty([0,\,T],\,A)\mid\,b\in B\}$ be a $L^1$-continuous family of $A$-valued functions. Then for each $\varepsilon>0$ one can construct a  family $\mathcal V^\varepsilon:=\{v^\varepsilon(t,\,b)\in L^\infty([0,\,T],\,\{p_1,\,\dots,\,p_r\})\mid\,b\in B\}$ of $\{p_1,\,p_2,\,\dots,\,p_r\}$-valued functions such that
\begin{itemize}
\item $\mathcal V^\varepsilon$ is $\delta$-continuous, i.e., $b\mapsto v^\varepsilon(\cdot,\,b)$ is $(B,\,\delta)$-continuous;
\item $\mathcal V^\varepsilon$ $\varepsilon$-approximates, uniformly w.r.t. $b$, the family $\mathcal V$ in relaxation metric, i.e., $\forall b\in B\quad\|v^\varepsilon(\cdot,\,b)-v(\cdot,\,b)\|_{rx}<\varepsilon$ and;
\item The elements of $\mathcal V^\varepsilon$ are piecewise constant and the number of intervals of constancy is the same for all $b\in B$.
\end{itemize}
\end{Lemma}
\begin{Remark}\label{delta.Lq,leng}
Note that $\delta$-metric and $L^q$-metric $(1\leq q<+\infty)$ are equivalent in the subset of piecewise constant functions taking values on a fixed finite set. Also, since the controls take values on a finite set, the $\delta$-continuity of $\mathcal V^\varepsilon$ is equivallent to the continuity of the lengths of the intervals of constancy of the controls.
\end{Remark}
In Lemma \ref{app.lemma} is said that the intervals of constancy can be taken the same for all $b\in B$ but, looking at [\cite{gam}; ch. 3] some of those intervals may degenerate to a single point. We claim that we can suppose non-degeneracy of the intervals,\footnote{Note that it is not enough to eliminate the degenerate intervals because the number of intervals would not be the same for all $b\in B$.}. We may even suppose that there exists a lower bound $\theta^\varepsilon$ for the lengths of the intervals of constancy of the family $\mathcal V^\varepsilon$, i.e., for all $b\in B$ none $v^\varepsilon(\cdot,\,b)$ has an interval of constancy with length less than $\theta^\varepsilon$.\par
Looking at [\cite{gam}; ch. 3] we see that the intervals of constancy are (or can be) constructed in the following way: First we subdivide the interval $[0,\,T]$ into $n^2$ intervals --- $L_i,\quad i=1,\,\dots,\,n^2$ --- with the same length --- $\frac{T}{n^2}$. Then subdivide each one of these intervals $L_i$ into $r$ subintervals --- $L_{ij}\quad j=1,\,\dots,\,r$ --- which lengths $length(L_{ij})=\int_{L_i}v_j(\tau,\,b)\,d\tau,\quad v(\cdot,\,b)=\sum_{j=1}^r v_j(\cdot,\,b)p_j,\,v_j(t,\,b)\in [0,\,1]$, depend (continuously) on $b$ and the interval $L_{ij_1}$ lays on the left of  $L_{ij_2}$ if $j_1<j_2$. To this partition is associated the piecewise constant control $v_{n^2}(\cdot,\,b)$ defined by:
$$
t\in L_{ij}\Rightarrow v_{n^2}(t,\,b)=p_j;\qquad i=1,\,\dots,\,n^2,\;j=1,\,\dots,\,r,
$$ 
i.e., in each interval $L_i$ we use all the controls from $\{p_1,\,\dots,\,p_r\}$ using $p_{j_1}$ before $p_{j_2}$ if $j_1<j_2$. Note that as we have said before some control may be used for time zero.\par Each $v^\varepsilon(\cdot,\,b)\in\mathcal V^\varepsilon$ have the form $v^\varepsilon=v_{n^2}$ for some $n$ depending only on $\varepsilon$.\par
We shall need the following strong result without non degeneracy of the intervals of constancy:
\begin{Corollary}[Approximation Corollary]\label{app.corol}
Let $A\subseteq \mathbb R^d$ be the convexification of a finite set of points:
$$
A:=Conv\{p_1,\,p_2,\,\dots,\,p_r\}
$$
and, $\mathcal V:=\{v(t,\,b)\in L^\infty([0,\,T],\,A)\mid\,b\in B\}$ be a $L^1$-continuous family of $A$-valued functions. Then for each $\varepsilon>0$ there is $\theta^\varepsilon>0$ and a  family $\mathcal Z^\varepsilon:=\{z^\varepsilon(t,\,b)\in L^\infty([0,\,T],\,\{p_1,\,p_2,\,\dots,\,p_r\})\mid\,b\in B\}$ of $\{p_1,\,p_2,\,\dots,\,p_r\}$-valued functions such that
\begin{itemize}
\item $\mathcal Z^\varepsilon$ is $\delta$-continuous;
\item $\mathcal Z^\varepsilon$ $\varepsilon$-approximates, uniformly w.r.t. $b$, the family $\mathcal V$ in relaxation metric, i.e., $\forall b\in B\quad\|z^\varepsilon(\cdot,\,b)-v(\cdot,\,b)\|_{rx}<\varepsilon$;
\item The elements of $\mathcal Z^\varepsilon$ are piecewise constant and the number of intervals of constancy is the same for all $b\in B$ and,
\item For all $b\in B$ all the intervals of constancy of $z^\varepsilon(\cdot,\,b)$ have a length not less than $\theta^\varepsilon>0$.
\end{itemize}
\end{Corollary}
Before the proof consider the following Lemma:
\begin{Lemma}\label{PtoP}
Given $K>0$, $\gamma>0$ and $L\in\mathbb N_0$. Define the sets
\begin{align*}
\mathcal P_0&:=\{(x_1,\,\dots,\,x_L)\in\mathbb R^L\mid\,x_i\geq 0,\,\sum_{i=1}^L x_i=K\} \\
\mathcal P_\theta&:=\{(x_1,\,\dots,\,x_L)\in\mathbb R^L\mid\,x_i\geq \theta,\,\sum_{i=1}^L x_i=K\}
\end{align*}
where $\theta>0$.\\
Choose $n\in\mathbb N_0$ such that $\frac{(L+1)K}{nL}<\gamma$ and put
\begin{equation}\label{theta.KL}
\theta=\frac{K}{nL}.\qquad\footnotemark
\end{equation}
\footnotetext{So, $\theta$ depends on both $K,\,L$ and $\gamma$.}
Then the map $P_{0,\theta}=(P_{0,\theta}^1,\,\dots,\,P_{0,\theta}^L)$, defined on $\mathcal P_0$ by:
$$
P_{0,\theta}^i(x):=\bigl(1-\frac{1}{n}\bigr)\bigl(x_i-\frac{K}{L}\bigr)+\frac{K}{L},
$$
is continuous, take its values on $\mathcal P_\theta$ and, satisfies $|P_{0,\theta}^i(x)-x_i|<\gamma$.
\end{Lemma}

The Approximation Corollary follows from the Approximation Lemma and from the last Lemma. Indeed let $\varepsilon>0$ be a positive real number, put $D:=\max_{i\in\{1,\,\dots,\,r\}}\{\|p_i\|\}$, take $K=1,\,L=r$, $\gamma<\frac{\varepsilon}{2TDr}$, $n\in\mathbb N_0$ such that $\frac{r+1}{nr}<\gamma$ and $\theta=\frac{1}{nr}$.\\
The continuity of $P_{0,\theta}$ implies the continuity of the map
$$
v(t,\,b)=\sum_{j=1}^rv_j(t,\,b)p_j\mapsto z(t,\,b)=\sum_{j=1}^rP_{0,\theta}^j(v(t,\,b))
$$
and then the family of the controls $z (\cdot,\,b)=\sum_{j=1}^rP_{0,\theta}^j(v(t,\,b))$ is parametrized in $L^2(0,\,T,\,Conv\{p_1,\,\dots,\,p_r\})$-norm.\\
The family $\mathcal Z$ $\frac{\varepsilon}{2}$-approximates the family $\mathcal V$ in relaxation metric. Indeed
\begin{align*}
&\|z (\cdot,\,b)-v (\cdot,\,b)\|_{rx}\\
\leq&\int_0^T|z (\tau,\,b)-v (\tau,\,b)|\,d\tau=\int_0^T|\sum_{j=1}^rz_j(t,\,b)p_j-\sum_{j=1}^rv_j(t,\,b)p_j|\,d\tau\\
\leq& \int_0^T|\sum_{j=1}^r|z_j(t,\,b)-v_j(t,\,b)|D\,d\tau< T\gamma Dr<\frac{\varepsilon}{2}
\end{align*}
Now apply Approximation lemma to the family  $\mathcal Z$ and find a family  $\mathcal Z^\frac{\varepsilon}{2}$ that $\frac{\varepsilon}{2}$-approximates $\mathcal Z$ in relaxation metric. Hence
$$
\|z^\frac{\varepsilon}{2}(\cdot,\,b)-v (\cdot,\,b)\|_{rx}\leq\|z^\frac{\varepsilon}{2}(\cdot,\,b)-z (\cdot,\,b)\|_{rx}+\|z(\cdot,\,b)-v (\cdot,\,b)\|_{rx}<\varepsilon.
$$
Therefore the family  $\mathcal Z^\frac{\varepsilon}{2}$  $\varepsilon$-approximates $\mathcal V$ in relaxation metric.\par
From the Approximation Lemma the number of intervals of constancy of $z^\frac{\varepsilon}{2} (\cdot,\,b)$ is the same for all $b\in B$ and, the $L_{ij}$ interval of $z^\frac{\varepsilon}{2} (\cdot,\,b)$ has length $\int_{L_i} z_j(\tau,\,b)\,d\tau>\frac{T}{n^2}\theta$ for that $n\in\mathbb N_0$ such that $z^\frac{\varepsilon}{2} (\cdot,\,b)=z_{n^2}(\cdot,\,b)$. Thus all intervals have a length not smaller than $\frac{T}{n^2}\theta>0$ --- a positive constant depending only in $\varepsilon$.

\begin{proof}[Proof of Lemma \ref{PtoP}]
The continuity of each $P_{0,\theta}^i$, and then of $P_{0,\theta}$, is clear.\\
From
$$
\sum_{i=1}^LP_{0,\theta}^i(x)=\sum_{i=1}^L\Bigr[\bigl(1-\frac{1}{n}\bigr)\bigl(x_i-\frac{K}{L}\bigr)+\frac{K}{L}\Bigl]=K
$$
and
$$
P_{0,\theta}^i(x)=\bigl(1-\frac{1}{n}\bigr)x_i+\frac{K}{nL}\geq\frac{K}{nL}=\theta
$$
we conclude that $P_{0,\theta}$ takes its values on $\mathcal P_\theta$.\\
It remais to estimate $|P_{0,\theta}^i(x)-x_i|$:\\
$$
\Bigl|P_{0,\theta}^i(x)-x_i\Bigr|=\Bigl|-\frac{x_i}{n}+\frac{K}{nL}\Bigr|\leq \frac{K}{n}+\frac{K}{nL}\leq \frac{(L+1)K}{nL}<\gamma.
$$
\end{proof}
\ \\
Now, we are ready to start the proof of the induction step. Fix $u_0\in V$. By ``back-induction'' hypothesis system [\eqref{sysN}.N] is time-T solid controllable in observed component so, there is a family $\mathcal V:=\{v(\cdot,\,b)\in L^\infty([0,\,T],\,\mathbb R^{\kappa_N})\mid\,b\in B\}$ parametrized continuously in $L^2$-metric such that
\begin{align}
&\mathbb E_T(u_0,\,\mathcal V)\; \text{covers}\;\overline{\mathcal O}_R(u_0^{\#\mathcal O})\;\text{solidly}.\label{ET.sol}\\
&\|v(t,\,b)\|_{l_1}\leq\Xi_0\quad\forall (t,\,b)\in[0,\,T]\times B.\label{v.unif.bdd.Xi0}
\end{align}
Note that, by \eqref{v.unif.bdd.Xi0} we obtain
$$
\forall (t,\,b)\in[0,\,T]\times B\quad v(t,\,b)\in Conv\{\pm\Xi_0 e_k\mid\,k\in\mathcal K^N\}\quad\footnotemark
$$
\footnotetext{Note that the norm used in \eqref{v.unif.bdd.Xi0} is the $l_1$ one.}
Then there exists $\Xi>0$ such that
\begin{multline}\label{v.unif.bdd.Xi}
\forall (t,\,b)\in[0,\,T]\times B\\ v(t,\,b)\in Conv\{\pm\Xi e_k,\;\pm\Xi\delta_{m,n}\mid\,k\in\mathcal K^{N-1},\;(m,\,n)\in S_{N-1}\}=:C_N
\end{multline}
Now fix $\varepsilon>0$ and $\varepsilon$-approximate the family $\mathcal V$ taking values on $C_N$ by a family $\mathcal Z^\varepsilon:=\{z^\varepsilon(t,\,b)\in\mid\,b\in B\}$ taking values on $\{\pm\Xi e_k,\;\pm\Xi\delta_{m,n}\mid\,k\in\mathcal K^{N-1},\;(m,\,n)\in S_{N-1}\}$ like in Corollary \ref{app.corol}.\quad\footnote{Recall that $S_{N-1}$ has been defined, in the proof of Proposition \ref{Kj.to.Kj+1} for $N>2$ and in the proof of Proposition \ref{k1sat} for $N=2$, as the set of ``extracted'' pairs of modes associated with the new directions.}\\
Choose a real number $\gamma$ such that
$$
0<\gamma<\|u_0^{\#\mathcal O}-\partial\mathbb E_T(u_0,\,\mathcal V)\|_{l_1}-R.
$$
\begin{Remark}
It is clear that $\|u_0^{\#\mathcal O}-\partial\mathbb E_T(u_0,\,\mathcal V)\|_{l_1}\geq R$ because of \eqref{ET.sol}. To see that the inequality is, in fact, strict we suppose not and then, there is $p\in\partial\mathbb E_T(u_0,\,\mathcal V)$ such that $\|u_0^{\#\mathcal O}-p\|_{l_1}=R$. For the sequence of continuous functions $\phi_n$ defined in $B$ by
$$
\phi_n(b):=\Bigl(1-\frac{1}{n}\Bigr)\Bigl(\mathbb E_T(u_0,\,v(\cdot,\,b))-u_0^{\#\mathcal O}\Bigr)+u_0^{\#\mathcal O}
$$
we have
\begin{align*}
\|\phi_n-\mathbb E_T\|_{C(B,\,\mathbb R_\infty^{\#\mathcal O})}&=\frac{1}{n}\|\mathbb E_T\|_{C(B,\,\mathbb R_\infty^{\#\mathcal O})}+\frac{1}{n}\|u_0^{\#\mathcal O}\|_{l_\infty}\to 0;\\
\phi_n(\mathcal V)&=\Bigl(1-\frac{1}{n}\Bigr)\mathbb E_T(u_0,\,\mathcal V)+\frac{1}{n}u_0^{\#\mathcal O}\\
\Bigl\|\Bigl(1-\frac{1}{n}\Bigr)p+\frac{1}{n}u_0^{\#\mathcal O}-u_0^{\#\mathcal O}\Bigr\|_{l_1}&=\Bigl(1-\frac{1}{n}\Bigr)\|p-u_0^{\#\mathcal O}\|_{l_1}\leq R.
\end{align*}
Therefore
$$
\Bigl(1-\frac{1}{n}\Bigr)p+\frac{1}{n}u_0^{\#\mathcal O}\in\partial\phi_n(\mathcal V)\cap\mathcal O_R(u_0^\#\mathcal O)
$$
and then, for all $n\in\mathbb N_0$, $\phi_n(\mathcal V)$ does not cover $\overline{\mathcal O}_R(u_0^{\#\mathcal O})$ which contradicts \eqref{ET.sol}.
\end{Remark}
Now note that for $x\in \overline{\mathcal O}_R(u_0^{\#\mathcal O})$ we obtain $\|u_0^{\#\mathcal O}-\partial\mathbb E_T(u_0,\,\mathcal V)\|_{l_1}\leq R+\|x-\partial\mathbb E_T(u_0,\,\mathcal V)\|_{l_1}$, i.e., $\gamma<\|x-\partial\mathbb E_T(u_0,\,\mathcal V)\|_{l_1}$ and then
\begin{equation}\label{x>gamma}
\frac{\gamma}{\#\mathcal O}<\|x-\partial\mathbb E_T(u_0,\,\mathcal V)\|_{l_\infty}.
\end{equation}
By the continuity of $\mathbb S_{srx}$ there is $\bar\varepsilon>0$ such that
\begin{multline*}
\|v-w\|_{rx}<\bar\varepsilon\Rightarrow\\
\|\Pi_{\mathcal O}\circ\mathbb S_s(u_0,\,F,\,w,\,\nu)(T)-\Pi_{\mathcal O}\circ\mathbb S_s(u_0,\,F,\,v,\,\nu)(T)\|_{C(B,\,\mathbb R_\infty^{\#\mathcal O}}<\frac{\gamma}{2\#\mathcal O}.
\end{multline*}
 Fix $\varepsilon<\bar\varepsilon$ and put
$$
\mathbb E_T^\varepsilon(u_0,\,z^\varepsilon(b)):=\Pi_{\mathcal O}\circ\mathbb S_s(u_0,\,F,\,z^\varepsilon(b),\,\nu)(T)
$$
and consider $\phi$ defined on $B$ such that
$$
\|\phi(\cdot)-\mathbb E_T^\varepsilon(u_0,\,z^\varepsilon(\cdot))\|_{C(B,\,\mathbb R_\infty^{\#\mathcal O}}<\frac{\gamma}{2\#\mathcal O}.
$$ 
Hence
$$
\|\phi-\mathbb E_T\|_{C(B,\,\mathbb R_\infty^{\#\mathcal O}}\leq \|\phi-\mathbb E_T^\varepsilon\|_{C(B,\,\mathbb R_\infty^{\#\mathcal O}}+\|\mathbb E_T^\varepsilon-\mathbb E_T\|_{C(B,\,\mathbb R_\infty^{\#\mathcal O}}<\frac{\gamma}{\#\mathcal O}
$$
and, by \eqref{x>gamma} we conclude
$$
\forall x\in \overline{\mathcal O}_R(u_0^\#\mathcal O)\Bigl[deg(\phi,\,B,\,x)=deg(\mathbb E_T,\,B,\,x)\ne 0\Bigr]
$$
so, $\phi(B)$ covers $\overline{\mathcal O}_R(u_0^\#\mathcal O)$. Therefore $E_T^\varepsilon(u_0,\,\mathcal Z^\varepsilon)$ covers $\overline{\mathcal O}_R(u_0^\#\mathcal O)$ solidly for all positive $\varepsilon<\bar\varepsilon$.
Moreover, for all pairs $(t,\,b)$ we have $\|z^\varepsilon(t,\,b)\|_{l_1}\leq\Xi$. Therefore for $\varepsilon<\bar\varepsilon$ system [\eqref{sysN}.N] is is solid controllable by means of the family $\mathcal Z^\varepsilon$.\par
Now we fix $\varepsilon\in]0,\,\bar\varepsilon[$ and define $\mathcal Z:=\mathcal Z^\varepsilon;\quad z(t,\,b):=z^\varepsilon(t,\,b)$. We shall compare the trajectories generated by a control $z(t,\,b)\in\mathcal Z$ with another generated by some ``appropriated'' control taking values on $\mathbb R^{\kappa_{N-1}}$.\par
To the control $z(\cdot,\,b)$ is associated a partition $X(b)$ of $[0,\,T]$ into $m$ non-degenerated intervals of lengths $x_i\geq\theta>0$:
$$
X(b)=(x_1,\,x_2,\,\dots,\,x_m)\in\mathbb R^m,\;m\in\mathbb N_0,\;\sum_{i=1}^mx_i=T,$$
where $m$ and $\theta$ are independent of the parameter $b$. We put
$$
A(b)=\{(0=\alpha_0,\,\alpha_1,\,\dots,\,\alpha_m=T)\}
$$
for the end points of the intervals in $X(b)$. So, 
\begin{multline*}
A(b)\in\mathcal A_\theta:=\{(\alpha_0,\,\alpha_1,\,\dots,\,\alpha_m)\in\mathbb R^{m+1}\mid\,\alpha_0=0,\,\alpha_m=T,\\
\alpha_i-\alpha_{i-1}\geq\theta,\,\sum_{i=1}^m(\alpha_i-\alpha_{i-1})=T\}.
\end{multline*}
Let $w\in\mathbb R,\quad w\geq 3$. Now in $[0,\,T]$ we want to define a function $\phi_w(\cdot,\,b)$ associated to $X(b)$ with the following properties\\
\begin{itemize}
\item $\phi_w(\cdot,\,b)$ vanishes at the points $\alpha_i,\quad i=0,\,\dots,\,m$;
\item $\phi_w(\cdot,\,b)\in W^{1,\infty}([0,\,T],\,\mathbb R)$ with
$$
\|\phi_w(\cdot,\,b)\|_{C([0,\,T],\,\mathbb R)}\leq 1;\qquad \|\dot{\phi}_w(\cdot,\,b)\|_{L^\infty([0,\,T],\,\mathbb R)}\leq\frac{w(1+\theta)}{\theta}
$$
\item $\delta(\phi_w(t,\,b),\,\sin(wt))\leq \frac{2T}{w}$ and;
\item For fixed $w$, the map $A(b)\mapsto \phi_w(\cdot,\,b)$ is $(\mathcal A_\theta,\,W^{1,2}(0,\,T,\,\mathbb R))$-continuous (where $\mathcal A_\theta$ is endowed with the topology induced by $\mathbb R^{m+1}$).
\end{itemize}
We proceed with the construction of $\phi_w(\cdot,\,b)$: For each $i\in\{1,\,\dots,\,m\}$ put $\rho_i=\frac{x_i}{w}$. Then subdivide each interval $[\alpha_{i-1},\,\alpha_i]$ into
$$
[\alpha_{i-1},\,\alpha_i]=[\alpha_{i-1},\,\alpha_{i-1}+\rho_i[\cup[\alpha_{i-1}+\rho_i,\,\alpha_i-\rho_i]\cup]\alpha_i-\rho_i,\,\alpha_i].\quad\footnotemark
$$ 
\footnotetext{Note that, since $w\geq 3$ we have $\rho_i\leq\frac{x_i}{3}$ and the subdivision is well defined.}
For each $i=1,\,\dots,\,m$ put
$$
\phi_w(t,\,b)=
\begin{cases}
\frac{\sin(w(\alpha_{i-1}+\rho_i))}{\rho_i}(t-\alpha_{i-1})&\text{if}\quad t\in[\alpha_{i-1},\,\alpha_{i-1}+\rho_i];\\
\sin(wt)&\text{if}\quad t\in[\alpha_{i-1}+\rho_i,\,\alpha_i-\rho_i];\\
\frac{\sin(w(\alpha_i-\rho_i))}{-\rho_i}(t-\alpha_i)\quad&\text{if}\quad t\in[\alpha_i-\rho_i,\,\alpha_i].
\end{cases}
$$
Then the graph of the restriction of $\phi_w(\cdot,\,b)$ to an interval $[\alpha_{i-1},\,\alpha_i]$ is a concatenation of a straight line, a piece of the graph of $\sin(wt)$ and another straight line. From the construction is clear that $\phi_w(\cdot,\,b)$ vanishes at the points $\alpha_i,\quad i=0,\,\dots,\,m$ and that $\phi_w(\cdot,\,b)$ is continuous with $\|\phi_w(\cdot,\,b)\|_{C([0,\,T],\,\mathbb R)}\leq 1$.\\
In the subintervals $]\alpha_{i-1}+\rho_i,\,\alpha_i-\rho_i[$ we have $\dot\phi_w(t,\,b)=w\cos(wt)$ so,$|\dot\phi_w(t,\,b)|\leq w\leq\frac{w(1+\theta)}{\theta}$. In the subintervals $]\alpha_{i-1},\,\alpha_{i-1}+\rho_i[$ and $]\alpha_i-\rho_i,\,\alpha_i[$ we have $|\dot{\phi}_w(t,\,b)|\leq \frac{1}{\rho_i}=\frac{w}{x_i}\leq\frac{w}{\theta}\leq\frac{w(1+\theta)}{\theta}$. Hence we have $\|\dot{\phi}_w(\cdot,\,b)\|_{L^\infty([0,\,T],\,\mathbb R)}\leq \frac{w(1+\theta)}{\theta}$. Therefore $\phi_w(\cdot,\,b)\in W^{1,\infty}([0,\,T],\,\mathbb R)$ and
$$
\|\phi_w(\cdot,\,b)\|_{W^{1,\infty}([0,\,T],\,\mathbb R)}\leq 1+\frac{w(1+\theta)}{\theta}.
$$ 
We see that $\phi_w(t,\,b)$ differs from $\sin(wt)$ only in the intervals 
$[\alpha_{i-1},\,\alpha_{i-1}+\rho_i[$ and $]\alpha_i-\rho_i,\,\alpha_i]$ so,
$$
\delta(\phi_w(t,\,b),\,\sin(wt))=\sum_{i=1}^m 2\rho_i=\sum_{i=1}^m 2\frac{x_i}{w}\leq\frac{2T}{w}.
$$
It remains to check the continuity property. For that fix $w\geq 3$ and $A\in\mathcal A_\theta$. Let $\gamma>0$ and $B\in\mathcal A_\theta$ such that $\|A-B\|_{l_1}<\gamma$. Let $A=(\alpha_0,\,\alpha_1,\,\dots,\,\alpha_m)$ and $B=(\beta_0,\,\beta_1,\,\dots,\,\beta_m)$ be the coordinates of $A$ and $B$. Then we have $|\alpha_i-\beta_i|<\gamma$ for $i=1,\,\dots,\,m-1$, $\alpha_0=\beta_0=0$ and $\alpha_m=\beta_m=T$. For small $\gamma\quad (\gamma<\frac{\theta}{w})$\footnote{Note that for $\gamma<\frac{\theta}{w}$ we have $\gamma<\rho_i^A$ and $\gamma<\rho_i^B$ for all $i\in\{1,\,\dots,\,m\}$.}, putting $\phi_A$ and $\phi_B$ for the functions $\phi_w$ associated with $A$ and $B$, we have that  $\phi_A$ and $\phi_B$ differ only in the the following union of subintervals 
\begin{align}
&[0,\,\min\{\rho_1^A,\,\rho_1^B\}]\cup[\min\{\rho_1^A,\,\rho_1^B\},\,\max\{\rho_1^A,\,\rho_1^B\}]\notag\\
\cup&[\min\{\alpha_i-\rho_i^A,\,\beta_i-\rho_i^B\},\,\max\{\alpha_i-\rho_i^A,\,\beta_i-\rho_i^B\}]\notag\\
\cup&[\max\{\alpha_i-\rho_i^A,\,\beta_i-\rho_i^B\},\,\min\{\alpha_i,\,\beta_i\}]\notag\\
\cup&[\min\{\alpha_i,\,\beta_i\},\,\max\{\alpha_i,\,\beta_i\}]\notag\\
\cup&[\max\{\alpha_i,\,\beta_i\},\,\min\{\alpha_i+\rho_{i+1}^A,\,\beta_i+\rho_{i+1}^B\}]\notag\\
\cup&[\min\{\alpha_i+\rho_{i+1}^A,\,\beta_i+\rho_{i+1}^B\},\,\max\{\alpha_i+\rho_{i+1}^A,\,\beta_i+\rho_{i+1}^B\}]\notag\\
\cup&[\min\{T-\rho_m^A,\,T-\rho_m^B\},\,\max\{T-\rho_m^A,\,T-\rho_m^B\}]\notag\\
\cup&[\max\{T-\rho_m^A,\,T-\rho_m^B\},\,T].\label{unionI}
\end{align}
Now we prove that $\|\phi_A-\phi_B\|_{C([0,\,T],\,\mathbb R)}$ goes to zero as $\gamma$ does. For that is enough to prove that  $\|\phi_A-\phi_B\|_{C(I,\,\mathbb R)}$ goes to zero as $\gamma$ does for every interval $I$ of the union \eqref{unionI}.\\
\begin{itemize}
\item For $I=[0,\,\min\{\rho_1^A,\,\rho_1^B\}]$:
\begin{align*}
&\|\phi_A-\phi_B\|_{C(I,\,\mathbb R)}=\max_{t\in I}\{\frac{\sin(w\rho_1^A)}{\rho_1^A}t-\frac{\sin(w\rho_1^B)}{\rho_1^B}t\}\\
=&\Bigl|\frac{\sin(w\rho_1^A)}{\rho_1^A}\min\{\rho_1^A,\,\rho_1^B\}-\frac{\sin(w\rho_1^B)}{\rho_1^B}\min\{\rho_1^A,\,\rho_1^B\}\Bigr|\quad\footnotemark
\end{align*}
\footnotetext{Note that the maximum is attained at one of the end points of the interval $I$, because the functions $\phi_A$ and $\phi_B$ are affine in $I$.}
\item For $I=[\min\{\rho_1^A,\,\rho_1^B\},\,\max\{\rho_1^A,\,\rho_1^B\}]$:
\begin{align*}
\|\phi_A-\phi_B\|_{C(I,\,\mathbb R)}\leq&(\max\{\rho_1^A,\,\rho_1^B\}-\min\{\rho_1^A,\,\rho_1^B\})\|\phi_A-\phi_B\|_{C([0,\,T],\,\mathbb R)}\\
\leq& 2|\rho_1^A-\rho_1^B|= 2\frac{1}{w}|x_1-y_1|\leq 2\frac{\gamma}{w}
\end{align*}
\item For $i=1,\,\dots,\,m-1$ and
\begin{itemize}
\item For $I=[\min\{\alpha_i-\rho_i^A,\,\beta_i-\rho_i^B\},\,\max\{\alpha_i-\rho_i^A,\,\beta_i-\rho_i^B\}]$:
\begin{align*}
&\|\phi_A-\phi_B\|_{C(I,\,\mathbb R)}\\
\leq&(\max\{\alpha_i-\rho_i^A,\,\beta_i-\rho_i^B\}-\min\{\alpha_i-\rho_i^A,\,\beta_i-\rho_i^B\})\|\phi_A-\phi_B\|_{C^0}\\
\leq& 2|(\alpha_i-\rho_i^A)-(\beta_i-\rho_i^B)|\leq 2|\alpha_i-\beta_i|+2|\rho_i^A-\rho_i^B|\\\leq & 2\gamma+2\frac{1}{w}|x_i-y_i|\leq 2\gamma+2\frac{2\gamma}{w}=2\gamma(1+\frac{2}{w}).
\end{align*}
\item For $I=[\max\{\alpha_i-\rho_i^A,\,\beta_i-\rho_i^B\},\,\min\{\alpha_i,\,\beta_i\}]$:
\begin{align*}
&\|\phi_A-&&\phi_B\|_{C(I,\,\mathbb R)}\\
=&\max_{t\in I}\Bigl\{&&\Bigl|\frac{\sin(w(\alpha_i-\rho_i^A)}{-\rho_i^A}(t-\alpha_i)-\frac{\sin(w(\beta_i-\rho_i^B))}{-\rho_i^B}(t-\beta_i)\Bigr|\Bigr\}\\
=&\max\biggl\{&&\Bigl|\frac{\sin(w(\alpha_i-\rho_i^A)}{-\rho_i^A}(\max\{\alpha_i-\rho_i^A,\,\beta_i-\rho_i^B\}-\alpha_i)\\
& &&\quad -\frac{\sin(w(\beta_i-\rho_i^B))}{-\rho_i^B}(\max\{\alpha_i-\rho_i^A,\,\beta_i-\rho_i^B\}-\beta_i)\Bigr|,\,\\
& &&\frac{\Bigl|\sin(w(\alpha_i-\rho_i^A)}{-\rho_i^A}(\min\{\alpha_i,\,\beta_i\}-\alpha_i)\\
& &&\quad -\frac{\sin(w(\beta_i-\rho_i^B))}{-\rho_i^B}(\min\{\alpha_i,\,\beta_i\}-\beta_i)\Bigr|\biggr\}
\end{align*}
\item For $I=[\min\{\alpha_i,\,\beta_i\},\,\max\{\alpha_i,\,\beta_i\}]$:
$$
\|\phi_A-\phi_B\|_{C(I,\,\mathbb R)}\leq(\max\{\alpha_i,\,\beta_i\}-\min\{\alpha_i,\,\beta_i\})\|\phi_A-\phi_B\|_{C^0}\leq 2\gamma.
$$
\item For $I=[\max\{\alpha_i,\,\beta_i\},\,\min\{\alpha_i+\rho_{i+1}^A,\,\beta_i+\rho_{i+1}^B\}]$:
\begin{align*}
&\|\phi_A-&&\phi_B\|_{C(I,\,\mathbb R)}\\
=&\max_{t\in I}\Bigl\{&&\Bigl|\frac{\sin(w(\alpha_i+\rho_{i+1}^A)}{\rho_i^A}(t-\alpha_i)-\frac{\sin(w(\beta_i+\rho_{i+1}^B))}{\rho_i^B}(t-\beta_i)\Bigr|\Bigr\}\\
=&\max\biggl\{&&\Bigl|\frac{\sin(w(\alpha_i+\rho_{i+1}^A)}{\rho_i^A}(\max\{\alpha_i,\,\beta_i\}-\alpha_i)\\
& &&\quad -\frac{\sin(w(\beta_i+\rho_{i+1}^B))}{\rho_i^B}(\max\{\alpha_i,\,\beta_i\}-\beta_i)\Bigr|,\,\\
& && \Bigl|\frac{\sin(w(\alpha_i+\rho_{i+1}^A)}{\rho_i^A}(\min\{\alpha_i+\rho_{i+1}^A,\,\beta_i+\rho_{i+1}^B\}-\alpha_i)\\
& &&\quad -\frac{\sin(w(\beta_i+\rho_{i+1}^B))}{\rho_i^B}(\min\{\alpha_i+\rho_{i+1}^A,\,\beta_i+\rho_{i+1}^B\}-\beta_i)\Bigr|\biggr\}
\end{align*}
\item For $I=[\min\{\alpha_i+\rho_{i+1}^A,\,\beta_i+\rho_{i+1}^B\},\,\max\{\alpha_i+\rho_{i+1}^A,\,\beta_i+\rho_{i+1}^B\}]$:
\begin{align*}
&\|\phi_A-\phi_B\|_{C(I,\,\mathbb R)}\\\leq&(\max\{\alpha_i+\rho_{i+1}^A,\,\beta_i+\rho_{i+1}^B\}-\min\{\alpha_i+\rho_{i+1}^A,\,\beta_i+\rho_{i+1}^B\})\|\phi_A-\phi_B\|_{C^0}\\
\leq& 2|(\alpha_i+\rho_{i+1}^A)-(\beta_i+\rho_{i+1}^B)|\leq 2|\alpha_i-\beta_i|+2|\rho_{i+1}^A-\rho_{i+1}^B|\\ \leq&2\gamma+2\frac{1}{w}|x_{i+1}-y_{i+1}|\leq 2\gamma+2\frac{2\gamma}{w}=2\gamma(1+\frac{2}{w}).
\end{align*}
\end{itemize}
For $I=[\min\{T-\rho_m^A,\,T-\rho_m^B\},\,\max\{T-\rho_m^A,\,T-\rho_m^B\}]$:
\begin{align*}
&\|\phi_A-\phi_B\|_{C(I,\,\mathbb R)}\\
\leq&(\max\{T-\rho_m^A,\,T-\rho_m^B\}-\min\{T-\rho_m^A,\,T-\rho_m^B\})\|\phi_A-\phi_B\|_{C([0,\,T],\,\mathbb R)}\\
\leq&2|\rho_m^A-\rho_m^B|=\frac{2}{w}|x_m-y_m|\leq\frac{2\gamma}{w}.
\end{align*}
For $I=[\max\{T-\rho_m^A,\,T-\rho_m^B\},\,T]$:
\begin{align*}
&\|\phi_A-\phi_B\|_{C(I,\,\mathbb R)}\\
=&\max_{t\in I}\Bigl\{\Bigl|\frac{\sin(w(T-\rho_m^A)}{\rho_m^A}(t-T)-\frac{\sin(w(T-\rho_m^B))}{\rho_i^B}(t-T)\Bigr|\Bigr\}\\
=&\Bigl|\frac{\sin(w(T-\rho_m^A)}{\rho_m^A}(\max\{T-\rho_m^A,\,T-\rho_m^B\}-T)\\
\quad&-\frac{\sin(w(T-\rho_m^B))}{\rho_i^B}(\max\{T-\rho_m^A,\,T-\rho_m^B\}-T)\Bigr|
\end{align*}
\end{itemize}
Therefore since $\beta_i\to\alpha_i$ and $\rho_m^B\to\rho_m^A$ as $\gamma\to 0$ we have that for every interval $I$ in the union \eqref{unionI}
$$
\|\phi_A-\phi_B\|_{C(I,\,\mathbb R)}\to 0\quad\text{as}\quad\gamma\to 0.
$$
Hence
\begin{equation}\label{cont.phi.C}
\|\phi_A-\phi_B\|_{C(0,\,T,\,\mathbb R)}\to 0\quad\text{as}\quad B\to A.
\end{equation}
For the derivatives we have
\begin{align}
&\int_0^T|\dot\phi_A(t)-\dot\phi_B(t)|^2\,dt\notag\\
\footnotemark=&\int_{0}^{\min\{\rho_1^A,\,\rho_1^B\}}|\dot\phi_A(t)-\dot\phi_B(t)|^2\,dt+\int_{\min\{\rho_1^A,\,\rho_1^B\}}^{\max\{\rho_1^A,\,\rho_1^B\}}|\dot\phi_A(t)-\dot\phi_B(t)|^2\,dt\notag\\
+&\sum_{i=1}^{m-1}\Biggl[\int_{\min\{\alpha_i-\rho_i^A,\,\beta_i-\rho_i^B\}}^{\max\{\alpha_i-\rho_i^A,\,\beta_i-\rho_i^B\}}|\dot\phi_A(t)-\dot\phi_B(t)|^2\,dt\notag\\
+&\int_{\max\{\alpha_i-\rho_i^A,\,\beta_i-\rho_i^B\}}^{\min\{\alpha_i,\,\beta_i\}}|\dot\phi_A(t)-\dot\phi_B(t)|^2\,dt\notag\\
+&\int_{\min\{\alpha_i,\,\beta_i\}}^{\max\{\alpha_i,\,\beta_i\}}|\dot\phi_A(t)-\dot\phi_B(t)|^2\,dt\notag\\
+&\int_{\max\{\alpha_i,\,\beta_i\}}^{\min\{\alpha_i+\rho_{i+1}^A,\,\beta_i+\rho_{i+1}^B\}}|\dot\phi_A(t)-\dot\phi_B(t)|^2\,dt\notag\\
+&\int_{\min\{\alpha_i+\rho_{i+1}^A,\,\beta_i+\rho_{i+1}^B\}}^{\max\{\alpha_i+\rho_{i+1}^A,\,\beta_i+\rho_{i+1}^B\}}|\dot\phi_A(t)-\dot\phi_B(t)|^2\,dt\Biggr]\notag\\
+&\int_{\min\{T-\rho_m^A,\,T-\rho_m^B\}}^{\max\{T-\rho_m^A,\,T-\rho_m^B\}}|\dot\phi_A(t)-\dot\phi_B(t)|^2\,dt\\
+&\int_{\max\{T-\rho_m^A,\,T-\rho_m^B\}}^T|\dot\phi_A(t)-\dot\phi_B(t)|^2\,dt.\label{int.phidot.L2}
\end{align}
\footnotetext{Recall that out of the intervals of \eqref{unionI} $\dot\phi_A$ and $\dot\phi_B$ coincide.}
Now we prove that each one of the intervals on the right-hand side of \eqref{int.phidot.L2} goes to zero when $\gamma$ does: 
\begin{align*}
&\int_{0}^{\min\{\rho_1^A,\,\rho_1^B\}}|\dot\phi_A(t)-\dot\phi_B(t)|^2\,dt\\
\leq&\biggl[\frac{\sin(w\rho_1^A)}{\rho_1^A}-\frac{\sin(w\rho_1^B)}{\rho_1^B}\biggr]^2\min\{\rho_1^A,\,\rho_1^B\}.
\end{align*}

\begin{align*}
&\int_{\min\{\rho_1^A,\,\rho_1^B\}}^{\max\{\rho_1^A,\,\rho_1^B\}}|\dot\phi_A(t)-\dot\phi_B(t)|^2\,dt\\
\leq&(\max\{\rho_1^A,\,\rho_1^B\}-\min\{\rho_1^A,\,\rho_1^B\})\|\dot\phi_A-\dot\phi_B\|_{L^\infty(0,\,T,\,\mathbb R)}^2\\
=&|\rho_1^A-\rho_1^B|\Bigl(2\frac{w(1+\theta)}{\theta}\Bigr)^2=\Bigl(2\frac{w(1+\theta)}{\theta}\Bigr)^2\frac{1}{w}|x_1-y_1|\leq 4w\gamma\Bigl(\frac{1+\theta}{\theta}\Bigr)^2,
\end{align*}

For $i=1,\,\dots,\,m-1$:
\begin{align*}
&\int_{\min\{\alpha_i-\rho_i^A,\,\beta_i-\rho_i^B\}}^{\max\{\alpha_i-\rho_i^A,\,\beta_i-\rho_i^B\}}|\dot\phi_A(t)-\dot\phi_B(t)|^2\,dt\\
\leq&(\max\{\alpha_i-\rho_i^A,\,\beta_i-\rho_i^B\}-\min\{\alpha_i-\rho_i^A,\,\beta_i-\rho_i^B\})\|\dot\phi_A-\dot\phi_B\|_{L^\infty(0,\,T,\,\mathbb R)}^2\\
=&4w^2\Bigl(\frac{1+\theta}{\theta}\Bigr)^2|(\alpha_i-\rho_i^A)-(\beta_i-\rho_i^B)|\leq 4w^2\Bigl(\frac{1+\theta}{\theta}\Bigr)^2(|\alpha_i-\beta_i|+|\rho_i^A-\rho_i^B|)\\
\leq&4w^2\Bigl(\frac{1+\theta}{\theta}\Bigr)^2(\gamma+\frac{1}{w}|x_i-y_i|)
\leq 4w^2\Bigl(\frac{1+\theta}{\theta}\Bigr)^2\Bigl(\gamma+\frac{2\gamma}{w}\Bigr)\\
=&4\gamma w^2\Bigl(\frac{1+\theta}{\theta}\Bigr)^2\Bigl(1+\frac{2}{w}\Bigr);
\end{align*}

\begin{align*}
&\int_{\max\{\alpha_i-\rho_i^A,\,\beta_i-\rho_i^B\}}^{\min\{\alpha_i,\,\beta_i\}}|\dot\phi_A(t)-\dot\phi_B(t)|^2\,dt\\
=&\biggl[\frac{\sin(w(\alpha_i-\rho_i^A))}{-\rho_i^A}-\frac{\sin(w(\beta_i-\rho_i^B))}{-\rho_i^B}\biggr]^2\\&\qquad\qquad\qquad\cdot(\min\{\alpha_i,\,\beta_i\}-\max\{\alpha_i-\rho_i^A,\,\beta_i-\rho_i^B\});
\end{align*}

\begin{align*}
&\int_{\min\{\alpha_i,\,\beta_i\}}^{\max\{\alpha_i,\,\beta_i\}}|\dot\phi_A(t)-\dot\phi_B(t)|^2\,dt\\
\leq&(\max\{\alpha_i,\,\beta_i\}-\min\{\alpha_i,\,\beta_i\})\|\dot\phi_A-\dot\phi_B\|_{L^\infty(0,\,T,\,\mathbb R)}^2\leq 4w^2\Bigl(\frac{1+\theta}{\theta}\Bigr)^2\gamma;
\end{align*}

\begin{align*}
&\int_{\max\{\alpha_i,\,\beta_i\}}^{\min\{\alpha_i+\rho_{i+1}^A,\,\beta_i+\rho_{i+1}^B\}}|\dot\phi_A(t)-\dot\phi_B(t)|^2\,dt\\
=&\biggl[\frac{\sin(w(\alpha_i+\rho_{i+1}^A))}{\rho_{i+1}^A}-\frac{\sin(w(\beta_i+\rho_{i+1}^B))}{\rho_{i+1}^B}\biggr]^2\\&\qquad\qquad\qquad\cdot(\min\{\alpha_i+\rho_{i+1}^A,\,\beta_i+\rho_{i+1}^B\}-\max\{\alpha_i,\,\beta_i\});
\end{align*}

\begin{align*}
&\int_{\min\{\alpha_i+\rho_{i+1}^A,\,\beta_i+\rho_{i+1}^B\}}^{\max\{\alpha_i+\rho_{i+1}^A,\,\beta_i+\rho_{i+1}^B\}}|\dot\phi_A(t)-\dot\phi_B(t)|^2\,dt\\
\leq&(\max\{\alpha_i+\rho_{i+1}^A,\,\beta_i+\rho_{i+1}^B\}-\min\{\alpha_i+\rho_{i+1}^A,\,\beta_i+\rho_{i+1}^B\})\|\dot\phi_A-\dot\phi_B\|_{L^\infty(0,\,T,\,\mathbb R)}^2\\
\leq& 4w^2\Bigl(\frac{1+\theta}{\theta}\Bigr)^2|(\alpha_i+\rho_{i+1}^A)-(\beta_i+\rho_{i+1}^B)|\leq 4\gamma w^2\Bigl(\frac{1+\theta}{\theta}\Bigr)^2\Bigl(1+\frac{2}{w}\Bigr);
\end{align*}

\begin{align*}
&\int_{\min\{T-\rho_m^A,\,T-\rho_m^B\}}^{\max\{T-\rho_m^A,\,T-\rho_m^B\}}|\dot\phi_A(t)-\dot\phi_B(t)|^2\,dt\\
\leq&(\max\{T-\rho_m^A,\,T-\rho_m^B\}-\min\{T-\rho_m^A,\,T-\rho_m^B\})\|\dot\phi_A-\dot\phi_B\|_{L^\infty(0,\,T,\,\mathbb R)}^2\\
\leq&4w^2\Bigl(\frac{1+\theta}{\theta}\Bigr)^2|\rho_m^A-\rho_m^B|=4w^2\Bigl(\frac{1+\theta}{\theta}\Bigr)^2\frac{1}{w}|x_m-y_m|\leq4\gamma w\Bigl(\frac{1+\theta}{\theta}\Bigr)^2;
\end{align*}

Finally for the last integral of \eqref{int.phidot.L2} we have
\begin{align*}
&\int_{\max\{T-\rho_m^A,\,T-\rho_m^B\}}^T|\dot\phi_A(t)-\dot\phi_B(t)|^2\,dt\\
=&\biggl[\frac{\sin(w(T-\rho_m^A))}{-\rho_m^A}-\frac{\sin(w(T-\rho_m^B))}{-\rho_m^B}\biggr]^2(T-\max\{T-\rho_m^A,\,T-\rho_m^B\}).
\end{align*}
Since $\beta_i\to\alpha_i$ and $\rho_m^B\to\rho_m^A$ as $\gamma\to 0$ we have
that all the intervals on the right-hand side of  \eqref{int.phidot.L2} go to zero as $\gamma$ does, i.e.,
\begin{equation}\label{cont.phidot.L2}
\|\dot\phi_A-\dot\phi_B\|_{L^2(0,\,T,\,\mathbb R)}\to 0\quad\text{as}\quad B\to A.
\end{equation}

Hence from \eqref{cont.phi.C} and \eqref{cont.phidot.L2} follows the $(\mathcal A_\theta,\,W^{1,2}(0,\,T,\,\mathbb R))$-continuity of $A\mapsto\phi_A$.\par
Now from the $(B,\,\mathcal A_\theta)$-continuity of the map $b\mapsto A(b)$ (which is equivalent to the $\delta$-continuity of the family $\mathcal Z$) and, from the $(\mathcal A_\theta,\,W^{1,2})$-continuity of $A\mapsto\phi_A$ we have the following corollary
\begin{Corollary}
The map $b\mapsto \phi_w(\cdot,\,b)$ is $(B,\,W^{1,2}(0,\,T,\,\mathbb R))$-continuous.
\end{Corollary}
\subsubsection{Imitation.}
Now we ``imitate'' the control $z(\cdot,\,b)\in\mathcal Z$ taking values in $\{\pm \Xi e_k,\,\pm\Xi\delta_{m,n}\mid\,k\in\mathcal K^{N-1},\,(m,\,n)\in S_{N-1}\}$ by a control $z^w(\cdot,\,b)$ taking values in $\mathbb R^{\kappa_{N-1}}$.\par Take the solution $u^\infty(\cdot,\,b)$ of the equation
$$
u_t^\infty(\cdot,\,b)=-\nu Au^\infty-Bu^\infty+F+z(\cdot,\,b),\quad u(0)=u_0
$$
and, consider its projection onto $\mathbb R^{\kappa_{N-1}}$:
$$
q^\infty(\cdot,\,b)=P^{\kappa_{N-1}}u^\infty(\cdot,\,b).
$$
Let $\{0=\alpha_0<\alpha_1<\dots<\alpha_m=T\}$ be the end-points of the intervals of constancy of $z(\cdot,\,b)$. For $w\geq 3$ define the control $z^w(\cdot,\,b)$ by recursion in the following way:\\
\begin{itemize}
\item In the first interval of constancy $[\alpha_0,\,\alpha_1]$:
$$
z^w(\cdot,\,b):=
\begin{cases}
&z(\cdot,\,b) \quad \text{if}\quad z(\cdot,\,b)\in\{\pm \Xi e_k\mid\,k\in\mathcal K^{N-1}\};\\
&v^{\kappa_{N-1}}(q_1^\infty(\cdot,\,b)+\sqrt{2\Xi}\phi^w(\cdot,\,b)(e_m\pm e_n),\,U_0)\quad\footnotemark\\
&\qquad\qquad\text{if}\quad z(\cdot,\,b)\in\{\pm\Xi\delta_{m,n}\mid\,(m,\,n)\in S_{N-1}\}.
\end{cases}
$$
\footnotetext{Here $v^{\kappa_{N-1}}$ is the control given by Lemma \ref{NS:q.to.v}, for $J=\mathcal K^{N-1}$.}
where $U_0$ is the projection of $u_0$ onto $\mathbb J_V^\bot=(\mathbb R^{\kappa_N})_V^\bot$ --- the orthogonal space to $\mathbb R^{\kappa_{N-1}}$ in $V$ and, $q_1^\infty(\cdot,\,b)$ is the restriction of $q^\infty(\cdot,\,b)$ to $[\alpha_0,\,\alpha_1]$;
\item If the control $z^w(\cdot,\,b)$ is already defined in the first $p-1$ intervals of constancy (up to $\alpha_{p-1}$), we define it in the $p^{th}$ interval $[\alpha_{p-1},\,\alpha_p]$ by:
$$
z^w(\cdot,\,b):=
\begin{cases}
&z(\cdot,\,b)\quad\text{if}\quad z(\cdot,\,b)\in\{\pm \Xi e_k\mid\,k\in\mathcal K^{N-1}\}\\
&v^{\kappa_{N-1}}(q_p^\infty(\cdot,\,b)+\sqrt{2\Xi}\phi^w(\cdot,\,b)(e_m\pm e_n),\,U^w(\alpha_{p-1}))\\
&\qquad\qquad\text{if}\quad z(\cdot,\,b)\in\{\pm\Xi\delta_{m,n}\mid\,(m,\,n)\in S_{N-1}\}.
\end{cases}
$$
where $U^w$ is the projection onto $(\mathbb R^{\kappa_N})_V^\bot$ of the solution of the equation
$$
u_t^w(\cdot,\,b)=-\nu Au^w-P^\nabla Bu^w+F+z^w(\cdot,\,b),\quad u(0)=u_0,\,t\in[0,\,\alpha_{p-1}]
$$
and, $q_p^\infty(\cdot,\,b)$ is the restriction of $q^\infty(\cdot,\,b)$ to $[\alpha_{p-1},\,\alpha_p]$.
\end{itemize}
We shall prove that at time $T$, $u^w(T)$ goes, uniformly w.r.t. $b$, to $u^\infty(T)$ in $\mathbf L^2$-norm as $w$ goes to $\infty$, i.e.,  
\begin{Lemma}\label{winf.Tclose}
For any $\varepsilon>0$ there exists $w_\varepsilon\geq 3$ such that
$$
\forall b\in B\forall w\geq w_\varepsilon\quad |u^w(T,\,b)-u^\infty(T,\,b)|<\varepsilon.
$$
\end{Lemma}
We claim that if the statement of the previous Lemma is true and, if we put $\gamma=\frac{1}{2}\Bigl(\|\partial\mathbb E_T(u_0,\,\mathcal Z),\,u_0\|_{l_1}-R\Bigr)$ then, $\mathbb E_T(u_0,\,\mathcal Z^{w_{\varepsilon}})$ with $\varepsilon=\frac{\gamma}{\#\mathcal O C}$ and $C$ is a constant such that $\|x\|_{l_1}\leq C\|x\|_H\quad (x\in\mathbb R^{\#\mathcal O})$, covers $\overline{\mathcal O}_R(u_0)$ solidly. Indeed, let $\phi$ be a continuous function defined on the closure of $B$ such that
$$
\|\phi-\mathbb E_T(u_0,\,\mathcal Z^{w_{\varepsilon}})\|_{C(B,\,\mathbb R_\infty^{\#\mathcal O})}<\frac{\gamma}{\#\mathcal O},
$$
then
\begin{align*}
&\|\phi-\mathbb E_T(u_0,\,\mathcal Z)\|_{C(B,\,\mathbb R_\infty^{\#\mathcal O})}\\
\leq &\|\phi-\mathbb E_T(u_0,\,\mathcal Z^{w_{\varepsilon}})\|_{C(B,\,\mathbb R_\infty^{\#\mathcal O})}+\|\mathbb E_T(u_0,\,\mathcal Z^{w_{\varepsilon}})-\mathbb E_T(u_0,\,\mathcal Z)\|_{C(B,\,\mathbb R_\infty^{\#\mathcal O})}\\
<&\frac{\gamma}{\#\mathcal O}+C\varepsilon<2\frac{\gamma}{\#\mathcal O}.
\end{align*}
For $x\in\overline{\mathcal O}_R(u_0)$ we obtain
\begin{align*}
\|\phi-\mathbb E(u_0,\,\mathcal Z)\|_{C(B,\,\mathbb R_\infty^{\#\mathcal O})}<&\frac{1}{\#\mathcal O}\|\mathbb E_T(u_0,\,\mathcal Z),\,u_0\|_{l_1}-R\leq \frac{1}{\#\mathcal O}\|\mathbb E_T(u_0,\,\mathcal Z),\,x\|_{l_1}\\
\leq& \|\mathbb E_T(u_0,\,\mathcal Z),\,x\|_{l_\infty}.
\end{align*}
Hence $deg(\phi,\,B,\,x)=deg(\mathbb E_T,\,B,\,x)\ne 0$ and so, $\phi(B)$ covers$\overline{\mathcal O}_R(u_0)$.\par
Therefore what remains to conclude the proof of the back-induction step and then, the solid controllability in Observed Component of system \eqref{infsysu} is the proof of Lemma \ref{winf.Tclose} that we present in the next section.
\section{Proof of Lemma \ref{winf.Tclose}.}
First we note that at the times $\alpha_i$ we have that the projections of $u^\infty(\alpha_i,\,b)$ and $u^w(\alpha_i,\,b)$ onto $\mathbb R^{\kappa_{N-1}}$ coincide. In fact, denoting the projection of $u^w(\cdot,\,b)$ by $q^w(\cdot,\,b)$ we have
$$
q^w(\cdot,\,b)=q^\infty(\cdot,\,b)+\sum_{(m,\,n)\in S_{N-1}}\phi_{m,n}^w(\cdot,\,b)
$$
where
$$
\phi_{m,n}^w(t,\,b):=
\begin{cases}
0 &\text{if}\quad v(t,\,b)\in\{\pm \Xi e_k\mid\,k\in\mathcal K^{N-1}\}\\
\sqrt{2\Xi}\phi^w(t,\,b)(e_m\pm e_n)&\text{if}\quad z(t,\,b)\in\{\pm\Xi\delta_{m,n}\mid\,(m,\,n)\in S_{N-1}\}.
\end{cases}
$$
Since $\phi^w(\cdot,\,b)$ vanishes at the points $\alpha_i$ also $\phi_{m,n}^w(t,\,b)$ does. Hence $q^\infty(\alpha_i,\,b)=q^w(\alpha_i,\,b)$ for all $i=0,\,1,\,\dots,\,m$. In particular for $i=m$ we obtain $q^\infty(T,\,b)=q^w(T,\,b)$ and so, we need only to compare the projections onto $(\mathbb R^{\kappa_{N-1}})_H^\bot$. For that we shall need some lemmas:
\begin{Lemma}\label{NS.dotbddL2}
The solution of the controlled NSE
$$
u_t=-\nu Au-P^\nabla Bu+F+v,\qquad u(0)=u_0,\quad t\in[0,\,T]
$$
satisfies 
$$
\|u_t\|_{L^2(0,\,T,\,H)}\leq C_v
$$
where $C_v$ is a constant depending only on the norm $\|v\|_{L^2(0,\,T,\,H)}$.\footnote{Recall that we have fixed $\nu$, $F$, $T$ and $u_0$. Otherwise the constant would depend on them.}\par
Moreover for a given family of controls $\{v_\beta\mid\,\beta\in\mathcal B\}$ bounded in $L^2(0,\,T,\,H)$-norm, i.e., if there is a constant $K>0$ such that for all $\beta \in\mathcal B$ $\|v_\beta\|_{L^2(0,\,T,\,H)}\leq K$, then we can find a constant $C_K$ depending only on $K$ such that
$$
\|u_t^\beta\|_{L^2(0,\,T,\,H)}\leq C_K
$$
where $u^\beta$ is the solution of
$$
u_t^\beta=-\nu Au^\beta-P^\nabla Bu^\beta+F+v_\beta,\qquad u(0)=u_0,\quad t\in[0,\,T].
$$
\end{Lemma}
\begin{proof}
Multiplying the equation by $u_t$ we obtain
\begin{align*}
|u_t|^2&=-\frac{\nu}{2}\frac{d}{dt}\|u\|^2-(P^\nabla Bu,\,u_t)+(F+v,\,u_t)\\
&\leq -\frac{\nu}{2}\frac{d}{dt}\|u\|^2+C\|u\||u|_{[2]}|u_t|+|F||u_t|+|v||u_t|
\end{align*}
so,
\begin{align}
\frac{1}{2}|u_t|^2&\leq-\frac{\nu}{2}\frac{d}{dt}\|u\|^2+C_1\|u\|^2|u|_{[2]}^2+C_1|F|^2+C_1|v|^2\notag\\
\int_0^T|u_t(t)|^2\,dt&\leq C_2\|u\|_{C([0,\,T],\,V)}^2(1+\int_0^T|u(t)|_{[2]}^2\,dt)\\
&+2C_1T|F|^2+2C_1\int_0^T|v|^2\,dt.\label{for.bd.dotu}
\end{align}
Multiplying the equation by $u$ and $Au$, analogously as we have obtained the ``a priori'' estimates in the proofs of existence of weak and strong solutions, we obtain the inequalities
\begin{align}
\|u\|_{C([0,\,T],\,H)}&\leq |u_0|^2+C_1\|F+v\|_{L^2(0,\,T,\,V^\prime)}^2\notag\\
\|u\|_{L^2(0,\,T,\,V)}&\leq 2\|u\|_{C([0,\,T],\,H)}+C_1{\nu}\|F+v\|_{L^2(0,\,T,\,V^\prime)}^2\notag\\
\|u\|_{C([0,\,T],\,V)}&\leq \exp(C_2\|u\|_{C([0,\,T],\,H)}\|u\|_{L^2(0,\,T,\,V)})\Bigl(\|u_0\|^2\notag\\
&+C_2\|F+v\|_{L^2(0,\,T,\,H)}^2\Bigr)\notag\\
\|u\|_{L^2(0,\,T,\,D(A))}&\leq 2\|u\|_{C([0,\,T],\,V)}+C_3\|u\|_{C([0,\,T],\,H)}^2\|u\|_{C([0,\,T],\,V)}^4\notag\\
&+C_3\|F+v\|_{L^2(0,\,T,\,H)}^2.\label{unifbounds}
\end{align}
from which, looking at \eqref{for.bd.dotu}, we conclude that $\|u_t\|_{L^2(0,\,T,\,H)}\leq C_v$ for some constant $C_v$ depending only on $\|v\|_{L^2(0,\,T,\,H)}$.\\
The statement relative to the family $v_\beta$ is also clearly true from the previous expressions:\\
\begin{align}
&\text{In each one of the previous estimates if we ``replace $\|v_\beta\|_{L^2(0,\,T,\,H)}$}\notag\\
&\text{by the bound $K$'' we obtain bounds independent of $\beta$}.\label{unifbd.uCL2}
\end{align}
\end{proof}
\begin{Corollary}\label{dotu.bddL2}
There is a constant $C_\infty>0$ such that
$$
\forall b\in B\quad \|u^\infty_t(\cdot,\,b)\|_{L^2(0,\,T,\,H)}\leq C_\infty.
$$
\end{Corollary}
\begin{proof}
For all $t\in [0,\,T],\,b\in B$ we have
$$
\|v(t,\,b)\|_{l_1}\leq\Xi\max\Bigl\{\|x\|_{l_1}\mid\,x\in\{e_k,\,\delta_{m,n}\mid\,k\in\mathcal K^{N-1},\,(m,\,n)\in S_{N-1}\}\,\Bigr\}
$$
so,  we have that for some constant $K_2$ 
$$
\|v(\cdot,\,b)\|_{L^\infty(0,\,T,\,H)}\leq K_2.
$$
The result follows by the last part of Lemma \ref{NS.dotbddL2} and the continuity of the inclusion $L^\infty(0,\,T,\,H)\to L^2(0,\,T,\,H)$.
\end{proof}
If $v(t,\,b)\in\{\pm \Xi e_k\mid\,k\in\mathcal K^{N-1}\}$ in the interval of constancy $[\alpha_{i-1},\,\alpha_i]$ let us call this interval of the {\bf first kind}. Otherwise, if $z(t,\,b)\in\{\pm\Xi\delta_{m,n}\mid\,(m,\,n)\in S_{N-1}\}$ on $[\alpha_{i-1},\,\alpha_i]$ we call the interval of the {\bf second kind}.\par
Now we note that for $u^w(\cdot,\,b)$ we can not find a bound for $\|u^w_t\|_{L^2(0,\,T,\,H)}$ independent of the parameter $w$ because the projection $q^w(\cdot,\,b)$ in an interval of constancy $I=[\alpha_{i-1},\,\alpha_i]$ of the second kind, say $v(t,b)=\delta_{m,n}$ on $I$, reads 
$$
q^w(t,\,b)=q^\infty(t,\,b)+\phi_{m,n}^w(t,\,b)
$$
and,
\begin{align*}
&\|\phi_t^w\|_{L^2(0,\,T,\,\mathbb R)}^2\\
\geq&\int_{[\alpha_{i-1}+\frac{L}{w},\,\alpha_i-\frac{L}{w}]}|\phi_t^w(t)|^2\,dt=w^2\int_{[\alpha_{i-1}+\frac{L}{w},\,\alpha_i-\frac{L}{w}]}\cos^2(wt)\,dt\\
=& w^2\biggl[\frac{1}{2}(L-2\frac{L}{w})+\frac{1}{4w}\Bigl(\sin(2w(\alpha_i-\frac{L}{w}))-\sin(2w(\alpha_{i-1}+\frac{L}{w}))\Bigr)\biggr].
\end{align*}
where $L$ is the length of $I$. Since $w\geq 3$ (and $L>0$) we obtain
$$
\|\phi_t^w\|_{L^2(0,\,T,\,\mathbb R)}^2\geq w^2\biggl[\frac{L}{6}+\frac{1}{4w}\Bigl(\sin(2w(\alpha_i-\frac{L}{w}))-\sin(2w(\alpha_{i-1}+\frac{L}{w}))\Bigr)\biggr].
$$
Then we see that when $w$ goes to $\infty$ also $\|\phi_t^w\|_{L^2(0,\,T,\,\mathbb R)}$ does and then, so does $\|u^w_t\|_{L^2(0,\,T,\,H)}$.\par
If we consider the projection $U^w(\cdot,b)$ of $u^w(\cdot,b)$ onto $(\mathbb R^{\kappa_{N-1}})_H^\bot$, there holds that we can find a bound for $\|U_t^w(\cdot,b)\|_{L^2(0,\,T,\,(\mathbb R^{\kappa_{N-1}})_H^\bot)}$ independent of both parameters $b$ and $w$. If fact multiplying by $U_t^w(\cdot,b)$ the equation
$$
U_t^w(\cdot,b)=-\nu AU^w(\cdot,b)-P^{-\kappa_{N-1}}P^\nabla Bu^w(\cdot,b)+P^{-\kappa_{N-1}}F
$$
that is satisfied by $U^w$ we obtain
\begin{multline}
|U_t^w(t,b)|^2\leq-\frac{\nu}{2}\frac{d}{dt}\|U^w(t,b)\|^2+|F||U^w_t(t,b)|\\
+C\|U^w(t,b)+q^w(t,b)\||U^w(t,b)+q^w(t,b)|_{[2]}|U^w_t(t,b)|.
\end{multline}
After using appropriate Young Inequalities and integrating
\begin{align}
&\int_0^T|U_t^w(t,b)|^2\,dt\notag\\
\leq& C_1\|U^w(\cdot,b)\|_{C([0,\,T],\,(\mathbb R^{\kappa_{N-1}})_V^\bot)}^2+C_1|F|^2\notag\\
+&C_1\|U^w(\cdot,b)+q^w(\cdot,b)\|_{C([0,\,T],\,V)}^2\|U^w(\cdot,b)+q^w(\cdot,b)\|_{L^2(0,\,T,\,D(A))}^2\notag\\
\leq& C\|U^w(\cdot,b)\|_{C([0,\,T],\,(\mathbb R^{\kappa_N})_V^\bot)}+T|F|^2\notag\\
+&C\biggl[\|U^w(\cdot,b)\|_{C([0,\,T],\,(\mathbb R^{\kappa_N})_V^\bot)}^2+\|q^w(\cdot,b)\|_{C([0,\,T],\,\mathbb R^{\kappa_{N-1}})}^2\biggr]\notag\\
&\qquad\qquad\times\biggl[\|U^w(\cdot,b)\|_{L^2(0,\,T,\,(\mathbb R^{\kappa_{N-1}})_{D(A)}^\bot)}^2+\|q^w(\cdot,b)\|_{L^2(0,\,T,\,\mathbb R^{\kappa_{N-1}})}^2\biggr].\label{for.bddotw}
\end{align}
Multiplying the equation by $U^w(\cdot,b)$ and by $AU^w(\cdot,b)$ we obtain that $U^w(\cdot,b)$ satisfy the estimates \eqref{qtoQ.bdH}, \eqref{qtoQ.bd2V}, \eqref{qtoQ.bdV} and \eqref{qtoQ.bd2D(A)} with $U^{J,L}$ replaced by $U^w(\cdot,b)$ and we can easily see that there is a constant $C_{w,b}$ depending only on $\|q^w(\cdot,b)\|_{C([0,\,T],\,\mathbb R^{\kappa_{N-1}})}$ such that
\begin{equation}\label{unifbd.UwCL2}
\|U^w(\cdot,b)\|_{C([0,\,T],\,(\mathbb R^{\kappa_N})_V^\bot)}\leq C_{w,b}\And \|U^w(\cdot,b)\|_{L^2(0,\,T,\,(\mathbb R^{\kappa_{N-1}})_{D(A)}^\bot)}\leq C_{w,b}.
\end{equation}
Moreover the family $\{U^w_t(\cdot,b)\}$ is uniformly bounded, w.r.t. $b$ and $w$, in the norm of $L^2(0,\,T,\,H)$ because, in the ``a priori''-like estimates above and in \eqref{for.bddotw} we can take bounds independent of the parameters depending only in the bound

\begin{align*}
&max_{(b,\,w)\in B\times[3,+\infty[}\{\|q^w(\cdot,b)\|_{C([0,\,T],\,\mathbb R^{\kappa_{N-1}})}\}\\
\leq&\|q^\infty(\cdot,b)\|_{C([0,\,T],\,\mathbb R^{\kappa_{N-1}})}+\Bigl\|\sum_{(m,\,n)\in S_{N-1}}\phi_{m,n}^w(\cdot,\,b)\Bigr\|_{C([0,\,T],\,\mathbb R^{\kappa_{N-1}})}\leq C
\end{align*}

where $C$ does not depend neither on $b$ nor on $w$. Indeed by \eqref{v.unif.bdd.Xi} the family $v(\cdot,\,b)$ is uniformly bounded in $L^2(0,\,T,\,\mathbb R^{\kappa_{N-1}})$-norm which implies that $q^\infty(\cdot,\,b)$ is uniformly bounded in $C([0,\,T],\,\mathbb R^{\kappa_{N-1}})$-norm. Indeed replacing $u$ by $u^\infty(\cdot,\,b)$ in the equations \eqref{unifbounds} we see that, since the family of controls $\{v(\cdot,b)\}$ is uniformly bounded in $L^2(0,\,T,\,H)$ norm, we can find a uniform bound $D$ for the solutions $u^\infty(\cdot,\,b)$ in anyone of the norms $C([0,\,T],\,H)$, $L^2(0,\,T,\,V)$, $C([0,\,T],\,V)$ and $L^2(0,\,T,\,D(A))$. In particular\\ $\|q^\infty(\cdot,\,b)\|_{C([0,\,T],\,\mathbb R^{\kappa_{N-1}})}\leq D_1\|u^\infty(\cdot,\,b)\|_{C([0,\,T],\,H)}\leq D_1D$.\\ On the other hand $\|\sum_{(m,\,n)\in S_{N-1}}\phi_{m,n}^w(\cdot,\,b)\|_{C([0,\,T],\,\mathbb R^{\kappa_{N-1}})}\}$ never exceeds the value $2\sqrt{2\Xi}$. The we can the write:
\begin{Corollary}\label{dotUwbddL2}
There is a constant $C>0$ such that for all $(w,\,b)$ in the product $[3,\,+\infty[\times B$:
$$
\|U_t^w(\cdot,\,b)\|_{L^2(0,\,T,\,\mathbb (R^{\kappa_{N-1}})_H^\bot)}\leq C.
$$
\end{Corollary}
Another Lemma we shall use is the following:
\begin{Lemma}\label{bd.W12.rx}
Let $\{z(\cdot,\,\sigma)\in W^{1,2}([t_i,\,t_f],\,\mathbb R)\mid\,\sigma\in\Sigma\}$ be a uniformly bounded, w.r.t. $\sigma$, family, i.e.,
$$
\exists C>0\,\forall\sigma\in\Sigma\quad\|z(\cdot,\,\sigma)\|_{W^{1,2}([t_i,\,t_f],\,\mathbb R)}\leq C.
$$
Then there is a constant $D_1$ depending only in $C$ and in the length $t_f-t_i$ of the interval $[t_i,\,t_f]$ and so, independent of $\sigma$ and $w$, such that 
$$
\|\sin(wt)z(t,\,\sigma)\|_{rx}\leq D_1w^{-1},\quad\text{and}\quad \|\cos(wt)z(t,\,\sigma)\|_{rx}\leq D_1w^{-1}.
$$
\end{Lemma}
\begin{proof}
The Lemma follows by direct computation: Let $s,\,r$ belong to $[t_i,\,t_f]$:
\begin{align*}
&\int_s^r \sin(wt)z(t,\,\sigma)\,dt\\
=&-w^{-1}\int_s^r -w\sin(wt)z(t,\,\sigma)\,dt=w^{-1}\int_s^r \bigl(\frac{d}{dt}\cos(wt)\bigr)z(t,\,\sigma)\,dt\\
=&-w^{-1}\biggl(\bigl[cos(wt)z(t,\,\sigma)\bigr]_s^r-\int_s^r cos(wt)\bigl(\frac{d}{dt}z(t,\,\sigma)\bigr)\,dt\biggr)\\
\leq&-w^{-1}\bigl[cos(wr)z(r,\,\sigma)-cos(ws)z(s,\,\sigma)\bigr]\\
+&w^{-1}\biggl(\int_s^r \cos^2(wt)\,dt\biggr)^\frac{1}{2}\biggl(\int_s^r \|\frac{d}{dt}z(t,\,\sigma)\|^2\,dt\biggr)^\frac{1}{2}\\
\leq& 2w^{-1}\|z(t,\,\sigma)\|_{C([t_i,\,t_f],\,\mathbb R)}+w^{-1}|r-s|^\frac{1}{2}C.
\end{align*}
By the continuity of the embedding $W^{1,2}\to C^0$ there is a constant $C_1$ such that $\|u\|_{C^0}\leq C_1\|u\|_{W^{1,2}}$. Putting $C_2=\max\{C,\,C_1C\}$ we obtain

$\int_s^r \sin(wt)z(t,\,\sigma)\,dt\leq w^{-1}C_2\bigl(2+(t_f-t_i)^\frac{1}{2}\bigr).$

Analogously we arrive to
$$
\int_s^r\cos(wt)z(t,\,\sigma)\,dt\leq w^{-1}C_2\bigl(2+(t_f-t_i)^\frac{1}{2}\bigr).
$$
Hence
$$
\forall\sigma\in\Sigma\,\forall s,\,r\in[t_i,\,t_f]
\begin{cases}
\int_s^r\sin(wt)z(t,\,\sigma)\,dt\leq  w^{-1}C_2\bigl(2+(t_f-t_i)^\frac{1}{2}\bigr)\\
\int_s^r\cos(wt)z(t,\,\sigma)\,dt\leq w^{-1}C_2\bigl(2+(t_f-t_i)^\frac{1}{2}\bigr).
\end{cases}
$$
Therefore
$$
\forall\sigma\in\Sigma
\begin{cases}
\|\sin(wt)z(t,\,\sigma)\|_{rx}\leq w^{-1}C_2\bigl(2+(t_f-t_i)^\frac{1}{2}\bigr)\\
\|\cos(wt)z(t,\,\sigma)\|_{rx}\leq w^{-1}C_2\bigl(2+(t_f-t_i)^\frac{1}{2}\bigr).
\end{cases}
$$
Choose $D_1=C_2\bigl(2+(t_f-t_i)^\frac{1}{2}\bigr)$.
\end{proof}
\begin{Corollary}\label{bd.W12.rx.phi}
Let $\{z(\cdot,\,\sigma)\in W^{1,2}([0,\,T],\,\mathbb R)\mid\,\sigma\in\Sigma\}$ be a uniformly bounded, w.r.t. $\sigma$, family:
$$
\exists C>0\,\forall\sigma\in\Sigma\quad\|z(\cdot,\,\sigma)\|_{W^{1,2}([0,\,T],\,\mathbb R)}\leq C.
$$
Then there is a constant $D_2$ depending only in $C$ and so, independent of $\sigma$ and of the parameter $b$ of our controls, such that 
$$
\|\phi_w(\cdot,\,b)z(\cdot,\,\sigma)\|_{rx}\leq D_2w^{-1}.
$$
\end{Corollary}
\begin{proof}
Let $s,\,r$ be in $[0,\,T]$ and suppose without loss of generality that $s<r$. Put $D:=\{t\in[0,\,T]\mid\,\phi_w(t,\,b)\ne\sin(wt)\}$. 
\begin{align*}
\int_s^r \phi_w(t,\,b)z(t,\,\sigma)\,dt&=\int_{D\cap[s,\,r]} \phi_w(t,\,b)z(t,\,\sigma)\,dt+\int_{[s,\,r]\setminus D} \sin(wt)z(t,\,\sigma)\,dt\\
&\leq 2\frac{T}{w}C_1C+mw^{-1}C_2\bigl(2+T^\frac{1}{2}\bigr)\\
&\leq w^{-1}C_2\bigl(2T+2m+T^\frac{1}{2}m\bigr).\quad\footnotemark
\end{align*}
\footnotetext{Note that $m$ being the number of intervals of constancy, then $[s,\,r]\setminus D$ is a union of at most $m$ intervals.}
Choose $D_2=C_2\bigl(2T+2m+T^\frac{1}{2}m\bigr)$.
\end{proof}
Now we compare the projections $U^w(\cdot,\,b)$ and $U^\infty(\cdot,\,b)$ of respectively $u^w(\cdot,\,b)$ and $u^\infty(\cdot,\,b)$ onto $(\mathbb R^{\kappa_{N-1}})_H^\bot$.  We claim that at time $\alpha_i,\quad(i=1,\,\dots,\,m)$ there holds
\begin{equation}\label{Pr.closealphai}
|U^w(\alpha_i,\,b)-U^\infty(\alpha_i,\,b)|\leq \overline {C_i}w^{-1}
\end{equation} 
where $\overline {C_i}$ is a constant independent of the parameters $w$ and $b$. In particular at time $T$ there holds
\begin{equation}\label{Pr.closeT}
|U^w(T,\,b)-U^\infty(T,\,b)|\leq \overline {C_m}w^{-1}.
\end{equation} 
Note that \eqref{Pr.closeT} implies Lemma \ref{winf.Tclose} because as we have seen at time $T$ we have $q^w(T,\,b)=q^\infty(T,\,b)$. To prove \eqref{Pr.closeT} we shall compare $U^w(\cdot,\,b)$ and $U^\infty(\cdot,\,b)$ in each interval of constancy:\\
\begin{itemize}
\item In an interval of the first kind  $U^w(\cdot,\,b)$ and $U^\infty(\cdot,\,b)$ satisfy the same equation
\begin{align*}
U_t^\infty(\cdot,\,b)&=-\nu AU^\infty(\cdot,\,b)-P^{-\kappa_{N-1}}P^\nabla \bigl(U^\infty(\cdot,\,b)+q^\infty(\cdot,\,b)\bigr)+P^{-\kappa_{N-1}} F;\\
U_t^w(\cdot,\,b)&=-\nu AU^w(\cdot,\,b)-P^{-\kappa_{N-1}}P^\nabla \bigl(U^w(\cdot,\,b)+q^\infty(\cdot,\,b)\bigr)+P^{-\kappa_{N-1}} F
\end{align*}
\item In an interval of the second kind, say $v(\cdot,\,b)=\pm\Xi\delta_{m,n}\quad (m,n)\in S_{N-1}$,  $U^w(\cdot,\,b)$ and $U^\infty(\cdot,\,b)$ satisfy the equations
\begin{align*}
U_t^\infty(\cdot,\,b)&=-\nu AU^\infty(\cdot,\,b)-P^{-\kappa_{N-1}}P^\nabla \bigl(U^\infty(\cdot,\,b)+q^\infty(\cdot,\,b)\bigr)+P^{-\kappa_{N-1}} F;\\
U_t^w(\cdot,\,b)&=-\nu AU^w(\cdot,\,b)-P^{-\kappa_{N-1}}P^\nabla \Bigl[U^w(\cdot,\,b)+q^\infty(\cdot,\,b)\\
&\qquad+\sqrt{2\Xi}\phi^w(\cdot,\,b)(e_m\pm e_n)\Bigr]+P^{-\kappa_{N-1}} F\pm\Xi\delta_{m,n}.
\end{align*}
\end{itemize}
Define the difference
$$
\eta^w(\cdot,\,b):=U^w(\cdot,\,b)-U^\infty(\cdot,\,b).
$$
For an interval $[\alpha_{i-1},\,\alpha_i]$  of second kind we find
\begin{align*}
\eta_t^w(s,\,b)=&-\nu A\eta^w(\cdot,\,b)+P^{-\kappa_{N-1}}P^\nabla \Bigl[U^\infty(\cdot,\,b)+q^\infty(\cdot,\,b)\Bigr]\\
&-P^{-\kappa_{N-1}}P^\nabla\Bigl[U^w(\cdot,\,b)+q^\infty(\cdot,\,b)+\sqrt{2\Xi}\phi^w(\cdot,\,b)(e_m\pm e_n)\Bigr)\\
&\pm\Xi P^{-\kappa_{N-1}}\delta_{m,n}.
\end{align*}
Multiplying by $\eta^w(\cdot,\,b)$ we obtain

\begin{align*}
&\frac{1}{2}\frac{d}{dt}|\eta_t^w(s,\,b)|^2\\
=&-\nu \|\eta^w(s,\,b)\|^2+\Bigl[ P^{-\kappa_{N-1}}P^\nabla \bigl(U^\infty(s,\,b)+q^\infty(s,\,b)\bigr),\,\eta^w(s,\,b)\Bigr]\\
&-\Bigl[ P^{-\kappa_{N-1}}P^\nabla \bigl(U^w(s,\,b)+q^\infty(s,\,b)+\sqrt{2\Xi}\phi^w(s,\,b)(e_m\pm e_n)\bigr),\,\eta^w(s,\,b)\Bigr]\\
&\pm\Bigl(\Xi P^{-\kappa_{N-1}}\delta_{m,n},\,\eta^w(s,\,b)\Bigr)
\end{align*}
so,
\begin{align}
&\frac{1}{2}\frac{d}{dt}|\eta^w(s,\,b)|^2+\nu \|\eta^w(s,\,b)\|^2\notag\\
\leq& -\Bigl( P^{-\kappa_{N-1}}P^\nabla \bigl(U^w(s,\,b)+q^\infty(s,\,b)\bigr),\,\eta^w(s,\,b)\Bigr)\notag\\&+\Bigl( P^{-\kappa_{N-1}}P^\nabla \bigl(U^\infty(s,\,b)+q^\infty(s,\,b)\bigr),\,\eta^w(s,\,b)\Bigr)\notag\\
&-\Bigl( P^{-\kappa_{N-1}}P^\nabla \bigl(U^w(s,\,b)+q^\infty(s,\,b),\,\sqrt{2\Xi}\phi^w(s,\,b)(e_m\pm e_n)\bigr),\,\eta^w(s,\,b)\Bigr)\notag\\
&-\Bigl( P^{-\kappa_{N-1}}P^\nabla \bigl(\sqrt{2\Xi}\phi^w(s,\,b)(e_m\pm e_n),\,U^w(s,\,b)+q^\infty(s,\,b)\bigr),\,\eta^w(s,\,b)\Bigr)\notag\\
&-\Bigl( P^{-\kappa_{N-1}}P^\nabla \bigl(\sqrt{2\Xi}\phi^w(s,\,b)(e_m\pm e_n)\bigr),\,\eta^w(s,\,b)\Bigr)\notag\\
&\pm\Bigl(\Xi P^{-\kappa_{N-1}}\delta_{m,n},\,\eta^w(s,\,b)\Bigr)\notag\\
=& A_1(s)+A_2(s)+A_3(s)+A_4(s)\label{eta2}
\end{align}
where we define
\begin{align*}
A_1(s):=&- \Bigl( P^{-\kappa_{N-1}}P^\nabla \bigl(U^w(s,\,b)+q^\infty(s,\,b),\,\eta^w(s,\,b)\Bigr)\\&+\Bigl( P^{-\kappa_{N-1}}P^\nabla \bigl(U^\infty(s,\,b)+q^\infty(s,\,b)\bigr),\,\eta^w(s,\,b)\Bigr)\\
A_2(s):=&\\
-\Bigl( P&^{-\kappa_{N-1}}P^\nabla \bigl(U^w(s,\,b)+q^\infty(s,\,b),\,\sqrt{2\Xi}\phi^w(s,\,b)(e_m\pm e_n)\bigr),\,\eta^w(s,\,b)\Bigr)\\
A_3(s):=&\\
-\Bigl( P&^{-\kappa_{N-1}}P^\nabla \bigl(\sqrt{2\Xi}\phi^w(s,\,b)(e_m\pm e_n),\,U^w(s,\,b)+q^\infty(s,\,b)\bigr),\,\eta^w(s,\,b)\Bigr)\\
A_4(s):=&-\Bigl( P^{-\kappa_{N-1}}P^\nabla \bigl(\sqrt{2\Xi}\phi^w(s,\,b)(e_m\pm e_n)\bigr),\,\eta^w(s,\,b)\Bigr)\\
&\pm\Bigl(\Xi P^{-\kappa_{N-1}}\delta_{m,n},\,\eta^w(s,\,b)\Bigr).
\end{align*}
To estimate $|\eta^w(z,\,b)|^2\quad z\in [\alpha_{i-1},\,\alpha_i]$ we will integrate \eqref{eta2}. In the second member, for $A_1$ we have:
\begin{equation}\label{bd.eta2A1}
|A_1|\leq C|\eta^w(z,\,b)|\|\eta^w(z,\,b)\|\|U^\infty(z,\,b)+q^\infty(z,\,b)\|
\end{equation}
and by Young Inequality we arrive to
\begin{align*}
&\frac{d}{dt}|\eta^w(z,\,b)|^2+\nu \|\eta^w(z,\,b)\|^2\\
\leq& D|\eta^w(z,\,b)|^2\|U^\infty(z,\,b)+q^\infty(z,\,b)\|^2+2\sum_{r=2}^4 A_r(z)
\end{align*}
and, by Gronwall Inequality we obtain
\begin{align}
|\eta^w(z,\,b)|^2\notag
\leq& |\eta^w(\alpha_{i-1},\,b)|^2\exp\Bigl(D\int_{\alpha_{i-1}}^{\alpha_i}\|U^\infty(s,\,b)+q^\infty(s,\,b)\|^2\,ds\Bigr)\\
+&2\sum_{r=2}^4 \int_{\alpha_{i-1}}^z A_r(s)E(s)\,ds\label{eta2.2}
\end{align}
with $E(s):=\exp\int_z^s -D\|u^\infty(t,\,b)\|^2\,dt$.
For $r=2$ we have
\begin{align}
&\int_{\alpha_{i-1}}^{z}A_2(s)E(s)\,ds=\notag\\
-&\sqrt{2\Xi}\int_{\alpha_{i-1}}^z\phi^w(s,\,b)\Bigl( P^{-\kappa_{N-1}}P^\nabla \bigl(U^w(s,\,b)+q^\infty(s,\,b),\,e_m\bigr),\,\eta^w(s,\,b)\Bigr)E(s)\,ds\notag\\
\mp& \sqrt{2\Xi}\int_{\alpha_{i-1}}^z\phi^w(s,\,b)\Bigl( P^{-\kappa_{N-1}}P^\nabla \bigl(U^w(s,\,b)+q^\infty(s,\,b),\,e_n\bigr),\,\eta^w(s,\,b)\Bigr)E(s)\,ds.\label{intA2}
\end{align}
Now we estimate the derivative of the product
$$
\Bigl( P^{-\kappa_{N-1}}P^\nabla \bigl(U^w(s,\,b)+q^\infty(s,\,b),\,e_m\bigr),\,\eta^w(s,\,b)\Bigr)E(s):
$$
\begin{align*}
&\frac{d}{ds}\Bigl( P^{-\kappa_{N-1}}P^\nabla \bigl(U^w(s,\,b)+q^\infty(s,\,b),\,e_m\bigr),\,\eta^w(s,\,b)\Bigr)\\
=&\Bigl( P^{-\kappa_{N-1}}P^\nabla \bigl(\dot U^w(s,\,b)+\dot q^\infty(s,\,b),\,e_m\bigr),\,\eta^w(s,\,b)\Bigr)\\
&+\Bigl( P^{-\kappa_{N-1}}P^\nabla \bigl(U^w(s,\,b)+q^\infty(s,\,b),\,e_m\bigr),\,\dot\eta^w(s,\,b)\Bigr)\\
\leq& C|\dot U^w(s,\,b)+\dot q^\infty(s,\,b)||e_m|_{[2]}\|\eta^w(s,\,b)\|\\
+&C\|U^w(s,\,b)+q^\infty(s,\,b)\||e_m|_{[2]}|\dot\eta^w(s,\,b)|\\
\leq& D|\dot U^w(s,\,b)+\dot q^\infty(s,\,b)|\|\eta^w(s,\,b)\|+D\|U^w(s,\,b)+q^\infty(s,\,b)\||\dot\eta^w(s,\,b)|
\end{align*}
so,
\begin{align*}
&\int_{\alpha_{i-1}}^z|\frac{d}{ds}\Bigl( P^{-\kappa_{N-1}}P^\nabla \bigl(U^w(s,\,b)+q^\infty(s,\,b),\,e_m\bigr),\,\eta^w(s,\,b)\Bigr)|^2\,ds\\
\leq&D_1\int_{\alpha_{i-1}}^{\alpha_i}|\dot U^w(s,\,b)+\dot q^\infty(s,\,b)|^2\|\eta^w(s,\,b)\|^2\,ds\\+&D_1\int_{\alpha_{i-1}}^{\alpha_i}\|U^w(s,\,b)+q^\infty(s,\,b)\|^2|\dot\eta^w(s,\,b)|^2\,ds
\end{align*}
and, for $\Bigl|\Bigl( P^{-\kappa_{N-1}}P^\nabla \bigl(U^w(s,\,b)+q^\infty(s,\,b),\,e_m\bigr),\,\eta^w(s,\,b)\Bigr)\Bigr|$ we obtain the bound
\begin{align*}
&\Bigl|\Bigl( P^{-\kappa_{N-1}}P^\nabla \bigl(U^w(s,\,b)+q^\infty(s,\,b),\,e_m\bigr),\,\eta^w(s,\,b)\Bigr)\Bigr|\\
\leq& C\|U^w(s,\,b)+q^\infty(s,\,b)\|\|e_m\|\|\eta^w(s,\,b)\|\leq D\|U^w(s,\,b)+q^\infty(s,\,b)\|\|\eta^w(s,\,b)\|.
\end{align*}
Hence, by \eqref{unifbd.uCL2} and \eqref{unifbd.UwCL2} and Corollaries \ref{dotu.bddL2}, \ref{dotUwbddL2} we conclude that  the family
$$
\Bigl\{\Bigl( P^{-\kappa_{N-1}}P^\nabla \bigl(U^w(s,\,b)+q^\infty(s,\,b),\,e_m\bigr),\,\eta^w(s,\,b)\Bigr)\Bigr\}
$$
is uniformly bounded, w.r.t. $w$ and $b$, in $W^{1,2}([\alpha_{i-1},\,\alpha_i])$-norm.\par
For $E(s)$ we find
\begin{align}
|E(s)|\leq \exp(TD\|u^\infty(\cdot,\,b)\|_{C([0,\,T],\,V)})\leq C_1\label{bd.EC}\\
|\frac{d}{ds}E(s)|\leq \|u^\infty(\cdot,\,b)\|_{C([0,\,T],\,V)}C_1\leq C_2\label{bd.EdotLinf}
\end{align}
where $C_1$ and $C_2$ do not depend neither on $b$ nor on $s$. In particular $E$ is uniformly bounded in $W^{1,\infty}([\alpha_{i-1},\,\alpha_i],\,\mathbb R)$-norm and then also the family
$$
\Bigl\{\Bigl( P^{-\kappa_{N-1}}P^\nabla \bigl(U^w(s,\,b)+q^\infty(s,\,b),\,e_m\bigr),\,\eta^w(s,\,b)\Bigr)E(s)\Bigr\}
$$
 is uniformly bounded in $W^{1,2}([\alpha_{i-1},\,\alpha_i],\,\mathbb R)$-norm.\par
Therefore by \eqref{intA2} and Corollary \ref{bd.W12.rx.phi} there is a constant $K_1$, independent of $z$ and of the parameters $w$ and $b$, such that
\begin{align*}
&\biggl|\int_{\alpha_{i-1}}^z\phi^w(s,\,b)\Bigl( P^{-\kappa_{N-1}}P^\nabla \bigl(U^w(s,\,b)+q^\infty(s,\,b),\,e_m\bigr),\,\eta^w(s,\,b)\Bigr)E(s)\,ds\biggr|\\
\leq& K_1w^{-1}.
\end{align*}
If we replace $m$ by $n$ we obtain a similar estimate and then, we conclude that for some constant $K_2$ independent of $z$, $w$ and $b$:
\begin{equation}\label{bd.eta2A2}
\biggl|\int_{\alpha_{i-1}}^z A_2(s)E(s)\,ds\biggr|\leq K_2w^{-1}.
\end{equation}
Analogously we conclude that for some constant $K_3$ independent of $z$, $w$ and $b$:
\begin{equation}\label{bd.eta2A3}
\biggl|\int_{\alpha_{i-1}}^z A_3(s)E(s)\,ds\biggr|\leq K_3w^{-1}.
\end{equation}
Recalling equation \eqref{PHunu} and supposing without loss of generality that $m<n$, we have
\begin{align*}
P^\nabla B\bigl[\sqrt{2\Xi}\phi^w(s,\,b)(e_m\pm e_n)\bigr]&=2\Xi\bigl(\phi^w(s,\,b)\bigr)^2 P^\nabla B\bigl[(e_m\pm e_n)\bigr]\\
&=\pm 2\Xi\bigl(\phi^w(s,\,b)\bigr)^2 \delta_{m,n}.
\end{align*}
Therefore
\begin{align*}
&\int_{\alpha_{i-1}}^z A_4(s)E(s)\,ds\\
=&\int_{\alpha_{i-1}}^z \bigl[\mp 2\Xi\bigl(\phi^w(s,\,b)\bigr)^2\pm\Xi\bigr]\Bigl(P^{-\kappa_{N-1}}\delta_{m,n},\,\eta^w(s,\,b)\Bigr)E(s)\,ds\\
 =&\int_{\alpha_{i-1}}^z \bigl[\mp 2\Xi\bigl(\phi^w(s,\,b)\bigr)^2\pm\Xi\bigr]\Bigl(\delta_{m,n},\,\eta^w(s,\,b)\Bigr)E(s)\,ds\\
=&\mp\Xi\int_{\alpha_{i-1}}^z\bigl[2\bigl(\phi^w(s,\,b)\bigr)^2-1\bigr]\Bigl(\delta_{m,n},\,\eta^w(s,\,b)\Bigr)E(s)\,ds\\
\leq&\mp\Xi\int_I\cos(2wt)\Bigl(\delta_{m,n},\,\eta^w(s,\,b)\Bigr)E(s)\,ds\\
&\mp\Xi\int_{[\alpha_{i-1},\,z]\setminus I}|(\delta_{m,n},\,\eta^w(s,\,b)||E(s)|\,ds.
\end{align*}
where $I=[\alpha_{i-1},\,z]\cap[\alpha_{i-1}+\frac{L_i}{w},\,\alpha_i-\frac{L_i}{w}]$ and $L_i=\alpha_i-\alpha_{i-1}$.\\
By \eqref{bd.EC} we have that
\begin{align*}
&\int_{[\alpha_{i-1},\,z]\setminus I}|(\delta_{m,n},\,\eta^w(s,\,b)||E(s)|\,ds\\
\leq& C_1|(\delta_{m,n},\,\eta^w(s,\,b))|_{C([0,\,T],\,\mathbb R)}\frac{2L_i}{w}\leq w^{-1}K_0\|\eta^w(s,\,b))\|_{C([0,\,T],\,(\mathbb R^{\kappa_{N-1}})_H^\bot)}
\end{align*}
where the constant $K_0$ can be taken independent of $(m,\,n)\in S_{N-1}$, because $S_{N-1}$ is finite.\\
For the derivative of the product $(\delta_{m,n},\,\eta^w(s,\,b))$ we find
$$
\frac{d}{dt}(\delta_{m,n},\,\eta^w(s,\,b))=(\delta_{m,n},\,\dot\eta^w(s,\,b))\leq K_1|\dot\eta^w(s,\,b)|
$$
so, by \eqref{unifbd.uCL2}, \eqref{unifbd.UwCL2} and corollaries \ref{dotu.bddL2} and \ref{dotUwbddL2} we conclude that for some constant $C_2$ independent of $w$ and $b$: $\|(\delta_{m,n},\,\eta^w(\cdot,\,b))\|_{W^{1,2}([0,\,T],\,\mathbb R)}\leq C_2$. So by the uniform boudedness of $E(\cdot)$ in $W^{1,\infty}([0,\,T],\,\mathbb R)$ we conclude the  uniform boudedness of $(\delta_{m,n},\,\eta^w(s,\,b))E(s)$ in $W^{1,2}([0,\,T],\,\mathbb R)$.\par
By Lemma \ref{bd.W12.rx} there is a constant $C_1$, independent of $w$ and $b$, such that
$$
\biggl|\int_I\cos(2wt)\Bigl(\delta_{m,n},\,\eta^w(s,\,b)\Bigr)E(s)\,ds\biggr|\leq C_1(2w)^{-1}.
$$
Hence
$$
\biggl|\int_{\alpha_{i-1}}^z A_4(s)E(s)\,ds\biggr|\leq \frac{C_1}{2}w^{-1}+w^{-1}K_0\|\eta^w(s,\,b)\|_{C([0,\,T],\,(\mathbb R^{\kappa_{N-1}})_H^\bot)}
$$
and, again by \eqref{unifbd.uCL2} and \eqref{unifbd.UwCL2} we arrive to
\begin{equation}\label{bd.eta2A4}
\biggl|\int_{\alpha_{i-1}}^z A_4(s)E(s)\,ds\biggr|\leq K_2w^{-1}.
\end{equation}
where $K_2$ is independent of $b$ and $w$.\par
By equation \eqref{eta2.2} and estimates \eqref{bd.eta2A2}, \eqref{bd.eta2A3} and \eqref{bd.eta2A4} we obtain
\begin{align*}
&|\eta^w(z,\,b)|^2\\
\leq& |\eta^w(\alpha_{i-1},\,b)|^2\exp\Bigl(D\int_0^T\|u^\infty(s,\,b)\|^2\,ds\Bigr)+2\sum_{r=2}^4\int_{\alpha_{i-1}}^z A_r(s)E(s)\,ds\\
\leq& D_1|\eta^w(\alpha_{i-1},\,b)|^2+D_2w^{-1}.
\end{align*}
In an interval of the first kind we just have $A_2=A_3=A_4=0$ on \eqref{eta2}
and then we obtain
\begin{align*}
&|\eta^w(z,\,b)|^2\\
\leq& |\eta^w(\alpha_{i-1},\,b)|^2\exp\Bigl(D\int_0^T\|U^\infty(s,\,b)+q^\infty(s,\,b)\|^2\,ds\Bigr)\leq D_1|\eta^w(\alpha_{i-1},\,b)|^2.
\end{align*}
Therefore in any interval of constancy (either of first or second kind) we have
\begin{equation}\label{incr}
|\eta^w(z,\,b)|^2\leq D_1|\eta^w(\alpha_{i-1},\,b)|^2+D_2w^{-1}
\end{equation}
where $D_1$ and $D_2$ are independent of $z\in[\alpha_{i-1},\,\alpha_i]$, of the parameters $w$ and $b$ and of $i\in\{1,\,2,\,\dots,\,m\}$, i.e., of the interval of constancy.\par
We prove \eqref{Pr.closealphai} by induction on $i$:
\begin{itemize}
\item $i=1$: At time $\alpha_1$, by \eqref{incr}, we have
$$
|U^w(\alpha_1,\,b)-U^\infty(\alpha_1,\,b)|^2\leq D_2w^{-1}\quad\footnotemark
$$
\footnotetext{Remember that both $u^\infty$ and $u^w$ start at $u_0$ at time $0$.}
Choose $\overline {C_1}=D_2$.
\item Induction Step: Suppose that at time $\alpha_{i-1}$ we have
$$
|U^w(\alpha_{i-1},\,b)-U^\infty(\alpha_{i-1},\,b)|^2\leq \overline{C_{i-1}}w^{-1}.$$
Then by \eqref{incr}, we have
\begin{multline*}
|U^w(\alpha_{i-1},\,b)-U^\infty(\alpha_{i-1},\,b)|^2\\
\leq D_1\bigl(\overline{C_{i-1}}w^{-1}\bigr)^2+D_2w^{-1}
=\bigl[D_1\overline{C_{i-1}}^2w^{-1}+D_2\bigr]w^{-1}\\
\leq \bigl[D_1\overline{C_{i-1}}^2\frac{1}{3}+D_2\bigr]w^{-1}.
\end{multline*}
Choose $\overline {C_i}=D_1\overline{C_{i-1}}^2\frac{1}{3}+D_2$.
\end{itemize}
 Therefore \eqref{Pr.closealphai} holds.

\chapter{$\mathbf L^2$-Approximate Controllability.}\label{L2AC}
The following Proposition says that for any $T>0$, system \eqref{infsysu} is {\bf time-$T$ approximately controllable} in $\mathbf L^2$-norm.
\begin{Proposition}
For any $u_0\in V$ and $T>0$, the attainable set at time $T$ from $u_0$ of system \eqref{infsysu} is dense in $H$.
\end{Proposition} 
\begin{proof}
Fix $\varepsilon>0$, $u_0\in V$, $x_1\in H$ and $T>0$. We prove that it is possible to drive the system from $u_0$ to the ball $\{y\in H\mid\,|y-x_1|\leq\varepsilon\}$. For that, first we set $M\in\mathbb N_0$ such that $|x_1-P^{\kappa_M}x_1|<\frac{\varepsilon}{2}$, where $P^{\kappa_M}$ is the projection map onto $\mathbb R^{\kappa_M}$ and, consider the system [\eqref{sysN}.M]. As we have seen in the first step of the proof of Proposition \ref{solcopstN} (section \ref{S:1step}), there is a control $v_M$ taking values on $\mathbb R^{\kappa_M}$ and driving the system from $u_0$ to some point $u_T^M$ such that $P^{\kappa_M}u_T^M=P^{\kappa_M}x_1$. Moreover by item \ref{rt.comp.tr.ini} of Corollary \ref{rt.comp.tr} we have
$$
|u_T^M-P^{\kappa_M}u_T^M|\leq K[T\exp(T)]^\frac{1}{2}.
$$ 
Now (see Lemma \ref{winf.Tclose})  we imitate $v_M$ by another control $v_{M-1}$ taking values in  $\mathbb R^{\kappa_{M-1}}$ and driving the system to a point $u_T^{M-1}\in V$ such that 
$$
|u_T^{M-1}-u_T^M|\leq \frac{\varepsilon}{2M}.
$$
Repeating this procedure of imitation, at each step we find a control $v_{i-1}$ taking values in  $\mathbb R^{\kappa_{i-1}}$ and driving the system to a point $u_T^{i-1}\in V$ such that 
$$
|u_T^{i-1}-u_T^i|\leq \frac{\varepsilon}{2M}.
$$
Let $T\leq T_1$ where $T_1$ equals the unique solution of $T_1\exp(T_1)=\Bigl(\frac{\varepsilon}{2KM}\Bigr)^2$. The control $v^1$ takes its values on $\mathbb R^{\kappa_1}$ and drives the system to a point $u_T^1\in V$ satisfying
$$
|u_T^1-x_1|\leq \sum_{i=2}^{M}|u_T^{i-1}-u_T^i|+|u_T^M-P^{\kappa_M}u_T^M|+|P^{\kappa_M}u_T^M-x_1|.
$$
Since $P^{\kappa_M}u_T^M=P^{\kappa_M}x_1$ we have
$$
|u_T^1-x_1|\leq \frac{\varepsilon}{2M}(M-1)+K[T\exp(T)]^\frac{1}{2}+\frac{\varepsilon}{2}\leq\varepsilon.
$$
Hence the Proposition is proved for $T\leq T_1$. If $T>T_1$ we apply an arbitrary control up to time $T-T_1$ arriving to some point $y_0\in V$ and then drive the system in time $T_1$ from $y_0$ to the ball $\{y\in H\mid\,|y-x_1|\leq\varepsilon\}$.
\end{proof}


\begin{thebibliography}{99}
\bibitem{agm}{S. Agmon,
            \emph{Lectures on Elliptic Boundary Value Problems}, van Nostrand 1965.}
\bibitem{bas}{A. A. Agrachev, Y. L. Sachkov
            \emph{Control Theory from the Geometric Viewpoint}, Encyclopaedia of Mathematical Sciences, 87, Springer,2004}
\bibitem{pas}{A. A. Agrachev, A. V. Sarychev,
            \emph{Navier-Stokes Equations: Controllability by Means
of Low Modes Forcing}, To appear on the Journal of Mathematical Fluid Mechanics.}
\bibitem{bre}{H. Brezis,\emph{ Analyse Fonctionnelle, Th\'{e}orie et Applications},
Masson, 1993.}
\bibitem{ency}{C. Foias, O. Manley, R. Rosa, R. Temam,\emph{ Navier-Stokes Equations and Turbulence}, Encyclopedia of Mathematics and its Applications, Cambridge university Press, 2001.}
\bibitem{deg.th}{I. Fonseca, W. Gangbo, \emph{Degree Theory in Analysis and Applications}; Oxford Lectures Series in Mathematics and its Applications, Oxford University Press, 1995.}
\bibitem{gam}{R.V. Gamkrelidze, \emph{Principles of Optimal Control Theory}, Plenum Press, 1978.}
\bibitem{jur}{V. Jurdjevic, \emph{Geometric Control Theory}, Cambridge Studies in Advanced Mathematics 51, Cambridge University Press, 1997.}
\bibitem{lio}{P.L. Lions, E. Magenes, \emph{Non-Homogeneous Boundary Value Problems and Applications}, vol.I, Die Grundlehren der Mathematischen Wissenschaften in Einzeldarstellungen, band 181, Springer-Verlag, 1972.}
\bibitem{pee}{J. Peetre, \emph{Espaces d'Interpolation et th\'{e}or\`{e}me de Soboleff}, Annales de l'Institut Fourier, {\bf 279-317 }, 1978.}
\bibitem{renrog}{M. Renardy, R.C. Rogers,\emph{ An Introduction to Partial Differential Equations}, Texts in Applied Mathematics, 13, Springer-Verlag, 1993.}
\bibitem{stric}{R.S. Strichartz,\emph{ A Guide to Distribution Theory and Fourier Transformations}, World Scientific, 2003.}
\bibitem{temsp68}{R. Temam,\emph{ Infinite-Dimensional dynamical Systems in Mechanics and Physics}, 2nd ed., Applied Mathematical Sciences, 68, Springer, 1997.}
\bibitem{temnsf}{R. Temam,\emph{ Navier-Stokes Equations and Nonlinear Functional Analysis}, 2nd Ed.,CBMS-NSF Regional Conference Series in Applied Mathematics, Society for Industrial and Applied Mathematics, 1995.}
\bibitem{temams}{R. Temam,\emph{ Navier-Stokes Equations: Theory and Numerical Analysis}, AMS Chelsea Publishing, 2001.}
\bibitem{wei.mat}{E. Weinam, J.C. Mattingly \emph{Ergodicity for the Navier-Stokes Equation with Degenerate Random Forcing: Finite Dimensional Approximation}, Comm. on Pure and Applied Math., Vol. 54, {\bf 1386-1402}, 2001.}

\end{thebibliography}
\end{document}